\documentclass[a4paper,twoside,leqno]{article}

\title{Non annulation des fonctions $L$ des formes modulaires de Hilbert en le point 
central}
\date{}
\author{Denis Trotabas}

\usepackage{amsmath,latexsym,amssymb,amsfonts,mathrsfs,bbm}
\usepackage[polutonikogreek,french, english]{babel}
\usepackage{fancyhdr}
\usepackage[colorlinks=true, linkcolor=black,citecolor=black,
urlcolor=black]{hyperref}
\newcommand{\HH}{\mathbbmss{H}}
\newcommand{\Q}{\mathbbm{Q}}
\newcommand{\q}{\mathfrak{q}}

\newcommand{\entier}{\mathcal{O}_F}
\newcommand{\unite}{\mathcal{O}_F^{\times+}}
\newcommand{\sym}{{\rm{sym}}^2}
\newcommand{\A}{\mathfrak{a}}
\newcommand{\B}{\mathfrak{b}}
\newcommand{\adele}{\mathbbmss{A}_F}
\newcommand{\gl}{{\rm{GL}_2}}
\newcommand{\pgl}{{\rm{PGL}_2}}
\newcommand{\GL}{{\rm{GL}_2^{+}}}
\newcommand{\cc}{\mathfrak{c}}
\newcommand{\kl}{\mathcal{K}\!\ell}
\newcommand{\N}{\mathbbm{N}}
\newcommand{\Z}{\mathbbm{Z}}
\newcommand{\n}{\mathfrak{n}}
\newcommand{\M}{{\rm{M}}}
\newcommand{\m}{\mathfrak{m}}
\newcommand{\p}{\mathfrak{p}}
\newcommand{\D}{\mathfrak{d}}
\newcommand{\de}{\mathfrak{D}}
\newcommand{\E}{\mathfrak{e}}
\newcommand{\res}{{\rm{res}}}
\newcommand{\e}{{\rm{e}}}
\newcommand{\disc}{\mathfrak{d}_F}
\newcommand{\diff}{\mathfrak{D}_F}
\newcommand{\R}{\mathbbm{R}}
\newcommand{\C}{\mathbbm{C}}
\newcommand{\f}{\mathfrak{f}}
\newcommand{\Whitt}{{\mathcal{W}}(\pi,\psi)}
\newcommand{\cond}{ {\rm{{q}}}_\pi}
\newcommand{\g}{\boldsymbol}
\newcommand{\PM}{\widehat{P_M}}
\newcommand{\Dir}{\mathscr{D}_\alpha({\rm{sym}}^2\pi_f)}
\newcounter{rem}
\def\rem{\addtocounter{rem}{1}{\noindent\bf{Remarque} \therem}\ }
\newcommand{\preuve}{{\sc{Preuve}}}
\newcommand{\Mo}{\mathscr{M}}
\def\build#1_#2^#3{\mathrel{\mathop{\kern 0pt#1}\limits_{#2}^{#3}}} 
\def\equi_#1{\build{\sim}_#1^{}}
\def\Equi#1#2{\equi_{\scriptscriptstyle#1\to #2}}

\newtheorem{theoreme}{\bf Th\'eor\`eme}
\newtheorem{definition}{\bf{D\'efinition}}
\newtheorem{lemme}{\bf{Lemme}}
\newtheorem{proposition}{\bf{Proposition}}
\newtheorem{corollaire}{\bf{Corollaire}}

\begin{document}
\maketitle
\thispagestyle{empty}
\selectlanguage{english}
\begin{abstract}
Birch and Swinnerton-Dyer conjecture allows for sharp estimates on the rank of certain abelian varieties defined over $ \Q$. in the case of the jacobian of the modular curves, this problem is equivalent to the estimation of the order of vanishing at $1/2$ of $L$-functions of classical modular forms, and was treated, without assuming the Riemann hypothesis, by Kowalski, Michel and VanderKam. The purpose of this paper is to extend this approach in the case of an arbitrary totally real field, which necessitates an appeal of Jacquet-Langlands' theory and the adelization of the problem. To show that the $L$-function (resp. its derivative) of a positive density of forms does not vanish at $1/2$, we follow Selberg's method of mollified moments (Iwaniec, Sarnak, Kowalski, Michel and VanderKam among others applied it successfully in the case of classical modular forms). We generalize the Petersson formula, and use it to estimate the first two harmonic moments, this then allows us to match the same unconditional densities as the ones proved over $\Q$ by Kowalski, Michel and VanderKam. In this setting, there is an additional term, coming from old forms, to control. Finally we convert our estimates for the harmonic moments into ones for the natural moments.

{\bf MSC:} 11F41, 11M41, 11F70.
\end{abstract}

\selectlanguage{french}

\tableofcontents

\section{Introduction et r\'esultats}

Soit $F/\,\Q$ une extension finie de degr\'e $d$,  totalement r\'eelle, d'anneau 
d'entiers $\entier$, et soit $\q$ un id\'eal premier de $\entier$. Les 
repr\'esentations automorphes cuspidales de caract\`ere central trivial de $\gl(\adele)$ sont les facteurs 
irr\'eductibles de l'action de $\gl(\adele)$ sur $L^2_0(\gl(F)Z({\adele})\backslash 
\gl(\adele))$. On notera $(\pi,V_\pi)$ ou simplement $\pi$ un tel constituant, et on sait que l'on 
a une factorisation: $\pi \cong \widehat{\bigotimes}_v \pi_v$, $v$ parcourant l'ensemble de toutes 
les places de $F$, chaque $\pi_v$ \'etant une repr\'esentation irr\'eductible de $\gl(F_v)$ 
uniquement d\'etermin\'ee par cet isomorphisme. En s\'eparant les places infinies et finies, 
on \'ecrit: $\pi\cong \pi_\infty\otimes \pi_f$, et on dit que $\pi$ est une forme modulaire 
de Hilbert de poids $\g{k}$ s'il existe $\g{k}=(k_j)_j\in 2\N_{\geq 1}^d$ tel que \begin{displaymath}\pi_\infty\cong \bigotimes_{j=1}^d \mathcal{D}
(k_j-1),\end{displaymath} produit de s\'eries discr\`etes de caract\`ere central trivial,
 et de param\`etres $k_j-1$.

Pour $F=\Q$, cela \'equivaut aux formes modulaires classiques (cf \cite{G}), et il est n\'ecessaire 
dans le cas d'un corps de nombres g\'en\'eral de travailler ad\'eliquement.

Soit $L(s,\pi_f)=\sum\lambda_\pi(\n) N(\n)^{-s}$ (avec la convention qu'une telle somme ne porte que sur les id\'eaux \emph{non nuls} de $\entier$) la fonction $L$ finie de $\pi$, convergente pour $\Re(s)>1$, $\rm{q}_\pi$ le conducteur de 
$\pi$ (c'est un id\'eal de $\entier$). 

Soit $$\Lambda(s,\pi):=N(\cond)^{s/2}L(s,\pi)=N(\cond)^{s/2}L(s,\pi_\infty)L(s,\pi_f)$$ la fonction $L$ compl\'et\'ee 
(i.e. tenant compte des places archim\'ediennes), qui se prolonge analytiquement au plan complexe, 
et satisfait \`a l'\'equation fonctionnelle:
\begin{eqnarray}\label{equation}\Lambda(s,\pi)=\varepsilon_\pi\Lambda(1-s,\pi)
\end{eqnarray}
pour $\varepsilon_\pi\in\{-1,1\}$ (car $\pi\cong \check{\pi}$).

Les valeurs $L(1/2,\pi_f)$ sont li\'ees \`a des probl\`emes arithm\'etico-g\'eom\'etriques 
(cf. la conjecture de Birch et Swinnerton-Dyer), et on s'int\'eresse ici \`a leur non 
annulation. Plus pr\'ecis\'ement, si $\g{k}\in \N_{\geq 1}^d$ est fix\'e, on consid\`ere 
$\Pi_\q^{\g{k}}$ l'ensemble (fini) des formes modulaires de Hilbert de poids $\g{k}$ et de conducteur  $\q$, dont on note 
le cardinal $|\Pi_\q^{\g{k}}|$, et on montre inconditionnellement les

\begin{theoreme}\label{1} Pour $\g{k}\geq \g{2}$ pair, et $\q$ parcourant les id\'eaux premiers de $\entier$,
$$\liminf_{N(\q)\to\infty} \frac{ | \{\pi\in\Pi_\q^{\g{k}};L(1/2,\pi)\neq 0\}| }{|\Pi_\q^{\g{k}}|}
 \geq\frac{1}{4}.$$
\end{theoreme}

\begin{theoreme}\label{2} Pour $\g{k}\geq \g{2}$ pair, et $\q$ parcourant les id\'eaux premiers de $\entier$,
$$\liminf_{N(\q)\to\infty} \frac{ | \{\pi\in\Pi_\q^{\g{k}}; \varepsilon_\pi=-1{\textrm{ et }}L'(1/2,\pi_f)\neq 0\}| }{|\Pi_\q^{\g{k}}|}
 \geq\frac{7}{16}.$$
\end{theoreme}
Selon la terminologie de Kowalski et Michel \cite{KM}, on dit qu'on a une densit\'e \emph{naturelle} positive de formes dont la fonction $L$ 
(resp. la d\'eriv\'ee de la fonction $L$) ne s'annule pas en $1/2$.\\

\rem\label{r}: si $\varepsilon_\pi=-1$, alors $\Lambda(1/2,
\pi)=0$ d'apr\`es l'\'equation fonctionnelle. Comme on a asymptotiquement une m\^eme proportion entre les formes  
de signe $1$ et $-1$, le r\'esultat prouv\'e se r\'e\'ecrit:
$$\liminf _{N(\q)\to\infty}\frac{| \{\pi\in\Pi_\q^{\g{k}}|L(1/2,\pi)\neq 0\}|}{|
\{\pi\in\Pi_\q^{\g{k}}|\varepsilon_\pi=+1\}|}\geq \frac{1}{2}.$$
Le travail de Iwaniec, Luo et Sarnak \cite{ILS}, dans le cadre des formes modulaires classiques, 
permet d'atteindre une proportion de $9/16$, sous l'hypoth\`ese de Riemann (GRH).
De m\^eme, le r\'esultat pour la d\'eriv\'ee s'\'ecrit:
$$\liminf _{N(\q)\to\infty}\frac{| \{\pi\in\Pi_\q^{\g{k}}|\varepsilon_\pi=-1\textrm{ et }L'(1/2,\pi_f)\neq 0\}|}{|
\{\pi\in\Pi_\q^{\g{k}}|\varepsilon_\pi=-1\}|}\geq \frac{7}{8}$$
et \cite{ILS} ont atteint, sous (GRH), $15/16$. Kowalski, Michel et VanderKam \cite{KMV} ont montr\'e (sur $\Q$) que cette meilleure pr\'ecision pour les d\'eriv\'ees 
n'est pas un hasard. On conjecture en fait que les proportions ci-dessus sont $1$, mais cela n'est pas atteignable 
avec les techniques d'analyse harmonique utilis\'ees ici, dans \cite{IS}, \cite{KMV} et \cite{Va1}.
\\

Dans ce travail, on prouve d'abord qu'il y a une densit\'e \emph{harmonique} positive de telles formes:
\begin{theoreme}\label{nonannulation}
Pour $\g{k}\geq \g{2}$ pair, et $\q$ parcourant les id\'eaux premiers de $\entier$,
$$\liminf_{N(\q)\to\infty}\sum_{\pi\in\Pi_\q^{\g{k}}}^h \mathbbmss{1}_{\Lambda(1/2,\pi)\neq 0}
\geq\frac{1}{4}.$$
\end{theoreme}
 avec la notation: $ \mathbbmss{1}_{\Lambda(1/2,\pi)\neq 0}$ vaut 1 si $\Lambda(
1/2,\pi)\neq 0$, 0 sinon. 

Le symbole $\sum^h$ indique que l'on pond\`ere la somme par des coefficients, introduits 
plus tard, provenant de l'extension de la formule de Petersson \`a $\Pi_\q^{\g{k}}$ (cf. 
section \ref{fonctionsL}, d\'efinition \ref{poids}).

La m\'ethode suivie ici est celle des moments amollis, initi\'ee par Selberg, et l'amollisseur choisi a \'et\'e introduit par Iwaniec et Sarnak \cite{IS}, g\'en\'eralis\'e par \cite{KMV}: par l'in\'egalit\'e de Cauchy-Schwarz, on 
peut en effet \'ecrire (tous les nombres sont r\'eels):
$$ \sum_{\pi\in\Pi_\q^{\g{k}}}^h \mathbbmss{1}_{\Lambda(1/2,\pi)\neq 0} \geq \frac{
\big(\sum_{\pi\in\Pi_\q^{\g{k}}}^h \Lambda(1/2,\pi)\big)^2}
{\sum_{\pi\in\Pi_\q^{\g{k}}}^h \Lambda(1/2,\pi)^2}. $$
Malheureusement, l'expression de droite tend vers 0 quand $N(\q)$ tend vers l'infini  
(elle est d'ordre $\log(N(\q))^{-1}$), 
ce qui a sugg\'er\'e d'\'ecrire:
$$ \sum_{\pi\in\Pi_\q^{\g{k}}}^h \mathbbmss{1}_{\Lambda(1/2,\pi)\neq 0} \geq \frac{
\big(\sum_{\pi\in\Pi_\q^{\g{k}}}^h \Lambda(1/2,\pi)\M(\pi)\big)^2}
{\sum_{\pi\in\Pi_\q^{\g{k}}}^h \Lambda(1/2,\pi)^2\M(\pi)^2}. $$
Si la suite de nombres $\Big\{\M(\pi)\Big\}_{\pi\in\Pi_\q^{\g{k}}}$ est bien choisie, on peut 
esp\'erer stabiliser le quotient, et obtenir une densit\'e positive: on nomme alors cette 
suite un \og amollisseur\fg. \cite{IS}, puis \cite{KMV} ont trouv\'e, sur $\Q$, une famille 
d'amollisseurs optimaux pour ce probl\`eme, parmi ceux de la forme:
$$\M(\pi)=\sum_{N(\m)\leq M}\lambda_\pi(\m)P_\m$$
avec $M=N(\q)^{\Delta/2}$, pour $\Delta$ dans $]0,1[$. Ici, ce m\^eme type d'amollisseurs 
est efficace, et en suivant \cite{KMV}, on prouve:

\begin{proposition}
Soit $P$ un polyn\^ome tel que $P(0)=P'(0)=0$, $\Delta\in ]0,1[$ tel que $M=N(\q)^{\Delta/2}\notin \N$, et 
$$\M(\pi)=\sum_{N(\m)\leq M}\frac{\mu(\m)P\left(\frac{\log(M/N(\m))}{\log(M)}\right)}{\psi(\m)N(\m)^{1/2}}\lambda_\pi(\m)$$
Quand $N(\q)\to\infty$, parmi les id\'eaux premiers de $\entier$, $\g{k}\geq \g{2}$ pair, on a:
\begin{multline}\label{moment1}
M_1(\q):=\sum_{\pi\in\Pi_\q^{\g{k}}}^h \Lambda(1/2,\pi) \M(\pi)\\ =\frac{\zeta_F(2) \Gamma\left(\frac{\g{k}}{2}\right)}
{(2\pi)^{\frac{\g{k}}{2}} \res_{s=1}(\zeta_F)}\times \frac{2N(\q)^{1/4}}{\Delta\log(N(\q))}
\Bigg(P'(1)+\mathcal{O}\bigg(\frac{1}{\log(N(\q))}\bigg)\Bigg)
\end{multline}
\begin{multline}\label{moment2}
M_2(\q):=\sum_{\pi\in\Pi_\q^{\g{k}}}^h \Lambda(1/2,\pi)^2\M(\pi)^2\\ = \frac{\zeta_F(2)^2\Gamma\left(\frac{\g{k}}{2}\right)^2}
{(2\pi)^{\g{k}}\res_{s=1}(\zeta_F)^2}\times\frac{8N(\q)^{1/2}}{\log(N(\q))^2}
\Bigg( \frac{\|P''\|_{L^2(0,1)}^2}
{\Delta^3}+\frac{P'(1)^2}{\Delta^2}+\mathcal{O}\bigg(\frac{1}
{\log(N(\q))}\bigg)\Bigg).
\end{multline}
\end{proposition}

Ce r\'esultat entra\^ine le th\'eor\`eme \ref{nonannulation} de fa\c con \'evidente (voir section 
\ref{Conclusion}). Remarquons les constantes (impliquant la g\'eom\'etrie de $F$) intervenant dans 
l'asymptotique des deux moments amollis (\`a comparer avec \cite{KMV}, propositions 4.1 et 5.1). Il est assez \'etonnant que 
la proportion des formes, elle, n'en soit pas affect\'ee, et qu'ainsi on puisse atteindre les m\^emes bornes que sur $\Q$ -- les
 meilleures connues inconditionnellement. 

Les deux expressions $M_1(\q)$ et $M_2(\q)$ sont les deux premiers moments amollis, et 
le terme principal des membres de droite proviennent de la \og diagonale\fg de la formule 
de Petersson. Un effort important doit \^etre fait pour montrer que le terme des sommes 
de Kloosterman a une contribution n\'egligeable: c'est ici que le choix de l'amollisseur (i.e. celui de \cite{KMV}) s'av\`ere crucial, puisqu'il permet d'\'eviter une contribution \og hors-diagonale\fg, qui \'etait pr\'esente lors d'un  travail ant\'erieur de Kowalski et Michel \cite{KM1}. Outre les difficult\'es techniques d\'ej\`a pr\'esentes sur $\Q$, 
il faut g\'erer les unit\'es de $F$, qui ont tendance \`a faire diverger les sommes. De plus, nous ne supposons pas que 
$\entier$ est principal, ce qui n\'ecessite l'intervention de la th\'eorie ad\'elique. Enfin, m\^eme pour $\g{k}=\g{2}$, 
il peut exister des formes non ramifi\'ees, et par cons\'equent un terme suppl\'ementaire \`a g\'erer: c'est un ph\'enom\`ene 
absent dans \cite{KMV}, mais une adaptation de la diagonalisation des formes anciennes selon \cite{ILS} nous permet de le 
contr\^oler.\\

Un corollaire de l'\'etude du second moment est le 
r\'esultat de grand crible suivant: 
\begin{theoreme} Soit $\{x_\n\}$ une suite de nombres complexes, $\q$ 
un id\'eal \emph{quelconque} de $\entier$. On pose:
$$\|x_X\|_2:=\Big(\sum_{N(\n)\leq X}|x_\n|^2\Big)^{1/2}$$
On a alors l'estimation:
\begin{eqnarray}\label{crible}
\sum_{\pi\in\Pi_\q^{\g{k}}}^h \Big|\sum_{N(\n)\leq X }  \lambda_\pi(\n)x_\n\Big|^2
 \ll _{F} \bigg(1+\frac{X}{N(\q)}\bigg)\|x_X\|_2^2.
\end{eqnarray}
\end{theoreme}

Cette estimation 
g\'en\'eralise l'in\'egalit\'e de grand crible classique pour les formes modulaires sur 
$\Q$, ainsi que celle de \cite{luo}, o\`u il est suppos\'e que $F$ a un groupe de 
classes \'etroit trivial. La preuve est donn\'ee en section \ref{Conclusion}.
\\

La d\'emarche donn\'ee ci-dessus vaut aussi pour l'\'etude de la non-annulation en moyenne harmonique 
de la d\'eriv\'ee , et permet de montrer le r\'esultat escompt\'e:

\begin{theoreme}\label{nonannulationderivee}
Soit $\g{k}\geq \g{2}$ pair. Quand $\q$ parcourt l'ensemble des id\'eaux maximaux de $\entier$, on a :
\begin{eqnarray}
\liminf_{N(\q)\to\infty}\sum_{\pi\in\Pi_{\q}^{\g{k}}}^h \mathbbmss{1}_{\varepsilon_\pi=-1,L'(1/2,\pi_f)\neq 0}\geq \frac{7}{16}.
\end{eqnarray}
\end{theoreme}

Enfin, pour achever la preuve des th\'eor\`emes \ref{1} et \ref{2}, on doit faire une \'etude de la fonction $L$ du 
carr\'e sym\'etrique, et adapter l'\'etude des moments: c'est ce qu'on appelle le passage de la moyenne harmonique \`a 
la moyenne naturelle.\\

\noindent{\bf Remerciements:}
Ce travail correspond pour l'essentiel \`a ma th\`ese, dirig\'ee par Philippe Michel, qui fut un excellent directeur, patient, passionn\'e et motivant. C'est lui qui m'a fait d\'ecouvrir la th\'eorie analytique des formes automorphes, et guid\'e dans le labyrinthe de D\'edale. Si ce qui suit a de l'int\'er\^et, c'est \` lui qu'il le doit. Je remercie A.Venkatesh et E.Kowalski pour avoir rapport\'e ma th\`ese. En particulier, c'est Emmanuel Kowalski qui m'a encourag\'e \`a r\'ediger enti\`erement la preuve dans ce contexte du th\'eor\`eme de densit\'e z\'ero, et qui a corrig\'e de nombreuses fautes.\\

\noindent{\bf{Organisation de ce travail}}\\
La section \ref{Notations} fixe les 
notations utilis\'ees de fa\c con r\'ecurrente dans le texte. La section \ref{Hilbert} d\'efinit les espaces des formes modulaires de Hilbert, la \ref{Automorphe} rappelle les fondements de la th\'eorie automorphe -- et la d\'efinition automorphe des formes de Hilbert. Le c\oe ur du travail commence en \ref{Petersson}, o\`u l'on prouve la formule de Petersson n\'ecessaire, que l'on utilise lors de toutes les sections ult\'erieures pour \'etudier les divers moments harmoniques, l'in\'egali\'e de grand crible en \ref{Conclusion}, jusqu'\`a la section \ref{Naturel} o\`u l'on traite les moments naturels. L'appendice final contient r\'esultats techniques concernant les formes anciennes surtout, et le carr\'e sym\'etrique: il s'agit de montrer que les termes issus de la formule de Petersson, param\'etr\'es par les formes anciennes, n'ont pas de contribution asymptotique.

\section{Notations et rappels}\label{Notations}

Dans toute la suite, on notera $F$ une extension totalement r\'eelle de $\Q$, de degr\'e 
$d$, d'anneau d'entiers $\entier$. On d\'esignera par $\q$  un id\'eal maximal, sauf dans la section 
\ref{Petersson}, plus g\'en\'erale. On notera usuellement $\adele$ l'anneau des ad\`eles 
de $F$. 
Les places de $F$ seront not\'ees $v$, $F_v$ d\'esignant le compl\'et\'e de $F$ en $v$, et 
$\mathcal{O}_{F_v}$ ou $\mathcal{O}_v$ l'anneau local des entiers quand $v\not| \infty$. L'\'ecriture $F_\infty$ est ici 
pour $\R^d$, et $F_\infty^{\times}=(\R^{\times})^d$ (respectivement: 
$F_\infty^{\times>0}=F_\infty^{\times+}=(\R^{\times}_+)^d$). $N_{F/\Q}$ d\'esigne la norme, ${\rm{Tr}}_{F/\Q}$ 
la trace, et $N=|N_{F/\Q}|$ (prolong\'ee aux id\'eaux fractionnaires).

La notation $\varpi_v$ (ou \'eventuellement $\varpi_\p$ si $v$ est la valuation associ\'ee \`a l'id\'eal maximal $\p$)
 d\'esigne une uniformisante de l'anneau $\mathcal{O}_v$. Si $\A=\prod \p^{n_\p(\A)}$ est un id\'eal fractionnaire, on 
notera $\textrm{id}(\A)$ l'id\`ele fini $(\varpi_\p^{n_\p(\A)})_\p$ (s'il est besoin valant 1 aux places infinies). C'est cet 
id\`ele que l'on dit correspondre \`a $\A$.

On note aussi $|X|$ le cardinal de l'ensemble fini $X$.

\subsection{G\'eom\'etrie de $F$}\label{geometrie}
\noindent$\bullet$ On note $\mathfrak{D}_F$ la diff\'erente de $F$: cet 
id\'eal de $\entier$ a pour norme $|\disc|$, o\`u $\disc$ est le discriminant de $F$. 
La norme d'un id\'eal $\A$ est \'egale \`a: $N(\A):=[ \entier : \A]$, d\'efinition 
que l'on peut prolonger par multiplicativit\'e au groupe des  id\'eaux fractionnaires 
de $F$, not\'e $\mathscr{I}(F)$. \\
Les $d$ plongements de $F$ dans $\R$ seront not\'es $\xi\mapsto \xi^{(j)}$ pour $j=1,\dots
,d$. Si $\xi$ v\'erifie: $\xi^{(j)}>0$ pour tout $j$, on notera $\xi\gg 0$ (on dit alors 
que $\xi$ est totalement positif), et
pour tout sous-ensemble $X$ de $F$, on pose: 
$$X^+=X^{\gg 0}:=\Big\{ x\in X ; x\gg 0\Big\}.$$
$\bullet$ L'ensemble $F^{\times\gg 0}$ est le sous-groupe de $\mathscr{I}(F)$ form\'e des id\'eaux 
principaux admettant un g\'en\'erateur totalement positif
. Le groupe des classes 
\'etroit est le quotient:
$$\mathscr{C}\ell^{+}(F):=\mathscr{I}(F)/F^{\times\gg 0}$$
Ce groupe admet la repr\'esentation ad\'elique:
$$\mathscr{C}\ell^+(F)=\adele^\times/F^\times F_\infty^{\times +}\widehat{\mathcal{O}}_F^\times$$
avec $\widehat{\mathcal{O}}_F=\prod_{v<\infty}\mathcal{O}_{F_v}$. C'est donc un groupe fini, de cardinal $h_F^+$, et on \emph{choisit} une fois pour toutes un syst\`eme de 
repr\'esentants $\big\{\A\big\}$ dans $\mathscr{I}_F$: ce choix est in\'el\'egant, mais reste n\'ecessaire 
dans la mesure o\`u il permet de traduire le probl\`eme pos\'e en terme d'analyse 
harmonique sur des espaces sym\'etriques r\'eels (probl\`eme de \og l'ad\'elisation\fg). 
De plus, la formule de Petersson 
pour un corps g\'en\'eral (i.e. dont le groupe des classes \'etroit n'est pas trivial) 
s'exprime avec des sommes de termes d\'ependant de ce choix, bien qu'invariante 
globalement. Techniquement, cela permet aussi de remplacer des sommes sur des id\'eaux par des 
sommes sur des entiers, et d'utiliser le lemme suivant (cf \cite{luo} page 131), cons\'equence du th\'eor\`eme de Dirichlet sur les unit\'es de $F$:
\begin{lemme}\label{sommation} 
Soit $F$ un corps de nombres totalement r\'eel. Il existe des constantes $C_1$, $C_2$ ne d\'ependant que de $F$ telles 
que 
\begin{multline}\forall\xi\in F,\exists\varepsilon\in\unite, \forall j\in\big\{1,\dots,d\big\}: \\ C_1|N(\xi)|^{1/d}\leq |(\varepsilon\xi)^{(j)}|\leq C_2 |N(\xi)|^{1/d}.
\end{multline}
\end{lemme} 
$\bullet$ Etant donn\'e $\A$ et $\B$ deux id\'eaux fractionnaires, on notera:
$$\A\sim \B\Leftrightarrow \textrm{$\A$ et $\B$ ont m\^eme image dans $\mathscr{I}_F$}
\Leftrightarrow \exists \xi\in F^{\times >0}; \A\B^{-1}=\xi\entier$$
et lorsque tel est le cas on notera $[\A\B^{-1}]$ le choix d'un $\xi$ satisfaisant \`a la 
relation pr\'ec\'edente.\\
$\bullet$ Nous noterons $\zeta_F$ la fonction de z\^eta de Dedekind du corps $F$. Cette s\'erie de Dirichlet 
permet de construire des fonctions arithm\'etiques, comme par exemple:\\
$\ast\,\mu$ la g\'en\'eralisation de la fonction de M\"obius, d\'efinie par la relation:
\begin{displaymath}
\mu(\n)= \left\{\begin{array}{ll}(-1)^r \textrm{ si $\n$ est produit de $r$ id\'eaux premiers distincts}\\ 0 \textrm{ sinon.}\end{array}\right.
\end{displaymath}
On v\'erifie ais\'ement l'identit\'e:
$$\zeta_F^{-1}(s)=\sum_{\n\subset \entier}\mu(\n)N(\n)^{-s},\forall \Re(s)>1$$ 
\\
$\ast\,\tau$, d\'efinie par $\tau(\n):=| \big\{\mathfrak{d}\subset \entier|\n\mathfrak{d}^{-1}
\subset \entier \big\}|$, donne le nombre de diviseurs, et v\'erifie toujours l'estimation, 
pour tout $\varepsilon >0$:
$$\tau(\n)\ll_{\varepsilon } N(\n)^\varepsilon $$
Elle peut se voir comme les coefficients de la s\'erie $\zeta_F^2$.\\
On utilisera aussi la fonction arithm\'etique $\psi$ 
$$\psi(\n):=\prod_{\p|\n}(1+N(\p)^{-1}).$$
Dans le produit, $\p$ d\'esigne un id\'eal maximal. 
Son introduction dans l'amollisseur $\M(\pi)$ permet de calculer explicitement les termes principaux 
des premiers et deuxi\`eme moments.
On se servira implicitement de l'estimation
$$| \big\{\n\subset \entier|N(\n)=n\big\}|\ll_{\varepsilon }n^\varepsilon, \,\forall \varepsilon >0 $$
pour la convergence de certaines sommes.

Enfin, nous noterons $\zeta_F^{(\q)}$ la fonction $\zeta_F\times (1-N(\q)^{-s})$, c'est-\`a-dire la 
fonction $\zeta$ \`a laquelle on a \^ot\'e le facteur en $\q$.

\subsection{Sommes de Kloosterman}\label{Kloosterman}
Par commodit\'e, on rappelle ici la d\'efinition donn\'ee par Venkatesh dans \cite{V} (d\'efinition 2) des sommes de Kloosterman.\\  

Soient $\A$, $\B$ deux id\'eaux fractionnaires de $F$, tels que $\B\subset\A$. On note $(\A/\B)^\times$ l'ensemble des 
$x\in\A/\B$ engendrant $\A/\B$ en tant que $\entier$-module. Pour un tel $x$, on note $\bar{x}$ l'unique \'el\'ement 
$y\in(\A^{-1}/\B\A^{-2})^\times$ tel que $xy\equiv 1\textrm{ mod } \B\A^{-1}$. Cela \'etant pos\'e, soient $\A_1,\A_2$ deux 
id\'eaux fractionnaires, et $\cc$ un id\'eal tel que $\cc^2\sim \A_1\A_2$. Soient aussi $\alpha_1\in\A_1^{-1}\diff^{-1}$, 
$\alpha_2\in\A_2^{-1}\diff^{-1}$ et $c\in\cc^{-1}\q$, $\q$ \'etant un id\'eal fix\'e de $\entier$. On pose 
\begin{multline}
KS(\alpha_1,\A_1;\alpha_2,\A_2;c,\cc)=\sum_{x\in(\A_1\cc^{-1}/\A_1c)^\times}\exp\left(2i\pi\textrm{Tr}_{F/\Q}
\left( \frac{\alpha_1x+\alpha_2\bar{x}}{c}\right)\right)
\end{multline}
Nous les renormaliserons en \ref{fonctionsL}, et donnerons alors la borne de Weil qu'elles satisfont.

\subsection{Rappels sur les groupes}\label{Groupes}
De fa\c con g\'en\'erale,  si $G$ est un groupe alg\'ebrique sur $\Z$, et $R$ un anneau quelconque, $G(R)$ d\'esigne le groupe 
des points de $G$ \`a valeurs dans $R$. Si $R$ est topologique, $G(R)$ peut \^etre muni d'une topologie \og forte\fg, issue de 
celle de $R$. En particulier, $G(F_\infty)$ d\'esigne bien s\^ur $G(\R)^d$ et 
$G^+(F_\infty)$ sa composante neutre. Nous travaillerons avec $G=\gl$, et aurons besoin de certains de ses sous-groupes: pour 
$R$ un anneau (\'eventuellement topologique), on notera:
\begin{eqnarray}
Z(R)&:=&\Big\{ \left(\begin{array}{cc}z & 0\\0 & z\end{array}\right); z\in R^\times \Big\}\nonumber\\
N(R)&:=&\Big\{ \left(\begin{array}{cc}1 & x\\0 & 1\end{array}\right); x\in R\Big\}\nonumber\\
A(R)&:=&\Big\{ \left(\begin{array}{cc}a & 0\\0 & 1\end{array}\right); a\in R^\times \Big\}\nonumber\\
P(R)&:=&\Big\{ \left(\begin{array}{cc}a & b\\0 & d\end{array}\right); a,d\in R^\times, b\in R \Big\}\nonumber
\end{eqnarray}
Le sous-groupes ${\rm{SL}}_2(R),{\rm{(S)O}}_2(R)$ interviendront avec $R=F_\infty$. On pourra noter $Z_\infty$ le centre 
de $\gl(F_\infty)$, et $Z_\infty^+\cong F_\infty^{\times +}$ sa composante neutre. 
Dans le cas o\`u $F$ est un corps de 
nombres, certains groupes compacts sont utiles: pour $v|\infty$, on a d\'ej\`a vu le sous-groupe compact maximal 
$K_v={\rm{SO}}_2(F_v)$ de $\GL(F_v)$, et en faisant le produit sur toutes les places infinies on note 
$K_\infty={\rm{SO}}_2(F_\infty)$ celui de $\GL(F_\infty)$. En $v$ finie, le compact maximal de $\gl(F_v)$ est $K_v=\gl(\mathcal{O}_{F_v})$. D'autres 
sous-groupes compacts dans le cas non-archim\'edien sont importants: si $\q_v\subset\mathcal{O}_{F_v}$ est un id\'eal, 
on notera 
$$K_0(\q_v)=\Big\{ \left(\begin{array}{cc}a & b\\c & d\end{array}\right)\in\gl(\mathcal{O}_{F_v}) , c_v\in\q_v\Big\}.$$
Ces sous-groupes locaux engendrent des sous-groupes compacts de $\gl(\adele)$. Soit $\q\subset\entier$ un id\'eal:
\begin{eqnarray}
K_f&:=&\prod_{v<\infty} \gl(\mathcal{O}_{F_v})\nonumber\\
K_0(\q)&:=&\Big\{ g\in K_f;\,\, \forall v\not|\infty, g_v\in K_0(\q\mathcal{O}_{F_v})\Big\}
\nonumber\\
K&:=&K_\infty K_f \nonumber
\end{eqnarray}

L'importance de ces groupes r\'eside notamment dans la d\'ecomposition d'Iwasawa:

\begin{proposition}\label{Iwasawa}
Soit $F$ totalement r\'eel, et $v$ une place de $F$. L'application:
\begin{displaymath}
\begin{array}{ccc}
Z(F_v)\times N(F_v)\times A(F_v)\times K_v &\longrightarrow & \gl(F_v)\\
(z,n,a,k) & \longmapsto & znak
\end{array}
\end{displaymath}
$\bullet$ est surjective, et sa restriction \`a $Z^+(F_v)\times N(F_v)\times A(F_v)\times K_v$ induit un hom\'eomorphisme si $v|\infty$.\\
$\bullet$ est surjective si $v\not|\infty$. De plus, si $(p_0,k_0)\in P(F_v)\times K_v$ a pour image $g\in\gl(F_v)$, toutes les 
d\'ecompositions d'Iwasawa de $g$ sont donn\'ees par $\{(p_0k^{-1},kk_0),k\in P\cap K_v\}$.
\end{proposition}
 Pour $v|\infty,v\not|\infty,v=\infty$, on notera pour $g_v\in\gl(F_v)$ une d\'ecomposition d'Iwasawa:
$$g_v=z(g_v)n(g_v)a(g_v)k(g_v)$$
On notera aussi que pour $v|\infty$ ou $v=\infty$, on a aussi une d\'ecomposition d'Iwasawa pour $\GL(F_v)$:
$$Z_\infty^+\times N(F_v)\times A^+(F_v)\times K_v $$
On pourra utiliser l'isomorphisme exceptionnel $\textrm{SO}_2(\R)\cong \R/2\pi\Z$ en notant $k\in\textrm{SO}_2(\R)$:
$$k=\left(\begin{array}{cc}\cos\theta & \sin\theta\\-\sin\theta 
& \cos\theta\end{array}\right).$$

\subsection{Mesures de Haar sur les corps locaux}

Soit $F$ un corps de nombres, que nous supposerons totalement r\'eel pour raccourcir notre propos. 

Dans ce cas, pour toute 
place archim\'edienne $v$, $F_v=\R$, et on munit $\R$ de la mesure de Lebesgue $dx$ (qui est la mesure Haar normalis\'ee dans 
ce cas). La mesure de Haar de $\R^\times$ est $d^\times x=\frac{dx}{|x|}$, et c'est aussi celle de $\R^{\times}_+$.

Pour $v$ finie, on utilise les mesures normalis\'ees de Tate (voir la th\`ese de Tate dans \cite{CF} pour plus de d\'etails)
 : si $\p$ est un id\'eal maximal correspondant \`a $v$, on 
note $n_\p(\diff)$ la $\p$-valuation de la diff\'erente globale (ou locale), et on choisit pour mesure de Haar normalis\'ee celle qui v\'erifie 
${\rm{vol}}(\mathcal{O}_v)=N(\p)^{-n_\p(\diff)/2}$, la notant $dx$. La mesure $\frac{dx}{|x|_v}$ est bien une mesure de Haar 
multiplicative, mais on note $d^\times x$ la mesure normalis\'ee de telle sorte que 
${\rm{vol}}^{\scriptscriptstyle{\times}}(\mathcal{O}_v^\times)=1$.

On peut alors mettre sur les groupes ad\'eliques $\adele$ et $\adele ^\times$ les mesures limite inductive (bien d\'efinies 
car on a pris la peine d'assurer ${\rm{vol}}(\mathcal{O}_v)=1$ p.p.($v$) et 
${\rm{vol}}^{\scriptscriptstyle{\times}}(\mathcal{O}_v^\times)=1$ p.p.($v$)). Ces mesures sont caract\'eris\'ees par leur 
valeur sur une base de la topologie: on pose ${\rm{mes}}(\prod_{v\in S} U_v\times \prod_{v\notin S}\mathcal{O}_v)=\prod _{v\in S}
{\rm{mes}}_v(U_v)$ (avec $S$ ensemble de places fini, $U_v$ ouvert dans $F_v$) dans le cas de la mesure additive par exemple, et on fait 
pareillement pour la mesure multiplicative. 
On aura juste besoin en fait de savoir 
que ${\rm{vol}}(\adele/F)=1$ pour la mesure additive, soit encore ${\rm{vol}}(F_\infty/\entier)=|\disc|^{1/2}$.\\

La d\'ecomposition d'Iwasawa permet \'egalement de normaliser la mesure de Haar de $\gl$. En fait, remarquant que 
$Z(F_v)\cong F_v^\times$, $N(F_v)\cong F_v$, $A(F_v)\cong F_v^\times$, on peut montrer que 
$$d\mu_{\gl(F_v)}(g):=d^\times z dx \frac{d^\times a}{|a|_v}d\mu_{K_v}$$ en notant 
$$g=\left(\begin{array}{cc}z & 0\\0 & z\end{array}\right)
\left(\begin{array}{cc}1 & x\\0 & 1\end{array}\right)\left(\begin{array}{cc}y & 0\\0 & 1\end{array}\right)k_v$$
avec $k_v\in K_v$. On normalise $d\mu_{K_v}$ pour en faire une mesure de probabilit\'e. Voir la preuve dans \cite{Bu}, 
proposition 2.1.5.

La mesure de Haar de $\gl(\R)$ induit celle de $\GL(\R)$ (on int\`egre sur $Z(\R)^+, N(\R), A(\R)^+$); celle du quotient 
$Z(\R)^+\backslash\GL(\R)$ est param\'etr\'ee par $dx \frac{d^\times a}{|a|_v}d\mu_{K_v}$.

On d\'eduit de ces mesures locales une mesure de Haar normalis\'ee sur $\gl(\adele)$  
d\'etermin\'ee par les valeurs ${\rm{mes}}\left(\prod_{v\in S}U_v\times\prod_{v\notin S}\gl(\mathcal{O}_v)\right):=
\prod_{v\in S}{\rm{mes}}_v(U_v)$. C'est dans un but global que l'on normalise les mesures locales.

\subsection{Notations des fonctions}\label{fonctions}
$F$ \'etant de degr\'e $d$, on aura \`a travailler avec des fonctions d\'efinies sur $\R^d$. On notera 
en caract\`eres gras les vecteurs. Si $n\in\N$, $\g{n}=(n,\dots,n)\in\N^d$.

Ainsi, si $\varphi=\otimes_j\varphi_j:\C^d\longrightarrow \C$, et $\g{z}=(z_1,\dots,z_d)$, 
on note:
$$\varphi(\g{z})=\prod_{j=1}^d\varphi_j(z_j).$$

De m\^eme, si $f:\C\longrightarrow \C$ est complexe, $f(\g{z})$ d\'esigne le nombre:
$$\prod_{j=1}^df(z_j)$$
en particulier: $2^{\g{z}}=2^{\sum_j z_j}$, $\Gamma(\g{k-1})=\prod_j\Gamma(k_j-1)$
, $2^{\g{1}}=2^d\neq 2^1$. Cela permet d'\'etendre la norme \`a $\C^d$ par: $N(\g{z}):=|\g{z}^{\g{1}}|.$\\

Si, pour $\nu\in\R$, $f_\nu$ d\'esigne une fonction \`a valeurs dans $\C$, alors pour $\g{\nu}=
(\nu_1,\dots,\nu_d)\in\R^d$ et $\g{z}\in\C^d$:
$$\g{f}_{\g{\nu}}(\g{z})=\prod_{j=1}^df_{\nu_j}(z_j).$$
Par exemple, pour $\g{k}\in\N^d$, $\g{x}\in\R^d$, on pose
$$\g{J}_{\g{k-1}}(\g{x}):=\prod_{j=1}^dJ_{k_j-1}(x_j).$$
Par respect de ces conventions, nous noterons un \'el\'ement du corps $F$ $\g{\xi}$ (et non $\xi$) quand il sera 
vu comme le vecteur $(\xi^{(j)})_{1\leq j \leq d}$.\\

\section{Formes modulaires de Hilbert}\label{Hilbert}

Supposons que $F$ est un corps totalement r\'eel, de degr\'e $d$ sur $\Q$. Soit $\HH^d$ le produit de $d$ copies du 
demi-plan de Poincar\'e. On d\'efinit une action de $\GL(F)=\{\gamma\in\gl(F);{\rm{det}}\gamma\gg 0\}$ sur $\HH^d$ 
de la mani\`ere suivante:

\begin{eqnarray}
\forall \g{z}=(z_j)_{1\leq j\leq d}\in\HH^d, \forall \gamma\in\GL(F), \gamma.\g{z}=(\gamma^{(j)}.z_j)_{1\leq j\leq d}
\end{eqnarray}
avec, si $\gamma=\left(\begin{array}{cc}a & b\\c & d\end{array}\right)$ et $1\leq  j \leq d$:
\begin{eqnarray}
\gamma^{(j)}.z_j=\frac{a^{(j)}z_j+b^{(j)}}{c^{(j)}z_j+d^{(j)}}
\end{eqnarray}

Par cons\'equent, pour tout sous-groupe $\Gamma\subset \GL(F)$ on a l'action induite $\Gamma\circlearrowleft\HH^d$. En particulier, les sous-groupes de congruence suivants, o\`u $ \A,\B$ d\'esignent des id\'eaux fractionnaires, sont centraux dans la th\'eorie:
\begin{eqnarray}
\Gamma_0(\A,\B)=\Big\{\left(\begin{array}{cc}a & b\\c & d\end{array}\right)\in\GL(F);\,\,
a,b\in\entier, c\in\A\B^{-1},b\in\B,ad-bc\in\mathcal{O}_F^{\times +}\Big\}\nonumber
\end{eqnarray}

\noindent{\bf{Formes modulaires de Hilbert classiques}}\\
Soit $\g{k}\in\N^d$ un entier pair, $ \Gamma=\Gamma_0(\A,\B)$. Soit 
$f:\HH^d\longrightarrow \C$ une fonction holomorphe (\`a plusieurs variables). Pour tout $\gamma\in\GL(F)$, on d\'efinit 
une nouvelle fonction $f|_\gamma$:

\begin{eqnarray}
f|_\gamma(\g{z})=\left( (\g{{\rm{det}}\gamma})^{-\frac{1}{2}}(\g{cz}+\g{d})\right)^{-\g{k}}f(\gamma.\g{z})
\end{eqnarray}
Rappelons que d'apr\`es nos conventions, avec $\gamma=\left(\begin{array}{cc}a & b\\c & d\end{array}\right)$, cela s'\'ecrit 
aussi:
$$f|_\gamma(\g{z})=\prod_{j=1}^d\left( ({\rm{det}}\gamma^{(j)})^{-\frac{1}{2}}(c^{(j)}z_j+d^{(j)})\right)^{-k_j}\times
f\left( (\gamma^{(j)}.z_j)_{1\leq j\leq d}\right)$$
\begin{definition}
On appelle forme modulaire classique de poids $\g{k}$, pour le groupe de congruence $\Gamma$, toute fonction holomorphe $f:\HH^d\longrightarrow \C$ 
telle que:
\begin{eqnarray}
f|_\gamma=f,\,\,\forall\gamma\in\Gamma
\end{eqnarray} 
\end{definition}
Le principe de Koechler donne alors l'holomorphie en les pointes, sous la seule hypoth\`ese que $ [F:\Q]>1$. Une forme modulaire de Hilbert classique admet un d\'eveloppement en s\'erie de Fourier:
$$f(\g{z})=\sum_{\substack{\xi\in\B^{\star}\\ \xi\geq 0  }}c(\xi,f)\exp(2i\pi {\rm{Tr}}(\xi \g{z}))$$
et l'on note $\mathscr{S}_{\g{k}}(\Gamma)$  l'espace des formes modulaires cuspidales classiques, dont le coefficient de Fourier constant $ c^\gamma(0,f)$ est nul pour chaque $\gamma\in\gl(F)$. La fonction $L$ que l'on peut d\'efinir \`a partir des formes classiques est une somme de Dirichlet index\'ee par les id\'eaux fractionnaires principaux de $F$: pour l'\'etude de probl\`emes g\'eom\'etriques -- comme les fonctions z\^eta des vari\'et\'es, il est n\'ecessaire de les ad\'eliser, afin d'avoir des sommes sur les id\'eaux quelconques.\\

Pour cela, commen\c cons par le th\'eor\`eme d'approximation forte:\\

\noindent{\sc{Th\'eor\`eme d'approximation forte pour $\rm SL_2$:}}\\
\emph{Pour tout corps totalement r\'eel $F$,  ${\rm{SL}}_2(F){\rm{SL}}_2(F_\infty)$ est dense dans ${\rm{SL}}_2(\adele)$.
}\\

Une version g\'en\'erale de ce th\'eor\`eme est prouv\'ee dans le livre de Platonov-Rapinchuk (th.7.12). On l'utilise seulement pour 
$\textrm{SL}_2$, pour lequel \cite{Ga} donne une preuve \'el\'ementaire (appendice A.3). 
\begin{proposition}[Approximation forte pour $\gl$]
Soit  $\{\A\}$ le syst\`eme  de la section \ref{geometrie}, et ${\rm{id}}(\A)$  un id\`ele fini 
correspondant \`a $\A$ 
$$\gl(\adele)=\coprod _{\overline{\A}\in\mathscr{C}\ell^+(F)}\gl(F)\GL(F_\infty)
\left(\begin{array}{cc}{\rm{id}}(\A) & 0\\0 & 1\end{array}\right)K_0(\q).$$
\end{proposition}
\preuve : On utilise d'abord le fait que (cf \ref{geometrie})
$$\adele^\times=\coprod_{\overline{\A}\in\mathscr{C}\ell^+(F)}F^\times F_\infty^{\times +} \textrm{id}(\A) \widehat{\mathcal{O}}_F^\times.$$
Le th\'eor\`eme d'approximation forte assure que $\textrm{SL}_2(F)\textrm{SL}_2(F_\infty)$ est dense dans $\textrm{SL}_2(\adele)$. 
En \'ecrivant $\textrm{det}d=\xi y_\infty \textrm{id}(\A)u_f$ selon la d\'ecomposition ci-dessus, on en d\'eduit donc qu'il existe 
deux suites $\gamma_n\in\textrm{SL}_2(F),g_{\infty,n}\in\textrm{SL}_2(F_\infty)$ telles que:
$${\gamma_n g_{\infty,n}}\,\, _{\overrightarrow{n\to\infty}}\,\, \left(\begin{array}{cc}\xi^{-1} & 0\\0 & 1\end{array}\right)
g\left(\begin{array}{cc}y_\infty^{-1} & 0\\0 & 1\end{array}\right)
\left(\begin{array}{cc}u_f^{-1} & 0\\0 & 1\end{array}\right)
 \left(\begin{array}{cc}\textrm{id}(\A)^{-1} & 0\\0 & 1\end{array}\right) 
$$
et donc $\coprod_\A \gl(F)\GL(F_\infty)\left(\begin{array}{cc}\textrm{id}(\A) & 0\\0 & 1\end{array}\right)
\left\{ \left(\begin{array}{cc}u & 0\\0 & 1\end{array}\right); u\in\widehat{\mathcal{O}}_F^\times\right\}$ est dense 
dans $\gl(\adele)$. Comme $K_0(\q)$ est ouvert, v\'erifiant:
$$\textrm{det}\left\{ \left(\begin{array}{cc}u & 0\\0 & 1\end{array}\right); u\in\widehat{\mathcal{O}}_F^\times\right\}=\textrm{det}(K_0(\q))$$
cela ach\`eve la preuve. $\blacksquare$\\
Soit:
$$\Gamma_0(\q,\A):=\GL(F)\cap \left(\begin{array}{cc}{\textrm{id}}(\A) & 0\\0 & 1\end{array}\right)
K_0(\q) \left(\begin{array}{cc}{\textrm{id}}(\A)^{-1} & 0\\0 & 1\end{array}\right)$$
ou encore:
$$\Gamma_0(\q,\A)=\Big\{\left(\begin{array}{cc}a & b\\c & d\end{array}\right)\in\gl(F)|\,\,
a,d\in\entier, c\in\q\A^{-1},b\in\A,ad-bc\in\mathcal{O}_F^{\times +}\Big\}$$
Le th\'eor\`eme d'approximation forte induit un hom\'eomorphisme:
\begin{eqnarray}
\gl(F)Z_\infty^+\backslash \gl(\adele)/K_0(\q)\cong \coprod _{\overline{\A}\in\mathscr{C}\ell(F)}
Z_\infty^+\Gamma_0(\q,\A)\backslash \GL(F_\infty).
\end{eqnarray}

\noindent{\sc{Notations}}: 
Soit $\Delta$ le Laplacien de $\GL(\R)$, d\'efini avec les coordonn\'ees d'Iwasawa par:
$$\Delta=-y^2\left( \frac{\partial ^2}{\partial x ^2}+\frac{\partial^2}{\partial y^2}\right)+
y\frac{\partial^2}{\partial x\partial \theta}$$

 $\Delta$ d\'efinit sur l'espace des fonctions lisses $\mathscr{C}^\infty(\GL(F_\infty))$
$d$ op\'erateurs $\{\Delta_j\}_{1\leq j\leq d}$ par:
$$\forall g_\infty, \Delta_j\varphi(g_\infty)=\Delta[h\mapsto \varphi(g_1,\dots,\underbrace{h}_{
{\textrm{place j}}},\dots,g_d,)]{\mid} _{h=g_j}.$$
Soit $\g{\Delta}=(\Delta_1,\dots,\Delta_d)$, 
et pour $\g{\lambda}\in\R_+^d$, $\varphi\in \mathscr{C}^\infty(\GL(F_\infty))$, 
on pose $\g{\Delta}\varphi=\g{\lambda}\varphi$ si, pour tout 
$j$ dans $\{1,\dots,d\}$:
$$\forall g_\infty\in \GL(F_\infty), \Delta_j\varphi(g_\infty)=\lambda_j\varphi(g_\infty).$$
Pour $k_\infty\in K_\infty$, $\e( \g{k}k_\infty)$ d\'esigne $\exp( i\g{k\theta})$ si 
$k_\infty=\left(\begin{array}{cc}\cos\g{\theta} & \sin\g{\theta}\\-\sin\g{\theta} 
& \cos\g{\theta}\end{array}\right)$ (avec les notations vectorielles de la section pr\'ec\'edente).\\

On peut maintenant d\'efinir un premier espace de formes ad\'eliques: c'est un espace auxiliaire, qui sera utile pour prouver 
la formule des traces de Petersson, car il permet une indexation pratique. 
\begin{multline}
\mathcal{S}_\q^{\g{k}}:=\Big\{ \varphi :\gl(\adele)\rightarrow \C \textrm{ born\'ee }   ; \varphi(\gamma 
z_\infty g k_\infty k_f)=\e(\g{k}k_\infty)\varphi(g)\\
 \forall (\gamma, z_\infty,g,k_\infty,k_f)\in 
 \gl(F)\times
Z_\infty^+ \times\gl(\adele)\times {\rm{SO}}_2(\R)^d\times K_0(\q);\\
\g{\Delta}\varphi=\g{\frac{k}{2}}\bigg(\g{1-\frac{k}{2}}\bigg)\varphi \textrm{ et }
\int_{F\backslash \adele}\varphi(n(x)g)dx=0\,\, \forall g\in\gl(\adele)  \Big\}
\end{multline} 
L'espace $\mathcal{S}_\q^{\g{k}}$ est muni du produit scalaire:
$$\left<\varphi,\psi\right>_{\mathcal{S}_\q^{\g{k}}}:=\int_{\gl(F)Z_\infty^+\backslash \gl(\adele)/K_0(\q)}\varphi(g)
\overline{\psi(g)}dg$$

Un corollaire de l'approximation forte pour $\gl$ est ce que l'on appelle \og l'ad\'elisation des formes modulaires\fg :

\begin{corollaire}
On a un isomorphisme
\begin{eqnarray}\label{adelisation}
\begin{array}{lll}\nonumber
\mathcal{S}_\q^{\g{k}}&\longrightarrow & \displaystyle{\bigoplus_{\overline{\A}\in\mathscr{C}\ell^+(F)}\mathscr{S}_{\g{k}}
(\Gamma_0(\q,\A))}\\
\varphi & \longmapsto & \bigg\{\g{z}=\g{x}+i\g{y}\mapsto \g{y}^{-\g{k}/2}\varphi\Big( \left(\begin{array}{cc}\g{y} & \g{x}\\0 & 1\end{array}\right)
\left(\begin{array}{cc}{\rm{id}}(\A) & 0\\0 & 1\end{array}\right)\Big)\bigg\}_{\overline{\A}
\in\mathscr{C}
\ell^+(F)}
\end{array}
\end{eqnarray}
\end{corollaire}
\rem En munissant l'espace $\oplus _{i}\mathcal{H}_i$, somme directe d'espaces de Hilbert $\mathcal{H}_i$, du produit scalaire 
$\left< \sum_i x_i,\sum_i y_i\right>_{\mathcal{H}_i}=\sum_i \left<x_i,y_i\right>$, l'isomorphisme ci-dessus est une isom\'etrie.\\ 

 On a sur $\mathcal{S}_\q^{\g{k}}$ une action $\rho$ non triviale de $\mathscr{C}\ell^+(F)$, quotient de l'action du centre: 
$$\rho(\B)\varphi(g)=
\varphi\left( \left(\begin{array}{cc}\textrm{id}(\B)& 0\\ 0 & \textrm{id}(\B)\end{array}\right) g\right)$$ 
et par cons\'equent on a la d\'ecomposition:
\begin{eqnarray}\label{decomposition}
\mathcal{S}_\q^{\g{k}}=\bigoplus_{\chi\in \widehat{\mathscr{C}\!\ell^+(F)}}\mathcal{S}_\q^{\g{k}}[\chi]
\end{eqnarray}
$ \widehat{\mathscr{C}\!\ell^+(F)}$ d\'esignant le dual du groupe fini commutatif  $\mathscr{C}\!\ell^+(F)$. \\
\begin{definition}
\noindent{\textsc{L'espace des formes modulaires de Hilbert ad\'eliques}}\\
L'\emph{espace des formes modulaires de Hilbert} de poids $\g{k}\in\N^d$, de 
niveau $\q$ est d\'efini par :
\begin{multline}
\mathcal{H}_\q^{\g{k}}=\Big\{ \varphi :\gl(\adele)\rightarrow \C \textrm{ born\'ee }   ; \varphi(\gamma 
z g k_\infty k_f)=\varphi(g)\e(\g{k}k_\infty)\\
 \forall (\gamma, z,g,k_\infty,k_f)\in 
 \gl(F)\times
Z(\adele) \times\gl(\adele)\times {\rm{SO}}_2(\R)^d\times K_0(\q);\\
\g{\Delta}\varphi=\g{\frac{k}{2}}\bigg(\g{1-\frac{k}{2}}\bigg)\varphi \textrm{ et }
\int_{F\backslash \adele}\varphi(n(x)g)dx=0\,\, \forall g\in\gl(\adele)  \Big\}
\end{multline} 
\end{definition}
Cet espace est nul si $\g{k}\notin 2\mathbb{Z}^d_{\geq 1}$. 
On a une injection \'evidente $\mathcal{H}_\q^{\g{k}}\hookrightarrow \mathcal{S}
_\q^{\g{k}}$ (et m\^eme $\mathcal{H}_\q^{\g{k}}=\mathcal{S}_\q^{\g{k}}[1]$).

L'espace $\mathcal{H}_\q^{\g{k}}$ admet une structure hermitienne, avec le produit scalaire:
$$\left<\varphi,\psi\right>_{\mathcal{H}_\q^{\g{k}}}= \int_{\gl(F)Z(\adele)\backslash \gl(\adele)/K_0(\q)} \varphi(g)\overline{\psi(g)}dg$$
et on peut remarquer d\`es \`a pr\'esent, ce sera utile dans la suite, que pour toute forme $\varphi$:
$$||\varphi||^2_{\mathcal{H}_\q^{\g{k}}}=[K_f:K_0(\q)]\int_{\gl(F)Z(\adele)\backslash \gl(\adele)} 
|\varphi(g)|^2dg$$
avec $[K_f:K_0(\q)]=N(\q)+1$ si $\q$ premier.

\section{Representations automorphes de $\gl$}\label{Automorphe}
\subsection{Formes automorphes pour $\gl$}
Soit $\adele$ l'anneau des ad\`eles de $F$, $\omega$ un caract\`ere multliplicatif de $\adele^\times/F^\times$. On se place dans $\gl(\adele)$, et on utilisera les 
notations de la section \ref{Groupes}. Outre les excellentes r\'ef\'erences usuelles \cite{G}, \cite{Bu} sur les formes automorphes, les notes de Godement \cite{Go}, moins populaires, m'ont \'et\'e d'une grande utilit\'e, pour tout ce qui regarde les fonctions $L$, ainsi que le livre r\'ecent de Bushnell et Henniart \cite{BH}.\\

On ne red\'efinira pas l'espace des formes automorphes du groupe $\gl$ not\'e $\mathcal{A}(\gl(F)\backslash\gl(\adele))$ 
(cf \cite{Bu}, chapitre 3): ce sont des fonctions \`a croissance mod\'er\'ee, finies (pour l'action par translation de 
$K_f$ et du centre de l'alg\`ebre enveloppante $Z(\mathcal{U}(\mathfrak{gl}(2)_\C))$). Le sous-espace des formes  
paraboliques, d\'efini par  $$\mathcal{A}_0(\gl(F)\backslash\gl(\adele))=\{\varphi\in\mathcal{A}(\gl(F)\backslash\gl(\adele));
\int_{F\backslash \adele}\varphi(n(x)g)dx=0\}$$ est constitu\'e de fonctions \`a d\'ecroissance rapide; on appelle 
repr\'esentation automorphe (resp. parabolique) irr\'eductible un sous-espace irr\'eductible de $\mathcal{A}$
(resp. $\mathcal{A}_0$). Il se trouve que $\mathcal{A}_0(\omega)$, 
sous-espace de $\mathcal{A}_0$ sur lequel $Z(\adele)$ agit selon $\omega$, est somme directe (alg\'ebrique) de 
repr\'esentations automorphes.\\

Soit $(\pi,V_\pi)$ une repr\'esentation automorphe irr\'eductible de $\gl(\adele)$. Elle se factorise
sous la forme d'un produit tensoriel restreint: 
\begin{eqnarray}\label{tensor}
\pi\cong \bigotimes_{v}\pi_v
\end{eqnarray}
$\pi_v$ \'etant une repr\'esentation admissible de $\gl(F_v)$
si $v$ finie (respectivement un $(\mathfrak{gl}_2(\R),{\rm{O}}_2(\R))$-module si $v$ infinie), 
sa restriction \`a $Z(\adele)$ induit un caract\`ere de $\adele^\times/F^\times$, not\'e $\omega_\pi$, dit caract\`ere central de $\pi$, 
dont chaque composante locale $(\omega_\pi)_v$ est le caract\`ere central $\omega_{\pi_v}$ de 
$\pi_v$. La repr\'esentation automorphe $(\pi,V_\pi)$ admet un mod\`ele de Whittaker (cf. \cite{Bu}, chapitre 3.3), pour tout choix d'un caract\`ere additif $\psi=\otimes\psi_v$,
 qui se factorise aussi:
$$\Whitt\cong \bigotimes_{v}\mathcal{W}(\pi_v,\psi_v)$$
avec la compatibilit\'e importante: si $\varphi=\otimes\varphi_v$ dans (\ref{tensor}), alors l'\'el\'ement de Whittaker correspondant $W\varphi$ v\'erifie: $W\varphi(g)=\prod_vW\varphi_v(g_v)$. De plus, on a par unicit\'e du mod\`ele de Whittaker la relation:
$$W\varphi(g)=\int_{F\backslash\adele}\varphi(n(x)g)\overline{\psi(x)}dx.$$

Dans ces factorisations, presque toute $\pi_v$ est non ramifi\'ee, i.e. admet un vecteur fixe par 
$\gl(\mathcal{O}_{F_v})$, et on sait alors que ce vecteur est unique (\`a homoth\'eties pr\`es). 
Si $\pi_v$ est ramifi\'ee, l'admisibilit\'e de  $\pi_v$ assure quand m\^eme l'existence 
d'un entier $n(\pi_v)$ minimal tel que :
$$\Big\{x\in V_v|\forall k\in K_0(\varpi_v^{n(\pi_v)})\,\,\, \pi_v(k)x=\omega_{\pi_v}(k)x\Big\} \neq \{0\}$$
et un r\'esultat c\'el\`ebre de Casselman dit que cet espace est unidimensionnel (pour $v$ finie). 
Ici on a pos\'e: 
$$\forall k=\left(\begin{array}{cc}a & b\\ c & d\end{array}\right) \in K_0(\varpi_v^{n(\pi_v)}),
\,\,\, \omega_{\pi_v}(k):=\omega_{\pi_v}(d).$$
Si $v$ correspond \`a l'id\'eal maximal $\p$ on note $n(\pi_v)=n(\pi_p)$ et l'id\'eal 
${\rm{q}}_\pi=\prod \p^{n(\pi_\p)}$ est nomm\'e conducteur de $\pi$. L'id\'eal local $\varpi_v^{n(\pi_v)}\mathcal{O}_v$ est appel\'e conducteur de $\pi_v$ ou conducteur en $v$ de $\pi$.\\
Le point important est que les vecteurs (presqu')invariants du th\'eor\`eme de Casselman, dits 
vecteurs sp\'eciaux ou essentiels, ont pour transform\'ee de Mellin les facteurs locaux 
$L(s,\pi_v)$ de la fonction $L$ de $\pi$ (voir \cite{popa}), et permettent de traduire des 
\'enonc\'es concernant des fonctions $L$ dans le langage de l'analyse harmonique (cf section \ref{Petersson}). En r\'esum\'e: 
\begin{proposition}Pour toute place $v$ finie, on note $W_v^0$ l'unique \'el\'ement sp\'ecial
 du mod\`ele de Whittaker local tel que $W_v^0(1)=1$. Si $\psi_v$ est non ramifi\'e, $W_v^0$ v\'erifie alors:
\begin{eqnarray}\label{mellin}L(s,\pi_v)=\int_{F_v^\times}W_v^0\Big(\left(\begin{array}{cc}
a & 0\\ 0 & 1\end{array}\right)\Big)|a|_v^{s-1/2}d^\times a
\end{eqnarray}
(avec convergence pour $\Re(s)>0$).
\end{proposition}
Modulo des adaptations techniques, on peut montrer une telle \'egalit\'e dans le cas archim\'edien aussi (m\^eme r\'ef\'erence). D'ailleurs, la fonction not\'ee $W_\infty^0$ au d\'ebut de la section \ref{Petersson} est pr\'ecis\'ement, dans le cas des s\'eries discr\`etes de $\gl(F_\infty)$, cet \'el\'ement sp\'ecial.\\

Supposons, ce qui sera notre cas dans la suite, que l'on parte d'une repr\'esentation globale parabolique $\pi$ dont le caract\`ere central est trivial, et dont le conducteur est sans facteur carr\'e: on peut alors donner la forme des fonctions $L$ locales. Le r\'esultat suivant r\'esume la situation (cf. \cite{Go} pour les preuves):

\begin{proposition}
Soit $v\leftrightarrow \p$ une place finie de $F$, avec ces hypoth\`eses\\
$\bullet$ si $\pi_v$ est non-ramifi\'ee, alors $\pi_v$ est une s\'erie principale et on peut \'ecrire 
$$L(s,\pi_v)=(1-\alpha_{\pi,1}(\p)N(\p)^{-s})^{-1}(1-\alpha_{\pi,2}(\p)N(\p)^{-s})^{-1}$$
On a l'\'equation fonctionnelle locale:
$$L(s,\pi_v)=L(1-s,\pi_v)$$
$\bullet$ si $\pi_v$ est ramifi\'ee de conducteur $\varpi_v\mathcal{O}_v$, alors $\pi_v$ est une repr\'esentation sp\'eciale et on peut \'ecrire:
$$L(s,\pi_v)=(1-\alpha_{\pi}(\p)N(\p)^{-s})^{-1}$$
On a l'\'equation fonctionnelle :
$$L(s,\pi_v)=\varepsilon_{\pi_v}N(\p)^{\frac{1}{2}-s}L(1-s,\pi_v)$$
avec $\varepsilon_{\pi_v}=-N(\p)^{1/2}\alpha_{\pi}(\p)\in\{1,-1\}.$
\end{proposition}
  
\rem: Si $k$ est un entier positif, et $\pi_\R$ une repr\'esentation irr\'eductible de $\gl(\R)$, isomorphe \`a une s\'erie 
discr\`ete $\mathcal{D}(k)$ de caract\`ere central trivial, alors: 
$$L(s,\pi_\R)=(2\pi)^{-\left(s+\frac{k}{2}\right)}\Gamma\left(s+\frac{k}{2}\right).$$
L'\'equation fonctionnelle est dans ce cas:
$$L(s,\pi_\R)=i^{k+1}L(1-s,\pi_\R)$$

Ces r\'esultats donnent une description compl\`ete 
 des fonctions $L$ d'une certaine classe de repr\'esentations automorphes: celle des formes modulaires de Hilbert de conducteur 
sans facteurs carr\'es 
(voir ci-dessous la justification de la confusion entre formes et repr\'esentations, qui sera exploit\'ee en \ref{fonctionsL} dans un cas simple).

\begin{definition} Soit $\g{k}\in 2\N^d$. Une forme modulaire de 
Hilbert de poids $\g{k}$ est une repr\'esentation automorphe parabolique $\pi$ de caract\`ere central trivial telle que $\pi_\infty\cong \mathcal{D}(\g{k}-\g{1})=\bigotimes_{1\leq j\leq d}\mathcal{D}(k_j-1)$.
\end{definition}
{\sc{Notation}}: On note $\Pi_\q^{\g{k}}$ l'ensemble des formes (repr\'esentations) modulaires de Hilbert de poids $\g{k}$ (r\'ef\'erant au param\`etre \`a l'infini), de conducteur $\q$.\\

Une premi\`ere chose \`a pr\'eciser est le lien entre formes (fonctions) automorphes, et les repr\'esentations, pour justifier la \og confusion\fg{} entre les deux. En effet $\mathcal{H}_\q^{\g{k}}$ se d\'ecompose \`a l'aide des $\Pi_\q^{\g{k}}$: pour $\pi$ 
dans cette famille, \'ecrivons 
\begin{multline}
(\pi,V_\pi)^{K_\infty K_0(\q)}=\Big\{\varphi\in (\pi,V_\pi)|\varphi(gk_\infty k_f)=\varphi(g)\e(\g{k}
k_\infty),\\ \forall (g,k_\infty,k_f)\in \gl(\adele)\times K_\infty\times K_0(\q)\Big\}\nonumber
\end{multline}
Cet ensemble est unidimensionnel si $\rm{q}_\pi=\q$, engendr\'e par un vecteur sp\'ecial global 
$\varphi_\pi=\otimes_v\varphi_v^0$, produit des vecteurs sp\'eciaux locaux. Si ${\rm{q}}_\pi\neq\q$, 
on sait calculer sa dimension (gr\^ace \`a Casselman). 
Ceci donne donc:
\begin{eqnarray}\label{newold}
\mathcal{H}_\q^{\g{k}}=\bigoplus_{\mathfrak{r}|\q}\bigoplus_{\pi\in\Pi_{\mathfrak{r}}^{\g{k}}}(\pi,V_\pi)^{K_\infty\times 
K_0(\q)}
\end{eqnarray}
On appelle espace des \emph{formes nouvelles} l'espace $\bigoplus_{\pi\in\Pi_\q^{\g{k}}}\C\varphi_\pi$, et l'espaces des 
\emph{formes anciennes} 
est $ \bigoplus_{\substack{\mathfrak{r}|\q\\\mathfrak{r}\not=\q}}\bigoplus_{\pi\in\Pi_{\mathfrak{r}}^{\g{k}}}(\pi,V_\pi)^{K_\infty\times 
K_0(\q)}$. \\

Revenant aux fonctions $L$, posons $L(s,\pi_\infty)=(2\pi)^{-(\g{s}+\frac{\g{k}-\g{1}}{2})}\Gamma\left(s+\frac{\g{k}-\g{1}}{2}\right)$. On pose donc:
\begin{eqnarray}
L(s,\pi)=\prod_vL(s,\pi_v)=L(s,\pi_\infty)L(s,\pi_f).
\end{eqnarray}
produit eul\'erien convergent pour $\Re(s)>1$, et la partie finie est une s\'erie de Dirichlet, une fois le produit 
d\'evelopp\'e, convergente pour $\Re(s)>1$ de la forme:
$$L(s,\pi_f)=\sum_{\n\subset\entier}\lambda_\pi(\n)N(\n)^{-s}$$
index\'ee par les id\'eaux \emph{non nuls} de $\entier$. Cette fonction admet un prolongement analytique. L'\'equation fonctionnelle de $L(s,\pi)$ est le produit des 
\'equations locales, soit dans le cas o\`u le conducteur global $\textrm{q}_\pi=\q$ est premier:
\begin{eqnarray}
N(\textrm{q}_\pi)^{\frac{s}{2}}L(s,\pi)=\varepsilon_\pi N(\textrm{q}_\pi)^{\frac{1-s}{2}}L(1-s,\pi)
\end{eqnarray}
avec $\varepsilon_\pi=-i^{\g{k}}\lambda_\pi(\q)N(\q)^{1/2}\in\{1,-1\}$ (il n'y a de s\'erie discr\`ete de param\`etre $k-1$ que si $k$ est un entier pair).\\

Vu l'\'equation fonctionnelle, on pose $\Lambda(s,\pi)=N(\textrm{q}_\pi)^{s/2}L(s,\pi)$. 
Nous aurons besoin des valeurs $\Lambda(1/2,\pi)$ et $\Lambda(1/2,\pi)^2$, qui sont en dehors de la 
zone de convergence de (\ref{mellin}). 
\begin{proposition}
Soit $\pi$ une forme modulaire de poids $\g{k}$, conducteur $\q$ (quelconque). On a:
\begin{eqnarray}\label{valeur1}
\Lambda(1/2,\pi)=(1+\varepsilon_\pi)N(\q)^{1/4}\sum_{\n\subset \entier}
\frac{\lambda_\pi(\n)}{\sqrt{N(\n)}}F(N(\n)/N(\q)^{1/2})
\end{eqnarray}
avec:
\begin{eqnarray}
F(y)=\frac{1}{2i\pi}\int_{(3/2)}y^{-s}L\left(s+\frac{1}{2},\pi_\infty\right)\frac{ds}{s}
\end{eqnarray}
\begin{eqnarray}\label{valeur2}
\Lambda(1/2,\pi)^2=2N(\q)^{1/2}\sum_{\n\subset \entier}
\frac{\lambda_\pi(\n)}{\sqrt{N(\n)}}\tau(\n) G(N(\n)/N(\q))
\end{eqnarray}
avec:
\begin{eqnarray}
G(y)=\frac{1}{2i\pi}\int_{(3/2)}y^{-s}\zeta_F^{(\q)}(1+2s)L\left(s+\frac{1}{2},\pi_\infty\right)^2\frac{ds}{s}.
\end{eqnarray}
\end{proposition}

Pour ce faire, on pose pour $\pi\in\Pi_\q^{\g{k}}$:
$$I(\pi)=\int_{(3/2)}\Lambda(s+1/2,\pi)\frac{ds}{s}.$$
L'equation fonctionnelle (\ref{equation}) assure que:
$$(1+\varepsilon_\pi)I(\pi)=\Lambda(1/2,\pi)$$
puis, en d\'eveloppant $I(\pi)$ en s\'erie, on trouve bien:
\begin{eqnarray}
\Lambda(1/2,\pi)=(1+\varepsilon_\pi)N(\q)^{1/4}\sum_{\n\subset \entier}
\frac{\lambda_\pi(\n)}{\sqrt{N(\n)}}F(N(\n)/N(\q)^{1/2}).
\end{eqnarray}

L'autre \'egalit\'e se montre de la m\^eme mani\`ere. 
Ces sommes sont absolument convergentes.

\subsection{Repr\'esentations et fonctions de Bessel}
Dans cette section, $\pi$ d\'esigne une repr\'esentation unitaire irr\'eductible de $\gl(\R)$, que 
nous supposerons de caract\`ere central trivial pour simplifier. Soit $\psi$ un caract\`ere additif de $\R$. \\
Les s\'eries de Poincar\'e que nous d\'efinirons dans la section \ref{Petersson} sont produites 
\`a partir d'\'el\'ements particuliers du mod\`ele de Whittaker $\Whitt$, et lors du 
calcul du produit scalaire apparaissent des int\'egrales d'\og entrelacement\fg. Dans la proposition 
suivante, $J_{\pi,\psi}$ d\'esigne la fonction de Bessel associ\'ee (cf \cite{CPS}, \cite{Sou}), 
et $\omega=\left(\begin{array}{cc}0 & -1\\1 & 0\end{array}\right)$.

\begin{proposition}\label{entrelacement}
Soit $\pi$ une repr\'esentation unitaire irr\'eductible de $\pgl(\R)$. Pour $b\in\R^\times$, 
on pose:
\begin{displaymath}
\mathscr{E}:\begin{array}{cll}\Whitt & \longrightarrow & \Whitt\\
W & \longmapsto & \left[g\mapsto \displaystyle{\int_{N(\R)}\psi(-n)W\Big( \left(\begin{array}{cc}b & 0\\0 & 1\end{array}
\right)\omega ng\Big)dn}\right]
\end{array}
\end{displaymath}
L'int\'egrale d\'efinissant $\mathscr{E}$ est convergente, et on a:
\begin{eqnarray} \forall g\in\pgl(\R),\,
\mathscr{E}(W)(g)=\frac{1}{|m|}J_{\pi,\psi}(b)W(g)
\end{eqnarray}
si $\psi=\exp(2im\pi.)$.
\end{proposition}
C'est ce r\'esultat, \'ecrit sous une autre forme, qui est au c\oe ur de la proposition 9.4 de 
\cite{BM}, et par cons\'equent de la formule de Kuznetsov. Dans notre situation, en section \ref{Petersson},  
on pourrait conclure plus \'el\'ementairement (pour $F=\Q$, la formule de Petersson se passe 
de celui-ci), mais il est important de remarquer qu'il est central dans toutes les formes 
\og difficiles\fg des formules de traces de Petersson/Kuznetsov.

Nous serons dans le cas o\`u $\pi\cong \mathcal{D}(k-1)$, auquel cas \cite{CPS} donne, pour 
$\psi(x)=\exp(2i\pi mx)$:
\begin{eqnarray} \forall x>0,\,\, J_{\pi,\psi}(x)=(-1)^{k/2}2\pi|m|\sqrt{x}J_{k-1}(4\pi|m|\sqrt{x})
\end{eqnarray}
$J_\nu$ \'etant la fonction de Bessel classique, 
dont la forme la plus utilis\'ee par la suite sera (cf \cite{luo}):
\begin{eqnarray}\label{bessel}
J_\nu(x)=\int_{(\sigma)}\frac{\Gamma\Big(\frac{\nu-s}{2}\Big)}{\Gamma\Big(\frac{\nu+s}{2}+1\Big)}
\Big(\frac{x}{2}\Big)^sds, \,\, \textrm{si } 0<\sigma<\Re(\nu).
\end{eqnarray}

\section{La formule de Petersson}\label{Petersson}

Soient $\g{k}\in 2\N^d$ un vecteur pair, $\q$ un id\'eal \emph{quelconque} de $\entier$. 
Soit $\Pi_\q^k$ la famille des formes modulaires de Hilbert de conducteur $\q$, de caract\`ere 
central trivial. \\
Pour $g\in\gl(\adele)$, on note $g=g_\infty g_f$ la factorisation en places infinies et finies, et 
on \'ecrira parfois $g_\infty=(g_1,\dots,g_d)$ les composantes infinies (bien s\^ur, si $v$ est une 
place, $g_v$ d\'enote la composante en $v$ de $g$).\\
Toute cette section est adapt\'ee du travail de Venkatesh (\cite{V}, section 6), qui a \'etendu la formule de Kuznetsov dans ce contexte ad\'elique. Le lecteur est invit\'e \`a se r\'ef\'erer \`a cet article, dont nous suivons les notations.

\subsection{Analyse harmonique r\'eelle}
Sur les espaces $\mathcal{S}_{\g{k}}(Z_\infty^+\Gamma_0(\q,\A)\backslash \GL(F_\infty))$, on 
dispose de l'analyse de Fourier classique. Plus pr\'ecis\'ement, soit $\psi_\infty$ le quasi-caract\`ere 
additif de $\C^d$ d\'efini par :
$$\forall \g{z}\in\C^d,\, \psi_\infty(\g{z})=\exp(2i\pi\g{z})=\exp(2i\pi\sum_{j=1}^dz_j)$$
Comme un \'el\'ement $g$ de $\GL(F_\infty)$ admet une unique repr\'esentation
\begin{eqnarray}
g=z(g)n(g)a(g)k(g)= \left(\begin{array}{cc}\g{z} & 0\\0 & \g{z}\end{array}\right)
\left(\begin{array}{cc}1 & \g{x}\\0 & 1\end{array}\right)\left(\begin{array}{cc}\g{y} & 0\\0 & 1\end{array}\right)
k_\infty
\end{eqnarray}
avec $\g{z}\in F_\infty^{\times>0}$, $\g{x}\in F_\infty$, $\g{y}\in F_\infty^{\times>0}$, $k_\infty\in K_\infty$ (d\'ecomposition 
d'Iwasawa), on peut d\'efinir:
\begin{multline}
W_\infty^0(g)=\g{y}^{\g{k}/\g{2}}\psi_\infty(\g{k}(\g{x}+i\g{y}))\e(\g{k}k_\infty)\\
=\prod_{j=1}^dy_j^{k_j/2}\exp(2i\pi\sum_{j=1}^d k_j(x_j+iy_j))\exp(i\g{k\theta})
\end{multline}
La notation est coh\'erente (avec (\ref{mellin})), car cet \'el\'ement est pr\'ecis\'ement le vecteur nouveau du mod\`ele de Whittaker des s\'eries 
discr\`etes $\bigotimes_j\mathcal{D}(k_j-1)$ de $\gl(F_\infty)$ (restreint \`a $\GL(F_\infty)$).\\
La th\'eorie de Fourier classique permet d'\'ecrire, pour $g\in\GL(F_\infty)$, 
$\varphi_{\A}\in\mathcal{S}_{\g{k}}(Z_\infty^+
\Gamma_0(\q,\A)\backslash \GL(F_\infty))$, sous forme de s\'erie normalement convergente:

\begin{eqnarray}\label{fourier}
\varphi_{\A}(g)=\sum_{\substack{\alpha\in\A^{-1}\diff^{-1}\\\alpha\gg 0}}c(\alpha,\varphi_{\A})W_\infty^0\bigg(
\left(\begin{array}{cc}\g{\alpha} & 0\\0 & 1\end{array}\right)g\bigg)
\end{eqnarray}

\subsection{S\'eries de Poincar\'e}
On va d\'efinir sur $\mathcal{S}_\q^{\g{k}}$ des fonctions, g\'en\'eralisant les s\'eries de Poincar\'e 
classiques (sur les espaces de formes modulaires), comme l'a fait Venkatesh dans le contexte des formes de Maass.

Soient $\A$ et $\B$ deux id\'eaux fractionnaires.
On se donne aussi $\alpha$ (resp. $\alpha'$) un \'el\'ement non nul de $\A^{-1}\diff^{-1}$ 
(resp. $(\A\B^2)^{-1}\diff^{-1}$).
On notera $\Gamma_N(\q,\A)=N(F_\infty)\cap \Gamma_0(\q,\A)$, $Z_\Gamma=Z_\infty^+\cap \Gamma
_0(\q,\A)$.
Pour $g_\infty\in\GL(F_\infty)$, on pose, avec les coordonn\'ees d'Iwasawa:
\begin{eqnarray}
W_\infty^{\alpha}(g_\infty):=\g{y}^{\g{k}/2}\psi_\infty(\g{\alpha}(\g{x}+i\g{y}))\e(\g{k}k_\infty)
\end{eqnarray}
Ce qui revient au m\^eme:
\begin{eqnarray}
W_\infty^{\alpha}(g_\infty)=\g{j}(g_\infty)^{-\g{k}}\psi_\infty(\g{\alpha}g_\infty . \g{i})
\end{eqnarray}
Sur l'espace $\mathcal{S}_{\g{k}}(Z_\infty^+\Gamma_0(\q,\A)\backslash \GL(F_\infty))$, on peut 
d\'efinir la s\'erie convergente:
$$P_\A^\alpha(g_\infty)=\sum_{\gamma\in Z_\Gamma\Gamma_N(\q,\A)\backslash \Gamma_0(\q,\A)}W_\infty^{\alpha}(\gamma
g_\infty)$$
La convergence de cette s\'erie est bien connue (cf \cite{Ga}), ainsi que ses propri\'et\'es de 
modularit\'e, gr\^ace \`a la derni\`ere \'ecriture de $W_\infty^\alpha$.
\begin{definition} Avec ces notations, on d\'efinit pour tout couple $(\alpha,\A)$ une s\'erie de 
Poincar\'e encore not\'ee $P_\A^\alpha\in\mathcal{S}_\q^{\g{k}}$ par: 
\begin{displaymath}
\left\{\begin{array}{ll}
\big(P_\A^\alpha\big)_{\overline{\B}} = 0 \textrm{ si } \overline{\B}\neq\overline{\A} \\
\big(P_\A^\alpha\big)_{\overline{\B}}  =  P_\A^\alpha \textrm{ si } \overline{\B}=\overline{\A}
\end{array}\right.
\end{displaymath}
les composantes \'etant celles de l'isomorphisme (\ref{adelisation}).
\end{definition}
On obtient traditionnellement la formule de Petersson en calculant de deux fa\c cons diff\'erentes 
le produit scalaire de deux s\'eries de Poincar\'e. Il sera crucial plus tard de disposer d'une formule du type \og $\sum_\pi \lambda_\pi(\n)\lambda_\pi(\m)$\fg avec $\m$ et $\n$ quelconques (et pas forc\'ement dans la m\^eme classe): pour cela, il faut permuter les composantes connexes de $\gl(F)Z_\infty^+\backslash \gl(\adele)/K_0(\q)$, et donc user de l'action naturelle \`a droite de $\gl(\adele)$ sur $\mathcal{S}_\q^{\g{k}}$, not\'ee $\rho$. 
On pose:
$$\rho(\B)P_{\A\B^2}^{\alpha'}(g):=P_{\A\B^2}^{\alpha'}\Big( g \left(\begin{array}{cc} {\rm{id}}(\B) 
& 0\\0 & \rm{id}(\B)\end{array}\right) \Big),\,\,\forall g\in\gl(\adele)$$
Puisqu'on peut \'ecrire, par le principe de l'approximation forte:
$$\left(\begin{array}{cc} {\rm{id}}(\B) 
& 0\\0 & \rm{id}(\B)\end{array}\right) \left(\begin{array}{cc} {\rm{id}}(\A) 
& 0\\0 & 1\end{array}\right)g_\infty =\gamma^{-1} \left(\begin{array}{cc} {\rm{id}}(\A\B^2) 
& 0\\0 & 1\end{array}\right) g_\infty' k_f$$
avec $\gamma\in\gl(F)$, $k_f\in K_0(\q)$ et $g_\infty,g'_\infty\in\GL(F_\infty)$, il vient d'une 
part: $\gamma\in\Gamma(\A\rightarrow \A\B^2):=\GL(F)\cap \left(\begin{array}{cc} {\rm{id}}(\A\B^2) 
& 0\\0 & 1\end{array}\right) K_0(\q)\left(\begin{array}{cc} {\rm{id}}(\A\B)^{-1} 
& 0\\0 & \rm{id}(\B)^{-1}\end{array}\right) $, d'autre part, en notant $\gamma_
{\A\rightarrow \A\B^2}$ la composante infinie d'un tel $\gamma$:
\begin{eqnarray}
\rho(\B)P_{\A\B^2}^{\alpha'}\Big(g_\infty\left(\begin{array}{cc} {\rm{id}}(\A) 
& 0\\0 & 1\end{array}\right) \Big) & = & P_{\A\B^2}^{\alpha'}\Big(\gamma_{\A\rightarrow \A\B^2} g_\infty\left(\begin{array}{cc}
 {\rm{id}}(\A\B^2) & 0\\0 & 1\end{array}\right) \Big)\nonumber\\
\rho(\B)P_{\A\B^2}^{\alpha'}\Big(g_\infty\left(\begin{array}{cc} {\rm{id}}(\A') 
& 0\\0 & 1\end{array}\right) \Big)&=& 0 \textrm{ si } \overline{\A}'\neq\overline{\A}\nonumber
\end{eqnarray}
On peut donc r\'e\'ecrire:
\begin{multline}
\rho(\B)P_{\A\B^2}^{\alpha'}\Big(g_\infty \left(\begin{array}{cc} {\rm{id}}(\A) 
& 0\\0 & 1\end{array}\right) \Big)=\sum_{\gamma\in Z_\Gamma\Gamma_N(\q,\A\B^2)\backslash \Gamma_0
(\q,\A\B^2)}W_\infty^{\alpha'}(\gamma \gamma_{\A\rightarrow \A\B^2}g_\infty) \\=
\sum_{\gamma\in Z_\Gamma\Gamma_N(\q,\A\B^2)\backslash \Gamma(\A\rightarrow\A\B^2)}W_\infty^
{\alpha'}(\gamma g_\infty)\nonumber
\end{multline}
En suivant \cite{V}, on peut \'enoncer:

\begin{lemme}\label{parametrisation} On a les propri\'et\'es suivantes:\\
$\bullet $ On a explicitement:
\begin{multline}
\Gamma(\A\rightarrow \A\B^2)=\Big\{ \left(\begin{array}{cc} a 
& b\\c & d\end{array}\right) \in \GL(F) \mid a \in\B,b\in\A\B,c\in\q\A^{-1}\B^{-1},\\
d\in\B^{-1},ad-bc\in\mathcal{O}_F^{\times +}\Big\}
\end{multline}
$\bullet $ De plus l'application $ \left(\begin{array}{cc} a 
& b\\c & d\end{array}\right) \in\Gamma(\A\rightarrow \A\B^2)\rightarrow (c,a,ad-bc)$ induit une bijection 
de $\Gamma_N(\q,\A\B^2)\backslash \Gamma(\A\rightarrow \A\B^2)/\Gamma_N(\q,\A)$ sur l'ensemble 
$\Big\{ (c,x,\varepsilon)\mid c\in\q\A^{-1}\B^{-1},\varepsilon\in\mathcal{O}_F^{\times +},x\in\B,
x \textrm{ engendre }\B/\A\B^2(c)\Big\}$.\\
$\bullet $ Enfin, si on note $B^+(F)=\Big\{ \left(\begin{array}{cc} a 
& b\\0 & d\end{array}\right) \in\GL(F)\Big\}$ alors $\Gamma_N(\q,\A\B^2)\backslash \Gamma(\A\rightarrow
\A\B^2)\cap B^+(F)$ est non vide si et seulement si $\B$ est principal, et admet le syst\`eme de  
repr\'esentants $\Big\{ \left(\begin{array}{cc} [\B]\varepsilon 
& 0\\0 & [\B]^{-1}\end{array}\right) ,\varepsilon\in\mathcal{O}_F^{\times+}\Big\}$, $[\B]$ \'etant un 
g\'en\'erateur fix\'e de $\B$.

\end{lemme}
Nous renvoyons \`a \cite{V} pour la preuve. 

\subsection{La formule de Petersson pour $\mathcal{S}_\q^{\g{k}}$}
Nous allons maintenant \'evaluer le produit scalaire $\left<P_\A^\alpha,P_{\A\B^2}^{\alpha'}
\right>_{\mathcal{S}_\q^{\g{k}}}$, ceci avec les m\^emes notations que pr\'ec\'edemment. Comme pour 
la formule de Kuznetsov, on fait ce calcul de deux fa\c cons: l'une, directe, utilise la 
param\'etrisation du dernier lemme, avec la propri\'et\'e d'\og entrelacement\fg du lemme \ref{entrelacement} (dont 
l'intervention est moralement justifi\'ee par la d\'efinition des s\'eries de Poincar\'e, via des 
\'el\'ements sp\'eciaux des mod\`eles de Whittaker aux places infinies); l'autre vient de la 
d\'ecomposition spectrale de l'espace $\mathcal{S}_\q^{\g{k}}$, qui ici est techniquement facile 
\`a g\'erer, car cet espace est de dimension finie.\\

\noindent $\bullet $ {\bf{Calcul direct du produit scalaire:}}\\
On doit donc calculer:
\begin{multline}
\left<P_\A^\alpha,\rho(\B)P_{\A\B^2}^{\alpha'}\right>_{\mathcal{S}_\q^{\g{k}}}\\=\int_{Z_\infty^+\gl(F)\backslash \gl(\adele)/
K_0(\q)}\overline{P_\A^\alpha(g)}P_{\A\B^2}^{\alpha'}\Big(g\left(\begin{array}{cc} {\rm{id}}(\B) 
& 0\\0 & \rm{id}(\B)\end{array}\right) \Big)dg\\=\sum_{\overline{\A}\in\mathscr{C}\!\ell^+(F)}\int_{Z_\infty^+
\Gamma_0(\q,\A)\backslash\GL(F_\infty)} \Big(\overline{P_\A^\alpha}\Big)_\A(g) \Big(\rho(\B)P_{\A\B^2}^{
\alpha'}\Big)_\A(g)dg\\=\int_{Z_\infty^+\Gamma_0(\q,\A)\backslash\GL(F_\infty)}
\overline{P_\A^\alpha\Big(g\left(\begin{array}{cc} {\rm{id}}(\A) 
& 0\\0 & 1\end{array}\right) \Big)}P_{\A\B^2}^{\alpha'}\Big(\gamma_{\A\rightarrow\A\B^2}
g\left(\begin{array}{cc} 
{\rm{id}}(\A\B^2) & 0\\0 & 1\end{array}\right) \Big)dg\\
=\int_{Z_\infty^+\Gamma_0(\q,\A)\backslash\GL(F_\infty)}\bigg(\sum_{\gamma\in Z_\Gamma
\Gamma_N(\q,\A)\backslash \Gamma_0(\q,\A)}\overline{W_\infty^\alpha(\gamma g)}\bigg)
\Big(P_{\A\B^2}^{\alpha'}\Big)_{\A\B^2}(\gamma_{\A\rightarrow\A\B^2}g)dg\\
=\int_{Z_\infty^+\Gamma_N(\q,\A)\backslash \GL(F_\infty)}\overline{W_\infty^\alpha(g)}
\bigg(\sum_{\gamma\in Z_\Gamma\Gamma_N(\q,\A\B^2)\backslash \Gamma_0(\q,\A\B^2)}
W_\infty^{\alpha'}(\gamma\gamma_{\A\rightarrow\A\B^2}g)\bigg)dg\\
=\int_{Z_\infty^+\Gamma_N(\q,\A)\backslash \GL(F_\infty)}\overline{W_\infty^\alpha(g)}
\bigg(\sum_{\gamma\in Z_\Gamma\Gamma_N(\q,\A\B^2)\backslash \Gamma(\A\rightarrow 
\A\B^2)}W_\infty^{\alpha'}(\gamma g)\bigg)dg.
\nonumber
\end{multline}
On peut donc \'ecrire: $\left<P_\A^\alpha,P_{\A\B^2}^{\alpha'}\right>=I_1+I_2$, avec:
\begin{multline}
I_1=\int_{Z_\infty^+\Gamma_N(\q,\A)\backslash \GL(F_\infty)}\overline{W_\infty^\alpha(g)}\bigg(
\sum_{\gamma\in Z_\Gamma\Gamma_N(\q,\A\B^2)\backslash (\Gamma(\A\rightarrow\A\B^2)\cap B^+(F))}
W_\infty^{\alpha'}(\gamma g)\bigg)dg
\end{multline}
\begin{multline}
I_2=\int_{Z_\infty^+\Gamma_N(\q,\A)\backslash \GL(F_\infty)}\overline{W_\infty^\alpha(g)}\bigg(
\sum_{\substack{\gamma\in Z_\Gamma\Gamma_N(\q,\A\B^2)\backslash \Gamma(\A\rightarrow\A\B^2)\\
\gamma\notin B^+(F)}}W_\infty^{\alpha'}(\gamma 
g)\bigg) dg.
\end{multline}
Calculons d'abord $I_1$: Par le lemme \ref{parametrisation}, on peut \'ecrire, en se rappelant 
que l'ensemble d'indices n'est pas vide que si $\B$ est principal:
\begin{multline}
I_1=\sum_{\varepsilon\in\mathcal{O}_F^{\times+}}\int_{Z_\infty^+\Gamma_N(\q,\A)\backslash \GL(F_\infty)}
\overline{W_\infty^\alpha(g)}W_\infty^{\alpha'}\Big(\left(\begin{array}{cc} \varepsilon [\B] 
& 0\\0 & [\B]^{-1}\end{array}\right) g\Big) dg\\=\sum_{\varepsilon\in\mathcal{O}_F^{\times +}}
\int_{F_\infty^{\times+}}\frac{d^\times\g{y}}{\g{y}}
\overline{W_\infty^\alpha\Big(\left(\begin{array}{cc} \g{y} 
& 0\\0 & 1\end{array}\right) \Big)}W_\infty^{\alpha'}\Big(\left(\begin{array}{cc} \g{\varepsilon 
y[\B]^2} & 0\\0 & 1\end{array}\right) \Big)\times \\
\int_{\Gamma_N(\q,\A)\backslash N(F_\infty)}
\overline{W_\infty^\alpha(n(\g{x}))}W_\infty^{\alpha'}(n(\g{\varepsilon [\B]^2x}))dx
\nonumber
\end{multline} 
cette derni\`ere ligne \'etant obtenue avec les coordonn\'ees d'Iwasawa. L'int\'egrale interne vaut :
$$\int_{\A\backslash F_\infty}\overline{W_\infty^\alpha\Big(\left(\begin{array}{cc} 1 
& \g{x}\\0 & 1\end{array}\right) \Big)}W_\infty^{\alpha'}\Big(\left(\begin{array}{cc} 1 
& \g{x\varepsilon[\B]^2}\\0 & 1\end{array}\right) \Big)dx={\rm{vol}}(\A\backslash F_\infty)
\mathbbmss{1}_{\alpha=\varepsilon[\B]^2\alpha'}.$$
En l'ins\'erant dans le calcul, et en utilisant la valeur explicite de $W_\infty^{\alpha}$, il vient:
\begin{multline}\label{I1}
I_1={\rm{vol}}(\A\backslash F_\infty)\sum_{\varepsilon\in\mathcal{O}_F^{\times+}}
(\g{[\B]\varepsilon})^{\g{k}/2}\mathbbmss{1}_{\alpha=\varepsilon[\B]^2\alpha'}\int_{F_\infty^{\times
+}}\frac{d\g{y}}{\g{y}}\g{y}^{\g{k-1}}e^{-2\pi\g{\alpha y}-2\pi\g{\alpha'\varepsilon[\B]^2y}}\\
=\frac{\Gamma(\g{k-1})}{(4\pi)^{\g{k-1}}} {\rm{vol}}(\A\backslash F_\infty)
\frac{N(\B)}{(\g{\alpha\alpha'})^{\frac{\g{k-1}}{2}}}
\mathbbmss{1}_{\alpha\alpha'^{-1}[\B]^{-2}\in\mathcal{O}_F^{\times+}}.
\end{multline}
Pr\'ecisons que ${\rm{vol}}(\A\backslash F_\infty)$ d\'esigne le volume de $\A\backslash F_\infty$ 
pour la mesure de Lebesgue de $F_\infty=\R^d$, $\A$, plong\'e diagonalement, \'etant un r\'eseau. \\

Maintenant, $I_2$:

\begin{multline}
I_2=\int_{Z_\infty^+\Gamma_N(\q,\A)\backslash \GL(F_\infty)}\sum_{\gamma\in Z_\Gamma\Gamma_N(\q,\A\B^2)
\backslash \big(\Gamma(\A\rightarrow\A\B^2)\setminus  B^+(F)\big)}\overline{W_\infty^\alpha(g)}W_\infty^{\alpha'}(\gamma g)
dg\nonumber
\end{multline}

\begin{multline}
=\int_{Z_\infty^+\backslash \GL(F_\infty)}\sum_{\gamma\in Z_\Gamma\Gamma_N(\q,\A\B^2)\backslash 
\big(\Gamma(\A\rightarrow\A\B^2)\setminus B^+(F)\big)/\Gamma_N(\q,\A)}\overline{W_\infty^\alpha(g)}
W_\infty^{\alpha'}(\gamma g)dg\\
=\sum_{\gamma\in Z_\Gamma\Gamma_N(\q,\A\B^2)\backslash 
\big(\Gamma(\A\rightarrow\A\B^2)\setminus B^+(F)\big)/\Gamma_N(\q,\A)}\int_{Z_\infty^+\backslash 
\GL(F_\infty)}\overline{W_\infty^\alpha(g)}W_\infty^{\alpha'}(\gamma g)dg
\nonumber
\end{multline}
apr\`es interversion de somme et int\'egrale, justifi\'ee a posteriori par les calculs \`a venir.\\
Prenons un syst\`eme de repr\'esentants $\gamma$, et \'ecrivons leur d\'ecomposition de Bruhat:
$$\gamma=\left(\begin{array}{cc}a & b\\c & d\end{array}\right)=\left(\begin{array}{cc}1
 & a/c\\0 & 1\end{array}\right)\left(\begin{array}{cc}0 & 1\\-1 & 0\end{array}\right)
\left(\begin{array}{cc}c & 0\\0 & (ad-bc)c^{-1}\end{array}\right)
\left(\begin{array}{cc}1 & d/c\\0 & 1\end{array}\right)$$
l\'egitime, car $c\neq 0$, par hypoth\`ese sur l'ensemble d'indices. En notant $n_1(\gamma)$, 
$n_2(\gamma)$ les composantes unipotentes, $\omega$ pour l'\'el\'ement de Weyl, et $a_\gamma$ pour 
la composante centrale, on a donc:
\begin{multline}
I_2=\sum_{\gamma}\int_{Z_\infty^+\backslash \GL(F_\infty)}\overline{W_\infty^\alpha(g)}
W_\infty^{\alpha'}\bigg(n_1(\gamma) \omega a(\gamma)n_2(\gamma)g\bigg)dg\\
=\sum_{\gamma}W_\infty^\alpha(n_2(\gamma))W_\infty^{\alpha'}(n_1(\gamma))\int_{Z_\infty^+\backslash 
\GL(F_\infty)}\overline{W_\infty^\alpha(g)}W_\infty^{\alpha'}\bigg(\omega a(\gamma)g\bigg)dg
\nonumber
\end{multline}
avec $W_\infty^\alpha(n_2(\gamma))=\psi_\infty(\frac{\g{\alpha d}}{\g{c}})$, et $W_\infty^{\alpha'}
(n_1(\gamma))=\psi_\infty(\frac{\g{\alpha' a}}{\g{c}})$. On utilise alors la description explicite 
d'un ensemble de repr\'esentants donn\'ee par le lemme \ref{parametrisation} et cela donne:
\begin{multline}
I_2=\frac{1}{2}\sum_{c\in\A^{-1}\B^{-1}\q\setminus \{0\}}\sum_{\varepsilon\in\mathcal{O}_F^{\times +}/
\mathcal{O}_F^{\times 2}}KS(\varepsilon\alpha,\A;\alpha',\A\B^2;c,\A\B)\times\\
\int_{Z_\infty^+\backslash 
\GL(F_\infty)}\overline{W_\infty^\alpha(g)}W_\infty^{\alpha'}\bigg(\omega\left(\begin{array}{cc}
\g{c^2\varepsilon^{-1}} & 0\\0 & 1\end{array}\right)g\bigg)dg\nonumber
\end{multline}
o\`u l'on a not\'e la somme de Kloosterman comme Venkatesh (voir section \ref{Kloosterman}):
\begin{multline}
KS(\varepsilon\alpha,\A;\alpha',\A\B^2;c,\A\B)=\sum_{x\in(\B^{-1}/\A (c))^\times}\exp\left(2i\pi\textrm{Tr}_{F/\Q}
\left( \frac{\varepsilon\alpha x+\alpha'\bar{x}}{c}\right)\right)
\end{multline}
Remarquer qu'il faut \'echanger les r\^oles de $d$ et $a$ pour appliquer le lemme \ref{parametrisation} de param\'etrisation: comme $a$ parcourt $(\B/\A\B^2(c))^\times$, $d$ d\'ecrit $(\B^{-1}/\A(c))^\times$. Le facteur $\frac{1}{2}$ 
appara\^it car l'ensemble d'indices ici n'est pas exactement l'ensemble des doubles classes 
\'etudi\'e au lemme \ref{parametrisation}: il y a un quotient suppl\'ementaire par $Z_\Gamma$, et 
une invariance par multiplication par $-1$. On a donc:

\begin{multline}
I_2=\frac{1}{2}\sum_{\substack{c\in\A^{-1}\B^{-1}\q\setminus \{0\}\\ \varepsilon\in\mathcal{O}_F^{\times +}/\mathcal{O}_F^{\times 2}}}
KS(\varepsilon\alpha,\A;\alpha',\A\B^2;c,\A\B)\int_{F_\infty^{\times +}}\frac{d^\times\g{y}}{\g{y}}
\overline{W_\infty^\alpha\bigg(\left(\begin{array}{cc}\g{y} & 0\\0 & 1\end{array}\right)\bigg)}\times\\
\int_{F_\infty}\psi(-\g{\alpha x})W_\infty^{\alpha'}\bigg(\omega\left(\begin{array}{cc}
\g{c^2\varepsilon^{-1}} & 0\\0 & 1\end{array}\right)\left(\begin{array}{cc}
1 & \g{x}\\0 & 1\end{array}\right)\left(\begin{array}{cc}
\g{y} & 0\\0 & 1\end{array}\right) \bigg)dx\\
=\frac{1}{2}\sum_{\substack{c\in\A^{-1}\B^{-1}\q\setminus \{0\}\\ \varepsilon\in\mathcal{O}_F^{\times +}/\mathcal{O}_F^{\times 2}}}
KS(\varepsilon\alpha,\A;\alpha',\A\B^2;c,\A\B)\int_{F_\infty^{\times +}}\frac{d^\times\g{y}}{\g{y}}
\overline{W_\infty^\alpha\bigg(\left(\begin{array}{cc}\g{y} & 0\\0 & 1\end{array}\right)\bigg)}
\times\\ \frac{N(\alpha')}{N(\alpha)}\int_{F_\infty}\psi_\infty(\g{\alpha'x})
W_\infty^{\alpha'}\bigg(\left(\begin{array}{cc}
\frac{\g{\varepsilon\alpha}}{\g{\alpha'c^2}} & 0\\0 & 1\end{array}\right)\omega\left(\begin{array}{cc}
1 & \g{x}\\0 & 1\end{array}\right)\left(\begin{array}{cc}
\frac{\g{\alpha y}}{\g{\alpha'}} & 0\\0 & 1\end{array}\right) \bigg)dx.
\nonumber
\end{multline}
L'int\'egrale interne est calculable gr\^ace au lemme \ref{entrelacement}, et pr\'ecis\'ement, on a:
\begin{multline}
\int_{F_\infty}\psi_\infty(-\g{\alpha'x})
W_\infty^{\alpha'}\bigg(\left(\begin{array}{cc}\frac{\g{\varepsilon\alpha}}{\g{\alpha'c^2}}& 0
\\0 & 1\end{array}\right)\omega\left(\begin{array}{cc}
1 & \g{x}\\0 & 1\end{array}\right)\left(\begin{array}{cc}
\frac{\g{\alpha y}}{\g{\alpha'}} & 0\\0 & 1\end{array}\right) \bigg)dx\\
=N(\alpha')^{-1}J_{\mathcal{D}(\g{k-1}),\psi_\infty(\g{\alpha'}.)}\left(\frac{\g{\varepsilon\alpha}}
{\g{\alpha'c^2}}\right)W_\infty^{\alpha'}\left(\left(\begin{array}{cc}\frac{\g{\alpha y}}{\g{\alpha'}} & 0\\
0 & 1\end{array}\right)\right)
\end{multline}
avec:
\begin{eqnarray}
J_{\mathcal{D}(\g{k-1}),\psi_\infty(\g{\alpha'}.)}(\g{z})=(2\pi)^d(-1)^{\g{k}/2}N(\alpha')
|\g{z}|^{1/2}J_{\g{k-1}}(4\pi|\g{\alpha'}|\sqrt{\g{ z}}).
\end{eqnarray}
En cons\'equence, en se rappelant que $\alpha,\alpha'$ sont totalement positifs:
\begin{multline}
I_2=\frac{1}{2}\sum_{\substack{c\in\A^{-1}\B^{-1}\q\setminus \{0\}\\ \varepsilon\in\mathcal{O}_F^{\times +}/\mathcal{O}_F^{\times 2}}}
\frac{KS(\varepsilon\alpha,\A;\alpha',\A\B^2;c,\A\B)}{N(c)}\times\frac{\sqrt{N(\alpha\alpha')}}{N(\alpha')}
(-1)^{\g{k}/2}(2\pi)^d\times\\J_{\g{k-1}}\left(4\pi\frac{\sqrt{\g{\varepsilon\alpha\alpha'}}}{|\g{c}|}\right)
\int_{F_\infty^{\times +}}\frac{d^\times\g{y}}{\g{y}}
\overline{W_\infty^\alpha\bigg(\left(\begin{array}{cc}\g{y} & 0\\0 & 1\end{array}\right)\bigg)}
W_\infty^{\alpha'}\bigg(\left(\begin{array}{cc}\frac{\g{\alpha y}}{\alpha'} & 0\\0 & 1\end{array}\right)\bigg)\\
=\frac{1}{2}\sum_{\substack{c\in\A^{-1}\B^{-1}\q\setminus \{0\}\\ \varepsilon\in\mathcal{O}_F^{\times +}/\mathcal{O}_F^{\times 2}}}
\frac{KS(\varepsilon\alpha,\A;\alpha',\A\B^2;c,\A\B)}{N(c)}\times\frac{\sqrt{N(\alpha\alpha')}}{N(\alpha')}
(-1)^{\g{k}/2}(2\pi)^d\times\\J_{\g{k-1}}\left(4\pi\frac{\sqrt{\g{\varepsilon\alpha\alpha'}}}{|\g{c}|}\right)
\left(\frac{\g{\alpha}}{\g{\alpha'}}\right)^{\g{k}/2}\int_{F_\infty^{\times +}}\frac{d^\times\g{y}}{\g{y}}
\left|W_\infty^\alpha\bigg(\left(\begin{array}{cc}\g{y} & 0\\0 & 1\end{array}\right)\bigg)\right|^2\nonumber
\end{multline}

\begin{multline}
=\frac{1}{2}\sum_{\substack{c\in\A^{-1}\B^{-1}\q\setminus \{0\}\\ \varepsilon\in\mathcal{O}_F^{\times +}/\mathcal{O}_F^{\times 2}}}
\frac{KS(\varepsilon\alpha,\A;\alpha',\A\B^2;c,\A\B)}{N(c)}\times\frac{\sqrt{N(\alpha\alpha')}}{N(\alpha')}
(-1)^{\g{k}/2}(2\pi)^d\times\\J_{\g{k-1}}\left(4\pi\frac{\sqrt{\g{\varepsilon\alpha\alpha'}}}{|\g{c}|}\right)
\left(\frac{\g{\alpha}}{\g{\alpha'}}\right)^{\g{k}/2}\times\frac{\Gamma(\g{k-1})}{(4\pi\g{\alpha})^
{\g{k}-1}}
\nonumber
\end{multline}
soit, apr\`es simplification:
\begin{multline}\label{I2}
I_2=\frac{\Gamma(\g{k-1})(-1)^{\g{k}/2}(2\pi)^d}{2(4\pi\sqrt{\g{\alpha\alpha'}})^{\g{k-1}}}\times\\
\sum_{\substack{c\in\A^{-1}\B^{-1}\q\setminus \{0\}\\ \varepsilon\in\mathcal{O}_F^{\times +}/\mathcal{O}_F^{\times 2}}}
\frac{KS(\varepsilon\alpha,\A;\alpha',\A\B^2;c,\A\B)}{N(c)}J_{\g{k-1}}\left(4\pi\frac{\sqrt{
\g{\varepsilon\alpha\alpha'}}}{|\g{c}|}\right).
\end{multline}
\\

\noindent $\bullet ${\bf{ Calcul spectral:}}\\
Soit $\{\phi\}$ une base orthonorm\'ee quelconque (finie) de $\mathcal{S}_\q^{\g{k}}$. On supposera cette base \emph{adapt\'ee} \`a la d\'ecomposition (\ref{decomposition}), c'est-\`a-dire r\'eunion de bases orthonorm\'ees de chacune des composantes isotypiques de (\ref{decomposition}). 
On peut \'ecrire:
$$\Big<P_\A^\alpha,\rho(\B)P_{\A\B^2}^{\alpha'}\Big>_{\mathcal{S}_\q^{\g{k}}}=\sum_{\phi}
\Big<P_\A^\alpha,\phi\Big>_{\mathcal{S}_\q^{\g{k}}}\overline{\Big<\rho(\B)P_{\A\B^2}^{\alpha'},
\phi\Big>}_{\mathcal{S}_\q^{\g{k}}}.$$
On a:
\begin{multline}
\Big< P_\A^\alpha,\phi\Big>_{\mathcal{S}_\q^{\g{k}}} =
\int_{Z_\infty^+\gl(F)\backslash \gl(\adele)/K_0(\q)}\overline{P_\A^\alpha(g)}\phi(g)dg\\=
\int_{Z_\infty^+\Gamma_0(\q,\A)\backslash \GL(F_\infty)}\overline{P_\A^\alpha(g)}
\phi\left(g\left(\begin{array}{cc}{\rm{id}}(\A) & 0\\0 & 1 \end{array}\right)\right)dg\\=
\int_{Z_\infty^+\Gamma_N(\q,\A)\backslash \GL(F_\infty)}\overline{W_\infty^\alpha(g)}\phi
\left(g\left(\begin{array}{cc}{\rm{id}}(\A) & 0\\0 & 1 \end{array}\right)\right)dg.\\
\nonumber
\end{multline}
Rappelons maintenant que $$\phi\left(g\left(\begin{array}{cc}{\rm{id}}(\A) & 0\\0 & 1 \end{array}
\right)\right)=\phi_\A(g)=\sum_{\substack{\xi\in\A^{-1}\diff^{-1}\\\xi\gg 0}}c(\xi,\A,\phi)W_\infty^0\left(
\left(\begin{array}{cc}\g{\xi} & 0\\0 & 1 \end{array}\right)g\right).$$
Il vient donc:
\begin{multline}
\Big<P_\A^\alpha,\phi\Big>_{\mathcal{S}_\q^{\g{k}}}=\int_{F_\infty^{\times+}}\frac{d^\times\g{y}}{\g{y}}
\overline{W_\infty^\alpha\left(\left(\begin{array}{cc}\g{y} & 0\\0 & 1 \end{array}\right)\right)}\times\\
\int_{\A\backslash F_\infty}\psi_\infty(-\g{\alpha x})\phi_\A\left(\left(\begin{array}{cc}
1 & \g{x}\\0 & 1 \end{array}\right)\left(\begin{array}{cc}\g{y} & 0\\0 & 1 \end{array}\right)
\right)\\
=\int_{F_\infty^{\times+}}\frac{d^\times\g{y}}{\g{y}}
\overline{W_\infty^\alpha\left(\left(\begin{array}{cc}\g{y} & 0\\0 & 1 \end{array}\right)\right)}\times\\
\frac{c(\alpha,\A,\phi)}{\g{\alpha}^{\g{k}/2}}{\rm{vol}}(\A\backslash F_\infty)\frac{\Gamma(\g{k-1})}
{(4\pi)^{\g{k-1}}} W_\infty^0\left(\left(\begin{array}{cc}\g{\alpha y} & 0\\0 & 1 \end{array}\right)
\right)
\nonumber
\end{multline}
soit finalement:
\begin{eqnarray}
\Big<P_\A^\alpha,\phi\Big> _{\mathcal{S}_\q^{\g{k}}} =\frac{c(\alpha,\A,\phi)}{\g{\alpha}^{\g{k}/2-\g{1}}}\times{\rm{vol}}
(\A\backslash F_\infty)\times\frac{\Gamma(\g{k-1})}{(4\pi)^{\g{k-1}}}.
\end{eqnarray}
De la m\^eme mani\`ere,
\begin{multline}\nonumber
\Big<\rho(\B)P_{\A\B^2}^{\alpha'},\phi\Big>_{\mathcal{S}_\q^{\g{k}}} \\=\int_{Z_\infty^+\gl(F)\backslash \gl(\adele)/K_0(\q)}
\overline{P_{\A\B^2}^{\alpha'}\left(g\left(\begin{array}{cc}{\rm{id}}(\B) & 0\\0 & {\rm{id}}(\B) \end{array}
\right)\right)}\phi(g)dg\\
=\omega_\phi(\B^{-1})\int_{Z_\infty^+\gl(F)\backslash \gl(\adele)/K_0(\q)}
\overline{P_{\A\B^2}^{\alpha'}(g)}\phi(g)dg
\end{multline}
car on a choisi la base $\{\phi\}$ de telle sorte que l'action du centre soit un twist par un Gr\"ossencharakter. Donc:
\begin{multline}
\Big<\rho(\B)P_{\A\B^2}^{\alpha'},\phi\Big>_{\mathcal{S}_\q^{\g{k}}} \\=\omega_\phi(\B^{-1})\int_{Z_\infty^+\Gamma_0(\q,\A\B^2)\backslash \GL(F_\infty)}
\overline{P_{\A\B^2}^{\alpha'}(g)}\phi\left(g\left(\begin{array}{cc}{\rm{id}}(\A\B^2) & 0\\0 & 1 \end{array}
\right)\right)dg\\
=\omega_\phi(\B^{-1})\times\frac{c(\alpha',\A\B^2,\phi)}{\g{\alpha'}^{\g{k}/2-\g{1}}}\times
{\rm{vol}}(\A\B^2\backslash F_\infty)\times\frac{\Gamma(\g{k-1})}{(4\pi)^{\g{k-1}}}\nonumber
\end{multline}
ceci par les m\^emes calculs que pr\'ec\'edemment.
\\

\noindent$\bullet ${\bf{ Bilan}}\\
Avec les normalisations des mesures de Haar choisies (${\rm{vol}}(\adele/F)=1$), on a 
$${\rm{vol }}(\entier\backslash F_\infty)=|\disc|^{1/2}$$
et par cons\'equent, pour deux id\'eaux fractionnaires $\A$ et $\B$:
$${\rm{vol}}(\A\backslash F_\infty)=N(\A)|\disc|^{1/2},{\rm{vol}}(\A\B^2\backslash F_\infty)=N(\A\B^2)|\disc|^{1/2}.$$
Pour $\varphi\in\mathcal{S}_\q^{\g{k}}$ et $g\in\GL(F_\infty)$, on \'ecrit comme pr\'ec\'edemment:
$$\varphi\left(g\left(\begin{array}{cc}{\rm{id}}(\A) & 0\\0 & 1 \end{array}
\right)\right)=\sum_{\substack{\alpha\in\A^{-1}\diff^{-1}\\ \alpha\gg 0}}c(\alpha,\A,\varphi)W_\infty^0\left(
\left(\begin{array}{cc}\g{\alpha} & 0\\0 & 1 \end{array}\right)g\right)$$
et posons pour all\'eger:
\begin{eqnarray}
\lambda(\alpha,\A,\varphi)=c(\alpha,\A,\varphi)N(\alpha\A\diff)^{1/2}.
\end{eqnarray}
Nous avons donc d\'emontr\'e:

\begin{proposition}[Formule de Petersson pour $\mathcal{S}_\q^{\g{k}}$]
Soit $\g{k}\in 2\N^d_{\geq 1}$ et $\{\phi\}$ une base orthonorm\'ee adapt\'ee de $\mathcal{S}_\q^{\g{k}}$. Soient $\A,\B$ des id\'eaux fractionnaires de $F$ et $\alpha\in\A^{-1}\diff^{-1}$ 
(respectivement $\alpha'\in(\A\B^2)^{-1}\diff^{-1}$). 
On a alors la relation:
\begin{multline}
\sum_{\phi}\frac{\Gamma(\g{k-1})}{(4\pi)^{\g{k-1}}|\disc|^{1/2}}\omega_\phi(\B)\lambda(\alpha,\A,\phi)
\overline{\lambda(\alpha',\A\B^2,\phi)}=\mathbbmss{1}_{\alpha\A=\alpha'\A\B^2}\\+
\frac{(-1)^{\g{k}/2}(2\pi)^d}{2|\disc|^{1/2}}\sum_{\substack{c\in\A^{-1}\B^{-1}\q\setminus \{0\}\\ \varepsilon\in\mathcal{O}_F^{\times +}/\mathcal{O}_F^{\times 2}}}
\frac{KS(\varepsilon\alpha,\A;\alpha',\A\B^2;c,\A\B)}{N(c\A\B)}J_{\g{k-1}}\left(4\pi\frac{\sqrt{
\g{\varepsilon\alpha\alpha'}}}{|\g{c}|}\right).
\end{multline}
\end{proposition}

\subsection{La formule de Petersson pour $\mathcal{H}_\q^{\g{k}}$}
On va pouvoir maintenant \'enoncer le r\'esultat principal de cette partie, \`a savoir la formule de Petersson 
pour l'espace $\mathcal{H}_\q^{\g{k}}$. Pour cela, soit $\{\varphi\}$ une base \emph{orthogonale} de 
$\mathcal{H}_\q^{\g{k}}$. Notons d'abord que :
\begin{multline}
||\varphi|| _{\mathcal{H}_\q^{\g{k}}}^2=\int_{Z(\adele)\gl(F)\backslash \gl(\adele)/K_0(\q)}|\varphi(g)|^2dg\\
=\frac{1}{h_F^+}\int_{Z_\infty^+\gl(F)\backslash \gl(\adele)/K_0(\q)}|\varphi(g)|^2dg=\frac{1}{h_F^+}||\varphi||_
{\mathcal{S}_\q^{\g{k}}}^2.
\end{multline}
Puis, $\{\varphi\}$ est une base orthogonale de $\mathcal{S}_\q^{\g{k}}[1]$ dans la d\'ecomposition (\ref{decomposition}). 
Soient maintenant $\n,\m$ deux id\'eaux fractionnaires de $F$, et $\alpha\in\n^{-1}\diff^{-1}$ (resp. $\beta\in\m^{-1}
\diff^{-1}$). On veut calculer:
\begin{multline}\sum_{\varphi}\frac{\Gamma(\g{k-1})}{(4\pi)^{\g{k-1}}|\disc|^{1/2}||\varphi||_{\mathcal{H}_\q^{\g{k}}}^2}\lambda(\alpha,\n,\varphi)
\overline{\lambda(\beta,\m,\varphi)}\\=h_F^+\sum_{\varphi}\frac{\Gamma(\g{k-1})}{(4\pi)^{\g{k-1}}|\disc |^{1/2}||\varphi||
_{\mathcal{S}_\q^{\g{k}}}^2}\lambda(\alpha,\n,\varphi)\overline{\lambda(\beta,\m,\varphi)}.\nonumber
\end{multline}
Or, en prenant $\{\phi\}$ une base orthonorm\'ee de $\mathcal{S}_\q^{\g{k}}$ comme pr\'ec\'edemment:
\begin{multline}\label{a}
\sum_{\varphi}\frac{\Gamma(\g{k-1})}{(4\pi)^{\g{k-1}}|\disc|^{1/2}||\varphi||
_{\mathcal{S}_\q^{\g{k}}}^2}\lambda(\alpha,\n,\varphi)\overline{\lambda(\beta,\m,\varphi)}\\=\frac{1}{h_F^+}
\sum_{\overline{\B}\in\mathscr{C}\!\ell^+(F)}\sum_{\phi}\frac{\Gamma(\g{k-1})}{(4\pi)^{\g{k-1}}|\disc|^{1/2}}
\lambda(\alpha,\n,\phi)\overline{\lambda(\beta,\m,\phi)}\omega_\phi(\B)\\
=\frac{1}{h_F^+}
\sum_{\overline{\B}^2=\overline{\m\n^{-1}}}\sum_{\phi}\frac{\Gamma(\g{k-1})}{(4\pi)^{\g{k-1}}|\disc|^{1/2}}
\lambda(\alpha,\n,\phi)\overline{\lambda(\beta,\m,\phi)}\omega_\phi(\B).
\end{multline}
En effet, soit $\widehat{\mathscr{C}\!\ell^+(F)}(2)$ le groupe des caract\`eres $\chi$ de $\mathscr{C}\!\ell^+(F)$ tels que 
$\chi(\A^2)=1$, quel que soit $\A$. Si $\{\phi\}$ est une base orthonorm\'ee de $\mathcal{S}_\q^{\g{k}}$, alors 
$\{\phi\otimes \chi\}$, $\chi$ fix\'e, en est une autre. Donc:
\begin{multline}
\sum_{\B}\sum_{\chi\in \widehat{\mathscr{C}\!\ell^+(F)}(2)}\sum_{\phi}\frac{\Gamma(\g{k-1})}{(4\pi)^{\g{k-1}}|\disc|^{1/2}}
\lambda(\alpha,\n,\phi\otimes\chi)\overline{\lambda(\beta,\m,\phi\otimes\chi)}\omega_{\phi\otimes\chi}(\B)\\
=\sum_{\B}\sum_{\phi}\frac{\Gamma(\g{k-1})}{(4\pi)^{\g{k-1}}|\disc|^{1/2}}
\lambda(\alpha,\n,\phi)\overline{\lambda(\beta,\m,\phi)}\omega_{\phi}(\B)\sum_{\chi\in \widehat{\mathscr{C}\!\ell^+(F)}
(2)}\chi(\n\m^{-1})\\=\sum_{\B}\sum_{\phi}\frac{\Gamma(\g{k-1})}{(4\pi)^{\g{k-1}}|\disc|^{1/2}}
\lambda(\alpha,\n,\phi)\overline{\lambda(\beta,\m,\phi)}\omega_{\phi}(\B)\mathbbmss{1}_{\m\n^{-1}\in\mathscr{C}\!\ell
^+(F)^2}\times|\widehat{\mathscr{C}\!\ell^+(F)}(2)|\\=|\widehat{\mathscr{C}\!\ell^+(F)}(2)| \sum_{\B}\sum_{\phi}\frac{\Gamma(\g{k-1})}{(4\pi)^{\g{k-1}}|\disc|^{1/2}}
\lambda(\alpha,\n,\phi)\overline{\lambda(\beta',\n\A^2,\phi)}\omega_{\phi}(\B) \\=|\widehat{\mathscr{C}\!\ell^+(F)}(2)| 
\sum_{\B}\sum_{\phi}\frac{\Gamma(\g{k-1})}{(4\pi)^{\g{k-1}}|\disc|^{1/2}}
\lambda(\alpha,\n,\phi)\overline{\lambda(\beta',\n\A^2,\phi)}\omega_{\phi}(\A)\omega_\phi(\B\A^{-1})\\
=|\widehat{\mathscr{C}\!\ell^+(F)}(2)| \sum_{\phi}\frac{\Gamma(\g{k-1})}{(4\pi)^{\g{k-1}}|\disc|^{1/2}}
\lambda(\alpha,\n,\phi)\overline{\lambda(\beta',\n\A^2,\phi)}\omega_{\phi}(\A)\sum_{\B}\omega_\phi(\B\A^{-1})\\
=\sum_{\overline{\B}^2=\overline{\m\n^{-1}}}\sum_{\chi\in \widehat{\mathscr{C}\!\ell^+(F)}(2)}\sum_{\phi}\frac{\Gamma(\g{k-1})}{(4\pi)^{\g{k-1}}|\disc|^{1/2}}\lambda(\alpha,\n,\phi\otimes\chi)\overline{\lambda(\beta,\m,\phi\otimes\chi)}\omega_{\phi\otimes\chi}(\B)\\\nonumber
\end{multline}
\begin{multline}
=\sum_{\overline{\B}^2=\overline{\m\n^{-1}}}|\widehat{\mathscr{C}\!\ell^+(F)}(2)|\sum_{\phi}\frac{\Gamma(\g{k-1})}{(4\pi)^{\g{k-1}}|\disc|^{1/2}}
\lambda(\alpha,\n,\phi)\overline{\lambda(\beta,\m,\phi)}\omega_{\phi}(\B)\nonumber
\end{multline}
(\ref{a}) \'etant maintenant justifi\'ee, on la r\'e\'ecrit, en notant $[\m\n^{-1}\B^{-2}]$ un g\'en\'erateur totalement 
positif de cet id\'eal:
\begin{multline}
\sum_{\varphi}\frac{\Gamma(\g{k-1})}{(4\pi)^{\g{k-1}}|\disc|^{1/2}||\varphi||
_{\mathcal{S}_\q^{\g{k}}}^2}\lambda(\alpha,\n,\varphi)\overline{\lambda(\beta,\m,\varphi)}\\
=\frac{1}{h_F^+}
\sum_{\overline{\B}^2=\overline{\m\n^{-1}}}\sum_{\phi}\frac{\Gamma(\g{k-1})}{(4\pi)^{\g{k-1}}|\disc|^{1/2}}
\lambda(\alpha,\n,\phi)\overline{\lambda(\beta[\m\n^{-1}\B^{-2}],\n\B^{-2},\phi)}\omega_\phi(\B)
\end{multline}
on applique maintenant la formule de Petersson pour $\mathcal{S}_\q^{\g{k}}$, et il vient apr\`es r\'eindexation (on 
remplace la somme sur $\B$ par une somme sur $\cc=\n\B$), qui donne:
\begin{multline}
\sum_{\varphi}\frac{\Gamma(\g{k-1})}{(4\pi)^{\g{k-1}}|\disc|^{1/2}||\varphi||
_{\mathcal{S}_\q^{\g{k}}}^2}\lambda(\alpha,\n,\varphi)\overline{\lambda(\beta,\m,\varphi)}=\mathbbmss{1}_{\alpha\n=\beta\m}\\
+\frac{(-1)^{\g{k}/2}(2\pi)^d}{2|\disc|^{1/2}h_F^+}\sum_{\substack{\overline{\cc}^2=\overline{\m\n}\\
c\in\cc^{-1}\q\setminus \{0\}\\ 
\varepsilon\in\mathcal{O}_F^{\times +}/\mathcal{O}_F^{\times 2}}}
\frac{KS(\varepsilon\alpha,\n;\beta[\m\n\cc^{-2}],\m;c,\cc)}{N(c\cc)}J_{\g{k-1}}\left(4\pi\frac{\sqrt{
\g{\varepsilon\alpha\alpha'}[\frac{\m\n}{\cc^{-2}}]}}{|\g{c}|}\right)\nonumber
\end{multline}
et, en prenant la norme de $\mathcal{H}_\q^{\g{k}}$, on a prouv\'e:

\begin{theoreme}[Formule de Petersson]\label{petersson}Soit $\q$ un id\'eal quelconque de $\entier$, $\g{k}\in 2\N^d_{\geq 1}$. 
Soit $\{\varphi\}$ une base orthogonale de $\mathcal{H}_\q^{\g{k}}$, $\m,\n$ des id\'eaux fractionnaires,
$\alpha\in\n^{-1}\diff^{-1}$ (respectivement $\beta\in\m^{-1}\diff^{-1}$). Soient $\{\lambda(\xi,\A,\varphi)\}$ 
les coefficients tels que pour tout $g\in\GL(F_\infty)$:
$$\varphi\left(g\left(\begin{array}{cc}{\rm{id}}(\A)&0\\0&1\end{array}\right)\right)=\sum_{\substack{\xi\in\A^{-1}\diff^{-1}\\
\xi\gg 0}}\frac{\lambda(\xi,\A,\varphi)}{N(\xi\A\diff)^{1/2}}W_\infty^0\left(\left(\begin{array}{cc}\xi&0\\0&1
\end{array}\right)g\right)$$
On a alors la relation:
\begin{multline}
\sum_{\varphi}\frac{\Gamma(\g{k-1})}{(4\pi)^{\g{k-1}}|\disc|^{1/2}||\varphi||
_{\mathcal{H}_\q^{\g{k}}}^2}\lambda(\alpha,\n,\varphi)\overline{\lambda(\beta,\m,\varphi)}=
\mathbbmss{1}_{\alpha\n=\beta\m}\\
+\frac{C}{|\disc|^{1/2}}\sum_{\substack{\overline{\cc}^2=\overline{\m\n}\\
c\in\cc^{-1}\q\setminus \{0\}\\ 
\varepsilon\in\mathcal{O}_F^{\times +}/\mathcal{O}_F^{\times 2}}}
\frac{KS(\varepsilon\alpha,\n;\beta[\m\n\cc^{-2}],\m;c,\cc)}{N(c\cc)}J_{\g{k-1}}\left(4\pi\frac{\sqrt{
\g{\varepsilon\alpha\alpha'}[\m\n\cc^{-2}]}}{|\g{c}|}\right)\nonumber
\end{multline}
avec $C=\frac{(-1)^{\g{k}/2}(2\pi)^d}{2}$.
\end{theoreme}
\rem Dans la formule de Petersson, les sommes sur $\cc$ et $\varepsilon$ sont des sommes finies, tenant compte du fait 
que le corps $F$ n'a pas forc\'ement un groupe des classes \'etroit trivial. Si cela est le cas, on retrouve alors le 
r\'esultat de \cite{luo}.\\

\rem On a travaill\'e par simplicit\'e avec des formes de caract\`ere central trivial. Il est clair cependant que l'on 
peut adapter les calculs pr\'ec\'edents \`a des formes dont le caract\`ere central est arithm\'etique, i.e. trivial 
sur les places infinies, comme l'a d'ailleurs fait \cite{V}. La formule est la m\^eme, \`a ceci pr\`es qu'il faut 
remplacer les sommes de Kloosterman \og $KS$\fg par \og $KS_{\chi}$\fg, sommes de Kloosterman twist\'ees.

\section{Lien avec les fonctions $L$}\label{fonctionsL}
Nous allons ici expliquer comment utiliser la formule de Petersson quand on manipule des fonctions $L$. Ce qui est dit ici 
est bien connu, mais il est plus prudent de l'expliciter, compte tenu des maintes normalisations de la litt\'erature. 
Soit $\pi$ une repr\'esentation parabolique de $\gl$, sur $F$. On a vu, dans la section \ref{Notations}, que $\pi\cong 
\bigotimes \pi_v$, et que dans chaque mod\`ele de Whittaker local $\mathcal{W}(\pi_v,\psi_v)$
, on peut trouver un unique \'el\'ement sp\'ecial $W_v^0$ tel 
que $$L(s,\pi_v)=\int_{F_v^\times}W_v^0\left(\left(\begin{array}{cc}a&0\\0&1\end{array}\right)\right)|a|^{s-1/2}
d^\times a$$
\`a condition que $\psi_v$ soit pris non ramifi\'e. Or, travaillant avec l'analyse harmonique, on est conduit \`a 
choisir comme caract\`ere additif global $\psi=\psi_\Q\circ {\rm Tr}_{F/\Q}$, qui est ramifi\'e en les places divisant la 
diff\'erente. Cependant, avec ce choix, prenons $W_v^0$ l'unique \'el\'ement sp\'ecial (i.e. invariant par $K_0(\q_v)$) 
tel que $$W_v^0\left(\left(\begin{array}{cc}\varpi_v^{-d_v}&0\\0&1\end{array}\right)\right)=1$$
o\`u $d_v$ est d\'etermin\'e comme suit: on \'ecrit $\diff=\prod_{\p}\p^{d_\p}$, et si $v$ est la place finie 
correspondant \`a $\p$, on pose $d_v=d_{\p}$. Cet \'el\'ement v\'erifie alors:
$$\int_{F_v^\times}W_v^0\left(\left(\begin{array}{cc}a&0\\0&1\end{array}\right)\right)|a|_v^{s-1/2}
d^\times a=|\varpi_v|^{d_v(1/2-s)}L(s,\pi_v)$$ et donc l'\'el\'ement $W_\pi$ de $\Whitt$, d\'efini par 
$W_\pi=W_\infty^0\otimes( \bigotimes W_v^0)=W_\infty^0\otimes W_\pi^0$ satisfait \`a :
$$\int_{\mathbbmss{A}_F^\times}W_\pi\left(\left(\begin{array}{cc}a&0\\0&1\end{array}\right)\right)|a|^{s-1/2}
d^\times a=|\disc|^{s-1/2}L(s,\pi)$$
et en notant $\mathbbmss{A}_f$ les ad\`eles finis:
$$\int_{\mathbbmss{A}_{f}^\times}W_\pi^0\left(\left(\begin{array}{cc}a&0\\0&1\end{array}\right)\right)|a|^{s-1/2}
d^\times a=|\disc|^{s-1/2}L(s,\pi_f)$$
en d\'eveloppant l'int\'egrale de gauche, et en \'ecrivant
$$L(s,\pi_f)=\sum_{\n\subset \entier}\lambda_\pi(\n)N(\n)^{-s}$$
on obtient la relation:
\begin{eqnarray}
\lambda_\pi(\n)=N(\n)^{1/2}W_\pi^0\left(\left(\begin{array}{cc}{\rm{id}}(\n\diff^{-1})&0\\0&1\end{array}\right)\right)
\end{eqnarray}
et par cons\'equent, si $\varphi_\pi$ est l'\'el\'ement sp\'ecial global de l'espace de $\pi$ correspondant \`a $W_\pi$, 
et $\alpha\in\n^{-1}$:
\begin{eqnarray}\label{coeff}
\lambda_\pi(\alpha\n)=\lambda(\alpha,\n\diff^{-1},\varphi_\pi)
\end{eqnarray}

On a vu en \ref{Automorphe} le d\'ecomposition de $\mathcal{H}_\q^{\g{k}}$ \`a l'aide des $\Pi_\q^{\g{k}}$. 
Dans le cas o\`u $\q$ est un id\'eal 
maximal (cas dans lequel on se placera dans les sections ult\'erieures), on a donc la d\'ecomposition plus simple:
$$\mathcal{H}_\q^{\g{k}}=\bigoplus_{\pi\in\Pi_\q^{\g{k}}}\C\varphi_\pi\oplus\bigoplus_{\pi\in\Pi_{\entier}^{\g{k}}}(\pi,V_\pi)^{K_\infty K_0(\q)}$$
 L'espace des formes anciennes peut \^etre nul (
si $F=\Q$, et $k\leq 10$ par exemple), mais si $\pi$ est non ramifi\'ee, alors 
l'espace $(\pi,V_\pi)^{K_\infty K_0(\q)}$ est de dimension deux.\\

Dans ce qui va suivre, la contribution des formes anciennes est n\'egligeable, et on a pr\'ef\'er\'e indiquer dans un 
appendice pourquoi (il s'agit de trouver une base orthogonale \og adapt\'ee\fg des formes anciennes, et d'\'evaluer 
leur contribution): la raison principale vient du fait que si $\varphi$ est non ramifi\'ee, alors $||\varphi||_{\mathcal{H}
_\q^{\g{k}}}^2=(N(\q)+1)||\varphi||_{\mathcal{H}_{\entier}^{\g{k}}}^2$, et cela fait baisser d'un facteur $N(\q)$ cette 
contribution.\\

Venons-en maintenant \`a la formule de Petersson. 
D'abord, on normalise les sommes de Kloosterman:
\begin{definition}Soient $\m,\n$ des id\'eaux fractionnaires de $F$, $\cc$ tel que $\cc^2\sim \m\n$, et 
$\alpha\in\n^{-1}$, $\beta\in\m^{-1}$. On pose:
\begin{eqnarray}
\kl(\alpha,\n;\beta,\m;c,\cc)=KS(\alpha,\n\diff^{-1};\beta[\m\n\cc^{-2}],\m\diff^{-1};c,\cc)
\end{eqnarray}
On a alors la borne de Weil:
\begin{eqnarray}
\left|\kl(\alpha,\n;\beta,\m;c,\cc)\right|\ll_F N\big((\alpha)\n,(\beta)\m,(c)\cc\big)^{1/2}\tau\big((c)\cc\big)N(c\cc)^{1/2}
\end{eqnarray}
\end{definition}
Dans la d\'efinition, on a not\'e $(\A,\B,\cc)$ le p.g.c.d. des id\'eaux $\A,\B,\cc$. La borne de Weil est \'enonc\'ee dans \cite{V}, section 2.6.\\

On peut \'enoncer la forme \og pratique\fg de la formule de Petersson:
\begin{definition}\label{poids}
Soit $\{x_\pi\}_{\pi\in\Pi_\q^{\g{k}}}$ une suite de nombres complexes.
Soit $\omega_\pi$ le poids 
\begin{eqnarray}
\omega_\pi=\frac{\Gamma(\g{k-1})}{(4\pi)^{\g{k-1}}|\disc|^{1/2}||\varphi_\pi||_{\mathcal{H}_\q^{\g{k}}}^2}
\end{eqnarray}
On note alors:
\begin{eqnarray}
\sum_{\pi\in\Pi_\q^{\g{k}}}^hx_\pi:=\sum_{\pi\in\Pi_\q^{\g{k}}}\omega_\pi x_\pi
\end{eqnarray}
\end{definition}

\begin{proposition}\label{peterssonbis} Soit $\q$ un id\'eal quelconque de $\entier$, $\g{k}\in 2\N^d_{\geq 1}$. 
Soient  $\m,\n$ des id\'eaux fractionnaires,
$\alpha\in\n^{-1}$ (respectivement $\beta\in\m^{-1}$).
\begin{multline}
\sum_{ \pi\in\Pi_\q^{\g{k}}}^h\lambda_\pi(\alpha\n)\lambda_\pi(\beta\m)+\textrm{ (Formes Anciennes) }
=\mathbbmss{1}_{\alpha\n=\beta\m}\\
+\frac{C}{|\disc|^{1/2}} \sum_{\substack{\overline{\cc}^2=\overline{\m\n}\\
c\in\cc^{-1}\q\setminus \{0\}\\ 
\varepsilon\in\mathcal{O}_F^{\times +}/\mathcal{O}_F^{\times 2}}}
\frac{\kl(\varepsilon\alpha,\n;\beta,\m;c,\cc)}{N(c\cc)}J_{\g{k-1}}\left(4\pi\frac{\sqrt{
\g{\varepsilon\alpha\beta}[\m\n\cc^{-2}]}}{|\g{c}|}\right)
\end{multline}
avec 
\begin{eqnarray}
\textrm{ (Formes Anciennes)}=\sum_{\varphi}\frac{\Gamma(\g{k-1})}{(4\pi)^{\g{k-1}}|\disc|^{1/2} ||\varphi||
_{\mathcal{H}_\q^{\g{k}}}^2}\lambda(\alpha,\n\diff^{-1},\varphi)\overline{\lambda(\beta,\m\diff^{-1},\varphi)}\nonumber
\end{eqnarray}
$\{\varphi\}$ parcourant une base orthogonale des formes anciennes, et
\begin{eqnarray}
C=\frac{(-1)^{\g{k}/2}(2\pi)^d}{2}.
\end{eqnarray}
\end{proposition}
\rem Les coefficients $\lambda_\pi(\n)$ sont r\'eels, dans le cas \'etudi\'e: c'est pour cela que la conjugaison complexe 
n'appara\^it pas dans la partie des formes nouvelles.\\

La formule de Petersson, ainsi qu'une \'etude des termes provenant des formes anciennes, permet de voir que:
$$\lim_{N(\q)\to\infty}\sum_{\pi\in\Pi_\q^{\g{k}}}^h1=1.$$
Pour une suite $\{x_\varphi\}$, index\'ee par une base orthogonale de $\mathcal{H}_\q^{\g{k}}$, on posera 
$$\sum_\varphi^hx_\varphi:=\sum_\varphi\frac{\Gamma(\g{k-1})}{(4\pi)^{\g{k-1}}|\disc|^{1/2} ||\varphi||_
{\mathcal{H}_\q^{\g{k}}}^2}x_\varphi.$$

La formule de Petersson est valable pour un conducteur quelconque. Cependant, on ne saura \'evaluer le terme venant des formes anciennes que dans le cas d'un id\'eal premier: cette difficult\'e est pr\'evisible car la d\'ecomposition (\ref{newold}) d\'epend de la complexit\'e arithm\'etique du conducteur.

\section{Le premier moment}\label{Moment1}
Dans cette partie, nous allons prouver l'estimation (\ref{moment1}), rappelons la:
\begin{multline}
M_1(\q):=\sum_{\pi\in\Pi_\q^{\g{k}}}^h \Lambda(1/2,\pi) \M(\pi)\\ =\frac{\zeta_F(2) \Gamma\left(\frac{\g{k}}{2}\right)}
{(2\pi)^{\frac{\g{k}}{2}} \res_{s=1}(\zeta_F)}\times \frac{2N(\q)^{1/4}}{\Delta\log(N(\q))}
\Bigg(P'(1)+\mathcal{O}\bigg(\frac{1}{\log(N(\q))}\bigg)\Bigg).\nonumber
\end{multline}
 Rappelons d'abord que l'on prend pour amollisseur:
$$\M(\pi)=\sum_{N(\m)\leq M}\frac{\mu(\m)P\left(\frac{\log(M/N(\m))}{\log(M)}\right)}{\psi(\m)N(\m)^{1/2}}\lambda_\pi(\m)$$
o\`u $M=N(\q)^{\Delta/2}$, $\Delta\in]0,1[$, $P$ est un polyn\^ome tel que $P(0)=P'(0)=0$. 

La repr\'esentation $\pi$ parcourt l'ensemble des formes modulaires de Hilbert de poids $\g{k}$, de 
conducteur l'id\'eal premier $\q$, et de caract\`ere central trivial. Il s'agit donc d'obtenir une asymptotique quand 
$N(\q)$ tend vers l'infini, avec $\Delta$ fixe.

Nous supposerons \'egalement dans cette partie et la suivante que $M$ n'est pas un entier: cette hypoth\`ese, non restrictive, 
permet d'utiliser une \'ecriture int\'egrale de $P$. En effet, Kowalski, Michel et VanderKam \cite{KMV} associent \`a
$$P(X)=\sum_{k}a_kX^k$$
la fraction rationnelle
$$\PM(X)=\sum_{k}a_kk!(X\log(M))^{-k}$$
de telle sorte que pour $M$ non entier, et $P$ s'annulant en $0$:
\begin{eqnarray}
P\left(\frac{\log(M/N(\m))}{\log(M)}\right)\mathbbmss{1}_{N(\m)< M}=\frac{1}{2i\pi}\int_{(1)}\frac{M^s}{N(\m)^s}\PM(s)
\frac{ds}{s}.
\end{eqnarray}

Rappelons (voir (\ref{valeur1})) que $\Lambda(1/2,\pi)$ est donn\'e par:
$$\Lambda(1/2,\pi)=(1+\varepsilon_\pi)N(\q)^{1/4}\sum_{\n\subset \entier}
\frac{\lambda_\pi(\n)}{\sqrt{N(\n)}}F(N(\n)/N(\q)^{1/2})$$
avec:
$$F(y)=\frac{1}{2i\pi}\int_{(3/2)}y^{-s} L\left(s+\frac{1}{2},\pi_\infty\right) \frac{ds}{s}$$
et la valeur explicite du signe de l'\'equation fonctionnelle:
\begin{eqnarray}
\varepsilon_\pi=-i^{\g{k}}N(\q)^{1/2}\lambda_\pi(\q).
\end{eqnarray}
Afin de gagner en lisibilit\'e, nous userons de la notation condens\'ee suivante:
\begin{eqnarray}
F_\n:=\frac{F(N(\n)/N(\q)^{1/2})}{N(\n)^{1/2}}\\
P_\m:=\frac{\mu(\m)P\left(\frac{
\log(M/N(\m))}{\log(M)}\right)}{\psi(\m)N(\m)^{1/2}}.
\end{eqnarray}

On consid\`ere donc $$M_1(\q)=\sum_{\pi\in\Pi_\q^{\g{k}}}^h\Lambda(1/2,\pi)M(\pi).$$
En d\'eveloppant l'expression de $\M(\pi)$, on voit que $M_1(\q)$ est donn\'ee par:
\begin{multline}
M_1(\q)=N(\q)^{1/4}\sum_{\n\subset \entier}\sum_{ N(\m)\leq M}F_\n P_\m
\sum_{\pi\in\Pi_\q^{\g{k}}}^h\Big(1-N(\q)^{1/2}\lambda_\pi(\q)\Big)\lambda_\pi(\n)\lambda_\pi(\m)\\
=N(\q)^{1/4}\sum_{\substack{\n\subset \entier\\ N(\m)\leq M}}F_\n P_\m
\sum_{\pi\in\Pi_\q^{\g{k}}}^h\lambda_\pi(\n)\lambda_\pi(\m)\\-N(\q)^{3/4}\sum_{\substack{\n\subset \entier\\ N(\m)\leq M}}F_\n P_\m
\sum_{\pi\in\Pi_\q^{\g{k}}}^h\lambda_\pi(\n\q)\lambda_\pi(\m)=A+B.\nonumber
\end{multline}
\\
Pour utiliser la formule de Petersson, on choisit un syst\`eme de repr\'esentants de $\mathscr{C}\!\ell^+(F)$, comme 
\`a la section \ref{Notations}, not\'e $\{\bar{\A}\}$ (resp. $\{\bar{\B}\}$). On effectue donc le changement de variables $\n=\nu\A,\m=\xi\B$, et $\nu$ (resp. $\xi$) parcourt $(\A^{-1})^{\gg 0}\mod \unite$ (resp. $(\B^{-1})^{\gg 0}\mod \unite$):
\begin{multline}
A=N(\q)^{1/4}\sum_{\bar{\A},\bar{\B}}\sum_{\nu\in(\A^{-1})^{\gg 0}/\unite}\sum_{\substack{ \xi\in(\B^{-1})^{\gg 0}/\unite\\ N(\xi)\leq MN(\B)^{-1}}}
F_{\nu\A} P_{\xi\B}
\sum_{\pi\in\Pi_\q^{\g{k}}}^h\lambda_\pi(\nu\A)\lambda_\pi(\xi\B)
\end{multline}

\begin{multline}
B=-N(\q)^{3/4}\sum_{\bar{\A},\bar{\B}}\sum_{\nu\in(\A^{-1})^{\gg 0}/\unite}\sum_{\substack{\xi\in(\B^{-1})^{\gg 0}/\unite\\N(\xi)\leq MN(\B)^{-1}}}
F_{\nu\A}P_{\xi\B}
\sum_{\pi\in\Pi_\q^{\g{k}}}^h\lambda_\pi(\nu\A\q)\lambda_\pi(\xi\B).
\end{multline}

On peut maintenant appliquer la formule de Petersson (proposition \ref{peterssonbis}):
\begin{multline}
A=N(\q)^{1/4}\sum_{\bar{\A},\bar{\B}}\sum_{\nu\in(\A^{-1})^{\gg 0}/\unite}\sum_{\substack{ \xi\in(\B^{-1})^{\gg 0}/\unite\\ N(\xi)\leq MN(\B)^{-1}}}
F_{\nu\A} P_{\xi\B}\times\\ \Bigg\{
\mathbbmss{1}_{\nu\A=\xi\B}
+\frac{C}{|\disc|^{1/2}} \sum_{\substack{\overline{\cc}^2=\overline{\A\B}\\
c\in\cc^{-1}\q\setminus \{0\}\\ 
\varepsilon\in\mathcal{O}_F^{\times +}/\mathcal{O}_F^{\times 2}}}
\frac{\kl(\varepsilon\nu,\A;\xi,\B;c,\cc)}{N(c\cc)}J_{\g{k-1}}\left(4\pi\frac{\sqrt{
\g{\varepsilon\nu\xi}[\A\B\cc^{-2}]}}{|\g{c}|}\right)\\+ (\textrm{Formes Anciennes}) \Bigg\}
 \end{multline}

\begin{multline}
B=-N(\q)^{3/4}\sum_{\bar{\A},\bar{\B}}\sum_{\nu\in(\A^{-1})^{\gg 0}/\unite}\sum_{\substack{\xi\in(\B^{-1})^{\gg 0}/\unite\\N(\xi)\leq MN(\B)^{-1}}}
F_{\nu\A}P_{\xi\B}\times\\\Bigg\{
\mathbbmss{1}_{\nu\A\q=\xi\B}
+\frac{C}{|\disc|^{1/2}} \sum_{\substack{\overline{\cc}^2=\overline{\A\B}\\
c\in\cc^{-1}\q\setminus \{0\}\\ 
\varepsilon\in\mathcal{O}_F^{\times +}/\mathcal{O}_F^{\times 2}}}
\frac{\kl(\varepsilon\nu,\A\q;\xi,\B;c,\cc)}{N(c\cc)}J_{\g{k-1}}\left(4\pi\frac{\sqrt{
\g{\varepsilon\nu\xi}[\A\q\B\cc^{-2}]}}{|\g{c}|}\right)\\+(\textrm{Formes Anciennes})\Bigg\}.
\end{multline}
On revient \`a $M_1(\q)=A+B$ en regroupant les termes diagonaux, ceux en sommes de Kloosterman, et ceux venant des formes anciennes. On \'ecrit donc:
$$M_1(\q)=M_1^{diag}(\q)+M_1^{Kloost}(\q)+M_1(\q,\textrm{Formes Anciennes}).$$

La d\'emarche sera donc la suivante: on \'etudie le terme diagonal, et en donne une asymptotique qui est celle de (\ref{moment1}); on montre que le terme venant des sommes de Kloosterman est n\'egligeable. Enfin, il faut estimer $M_1(\q,\textrm{Formes Anciennes})$: ce sera expliqu\'e dans l'appendice.\\

\noindent $\bullet $ {\bf{ Le terme diagonal}}\\
En appliquant la formule de Petersson \`a $M_1(\q)$, les termes diagonaux ont la forme: $\mathbbmss{1}_{\nu\A=\xi\B}$ 
et $\mathbbmss{1}_{\nu\A\q=\xi\B}$. Le dernier est nul, car $N(\xi\B)\leq M < N(\q)^{1/2}$. Il n'y a donc qu'une 
diagonale, et l'on peut r\'eindexer avec des sommes d'id\'eaux, en posant $\m=\nu\A=\xi\B$:
\begin{multline}
M_1^{diag}(\q)=N(\q)^{1/4}\sum_{ N(\m)\leq M}\frac{\mu(\m)F(N(\m)/N(\q)^{1/2})P\left(\frac{\log(M/N(\m))}{\log(M)}\right)}
{\psi(\m)N(\m)}=\\
\frac{N(\q)^{1/4}}{(2i\pi)^2}\iint_{(3),(1)}\sum_{\m\subset \entier}\frac{\mu(\m)}{\psi(\m)N(\m)^{1+s+t}}
N(\q)^{s/2}M^t L\left(s+\frac{1}{2},\pi_\infty\right) \PM(t)\frac{ds}{s}\frac{dt}{t}\nonumber
\end{multline}
en utilisant les expressions int\'egrales de $P$ et $F$ donn\'ees pr\'ec\'edemment.\\
Il faut noter que la s\'erie 
$t\mapsto \sum_{\m}\frac{\mu(\m)}{\psi(\m)N(\m)^{1+t+s}} $ se prolonge analytiquement au voisinage de $\Re (t+s)=0$. En effet, pour $\Re(t+s)>0$, on a:
\begin{multline}
\sum_{\m\subset \entier}\frac{\mu(\m)}{\psi(\m)N(\m)^{1+t+s}} =\prod_{\p}\left(1-\frac{N(\p)^{-(1+t+s)}}{1+N(\p)^{-1}}\right)
=\zeta_F(1+t+s)^{-1}\times\\
\prod_{\p}\left(\frac{1+N(\p)^{-1}-N(\p)^{-(1+t+s)}}{(1+N(\p)^{-1})(1-N(\p)^{-(1+t+s)})}\right)\nonumber
\end{multline}
et un d\'eveloppement limit\'e assure que le produit de droite converge pour $\Re(s+t)> -1$. Notons $D_1(s,t)$ ce prolongement holomorphe.\\

Pour aller plus loin, nous devons effectuer des changements de droite d'int\'egration: Kowalski, Michel et VamderKam ont exploit\'e cette technique avec beaucoup d'efficacit\'e dans \cite{KMV}. Le lecteur est invit\'e en cas de besoin \`a se r\'ef\'erer \`a la section \ref{fonctions} o\`u en est rappel\'e le principe.\\
On change d'abord de droite d'int\'egration en $s$, ce qui donne:
\begin{multline}
\frac{N(\q)^{1/4}}{2i\pi}\int_{(1)}{\textrm{res}}_{s=0}\bigg(
N(\q)^{s/2}L\left(s+\frac{1}{2},\pi_\infty\right)D_1(s,t) s^{-1}
\bigg) M^t\PM(t)\frac{dt}{t}\\+
\frac{N(\q)^{1/4}}{(2i\pi)^2}\iint_{(-1),(1)}D_1(s,t)
N(\q)^{s/2}M^tL\left(s+\frac{1}{2},\pi_\infty\right)  \PM(t)\frac{ds}{s}\frac{dt}{t}.\nonumber
\end{multline}
Le second terme est en $\mathcal{O}\left(N(\q)^{1/4-\delta}\right)$ pour un $\delta>0$, puisque $|N(\q)^{s/2}M^t|=N(\q)^{\frac{-1+\Delta}{2}}$. On continue avec le premier:

\begin{multline}
M_1^{diag}(\q)=\frac{N(\q)^{1/4}L\left(\frac{1}{2},\pi_\infty\right)  }{2i\pi} 
\int_{(1)}D_1(0,t) M^t\PM(t)\frac{dt}{t}+\mathcal{O}(N(\q)^{1/4-\delta}).
\nonumber
\end{multline}
De la m\^eme mani\`ere on change de ligne d'int\'egration en $t$. On \'ecrit donc:
\begin{multline}
M_1^{diag}(\q)=N(\q)^{1/4} L\left(\frac{1}{2},\pi_\infty\right)   
{\textrm{res}}_{t=0}\bigg(D_1(0,t) M^t\PM(t)t^{-1}\bigg)\\
+\frac{N(\q)^{1/4}L\left(\frac{1}{2},\pi_\infty\right) }{2i\pi} 
\int_{(-1/4)}D_1(0,t) M^t\PM(t)\frac{dt}{t}+\mathcal{O}(N(\q)^{1/4-\delta}).
\nonumber
\end{multline}
L\`a encore, le terme int\'egral est un $\mathcal{O}(N(\q)^{1/4-\delta})$; il faut \'evaluer le r\'esidu. Pour cela, 
on d\'eveloppe en s\'erie enti\`ere $M^t$, $\PM(t)$, et on remarque que :
$$\sum_{\m\subset \entier}\frac{\mu(\m)}{\psi(\m)N(\m)^{1+t}} \Equi{t}{0} \zeta_F(1+t)^{-1}
\prod_{\p}\left(\frac{1}{1+N(\p)^{-2}}\right)\Equi{t}{0} \frac{\zeta_F(2)t}{{\textrm{res}}_{t=1}(\zeta_F)}.$$
Ceci montre que l'on peut \'ecrire:
$$D_1(0,t)=\sum_{\ell\geq 1}\alpha_\ell t^\ell.$$
La contribution d'un facteur $\ell\geq 1$ au r\'esidu est de la forme:
$$\sum_{k-j+\ell=0}\frac{\log(M)^k}{k!}\frac{a_jj!}{\log(M)^j}\alpha_\ell=\mathcal{O}(\log(M)^{-\ell})$$
ainsi, seul le premier terme ($\ell=1$) contribue. On fait le calcul, et on trouve:
\begin{multline}
M_1^{diag}(\q)=\frac{N(\q)^{1/4}L\left(\frac{1}{2},\pi_\infty\right)   \zeta_F(2) }{{\textrm{res}}_1(\zeta_F)}
\sum_{k-j=-1}\frac{\log(M)^k}{k!}\frac{a_jj!}{\log(M)^j} \\+\mathcal{O}(N(\q)^{1/4}/\log(M)^2)+\mathcal{O}(N(\q)^{1/4-\delta})\\
=\frac{N(\q)^{1/4}L\left(\frac{1}{2},\pi_\infty\right)   \zeta_F(2) }{{\textrm{res}}_1(\zeta_F)\log(M)}
\bigg(P'(1)+\mathcal{O}\Big(1/\log(N(\q))\Big)\bigg).
\end{multline}

\noindent$\bullet$ {\textbf{Le terme en sommes de Kloosterman}}\\
Nous allons montrer que la contribution des termes venant des sommes de Kloosterman est n\'egligeable devant le terme 
diagonal. 
On posera donc:

\begin{multline}
M_{1,\A,\B}^{Kloost}(\q)=N(\q)^{1/4}\times\\\sum_{\nu\in(\A^{-1})^{\gg 0}/\unite} \sum_{\substack{\xi\in(\B^{-1})^{\gg 0}/\unite\\N(\xi)\leq MN(\B)^{-1}}}
\frac{F(N(\nu\A)/N(\q)^{1/2})}{N(\nu\A)^{1/2}}\cdot 
\frac{\mu(\xi\B)P\left(\frac{
\log(M/N(\xi\B))}{\log(M)}\right)}{\psi(\xi\B)N(\xi\B)^{1/2}}\times\\
\bigg\{ \frac{C}{|\disc|^{1/2}}\sum_{\substack{\overline{\cc}^2=\overline{\A\B}\\
c\in\cc^{-1}\q\setminus \{0\}\\ 
\varepsilon\in\mathcal{O}_F^{\times +}/\mathcal{O}_F^{\times 2}}}
\frac{\kl(\varepsilon\nu,\A;\xi,\B;c,\cc)}{N(c\cc)}J_{\g{k-1}}\left(4\pi\frac{\sqrt{
\g{\varepsilon\nu\xi}[\A\B\cc^{-2}]}}{|\g{c}|}\right)\\+N(\q)^{1/2} \frac{C}{|\disc|^{1/2}}\sum_{\substack{\overline{\cc}^2=\overline{\A\q\B}\\
c\in\cc^{-1}\q\setminus \{0\}\\ 
\varepsilon\in\mathcal{O}_F^{\times +}/\mathcal{O}_F^{\times 2}}}
\frac{\kl(\varepsilon\nu,\A\q;\xi,\B;c,\cc)}{N(c\cc)}J_{\g{k-1}}\left(4\pi\frac{\sqrt{
\g{\varepsilon\nu\xi}[\A\B\cc^{-2}]}}{|\g{c}|}\right)\bigg\}\nonumber
\end{multline}
et on va montrer que $M_{1,\A,\B}^{Kloost}(\q)\ll N(\q)^{1/4-\eta}$, avec $\eta>0$. Ceci montrera bien que 
$$M_1^{Kloost}(\q)=\sum_{\A,\B}M_{1,\A,\B}^{Kloost}(\q)\ll N(\q)^{1/4-\eta}.$$ 
De plus, les quelques sommes d'ensembles d'indices finis ne contribuant pas, on fixe un repr\'esentant $\cc$ tel que $\cc^2\sim \A\B$ et on se ram\`ene \`a l'\'etude de :

\begin{multline}
Err(\q)=N(\q)^{1/4}\times\\ \sum_{\nu\in(\A^{-1})^{\gg 0}/\unite}\sum_{\substack{\xi\in(\B^{-1})^{\gg 0}/\unite\\N(\xi)\leq MN(\B)^{-1}}}
\frac{F(N(\nu\A)/N(\q)^{1/2})}{N(\nu\A)^{1/2}}\cdot 
\frac{\mu(\xi\B)P\left(\frac{
\log(M/N(\xi\B))}{\log(M)}\right)}{\psi(\xi\B)N(\xi\B)^{1/2}}\times\\
\bigg\{ \frac{C}{|\disc|^{1/2}}\sum_{c\in\cc^{-1}\q\setminus \{0\}}
\frac{\kl(\nu,\A;\xi,\B;c,\cc)}{N(c)}J_{\g{k-1}}\left(4\pi\frac{\sqrt{
\g{\nu\xi}[\A\B\cc^{-2}]}}{|\g{c}|}\right)\\+N(\q)^{1/2} \frac{C}{|\disc|^{1/2}}\sum_{c\in\cc^{-1}\q\setminus \{0\}}
\frac{\kl(\varepsilon\nu,\A\q;\xi,\B;c,\cc)}{N(c)}J_{\g{k-1}}\left(4\pi\frac{\sqrt{
\g{\nu\xi}[\A\B\cc^{-2}]}}{|\g{c}|}\right)\bigg\}.\nonumber
\end{multline}

Pour majorer ce terme d'erreur, on \'ecrit d'abord:
\begin{multline}
Err(\q)=N(\q)^{1/4}\times\\ \sum_{\nu\in(\A^{-1})^{\gg 0}/\unite}\sum_ {\substack{\xi\in(\B^{-1})^{\gg 0}/\unite\\N(\xi)\leq MN(\B)^{-1}}}
\frac{F(N(\nu\A)/N(\q)^{1/2})}{N(\nu\A)^{1/2}}\cdot 
\frac{\mu(\xi\B)P\left(\frac{
\log(M/N(\xi\B))}{\log(M)}\right)}{\psi(\xi\B)N(\xi\B)^{1/2}}\times\\
\bigg\{ \frac{C}{|\disc|^{1/2}}\sum_{c\in\cc^{-1}\q\setminus \{0\}/\unite}\sum_{\eta\in\unite}
\frac{\kl(\nu,\A;\xi,\B;c\eta^{-1},\cc)}{N(c)}J_{\g{k-1}}\left(4\pi\frac{\sqrt{
\g{\nu\xi}[\A\B\cc^{-2}]}\g{\eta}}{|\g{c}|}\right)\\+  \frac{CN(\q)^{1/2}}{|\disc|^{1/2}}\sum_{\substack{c\in\cc^{-1}\q\setminus \{0\}/\unite\\\eta\in\unite}}
\frac{\kl(\varepsilon\nu,\A\q;\xi,\B;c\eta,\cc)}{N(c)}J_{\g{k-1}}\left(4\pi\frac{\sqrt{
\g{\nu\xi}[\A\B\cc^{-2}]}\g{\eta}}{|\g{c}|}\right)\bigg\}.\nonumber
\end{multline}

On choisit les repr\'esentants $c\in\q/\unite$, $\nu\in\A^{-1}/\unite$ et $\xi\in\B^{-1}/\unite$ satisfaisant au lemme \ref{sommation} de la section \ref{geometrie}. Concr\`etement, on peut supposer que $N(c)^{1/d}\ll c^{(j)}\ll N(c)^{1/d}$ pour tout $j$, ainsi que pour $\nu$ et $\xi$.\\
On note aussi que l'on a l'estimation pour tout $\varepsilon_j\geq 0$:
\begin{eqnarray}
J_{k_j-1}(z)\ll_{\delta}|z|^{1-\varepsilon_j}
\end{eqnarray}
et l'on ins\`ere cette repr\'esentation dans l'expression de $Err(\q)$, en choisissant pour chaque $\eta
\in\unite$: $\varepsilon_j=0$ si $|\eta^{(j)}|\leq 1$, et $\varepsilon_j=\delta$ ($\delta>0$) si $|\eta^{(j)}|> 1$. 
Ce choix est fait, car on va utiliser le fait suivant, prouv\'e dans \cite{luo} page 136:
\begin{eqnarray}\label{unites}
\sum_{\eta\in\unite}\prod_{|\eta^{(j)}|>1}|\eta^{(j)}|^{-\delta}<\infty.
\end{eqnarray}
En utilisant la borne de Weil pour les sommes de Kloosterman, il vient donc:

\begin{multline}
Err(\q)\ll_{\delta}N(\q)^{1/4}\times\\\sum_{\nu\in(\A^{-1})^{\gg 0}/\unite}\sum_{\substack{\xi\in(\B^{-1})^{\gg 0}/\unite\\N(\xi)\leq MN(\B)^{-1}}}
\frac{|F(N(\nu\A)/N(\q)^{1/2})|}{N(\nu)^{\delta/2}}\cdot 
\left|\frac{\mu(\xi\B)P\left(\frac{\log(M/N(\xi\B))}{\log(M)}\right)}{\psi(\xi\B)N(\xi)^{\delta/2}}\right|\times\\
\bigg(\sum_{c\in\q\setminus\{0\}/\unite}\frac{N(\nu\A,\xi\B,c\cc)^{1/2}}{N(c)^{3/2-\delta}}+N(\q)^{1/2}
\sum_{c\in\q\setminus\{0\}/\unite}\frac{N(\nu\A\q,\xi\B,c\cc)^{1/2}}{N(c)^{3/2-\delta}}\bigg)\nonumber
\end{multline}

\begin{multline}
\ll_{\delta}N(\q)^{1/4}\sum_ {\substack{\nu,\xi\\N(\xi)\leq MN(\B)^{-1}}}
\frac{|F(N(\nu\A)/N(\q)^{1/2})|}{N(\nu)^{\delta/2}}\cdot 
\left|\frac{\mu(\xi\B)P\left(\frac{\log(M/N(\xi\B))}{\log(M)}\right)}{\psi(\xi\B)N(\xi)^{\delta/2}}\right|\times\\
\bigg(\sum_{\cc\subset\q}\frac{N(\nu\A,\xi\B,\cc)^{1/2}}{N(\cc)^{3/2-\delta}}+N(\q)^{1/2}
\sum_{\cc\subset\q}\frac{N(\nu\A\q,\xi\B,\cc)^{1/2}}{N(\cc)^{3/2-\delta}}\bigg).
\nonumber
\end{multline}
Or (cf. Kowalski, th\`ese, lemme 9):

\begin{multline}
\sum_{\cc\subset\q}\frac{N(\n,\m,\cc)^{1/2}}{N(\cc)^{3/2-\delta}}\leq \frac{N(\m,\n,\q)^{1/2}}{N(\q)^{3/2-\delta}}\sum_{\A\subset\entier}
\frac{N(\m,\n,\A)^{1/2}}{N(\A)^{3/2-\delta}}\\ \ll \frac{N(\m,\n,\q)^{1/2}}{N(\q)^{3/2-\delta}}\sum_{\mathfrak{d}|(\n,\m)}
\frac{N(\mathfrak{d})^{1/2}}{N(\mathfrak{d})^{3/2-\delta}}
\ll_{\delta} \frac{N(\m,\n,\q)^{1/2}}{N(\q)^{3/2-\delta}}N(\n,\m)^{\delta}.
\end{multline}

Cette derni\`ere estimation donne:
\begin{multline}
Err(\q) 
\ll_{\delta}N(\q)^{1/4}\sum_{\nu\in\A^{-1}/\unite}\sum_{\substack{\xi\in\B^{-1}/\unite\\N(\xi)\leq MN(\B)^{-1}}}
\frac{|F(N(\nu\A)/N(\q)^{1/2})|}{N(\nu)^{\delta/2}}\times\\ 
\left|\frac{\mu(\xi\B)P\left(\frac{\log(M/N(\xi\B))}{\log(M)}\right)}{\psi(\xi\B)N(\xi)^{\delta/2}}\right|
N(\nu\A,\mu\B)^{\delta}\Big(\frac{N(\nu\A,\xi\B,\q)^{1/2}}{N(\q)^{3/2-\delta}}+\frac{N(\nu\A\q,\xi\B,\q)^{1/2}}{N(\q)^{1-\delta}}\Big).
\nonumber
\end{multline}
On a par hypoth\`ese $N(\xi\B)\leq M$ donc $$(\nu\q,\xi,\q)=(\nu,\xi,\q)=\entier.$$
Un changement de droite d'int\'egration dans l'\'ecriture de $F$ montre que:
\begin{displaymath}
\begin{array}{ccc}
F(y)&\ll &1 \textrm{ si $y\leq 1$}\\
F(y)&\ll_{A}& y^{-A} \textrm{ si $y>1$.}
\end{array}
\end{displaymath}
Ceci donne:
\begin{eqnarray}
Err(\q)\ll_{\delta}N(\q)^{1/4}\times \frac{N(\q)^{1/2+\delta/4}M^{1+\delta/2}}{N(\q)^{1-\delta}}
\ll_\delta N(\q)^{1/4} \times N(\q)^{\frac{\Delta-1}{2}+2\delta}\nonumber
\end{eqnarray}
et $\Delta$ \'etant fix\'e dans $]0,1[$, cela termine la preuve de l'estimation asymptotique du premier moment.

\section{Le deuxi\`eme moment}\label{Moment2}

Cette partie est analogue \`a la pr\'ec\'edente, et proc\`edera de la m\^eme fa\c con, afin de prouver (\ref{moment2}):
\begin{multline}\nonumber
M_2(\q):=\sum_{\pi\in\Pi_\q^{\g{k}}}^h \Lambda(1/2,\pi)^2\M(\pi)^2\\ = \frac{\zeta_F(2)^2\Gamma\left(\frac{\g{k}}{2}\right)^2}
{(2\pi)^{\g{k}}\res_{s=1}(\zeta_F)^2}\times\frac{8N(\q)^{1/2}}{\log(N(\q))^2}
\Bigg( \frac{\|P''\|_{L^2(0,1)}^2}
{\Delta^3}+\frac{P'(1)^2}{\Delta^2}+\mathcal{O}\bigg(\frac{1}
{\log(N(\q))}\bigg)\Bigg).
\end{multline}
 
On cherche donc un \'equivalent de l'expression:
$$M_2(\q)=\sum^h_{\pi\in\Pi_\q^{\g{k}}}\Lambda(1/2,\pi)^2M(\pi)^2$$
avec (voir (\ref{valeur2})):
$$\Lambda(\pi,1/2)^2=2N(\q)^{1/2}\sum_{\n\subset \entier}
\frac{\lambda_\pi(\n)}{\sqrt{N(\n)}}\tau(\n) G(N(\n)/N(\q))$$
o\`u:
$$G(y)=\frac{1}{2i\pi}\int_{(3/2)}y^{-s}\zeta_F^{(\q)}(1+2s)L\left(s+\frac{1}{2},\pi_\infty\right)^2\frac{ds}{s}.$$
En d\'eveloppant, on arrive \`a:
\begin{multline}
M_2(\q)=2N(\q)^{1/2}\sum_{\n\subset \entier}\sum_{\substack{ N(\m_1)\leq M\\N(\m_2)\leq M}}\frac{G(N(\n)/N(\q))\tau(\n)}{N(\n)^{1/2}}\cdot
\frac{\mu(\m_1)P\left(\frac{\log(M/N(\m_1))}{\log(M)}\right)}{\psi(\m_1)N(\m_1)^{1/2}}\times\\
\frac{\mu(\m_2)P\left(\frac{\log(M/N(\m_2))}{\log(M)}\right)}{\psi(\m_2)N(\m_2)^{1/2}}
\sum_{\pi\in\Pi_\q^{\g{k}}}^h\lambda_\pi(\n)\lambda_\pi(\m_1)\lambda_\pi(\m_2).
\end{multline}
La formule de Petersson donne une expression pour les formes bilin\'eaires en les coefficients des fonctions $L$. Pour 
l'appliquer, il faut transformer la forme trilin\'eaire ci-dessus, ce qui se fait en utilisant les propri\'et\'es multiplicatives de ces coefficients:
$$\lambda_\pi(\m_1)\lambda_\pi(\m_2)=\sum_{\D|(\m_1,\m_2)}\lambda_\pi(\m_1\m_2\D^{-2})\chi_\q(\D)$$
avec $\chi_\q(\D)=0$ si $\q|\D$, 1 sinon. Soit, apr\`es un changement de variables ($\m_1\leftarrow \m_1\D^{-1},
\m_2\leftarrow \m_2\D^{-1}$) en notant qu'ici $\chi_\q(\D)=1$ (car $N(\D)\leq M<N(\q)$):

\begin{multline}
M_2(\q)=2N(\q)^{1/2}\sum_{\n\subset\entier}\sum_{N(\D)\leq M }\sum_{\substack{  N(\m_1)\leq M/N(\D)\\N(\m_2)\leq M/N(\D)}}
\frac{G(N(\n)/N(\q))\tau(\n)}{N(\n)^{1/2}}\times\\
\frac{\mu(\D\m_1)P\left(\frac{\log(M/N(\D\m_1))}{\log(M)}\right)}{\psi(\D\m_1)N(\D\m_1)^{1/2}}\cdot
\frac{\mu(\D\m_2)P\left(\frac{\log(M/N(\D\m_2))}{\log(M)}\right)}{\psi(\D\m_2)N(\D\m_2)^{1/2}}
\sum_{\pi\in\Pi_\q^{\g{k}}}^h\lambda_\pi(\n)\lambda_\pi(\m_1\m_2).
\end{multline}
Pareillement \`a la pr\'ec\'edente section, on pose $\n=\nu\A,\m_1=\xi_1\B_1,\m_2=\xi_2\B_2$: $\A,\B_1,\B_2$ parcourent un ensemble de repr\'esentants du groupe restreint $\mathscr{C}\ell^+(F)$, $\nu$ parcourt $(\A^{-1})^{\gg 0}\mod\unite$, $\xi_j$ d\'ecrit $(\B_j^{-1})^{\gg 0}\mod\unite$ et 
on applique la formule de Petersson. \\
Il s'ensuit un terme diagonal, qui est dominant; un terme en 
sommes de Kloosterman, qui est n\'egligeable; un terme en les formes anciennes, qui l'est aussi comme on le verra dans  l'appendice.\\
Nous montrerons donc ici les deux premiers points.\\
On aura recours \`a la notation:
\begin{eqnarray}
G_\n:=\frac{G(N(\n)/N(\q))\tau(\n)}{N(\n)^{1/2}}\\
P_\m:=\frac{\mu(\m)P\left(\frac{
\log(M/N(\m))}{\log(M)}\right)}{\psi(\m)N(\m)^{1/2}}.
\end{eqnarray}

\noindent$\bullet${\textbf{ Le terme diagonal}}\\
Le terme diagonal s'\'ecrit:
\begin{multline}
M_2^{diag}(\q)=2N(\q)^{1/2}\sum_{\bar{\A},\bar{\B_1},\bar{\B_2}}\sum_{N(\D)\leq M} \sum_{\nu\in(\A^{-1})^{\gg 0}/\unite}G_{\nu\A}\times\\
\sum_{\substack{\xi_1\in(\B_1^{-1})^{\gg 0}/\unite\\ N(\xi_1\B_1)\leq M/N(\D)}}\sum_{\substack{\xi_2\in(\B_2^{-1})^{\gg 0}/\unite\\N(\xi_2\B_2)\leq M/N(\D)}}
P_{\xi_1\B_1\D}P_{\xi_2\B_2\D}\times\mathbbmss{1}_{\nu\A=\xi_1\xi_2\B_1\B_2}.
\end{multline}
On peut donc reparam\'etrer ces sommes par les id\'eaux $\m_1=\xi_1\B_1,\m_2=\xi_2\B_2$ ce qui donne:
\begin{multline}
M_2^{diag}(\q)=2N(\q)^{1/2}\sum_{N(\D)\leq M}
\sum_{N(\m_1)\leq M/N(\D)}\sum_{N(\m_2)\leq M/N(\D)}
G_{\m_1\m_2}P_{\m_1\D}P_{\m_2\D}\end{multline}
et donc, avec les expressions int\'egrales des fonctions $G$ et $P$ rappel\'ees pr\'ec\'edemment:
\begin{multline}
M_2^{diag}(\q)=\frac{2N(\q)^{1/2}}{(2i\pi)^3}\iiint_{\left(\frac{1}{2}\right),\left(\frac{1}{2}\right),\left(\frac{1}{2}\right)}
M^{t_1+t_2}\PM(t_1)\PM(t_2)\times\\ N(\q)^{s}\zeta_F^{(\q)}(1+2s)L\left(s+\frac{1}{2},\pi_\infty\right)^2   D_2(s,t_1,t_2)\frac{ds}{s}\frac{dt_1}{t_1}
\frac{dt_2}{t_2}
\end{multline}
o\`u l'on a pos\'e pour $s,t_1,t_2$ de parties r\'eelles positives
\begin{eqnarray}
D_2(s,t_1,t_2)=\sum_{\D,\m_1,\m_2}\frac{\mu(\m_1\D)\mu(\m_2\D)\tau(\m_1\m_2)}{\psi(\m_1\D)\psi(m_2\D)N(\D)^{1+t_1+t_2}N(\m_1)^{1+s+t_1}
N(\m_2)^{1+s+t_2}}.\nonumber
\end{eqnarray}
Pour simplifier l'exposition qui suit, r\'esumons ici quelques faits \`a propos de $D_2$:

\begin{lemme}
On a le d\'eveloppement eul\'erien, convergent pour $\Re(t_1+t_2)>0,\Re(s+t_1)>0,\Re(s+t_2)>0$:
\begin{multline}
D_2(s,t_1,t_2)=\\\prod_{\p}\Big(1+\frac{N(\p)^{-(1+t_1+t_2)}}{(1+N(\p)^{-1})^2}-\frac{2N(\p)^{-(1+s+t_1)}}{1+N(\p)^{-1}}-
\frac{2N(\p)^{-(1+s+t_2)}}{1+N(\p)^{-1}}+\frac{3N(\p)^{-(2+2s+t_1+t_2)}}{(1+N(\p)^{-1})^2}\Big)\nonumber
\end{multline}
De plus, on peut \'ecrire:
\begin{eqnarray}
D_2(s,t_1,t_2)=\frac{\zeta_F(1+t_1+t_2)}{\zeta_F(1+s+t_1)^2\zeta_F(1+s+t_2)^2}\eta(s,t_1,t_2)
\end{eqnarray}
o\`u $\eta$ est homolomorphe sur (au moins) $\{\Re(t_1+t_2)>-1/2,\Re(s+t_1)>-1/2,\Re(s+t_2)>-1/2\}$, et vaut en z\'ero:
\begin{eqnarray}
\eta(0,0,0)=\zeta_F(2)^2.
\end{eqnarray}
\end{lemme}
La preuve est un calcul facile. L'expression de $D_2$ permet de voir qu'elle est prolongeable m\'eromorphiquement sur un 
voisinage de z\'ero, ce qui permettra d'effectuer des changements de ligne d'int\'egration comme dans la section 
pr\'ec\'edente. C'est dans le calcul explicite de $\eta(0,0,0)$, n\'ecessaire pour un \'equivalent explicite du second 
moment, que l'introduction du terme $\psi(\m)$ dans la d\'efinition de l'amollisseur est important, car il permet d'achever les 
calculs. Pour le v\'erifier, on regarde le facteur en $\p$ du d\'eveloppement eul\'erien de $\eta$, en prenant $s=t_1=t_2=0$:
\begin{multline}
\Big(1+\frac{N(\p)^{-1}}{(1+N(\p)^{-1})^2}-\frac{2N(\p)^{-1}}{1+N(\p)^{-1}}-
\frac{2N(\p)^{-1}}{1+N(\p)^{-1}}+\frac{3N(\p)^{-2}}{(1+N(\p)^{-1})^2}\Big)\times\\ \frac{1-N(\p)^{-1}}{(1-N(\p)^{-1})^4}
=\frac{1-N(\p)^{-1}}{(1+N(\p)^{-1})^2}\times\frac{1}{(1-N(\p)^{-1})^3}=(1-N(\p)^{-2})^{-2}.\nonumber
\end{multline}
En faisant le produit sur tous les id\'eaux premiers $\p$, on trouve bien $\eta(0,0,0)=\zeta_F(2)^2$. $\blacksquare$\\

Revenant au second moment amolli, on change de droite d'int\'egration en $s$ de $(\frac{1}{2})$ \`a $(-\frac{1}{2}+\delta)$:
\begin{multline}
M_2^{diag}(\q)=\frac{2N(\q)^{1/2}}{(2i\pi)^2}\iint_{\left(\frac{1}{2}\right),\left(\frac{1}{2}\right)}
M^{t_1+t_2}\PM(t_1)\PM(t_2)\times\\ {\textrm{res}}_{s=0}\Big(s^{-1}N(\q)^{s}\zeta_F^{(\q)}(1+2s)L\left(s+\frac{1}{2},\pi_\infty\right)^2
D_2(s,t_1,t_2)\Big)\frac{dt_1}{t_1}
\frac{dt_2}{t_2}\\+\frac{2N(\q)^{1/2}}{(2i\pi)^3}\iiint_{\left(-\frac{1}{2}+\delta\right),\left(\frac{1}{2}\right),\left(\frac{1}{2}\right)}
M^{t_1+t_2}\PM(t_1)\PM(t_2) N(\q)^{s}\zeta_F^{(\q)}(1+2s)\times\\ L\left(s+\frac{1}{2},\pi_\infty\right)^2D_2(s,t_1,t_2)\frac{ds}{s}\frac{dt_1}{t_1}
\frac{dt_2}{t_2}
\end{multline}
$\delta$ \'etant choisi tel que $0<\delta<\frac{1-\Delta}{2}$, auquel cas on a:
\begin{multline}
M_2^{diag}(\q)=\frac{2N(\q)^{1/2}}{(2i\pi)^2}\iint_{\left(\frac{1}{2}\right),\left(\frac{1}{2}\right)}
M^{t_1+t_2}\PM(t_1)\PM(t_2)\times\\ {\textrm{res}}_{s=0}\Big(s^{-1}N(\q)^{s}\zeta_F^{(\q)}(1+2s)L\left(s+\frac{1}{2},\pi_\infty\right)^2   
D_2(s,t_1,t_2)\Big)\frac{dt_1}{t_1}
\frac{dt_2}{t_2}\\+\mathcal{O}(N(\q)^{1/2-\delta'})
\end{multline}
avec $\delta'=\frac{1-\Delta}{2}-\delta>0$, not\'e $\delta$ dans la suite.

Il faut donc calculer le terme constant du d\'eveloppement en s\'erie enti\`ere (en $s$) de $N(\q)^{s}\zeta_F^{(\q)}(1+2s)
L\left(s+\frac{1}{2},\pi_\infty\right)^2   D_2(s,t_1,t_2)$; en fait, on voit que ce terme s'\'ecrit:
\begin{multline}
\frac{{\textrm{res}}_{1}(\zeta_F^{(\q)})}{2}\log(N(\q))D_2(0,t_1,t_2)L\left(\frac{1}{2},\pi_\infty\right)^2 
+{\textrm{res}}_{1}(\zeta_F^{(\q)})L\left(\frac{1}{2},\pi_\infty\right)^2 \times\\  \frac{\partial D_2}{\partial s}(0,t_1,t_2)
+\ldots\nonumber
\end{multline}
o\`u les autres termes n'ont pas de contribution au terme principal.  On a donc:

\begin{multline}
M_2^{diag}(\q)= \frac{2N(\q)^{1/2}}{(2i\pi)^2}\iint_{\left(\frac{1}{2}\right),\left(\frac{1}{2}\right)}
\frac{dt_1}{t_1}\frac{dt_2}{t_2} M^{t_1+t_2}\PM(t_1)\PM(t_2)\times\\
\Big( \frac{{\textrm{res}}_{1}(\zeta_F)}{2}\log(N(\q))L\left(\frac{1}{2},\pi_\infty\right)^2  D_2(0,t_1,t_2)+
{\textrm{res}}_{1}(\zeta_F)L\left(\frac{1}{2},\pi_\infty\right)^2\times\\   \frac{\partial D_2}{\partial s}(0,t_1,t_2)\Big)
+\mathcal{O}(N(\q)^{1/2-\delta}).
\end{multline}
On peut donc \'ecrire:
$$M_2^{diag}(\q)=M_{21}(\q)+M_{22}(\q)+\mathcal{O}(N(\q)^{1/2-\delta})$$
avec:
\begin{multline}
M_{21}(\q)= \frac{N(\q)^{1/2}\log(N(\q)){\textrm{res}}_1(\zeta_F)L\left(\frac{1}{2},\pi_\infty\right)^2}{(2i\pi)^2}\times\\
\iint_{\left(\frac{1}{2}\right),\left(\frac{1}{2}\right)}M^{t_1+t_2}\PM(t_1)\PM(t_2)D_2(0,t_1,t_2)
\frac{dt_1}{t_1}\frac{dt_2}{t_2} 
\end{multline}
\begin{multline}
M_{22}(\q)= \frac{2N(\q)^{1/2}{\textrm{ res}}_1(\zeta_F)L\left(\frac{1}{2},\pi_\infty\right)^2 }{(2i\pi)^2}\times\\
\iint_{\left(\frac{1}{2}\right),\left(\frac{1}{2}\right)}M^{t_1+t_2}\PM(t_1)\PM(t_2)\frac{\partial D_2}{\partial s} (0,t_1,t_2)
\frac{dt_1}{t_1}\frac{dt_2}{t_2}.
\end{multline}
En d\'epla\c cant les lignes d'int\'egration, on trouve donc ($\delta$ d\'esigne un r\'eel positif):
\begin{multline}
M_{21}(\q)= N(\q)^{1/2}\log(N(\q)){\textrm{res}}_1(\zeta_F)L\left(s+\frac{1}{2},\pi_\infty\right)^2\times\\
{\textrm{res}}_{(t_1,t_2)=(0,0)}\left(t_1^{-1}t_2^{-1}M^{t_1+t_2}\PM(t_1)\PM(t_2) D_2(0,t_1,t_2)\right)
+\mathcal{O}(N(\q)^{1/2-\delta}) 
\end{multline}
et
\begin{multline}
M_{22}(\q)= 2N(\q)^{1/2}{\textrm{ res}}_1(\zeta_F)L\left(\frac{1}{2},\pi_\infty\right)^2\times\\
{\textrm{res}}_{(t_1,t_2)=(0,0)}\left(t_1^{-1}t_2^{-1}M^{t_1+t_2}\PM(t_1)\PM(t_2)\frac{\partial D_2}{\partial s} (0,t_1,t_2)
\right)+\mathcal{O}(N(\q)^{1/2-\delta}).
\end{multline}
En se r\'ef\'erant au lemme pr\'ec\'edent, on note que dans le terme 
$\frac{\partial D_2}{\partial s}(0,t_1,t_2)$, seul $\zeta_F(1+t_1+t_2)\eta(0,t_1,t_2)
\frac{\partial }{\partial s}(\zeta_F(1+s+t_1)^{-2}\zeta_F(1+s+t_2)^{-2}) |_{s=0}$ 
produit un terme principal pour $M_{22}(\q)$. De plus, il est facile de voir que pour trouver les termes principaux de 
$M_{21}(\q)$ et $M_{22}(\q)$ on peut remplacer toutes les expressions en $\zeta_F(1+x)$ par ${\rm{res}}(\zeta_F)x^{-1}$: les 
autres termes du d\'eveloppement de $\zeta_F$ ont une contribution inf\'erieure d'une puissance en $\log(M)$:
\begin{multline}
M_{21}(\q)= N(\q)^{1/2}\log(N(\q)){\textrm{res}}_1(\zeta_F)^{-2}L\left(\frac{1}{2},\pi_\infty\right)^2\times\\\Big\{
{\textrm{res}}_{(t_1,t_2)=(0,0)}\left(M^{t_1+t_2}\PM(t_1)\PM(t_2)\eta(0,0,0) \frac{t_1t_2}{t_1+t_2}\right)+\mathcal{O}(\log(M)^{-4})
\Big\}\\+\mathcal{O}(N(\q)^{1/2-\delta}) \nonumber
\end{multline}
et
\begin{multline}
M_{22}(\q)= 2N(\q)^{1/2}{\textrm{ res}}_1(\zeta_F)^{-2}L\left(\frac{1}{2},\pi_\infty\right)^2\times\\\Big\{
{\textrm{res}}_{(t_1,t_2)=(0,0)}\left(M^{t_1+t_2}\PM(t_1)\PM(t_2)\eta(0,0,0)
\right)+\mathcal{O}(\log(M)^{-3})\Big\}\\+\mathcal{O}(N(\q)^{1/2-\delta}).\nonumber
\end{multline}
Enfin:
\begin{multline}
{\textrm{res}}_{(t_1,t_2)=(0,0)}\left( M^{t_1+t_2}\PM(t_1)\PM(t_2) \frac{t_1t_2}{t_1+t_2}\right)\\=
{\textrm{res}}_{t_2=0}\left\{ M^{t_2}\PM(t_2){\textrm{res}}_{t_1=0}\left(M^{t_1}\PM(t_1)\frac{t_1}{1+t_1/t_2}\right)\right\}
\\={\textrm{res}}_{t_2=0}\left\{M^{t_2}\PM(t_2){\textrm{res}}_{t_1=0}\left(M^{t_1}\PM(t_1)
\sum_{k\geq 0}(-t_2)^{-k}t_1^{k+1}\right)\right\}\\={\textrm{res}}_{t_2=0}\left\{M^{t_2}\PM(t_2)
\sum_{\ell+k+1-j=-1}\frac{\log(M)^\ell}{\ell !}(-t_2)^{-k}\frac{a_jj!}{\log(M)^j}\right\}\\=
\sum_{\ell+k+1-j=-1}\log(M)^{\ell-j+k-1}\frac{a_jj!}{\ell !}(-1)^k\sum_{r-s=k-1}\frac{a_ss!}{r!}\nonumber
\end{multline}
en notant alors $P^{(n)}$ la $n$-i\`eme d\'eriv\'ee de $P$, et $${}^{(n)}P(X)=\sum_{k}\frac{a_{n+k}}{(n+k)(n+k-1)\ldots(k+1)}
X^{n+k}$$ on a alors, en se r\'ef\'erant \`a l'appendice de \cite{KMV} o\`u ces r\'esidus sont calcul\'es:
\begin{multline}
{\textrm{res}}_{(t_1,t_2)=(0,0)}\left( M^{t_1+t_2}\PM(t_1)\PM(t_2) \frac{t_1t_2}{t_1+t_2}\right)\\= \log(M)^{-3}
\sum_{k\geq 0}(-1)^kP^{(k+2)}(1) {}^{(k-1)}P(1)= \log(M)^{-3}\int_0^1P''(t)^2dt.
\end{multline}
On trouve plus facilement:
\begin{eqnarray}
\res_{(t_1,t_2)=(0,0)}\left(M^{t_1+t_2}\PM(t_1)\PM(t_2)\right)
= \log(M)^{-2}P'(1)^2
\end{eqnarray}
ce qui, tout compte fait, donne:
\begin{multline}
M_2^{diag}(\q)=8N(\q)^{1/2} {\textrm{ res}}_1(\zeta_F)^{-2}L\left(\frac{1}{2},\pi_\infty\right)^2\zeta_F(2)^2\times\\
\Big\{\frac{\log(N(\q))}{8}\log(M)^{-3}\int_0^1P''(t)^2dt+
\frac{\log(M)^{-2}}{4}P'(1)^2+\mathcal{O}(\log(N(\q))^{-3})\Big\}\\
+\mathcal{O}(N(\q)^{1/2-\delta})
\end{multline}
qui est bien le terme annonc\'e dans (\ref{moment1}).\\

\noindent$\bullet${\textbf{ Le terme de Kloosterman}}\\
Comme dans la section \ref{Moment1}, on peut \'ecrire la contribution des termes en sommes de Kloosterman du deuxi\`eme 
moment sous la forme:
\begin{eqnarray}
M_{2}^{Kloost}(\q)=\sum_{\A,\B_1,\B_2}M_{2,\A,\B_1,\B_2}^{Kloost}(\q)
\end{eqnarray}
o\`u $\A,\B_1,\B_2$ d\'esignent une famille de repr\'esentants de $\mathscr{C}\!\ell^+(F)$ telle que celle choisie en section 
\ref{Notations}, et avec:

\begin{multline}
M_{2,\A,\B_1,\B_2}^{Kloost}(\q)=2N(\q)^{1/2} \sum_{\nu\in(\A^{-1})^{\gg 0}/\unite}G_{\nu\A}\times\\
\sum_{N(\D)\leq M}
\sum_{\substack{\xi_1\in(\B_1^{-1})^{\gg 0}/\unite\\ N(\xi_1\B_1)\leq M/N(\D)}}\sum_{\substack{\xi_2\in(\B_2^{-1})^{\gg 0}/\unite\\N(\xi_2\B_2)\leq M/N(\D)}}
P_{\xi_1\B_1\D}P_{\xi_2\B_2\D}\\
\sum_{\substack{\overline{\cc}^2=\overline{\A\B_1\B_2}\\
c\in\cc^{-1}\q\setminus \{0\}\\ 
\varepsilon\in\mathcal{O}_F^{\times +}/\mathcal{O}_F^{\times 2}}}
\frac{\kl(\varepsilon\nu,\A;\xi_1\xi_2,\B_1\B_2;c,\cc)}{N(c\cc)}J_{\g{k-1}}\left(4\pi\frac{\sqrt{
\g{\varepsilon\nu\xi_1\xi_2}[\A\B_1\B_2\cc^{-2}]}}{|\g{c}|}\right).
\end{multline}

Fixons donc $\A,\B_1,\B_2$ et $\cc$ tel que $\cc^2\sim\A\B_1\B_2$. La somme finie en $\varepsilon$ n'a pas d'incidence analytique; on doit donc traiter:
\begin{multline}
Err_1(\q)=N(\q)^{1/2} \sum_{\nu\in(\A^{-1})^{\gg 0}/\unite}G_{\nu\A}\times\\
\sum_{N(\D)\leq M}
\sum_{\substack{\xi_1\in(\B_1^{-1})^{\gg 0}/\unite\\ N(\xi_1\B_1)\leq M/N(\D)}}\sum_{\substack{\xi_2\in(\B_2^{-1})^{\gg 0}/\unite\\N(\xi_2\B_2)\leq M/N(\D)}}
P_{\xi_1\B_1\D}P_{\xi_2\B_2\D}\\
\sum_{c\in\cc^{-1}\q\setminus \{0 \} }
\frac{\kl(\nu,\A;\xi_1\xi_2,\B_1\B_2;c,\cc)}{N(c)}J_{\g{k-1}}\left(4\pi\frac{\sqrt{
\g{\nu\xi_1\xi_2}[\A\B_1\B_2\cc^{-2}]}}{|\g{c}|}\right).
\end{multline}
Bien entendu, on supposera que le syst\`eme $\nu,\xi_1,\xi_2$ satisfait au lemme \ref{sommation}, et on utilisera aussi 
(\ref{unites}) dans les m\^emes circonstances.\\

Dans un premier temps, on r\'eduit la somme en $\nu$, gr\^ace \`a l'estimation suivante, qui se prouve comme dans \cite{Va1}, 
lemme 2.3:
\begin{eqnarray}
G(y)\ll\exp(-Cy^{1/2d})
\end{eqnarray}
pour $y>1$, et $C>0$ une constante.\\
On a alors, en utilisant la borne de Weil pour les sommes de Kloosterman, comme le fait \cite{Va2} dans le lemme 3.2:

\begin{multline}
N(\q)^{1/2}\sum_{N(\nu\A)\geq N(\q)^{1+\varepsilon}}\sum_{\D}
\sum_{\xi_1,\xi_2}
\frac{G(N(\nu)/N(\q))\tau(\nu)}{N(\nu)^{1/2}}
P_{\xi_1\B_1\D}P_{\xi_2\B_2\D}\times\\
\sum_{c\in\cc^{-1}\q\setminus \{0 \} }
\frac{\kl(\nu,\A;\xi_1\xi_2,\B_1\B_2;c,\cc)}{N(c)}J_{\g{k-1}}\left(4\pi\frac{\sqrt{
\g{\nu\xi_1\xi_2}[\A\B_1\B_2\cc^{-2}]}}{|\g{c}|}\right)
\\ \ll_{\varepsilon}\exp(-CN(\q)^{\varepsilon})
\end{multline}
pour tout $\varepsilon>0$ fix\'e, ce qui r\'eduit l'\'etude \`a :
\begin{multline}
Err_2(\q)=
N(\q)^{1/2}\sum_{N(\nu\A)\leq N(\q)^{1+\varepsilon}}\sum_{\D}
\sum_{\xi_1,\xi_2}
\frac{G(N(\nu)/N(\q))\tau(\nu)}{N(\nu)^{1/2}}
P_{\xi_1\B_1\D}P_{\xi_2\B_2\D}\times\\
\sum_{\eta\in\unite}\sum_{c\in\cc^{-1}\q\setminus \{0 \}/\unite }
\frac{\kl(\nu,\A;\xi_1\xi_2,\B_1\B_2;c\eta,\cc)}{N(c)}J_{\g{k-1}}\left(4\pi\frac{\sqrt{
\g{\nu\xi_1\xi_2}[\A\B_1\B_2\cc^{-2}]}}{|\g{c\eta}|}\right).
\end{multline}
Pour \'etudier cette somme, on utilise l'expression int\'egrale suivante:
\begin{multline}\label{abc}
J_{\g{k-1}}(\g{x})=\int_{(\g{\sigma})}\frac{\Gamma\Big(\frac{\g{k-1-s}}{2}\Big)}{\Gamma\Big(\frac{\g{k-1+s}}{2}+\g{1}\Big)}
\left(\frac{\g{x}}{2}\right)^{\g{s}}d\g{s}\\=\prod_{j=1}^d\int_{(\sigma_j)}\frac{\Gamma\Big(\frac{k_j-1-s_j}{2}\Big)}
{\Gamma\Big(\frac{k_j-1+s_j}{2}+1\Big)}
\left(\frac{x_j}{2}\right)^{s_j}ds_j
\end{multline}
o\`u l'on choisit $\sigma_j=1$ si $|\eta_j|\leq 1$ et $\sigma_j=1+\delta$ sinon, avec $\delta>0$ , de  sorte que la somme 
sur $\eta$ soit convergente (toujours gr\^ace \`a \cite{luo}, page 136). Cela donne:
\begin{multline}
Err_2(\q)=N(\q)^{1/2}\sum_{\eta\in\unite}\int_{(\g{\sigma})}
\frac{\Gamma\Big(\frac{\g{k-1-s}}{2}\Big)}{\Gamma\Big(\frac{\g{k-1+s}}{2}+\g{1}\Big)}\left(\frac{4\pi[\A\B_1\B_2\cc^{-2}]}{\g{\eta}}\right)^{\g{s}}\\
 \sum_{\substack{\nu\in(\A^{-1})^{\gg 0}/\unite\\N(\nu)\leq N(\q)^{1+\varepsilon}}}G_{\nu\A}\g{\nu}^{\g{s}/2}
\sum_{N(\D)\leq M}
\sum_{\substack{\xi_1\in(\B_1^{-1})^{\gg 0}/\unite\\ N(\xi_1\B_1)\leq M/N(\D)}}\sum_{\substack{\xi_2\in(\B_2^{-1})^{\gg 0}/\unite\\N(\xi_2\B_2)\leq M/N(\D)}}
P_{\xi_1\B_1\D}P_{\xi_2\B_2\D}(\g{\xi_1\xi_2})^{\g{s}/2}\\
\sum_ {c\in\cc^{-1}\q\setminus \{ 0 \} }\frac{\kl(\nu,\A;\xi_1\xi_2,\B_1\B_2;c,\cc)}{N(c)\g{c}^{\g{s}/2}}
\frac{d\g{s}}{\g{s}}.
\end{multline}
Consid\'erons l'expression interne, form\'ee par les sommes en $\D,\nu,\xi_1,\xi_2,c$. En ouvrant les sommes de Kloosterman, elle s'\'ecrit:
\begin{multline}
 R(s)=\sum_{ c\in\cc^{-1}\q\setminus \{ 0\}/\unite }\frac{1}{N(c)\g{c}^{\g{s}}} \sum_{N(\D)\leq M}\sum_{ x\in ( \A\diff^{-1}\cc^{-1}/\A\diff^{-1} c )^\times }
 \sum_{\substack{\nu\in(\A^{-1})^{\gg 0}/\unite\\N(\nu)\leq N(\q)^{1+\varepsilon}}}G_{\nu\A}\times\\ \g{\nu}^{\g{s}/2}
\psi_\infty\left(\frac{\g{x\nu}}{\g{c}}\right)
\sum_{\substack{\xi_1\in(\B_1^{-1})^{\gg 0}/\unite\\ N(\xi_1\B_1)\leq M/N(\D)}}P_{\xi_1\B_1\D}
P_{\xi_2\B_2\D}(\g{\xi_1\xi_2})^{\g{s}/2}
\psi_\infty\left(\frac{\g{\bar{x}}[\A\B_1\B_2\cc^{-2}]\g{\xi_1\xi_2}}{\g{c}}\right).\nonumber
\end{multline}
Posons momentan\'ement:

\begin{eqnarray}
X_{\nu\A} & = & G_{\nu\A}\g{\nu}^{\g{s}/2}\psi_\infty\left(\frac{\g{x\nu}}{\g{c}}\right) \nonumber\\
 Y_{\xi_1,\xi_2,\B_1,\B_2,\D} & = &  P_{\xi_1\B_1\D}P_{\xi_2\B_2\D}(\g{\xi_1\xi_2})^{\g{s}/2}
\psi_\infty\left(\frac{\g{\bar{x}}[\A\B_1\B_2\cc^{-2}]\g{\xi_1\xi_2}}{\g{c}}\right).\nonumber
\end{eqnarray}
L'in\'egalit\'e de Cauchy-Schwarz appliqu\'ee \`a la somme en $x$ donne:

\begin{multline}
R(s)\leq \sum_{N(\D)\leq M}\sum_{ c\in\cc^{-1}\q\setminus \{ 0\}/\unite }\frac{1}{N(c)|\g{c}|^{\Re ( \g{s} ) }}\times\\
\left( \sum_{ x\in ( \A\diff^{-1}\cc^{-1}/\A\diff^{-1} c )^\times }
\Big|\sum_{\substack{\nu\in(\A^{-1})^{\gg 0}/\unite\\N(\nu)\leq N(\q)^{1+\varepsilon}}}X_{\nu\A}\Big|^2\right)^{1/2} \times\\
\left( \sum_{ x\in ( \A\diff^{-1}\cc^{-1}/\A\diff^{-1} c )^\times }
\Big|\sum_{\substack{\xi_1\in(\B_1^{-1})^{\gg 0}/\unite\\ N(\xi_1\B_1)\leq M/N(\D)}}Y_{\xi_1,\xi_2,\B_1,\B_2,\D}\Big|^2 \right)^{1/2}. \nonumber
\end{multline}
On fait maintenant usage du grand crible additif suivant (cf. \cite{IK}, p.178, (7.30) ou \cite{Gal}):
\begin{eqnarray}
\sum_{\g{d}\in\mathbb{Z}^d/\g{c}\mathbb{Z}^d}\Big|\sum_{\substack{\g{n}\in\mathbb{Z}^d\\ n_j\leq X}}
y_{\g{n}}\e\left(\frac{\g{d.n}}{\g{c}}\right)\Big|^2\ll \left(N(\g{c})+X^d\right)\sum_{\substack{\g{n}\in\mathbb{Z}^d\\ n_j\leq X}}|y_{\g{n}}|^2.
\end{eqnarray}
L'hypoth\`ese faite sur les choix de $\nu,\xi_1,\xi_2$ (\`a savoir que leurs composantes sont born\'ees par la norme, selon le lemme \ref{sommation} de la section \ref{geometrie} ) permet de l'appliquer ici et donne, en se souvenant que $\A,\B_1,\B_2,\cc$ ont \'et\'e fix\'es au d\'ebut:
\begin{multline}\nonumber
\sum_{ x\in ( \A\diff^{-1}\cc^{-1}/\A\diff^{-1} c )^\times }
\Big|\sum_{\substack{\nu\in(\A^{-1})^{\gg 0}/\unite\\N(\nu)\leq N(\q)^{1+\varepsilon}}}X_{\nu\A}\Big|^2\ll\\
(N(c)+N(\q)^{1+\varepsilon})\sum_{N(\nu)\leq N(\q)^{1+\varepsilon}}
|G_{\nu\A}\g{\nu}^{\Re(\g{s})}|^2
\end{multline}
\begin{multline}\nonumber
\sum_{ x\in ( \A\diff^{-1}\cc^{-1}/\A\diff^{-1} c )^\times }
\Big|\sum_{\substack{\xi_j\in(\B_j^{-1})^{\gg 0}/\unite\\ N(\xi_j\B_j)\leq M/N(\D)}}Y_{\xi_1,\xi_2,\B_1,\B_2,\D}\Big|^2 
\ll\\ \sum_{\substack{\xi_j\in(\B_j^{-1})^{\gg 0}/\unite\\ N(\xi_j\B_j)\leq M/N(\D)}}
\tau(\xi_1\xi_2)|P_{\xi_1\B_1\D}P_{\xi_2\B_2\D}|^2|\g{\xi_1\xi_2}|^{\Re(\g{s})}.
\end{multline}
En rempla\c cant $G_{\nu\A}$ et $P_{\xi_j\B_j}$ par leur valeur, il vient:
\begin{multline}
R(s)\ll \sum_{N(\D)\leq M}N(\D)^{-1}\sum_{c\in\cc^{-1}\q\setminus\{0\} /\unite}\left(N(c)\g{c}^{\Re(\g{s})}\right)^{-1}
(N(c)+N(\q)^{1+\varepsilon})^{\frac{1}{2}}\times\\ \left( \sum_{\substack{\nu\in(\A^{-1})^{\gg 0}/\unite\\N(\nu)\leq N(\q)^{1+\varepsilon}}}
|G(N(\nu\A)/N(\q))|^2\tau(\nu\A)^2|\nu|^{\Re(\g{s})-1}\right)^{\frac{1}{2}}
\left(N(c)+\left(\frac{M}{N(\D)}\right)^2\right)^{\frac{1}{2}}\times \\ \left(\sum_{\substack{\xi_j\in(\B_j^{-1})^{\gg 0}/\unite\\ N(\xi_j\B_j)\leq M/N(\D)}}
\frac{\tau(\xi_1\xi_2\B_1\B_2)P\left(\frac{\log(M/N(\D\B_1\xi_1))}{\log(M)}\right)^2P\left(\frac{\log(M/N(\D\B_2\xi_2))}{\log(M)}\right)^2|\xi_1\xi_2|^{\Re(\g{s})}}{\psi(\D\xi_1\B_1)^2\psi(\D\B_2\xi_2)^2|\xi_1\xi_2|}
\right)^{\frac{1}{2}}\\
\ll_{\delta} \sum_{c\in\q\setminus \{0\}/\unite}\left(N(c)\g{c}^{\Re(\g{s})}\right)^{-1}
(N(c)+N(\q)^{1+\varepsilon})^{1/2}N(\q)^{(1+\varepsilon)(1+d\delta)/2}\times\\
\left(N(c)+\left(\frac{M}{N(\D)}\right)^2\right)^{1/2}M^{1+d\delta}
\ll N(\q)^{-1-\delta}N(\q)^{(1/2+\varepsilon/2)(1+d\delta)}N(\q)^{\frac{\Delta}{2}(1+d\delta)}\\ \ll N(\q)^{-\alpha}
\end{multline}
pour un $\alpha>0$ d\'ependant de $\Delta,\varepsilon,\delta$, d\`es que $\varepsilon,\delta$ sont choisis assez petits. Cette 
derni\`ere in\'egalit\'e implique que:
\begin{eqnarray}
Err_2(\q)\ll N(\q)^{1/2-\alpha}
\end{eqnarray}
et ach\`eve donc cette section.

\section{Conclusion}\label{Conclusion}
\subsection{Fin de la preuve du th\'eor\`eme \ref{nonannulation}}
Les estimations (\ref{moment1}) et (\ref{moment2}) \'etant maintenant prouv\'ees, on a donc, comme dit en introduction:
\begin{eqnarray}
\liminf _{N(\q)\to\infty}\sum_{\pi\in\Pi_\q^{\g{k}}}^h\mathbbmss{1}_{\Lambda(1/2,\pi)\neq 0}\geq
\frac{\Delta P'(1)^2}{2(||P''||_{L^2(0,1)}+\Delta P'(1)^2)}
\end{eqnarray}
et ce pour tous $\Delta\in]0,1[$ tel que $N(\q)^{\Delta/2}\notin \N$ et $P$ polyn\^ome d'ordre au moins deux en z\'ero. Ceci 
entra\^ine:
\begin{eqnarray}
\liminf _{N(\q)\to\infty}\sum_{\pi\in\Pi_\q^{\g{k}}}^h\mathbbmss{1}_{\Lambda(1/2,\pi)\neq 0}\geq
\sup_{\substack{P\\P(0)=P'(0)=0}} \frac{1}{2\left(1+\frac{\int_0^1P''(t)^2dt}{P'(1)^2}\right)}
\end{eqnarray}
or, d'apr\`es l'in\'egalit\'e de Cauchy Schwarz:
$$P'(1)^2=\left(\int_0^1P''(t)dt\right)^2\leq\int_0^1P''(t)^2dt$$
avec \'egalit\'e si et seulement si $P''$ est constante. Ceci donne donc l'in\'egalit\'e
\begin{eqnarray}
\liminf _{N(\q)\to\infty}\sum_{\pi\in\Pi_\q^{\g{k}}}^h\mathbbmss{1}_{\Lambda(1/2,\pi)\neq 0}\geq
\frac{1}{4}
\end{eqnarray}
obtenue pour le polyn\^ome optimal $P(X)=X^2$.\\

\subsection{L'in\'egalit\'e de grand crible}
Venons-en maintenant au grand crible, dont la preuve suit le m\^eme mod\`ele que ce qui pr\'ec\`ede. Soit $x=\{x_\n\}$ 
une suite de nombres complexes, $X$ un r\'eel positif, et $\q$ un id\'eal quelconque. On \'ecrit, par positivit\'e:
\begin{multline}
\sum_{\pi\in\Pi_\q^{\g{k}}}^h\Big|\sum_{\n\leq X}\lambda_\pi(\n)x_\n\Big|^2\leq \sum_{\varphi}^h\Big|\sum_{\substack{
\overline{\A}\\\nu\in(\A^\times)^{\gg 0}/\unite\\ N(\nu\A)\leq X}}
\lambda(\nu,\A,\varphi)x_{\nu\A}\Big|^2=\\ \sum_{\varphi}^h\sum_{\substack{\overline{\A},\overline{\B}\\ \nu,\mu}}\lambda(\nu,\A,\varphi)\overline{\lambda(\mu,\B,\varphi)}
x_{\nu,\A}\overline{x_{\mu\B}}
\end{multline}
$\varphi$ parcourant une base orthogonale de $\mathcal{H}_\q^{\g{k}}$. Comme pr\'ec\'edemment, $\{\A,\B\}$ d\'esignent une famille de repr\'esentants du 
groupe des classes \'etroit et on a pos\'e $\n=\nu\A,\m=\mu\B$, $\nu$ et $\mu$ satisfaisant au lemme \ref{sommation}. On  applique la formule de Petersson (th\'eor\`eme \ref{petersson}):
\begin{multline}
\sum_{\pi\in\Pi_\q^{\g{k}}}^h\Big|\sum_{\n\leq X}\lambda_\pi(\n)x_\n\Big|^2\leq 
\sum_{\substack{N(\n)\leq X\\N(\m)\leq X}}x_\n\overline{x_\m}\mathbbmss{1}_{\n=\m}+\frac{C}{|\disc|^{1/2}}\times\\ \sum_{\bar{\A},\bar{\B}}\sum_{\substack{\nu\in(\A^{-1})^{\gg 0}/\unite
\\\mu\in(\B^{-1})^{\gg 0}/\unite\\N(\nu\A),N(\mu\B)\leq X}}
\sum_{\cc,c,\varepsilon}\frac{\kl(\nu,\A;\mu,\B;c,\cc)}{N(c\cc)}J_{\g{k-1}}\left(4\pi\frac{\sqrt{\varepsilon\mu\nu[\A\B\cc^{-1}]}}
{|c|}\right)x_{\nu\A}\overline{x_{\mu\B}}.\nonumber
\end{multline}
 Le premier terme est exactement $||x_X||_2^2$. On effectue le calcul du second terme pour chaque valeur de $\A,\B,\cc$ car il n'y en a qu'un nombre fini. D'autre part, les difficult\'es \'etant similaires \`a celles rencontr\'ees lors de l'estimation du deuxi\`eme moment, on supposera que $\A=\B=\cc=\entier$ pour simplifier les notations. On traite donc:
\begin{multline}
K(X)=\sum_{\substack{\nu\in\entier^{\gg 0}/\unite
\\ N(\nu)\leq X}} \sum_{\substack{ \mu\in\entier^{\gg 0}/\unite\\ N(\mu)\leq X}}
\sum_{c\in\q\setminus \{0\}}\frac{\kl(\nu;\mu;c)}{N(c)}J_{\g{k-1}}\left(4\pi\frac{\sqrt{\g{\varepsilon\mu\nu}}}
{|\g{c}|}\right)x_{\nu}\overline{x_{\mu}}=\\
\sum_{c\in\q\setminus \{0\}/\unite}\sum_{\eta\in\unite}\int_{\Re(\g{s})=\g{\sigma}}\frac{(4\pi)^{\g{s}}\Gamma\left(\frac{\g{k-1-s}}{2}\right)}
{|\g{c\eta}|^{\g{s}}\Gamma\left(\frac{\g{k
+s+1}}{2}\right)}\sum_{\nu,\mu}\frac{\kl(\nu;\mu;c\eta)}{N(c)}x_{\nu}\g{\nu}^{\g{s}/2}\overline{
x_{\mu}}\g{\mu}^{\g{s}/2}d\g{s}\nonumber
\end{multline}
en choisissant $\g{\sigma}$ comme lors de (\ref{abc}), de sorte que la somme en $\eta$ soit convergente, et on ouvre les sommes de Kloosterman (la d\'emarche est la m\^eme que celle qui nous a permis de traiter le terme de Kloosterman du deuxi\`eme moment):
\begin{multline}
|K(X)|=\Big|\sum_{c\in\q\setminus \{0\}/\unite}\int_{\Re(\g{s})=\g{\sigma}}\frac{(4\pi)^{\g{s}}
\Gamma\left(\frac{\g{k-1-s}}{2}\right)}{N(c)|\g{c}|^{\g{s}}\Gamma\left(\frac{\g{k+s+1}}{2}\right)}\times\\
\sum_{x}\sum_{N(\nu),N(\mu)\leq X}\psi_\infty\left(\frac{x\nu}{c}\right)x_{\nu}\g{\nu}^{\g{s}/2}\psi_\infty\left(\frac
{\overline{x}\mu}{c}\right) \overline{x_{\mu}}\g{\mu}^{\g{s}/2}d\g{s}\Big|\nonumber
\end{multline}

\begin{multline}
\leq 
\sum_{c\in\q\setminus \{0\}/\unite}\int_{\Re(\g{s})=\g{\sigma}}\frac{(4\pi)^{\g{\sigma}}|\Gamma\left(\frac{\g{k-1-s}}{2}\right)|}
{N(c)|\g{c}|^{\g{\sigma}}|\Gamma\left(\frac{\g{k
+s+1}}{2}\right)|}\times\\
\left(\sum_{x}\Big|\sum_{N(\nu )\leq X}\psi_\infty\left(\frac{x\nu}{c}\right)x_{\nu}\g{\nu}^{\g{s}/2}\Big|^2\times
\sum_{x}\Big|\sum_{N(\mu )\leq X}\psi_\infty\left(\frac{\overline{x}\mu}{c}\right)\overline{x_{\mu}}\g{\mu}^{\g{s}/2}
\Big|^2\right)^{\frac{1}{2}}d\g{s}\\
\ll \sum_{c\in\q\setminus \{0\}/\unite}\int_{\Re(\g{s})=\g{\sigma}}\frac{|\Gamma\left(\frac{\g{k-1-s}}{2}\right)|}
{N(c)|\g{c}|^{\g{\sigma}}|\Gamma\left(\frac{\g{k
+s+1}}{2}\right)|}(N(c)+X)\sum_{N(\nu)\leq X}|x_\nu|^2\g{\nu}^{\g{s}}
d\g{s}\\
\ll_{F,\sigma} \frac{X^\sigma}{N(\q)^{\sigma}}\left(1+\frac{X}{N(\q)}\right)||x_X||_2^2.
\nonumber
\end{multline}
Avec le choix de $\sigma$ pr\'ec\'edent, et si $X\leq N(\q)$, on trouve bien:
$$ \sum_{\varphi}^h\Big|\sum_{\substack{\overline{\A}\\\nu\in(\A^\times)^{\gg 0}/\unite\\ N(\nu\A)\leq X}}
\lambda(\nu,\A,\varphi)x_{\nu\A}\Big|^2\ll_{F}\left(1+\frac{X}{N(\q)}\right)||x_X||_2^2.$$
Pour $X>N(\q)$, l'astuce de Iwaniec fonctionne aussi:\\
cela consiste en le choix d'un premier $\p$ tel que $X\leq N(\q\p)$, tout en ayant $N(\p\q)\ll X$. On utilise alors 
l'injection $\mathcal{H}_\q^{\g{k}}\hookrightarrow \mathcal{H}_{\q\p}^{\g{k}}$ (non isom\'etrique!): par cette injection, 
on a $ ||\varphi||_{\mathcal{H}_{\q\p}^{\g{k}}}^2=[K_0(\q):K_0(\q\p)].||\varphi||_{ \mathcal{H}_{\q}^{\g{k}}} ^2$ avec 
$$[K_0(\q):K_0(\q\p)]\leq N(\p)+1$$ ce qui donne:
\begin{multline}
\sum_{\varphi}^h\Big|\sum_
{\substack{\overline{\A}\\\nu\in(\A^\times)^{\gg 0}/\unite\\ N(\nu\A)\leq X}}
\lambda(\nu,\A,\varphi)x_{\nu\A}\Big|^2 \\ \leq (N(\p)+1)
\sum_{\varphi\in{\rm{BO}}(\mathcal{H}_{\p\q}^{\g{k}})}^h\Big|\sum_
{\substack{\overline{\A}\\\nu\in(\A^\times)^{\gg 0}/\unite\\ N(\nu\A)\leq X}}
\lambda(\nu,\A,\varphi)x_{\nu\A}\Big|^2\\ \ll (N(\p)+1)\left(1+\frac{X}{N(\q\p)}\right)||x_X||_2^2 \ll
\left(1+\frac{X}{N(\q)}\right)||x_X||_2^2.\nonumber
\end{multline}
\rem : Nous avons en fait d\'emontr\'e une in\'egalit\'e de grand crible pour une base orthogonale de $\mathcal{H}_\q^{\g{k}}$. 
La m\^eme m\'ethode fonctionne dans le cas d'un caract\`ere central non trivial (avec la modification qui s'impose \`a la 
formule de Petersson). On peut \'etudier des in\'egalit\'es de grand crible, pour des familles de repr\'esentations 
plus g\'en\'erales, avec la formule de Kuznetsov. Typiquement, ces in\'egalit\'es donnent des majorations fines du quatri\`eme 
moment non amolli (cf \cite{luo}). En particulier, la m\'ethode pr\'esent\'ee ici nous permet d'obtenir:
$$\sum_{\pi\in\Pi_\q^{\g{k}}}^hL(1/2,\pi)^4\ll\left(\log(N(\q)\right)^6$$
ce qui est le bon ordre de grandeur; pour avoir un \'equivalent, il faudrait avoir les moyens d'\'etudier des termes 
\og non-diagonaux\fg{} (ce qui est a \'et\'e fait pour $F=\Q$ par Kowalski, Michel et VanderKam \cite{KMV2}, inspir\'es par des travaux de Duke, Friedlander et Iwaniec dans des probl\`emes similaires).

\section{Non-annulation de la d\'eriv\'ee}\label{Derivee}
L'objet de cette partie est d'expliquer la preuve du th\'eor\`eme \ref{nonannulationderivee}, 
puisqu'on expliquera dans la section suivante le passage de la moyenne harmonique \`a la moyenne naturelle pour la non-annulation 
des fonctions $L$. On a donc la proposition suivante:

\begin{proposition}
Soit $\g{k}\in 2\N^d_{\geq 1}$. Quand $\q$ parcourt l'ensemble des id\'eaux premiers de $\entier$, on a :
\begin{eqnarray}
\liminf_{N(\q)\to\infty}\sum_{\pi\in\Pi_{\q}^{\g{k}}}^h \mathbbmss{1}_{\varepsilon_\pi=-1,L'(1/2,\pi_f)\neq 0}\geq \frac{7}{16}.
\end{eqnarray}
\end{proposition}

On rappelle que $L(s,\pi_f)$ d\'esigne la fonction $L$ finie, s\'erie de Dirichlet quand $\Re(s)>1$, $\varepsilon_\pi\in
\{1,-1\}$ est le signe de l'\'equation fonctionnelle. En se souvenant que 
$\Lambda(s,\pi)=N(\cond)^{s/2}L(s,\pi_\infty)L(s,\pi_f)$ satisfait l'\'equation fonctionnelle (\ref{equation}), on en d\'eduit que
\begin{eqnarray}
\{\pi\in\Pi_\q^{\g{k}},\varepsilon_\pi=-1{\textrm{ et }}L'(1/2,\pi_f)\neq 0\}=\{\pi\in\Pi_\q^{\g{k}},\Lambda'(1/2,\pi)\neq 0\}.\nonumber
\end{eqnarray}
Les valeurs sp\'eciales $\Lambda'(1/2,\pi)$ et $\Lambda'(1/2,\pi)^2$ s'\'etablissent pareillement \`a (\ref{valeur1}) et 
(\ref{valeur2}). On trouve:
\begin{lemme} Soit $\pi$ une forme modulaire de poids $\g{k}\in 2\N^d_{\geq 1}$, conducteur $\q$. On a :
\begin{multline}
\Lambda'(1/2,\pi)=\frac{(1-\varepsilon_\pi)N(\q)^{1/4}}{2i\pi}\times\\\sum_{\n\subset\entier} \frac{\lambda_\pi(\n)}{\sqrt{N(\n)}} 
\int_{(1)} \frac{d}{ds} \left( \frac{N(\q)^{s/2}}{N(\n)^s}L(s+1/2,\pi_\infty)\right)\frac{ds}{s}\nonumber\\
\Lambda'(1/2,\pi)^2=\frac{2N(\q)^{1/2}}{2i\pi}\sum_{\n,\m\subset\entier}\frac{\lambda_\pi(\n)\lambda_\pi(\m)}{\sqrt{N(\n\m)}}
\int_{(1)}\frac{d}{ds} \left( \frac{N(\q)^{s/2}}{N(\n)^s}L(s+1/2,\pi_\infty)\right)\times \\ 
\frac{d}{ds}\left( \frac{N(\q)^{s/2}}{N(\m)^s}L(s+1/2,\pi_\infty)\right)\frac{ds}{s}.\nonumber
\end{multline}
\end{lemme}
On se donne un amollisseur $\{\M(\pi)\}_{\pi\in\Pi_\q^{\g{k}}}$ de la m\^eme forme que lors des derni\`eres sections, et on \'etudie 
les moments amollis, pour $\g{k}\in\N^d$ pair fix\'e:
\begin{proposition}Soit $\g{k}\in 2\N^d_{\geq 1}$. Pour $\q$ tendant vers l'infini parmi les id\'eaux premiers, 
$\Delta\in]0,1[$ fix\'e tel que $M=N(\q)^{\Delta/2}\notin\N$, on a les asymptotiques:
\begin{multline}
\mathcal{M}_1(\q):=\sum_{\pi\in\Pi_\q^{\g{k}}}^h\Lambda'(1/2,\pi)\M(\pi)\\=\frac{\zeta_F(2)L(1/2,\pi_\infty)N(\q)^{1/4}}
{{\rm{res}}_1\zeta_F}\left(P(1)+\frac{P'(1)}{\Delta}+\mathcal{O}(\log(N(\q)^{-1})\right)
\end{multline}
\begin{multline}
\mathcal{M}_2(\q):=\sum_{\pi\in\Pi_\q^{\g{k}}}^h\Lambda'(1/2,\pi)^2\M(\pi)^2=\frac{2\zeta_F(2)^2L(1/2,\pi_\infty)^2N(\q)^{1/2}}
{({\rm{res}}_1\zeta_F)^2}\times \\ \left(\frac{1}{3\Delta^3}\int_0^1P''(t)^2dt+\frac{P'(1)^2}{\Delta^2}+\frac{2}{\Delta}P'(1)P(1)+P(1)^2
+\mathcal{O}(\log(N(\q)^{-1})\right).
\end{multline}
\end{proposition}

La preuve est  tr\`es proche de l'article \cite{KMV} modulo le passage aux sommes d'id\'eaux et le 
traitement des sommes de Kloosterman.

\subsection{Le premier moment}
Montrons l'estimation concernant $\mathcal{M}_1(\q)$. On \'ecrit donc:
\begin{multline}
\mathcal{M}_1(\q)= N(\q)^{1/4}\sum_{\substack{\n\subset\entier\\N(\m)\leq M}} 
\frac{\mathcal{F}\left(N(\n)/N(\q)^{1/2}\right)}{\sqrt{N(\n)}} 
\frac{\mu(\m)P\left(\frac{\log(M/\m)}{\log(M)}\right)}{\psi(\m)\sqrt{N(\m)}}\times\\
\sum_{\pi\in\Pi_\q^{\g{k}}}^h(1-\varepsilon_\pi)  \lambda_\pi(\n)\lambda_\pi(\m)
\end{multline}
avec $M=N(\q)^{\frac{\Delta}{2}}$, $\Delta\in]0,1[$, et pour $y>0$:
\begin{eqnarray}
\mathcal{F}(y):=\frac{1}{2i\pi}\int_{(1)} \frac{d}{ds} \left( y^{-s}L(s+1/2,\pi_\infty)\right)\frac{ds}{s}.
\end{eqnarray}

Comme dans la section \ref{Moment1}, ce terme se coupe en trois morceaux: un terme diagonal, un terme en sommes de Kloosterman 
(qui se majore comme dans cette section), et un terme en les formes anciennes (qui se traite comme dans l'appendice).\\

 \'Etudions 
donc le terme principal. Il s'\'ecrit:
\begin{multline}
\mathcal{M}_1^{diag}(\q)=N(\q)^{1/4}\sum_{N(\m)\leq M} \frac{ \mathcal{F}\left(N(\m)/N(\q)^{1/2}\right)\mu(\m)
P\left(\frac{\log(M/\m)}{\log(M)}\right)}
{\psi(\m)N(\m)}. 
\end{multline}
Il est facile de v\'erifier que de l'expression d\'efinissant $\mathcal{F}$ seule contribue au terme principal la partie 
$$\int_{(1)}\log(y^{-1})y^{-s}L\left(s+1/2,\pi_\infty\right)\frac{ds}{s}.$$ De plus, on emprunte \`a \cite{KMV} l'\'ecriture:
$$\log(y^{-1})=\int_{\mathscr{C}(\delta)}y^{-z}\frac{dz}{z^2}$$
$\mathscr{C}(\delta)$ \'etant un cercle de centre l'origine, de rayon $\delta>0$ (petit). On a donc:
\begin{multline}
\mathcal{M}_1^{diag}(\q)=\frac{N(\q)^{1/4}}{(2i\pi)^3}\iiint_{(1)\times(1)\times\mathscr{C}(\delta)} N(\q)^{\frac{s+z}{2}}M^tL(s+1/2,\pi_\infty)
\PM(t)\times\\  \left(\sum_{\m\subset\entier}\frac{\mu(\m)}{\psi(\m)N(\m)^{1+s+t+z}}\right) \frac{ds}{s}\frac{dt}{t}\frac{dz}{z^2}
+\mathcal{O}\left(\frac{N(\q)^{1/4}}{\log(N(\q))}\right).\nonumber
\end{multline}
On peut \'ecrire $$\sum_{\m}\frac{\mu(\m)}{\psi(\m)N(\m)^{1+s+t+z}}=\zeta_F(1+s+t+z)^{-1}\eta(s,t,z)$$ avec 
$\eta$ d\'efinie sur $\Re(s+t+z)\geq -1$ et $\eta(0,0,0)=\zeta_F(2)$. En changeant de droite en $s$ de $(1)$ \`a $(-1/2)$, on a :
\begin{multline}
\mathcal{M}_1^{diag}(\q)=\frac{N(\q)^{1/4}L(1/2,\pi_\infty)}{(2i\pi)^2}\iint_{(1)\times\mathscr{C}(\delta)} N(\q)^{\frac{z}{2}}M^t
\PM(t)\times\\  
\zeta_F(1+t+z)^{-1}\eta(0,t,z)\frac{dt}{t}\frac{dz}{z^2}+\mathcal{O}\left(\frac{N(\q)^{1/4}}{\log(N(\q))}\right)
+\mathcal{O}(N(\q)^{1/4-\eta})\nonumber
\end{multline}
pour $\eta>0$ d\'ependant de $\Delta$ et $\delta$. Si $\delta$ est choisi suffisamment petit, l'int\'egrale en $z$ est le r\'esidu en z\'ero:
\begin{multline}
\mathcal{M}_1^{diag}(\q)=\frac{N(\q)^{1/4}L(1/2,\pi_\infty)}{(2i\pi)}\int_{(1)}M^t\PM(t)\times\\{\rm{res}}_{z=0}\left(
 N(\q)^{\frac{z}{2}}  
\frac{\eta(0,t,z)}{\zeta_F(1+t+z)}\right)\frac{dt}{t}+\mathcal{O}\left(N(\q)^{1/4}\log(N(\q))^{-1}\right).\nonumber
\end{multline}
Le r\'esidu vaut: $$\log(N(\q)) \frac{\eta(0,t,0)}{\zeta_F(1+t)}+ \frac{\eta'(O,t,0)}{\zeta_F(1+t)}- \frac{\zeta'_F(1+t)}
{\zeta_F(1+t)^2}\eta(0,t,0)$$ et il est donc clair que le second terme est domin\'e par le premier. Soit:
\begin{multline}
\mathcal{M}_1^{diag}(\q)=\frac{N(\q)^{1/4}\log(N(\q)) L(1/2,\pi_\infty)}{(2i\pi)}\int_{(1)}M^t\PM(t) \frac{\eta(0,t,0)}{\zeta_F(1+t)}
\frac{dt}{t}\\
-\frac{N(\q)^{1/4}L(1/2,\pi_\infty)}{(2i\pi)}\int_{(1)}M^t\PM(t) \frac{\zeta'_F(1+t)}
{\zeta_F(1+t)^2}\eta(0,t,0)\frac{dt}{t}
+\mathcal{O}\left(\frac{N(\q)^{1/4}}{\log(N(\q))}\right)\\
=\frac{N(\q)^{1/4}\zeta_F(2)L(1/2,\pi_\infty)}{{\rm{res}}_1\zeta_F}\left( \frac{\log(N(\q))}{\log(M)}P'(1)+P(1)\right)
+\mathcal{O}\left(\frac{N(\q)^{1/4}}{\log(N(\q))}\right).\nonumber
\end{multline}

\subsection{Le deuxi\`eme moment}
Montrons ici l'asymptotique de $\mathcal{M}_2(\q)$. On a pour $\pi\in\Pi_\q^{\g{k}}$:
\begin{multline}
\Lambda'(1/2,\pi)^2=\frac{2N(\q)^{1/2}}{2i\pi}\sum_{\D,\n}\frac{\chi_\q(\D)}{N(\D)}\frac{\lambda_\pi(\n)}{\sqrt{N(\n)}}\times\\
\sum_{\n_1\n_2=\n}\int_{(1)}\frac{d}{ds} \left( \frac{N(\q)^{s/2}}{N(\D\n_1)^s}L(s+1/2,\pi_\infty)\right) 
\frac{d}{ds}\left( \frac{N(\q)^{s/2}}{N(\D\n_2)^s}L(s+1/2,\pi_\infty)\right)\frac{ds}{s}.\nonumber
\end{multline}
Ici $\chi_\q(\D)$ d\'esigne le caract\`ere principal en $\q$, qui vaut 1 si $(\q,\D)=1$ et 0 si $\q|\D$. 
Avec l'expression de $\M(\pi)^2$ d\'ej\`a vue en section \ref{Moment2}, on a :
\begin{multline}
\mathcal{M}_2(\q)=2N(\q)^{1/2}\sum_{\D,\n}\sum_{\substack{ N(\E)\leq M \\ N(\m_1)\leq M/N(\E)\\N(\m_2)\leq M/N(\E)}}
 \frac{\chi_\q(\D)}{N(\D)N(\n)^{1/2}}\times\\
\sum_{\n_1\n_2=\n}\int_{(1)}\frac{d}{ds} \left( \frac{N(\q)^{s/2}}{N(\D\n_1)^s}L(s+1/2,\pi_\infty)\right) 
\frac{d}{ds}\left( \frac{N(\q)^{s/2}}{N(\D\n_2)^s}L(s+1/2,\pi_\infty)\right)\frac{ds}{s}\\
\frac{\mu(\E\m_1)P\left(\frac{\log(M/N(\E\m_1))}{\log(M)}\right)}{\psi(\E\m_1)N(\E\m_1)^{1/2}}
\frac{\mu(\E\m_2)P\left(\frac{\log(M/N(\E\m_2))}{\log(M)}\right)}{\psi(\E\m_2)N(\E\m_2)^{1/2}}
\sum_{\pi\in\Pi_\q^{\g{k}}}^h\lambda_\pi(\n)\lambda_\pi(\m_1\m_2).\nonumber
\end{multline}
On pose $\m=\m_1\m_2$ pour appliquer la formule de Petersson, dont il sort un terme diagonal que nous allons traiter, et les 
termes en sommes de Kloosterman et en les formes anciennes que nous laisserons de c\^ot\'e pour les m\^emes raisons que ci-dessus.\\

De l'int\'egrale intervenant dans le d\'eveloppement de $\Lambda'(1/2,\pi)^2$, nous ne gardons que la partie qui contribue au terme principal, soit:
\begin{multline}
\mathcal{M}_2^{diag}(\q)=2N(\q)^{1/2}\sum_{\D,\E,\m}
 \frac{\chi_\q(\D)}{N(\D)N(\E)N(\m)}\int_{(1)}\frac{N(\q)^{s}}{N(\D^2\m)^s}L(s+1/2,\pi_\infty)^2 \frac{ds}{s}\\
\sum_{\n_1\n_2=\m} \log\left(\frac{N(\q)^{1/2}}{N(\m_1\D)}\right)\log\left(\frac{N(\q)^{1/2}}{N(\m_2\D)}\right) \times\\
\sum_{\m_1\m_2=\m}\frac{\mu(\E\m_1)\mu(\E\m_2) 
P\left(\frac{\log(M/N(\E\m_1))}{\log(M)}\right) P\left(\frac{\log(M/N(\E\m_2))}{\log(M)}\right)} 
 { \psi(\E\m_1)\psi(\E\m_2) } \nonumber
\end{multline}
et donc, avec les \'ecritures int\'egrales idoines:
\begin{multline}
\frac{2N(\q)^{1/2}}{(2i\pi)^5}\int_{( \frac{1}{2} )\times (\frac{1}{2} )\times (\frac{1}{2} )\times\mathscr{C}(\delta)\times\mathscr{C}(\delta)}
N(\q)^{s+\frac{z_1+z_2}{2}}M^{t_1+t_2}L(s+1/2,\pi_\infty)^2\PM(t_1) \\ \PM(t_2)
\sum_{\D,\E,\m}\frac{\chi_\q(\D)}{N(\D)^{1+2s+z_1+z_2}}\frac{1}{N(\E)^{1+t_1+t_2}N(\m)^{1+s}}
\left(\sum_{\n_1\n_2=\m}\frac{1}{N(\n_1)^{z_1}N(\n_2)^{z_2}}\right)\\
\left(\sum_{\m_1\m_2=\m}\frac{\mu(\m_1\E)}{\psi(\E\m_1)N(\m_1)^{t_1}}\frac{\mu(\m_2\E)}{\psi(\E\m_2)N(\m_2)^{t_2}}\right)
\frac{dsdt_1dt_2dz_1dz_2}{st_1t_2z_1^2z_2^2}.\nonumber
\end{multline}
La somme de Dirichlet multiple portant sur les id\'eaux entiers de $\entier$ peut se factoriser simplement, et un calcul donne :
\begin{multline}
\sum_{\D,\E,\m}\frac{\chi_\q(\D)}{N(\D)^{1+2s+z_1+z_2}}\frac{1}{N(\E)^{1+t_1+t_2}N(\m)^{1+s}}\times\\
\left(\sum_{\m_1\m_2=\m}\frac{\mu(\m_1\E)}{\psi(\E\m_1)N(\m_1)^{t_1}}\frac{\mu(\m_2\E)}{\psi(\E\m_2)N(\m_2)^{t_2}}\right)
\times\left(\sum_{\n_1\n_2=\m}\frac{1}{N(\n_1)^{z_1}N(\n_2)^{z_2}}\right)\\
= \frac{\zeta_F^{(\q)}(1+2s+z_1+z_2)\zeta_F(1+t_1+t_2)}{\prod_{i,j}\zeta_F(1+s+z_i+t_j)}\eta(s,t_1,t_2,z_1,z_2)
\end{multline}
les valeurs 
de $i,j$ \'etant toutes les combinaisons possibles de $1,2$, et  $\eta$ se prolongeant analytiquement sur un voisinage de z\'ero avec $$\eta(0,0,0,0,0)=\zeta_F(2)^2.$$ On change d'abord la droite d'int\'egration de $s$ de $\Re(s)=1/2$ en $\Re(s)=-1/2$ ce qui donne:
\begin{multline}
\mathcal{M}_2^{diag}(\q)=\frac{2L(1/2,\pi_\infty)^2N(\q)^{1/2}}{(2i\pi)^4}\int_{(\frac{1}{2})\times (\frac{1}{2})\times\mathscr{C}(\delta)\times\mathscr{C}(\delta)}
N(\q)^{\frac{z_1+z_2}{2}}M^{t_1+t_2}\times\\ \PM(t_1)\PM(t_2)
\frac{\zeta_F(1+t_1+t_2)\zeta_F^{(\q)}(1+z_1+z_2)}{\prod_{i,j}\zeta_F(1+z_i+t_j)}\eta(0,t_1,t_2,z_1,z_2)
\frac{dt_1dt_2dz_1dz_2}{t_1t_2z_1^2z_2^2}\\+\mathcal{O}(N(\q)^{1/2-\delta}).\nonumber
\end{multline}
Ici comme pr\'ec\'edemment $\mathscr{C}(\delta)$ est un cercle de petit rayon $\delta>0$: on peut alors calculer les int\'egrales 
en $z_1,z_2$, car il n'y a que des r\'esidus en 0. On doit donc calculer:
$$F(t_1,t_2):={\rm{res}}_{(z_1,z_2)=(0,0)}\left(\frac{N(\q)^{\frac{z_1+z_2}{2}}}{z_1^2z_2^2}\times\frac{\zeta_F(1+z_1+z_2)}{\prod_{i,j}\zeta_F(1+z_i+t_j)}
\right).$$
D'abord, le r\'esidu en $z_1=0$ vaut:
\begin{multline}
\frac{\log(N(\q))}{2}\frac{\zeta_F(1+z_2)}{\prod_j\zeta_F(1+t_j)}+\frac{\zeta'_F(1+z_2)}{\prod_j\zeta_F(1+t_j)}\\
-\zeta_F(1+z_2)\frac{\zeta'_F(1+t_1)\zeta_F(1+t_2)+
\zeta'_F(1+t_2)\zeta_F(1+t_1)}{\prod_j\zeta_F(1+t_j)^2}.
\end{multline}
Ce qui donne 
\begin{multline}
F(t_1,t_2)= \frac{\log(N(\q))}{2\prod_j\zeta_F(1+t_j)}\times{\rm{res}}_{z_2=0}\left( \frac{N(\q)^{\frac{z_2}{2}}\zeta_F(1+z_2)}{z_2^2\prod_j\zeta_F(1+z_2+t_j)} \right)\\
+\frac{1}{\prod_j\zeta_F(1+t_j)}\times{\rm{res}}_{z_2=0}\left( \frac{N(\q)^{\frac{z_2}{2}}\zeta'_F(1+z_2)}{z_2^2\prod_j\zeta_F(1+z_2+t_j)} \right)\\
-\frac{\zeta'_F(1+t_1)\zeta_F(1+t_2)+
\zeta'_F(1+t_2)\zeta_F(1+t_1)} {\prod_j\zeta_F(1+t_j)^2}\times\\{\rm{res}}_{z_2=0}\left( \frac{N(\q)^{\frac{z_2}{2}}
\zeta_F(1+z_2)}{z_2^2\prod_j\zeta_F(1+z_2+t_j)} \right).\nonumber
\end{multline}
Les r\'esidus en $z_2$ sont  p\'enibles \`a calculer; cependant une observation simple permet d'\'evacuer 
un bon nombre de termes ne participant pas au terme principal. \'Ecrivons les d\'eveloppements en s\'erie 
enti\`ere $N(\q)^{z/2}=\sum_{k\geq 0}a_kz^k$, $\prod_j\zeta_F(1+t_j+z)^{-1}=\sum_{\ell\geq 0}b_\ell(t_1,t_2) z^\ell$, et 
$\zeta_F(1+z)=\sum_{m\geq -1}c_mz^m$ et regardons par exemple: 
\begin{multline}
{\rm{res}}_{z=0}\left(\frac{N(\q)^{\frac{z}{2}}\zeta_F(1+z)}
{z^2\prod_j\zeta_F(1+t_j+z)}\right)=c_{-1}(a_2b_0+a_1b_1+a_0b_2)+c_0(a_1b_0+a_0b_1)\\ +c_1a_0b_0.\nonumber
\end{multline}
 Le terme $a_k$ vaut 
$2^{-k}\log(N(\q))^{k}$, en arrangeant les termes en: $b_0(c_{-1}a_2+c_0a_1+c_1a_0)+b_1\times\ldots$, on voit que seul le terme 
produit par $c_{-1}$ contribue. Ainsi, on peut remplacer $\zeta_F(1+z)$ (resp. $\zeta_F'(1+z)$) par ${\rm{res}}_1(\zeta_F)z^{-1}$
 (resp. $-{\rm{res}}_1(\zeta_F)z^{-2}$). \`A cause de l'estimation $$\int_{(1)}M^t\PM(t)t^n\frac{dt}{t}\sim (\log(M))^{-n}$$ on peut aussi remplacer $\zeta_F(1+t_1+t_2)$ par ${\rm{res}}_1(\zeta_F)(t_1+t_2)^{-1}$ (resp. 
 $\zeta_F(1+t_j+z_2)$ par ${\rm{res}}_1(\zeta_F)(t_j+z_2)^{-1}$ et $\zeta_F(1+t_j)$ par ${\rm{res}}_1(\zeta_F)t_j^{-1}$). 
Cela se r\'esume par la relation:
\begin{multline}
\mathcal{M}_2^{diag}(\q)=\frac{2L(1/2,\pi_\infty)^2N(\q)^{1/2}}{(2i\pi)^2{\rm{res}}_1(\zeta_F)^2}\int_{(\frac{1}{2})\times (\frac{1}{2})}
\frac{M^{t_1+t_2} \PM(t_1)\PM(t_2)}{t_1+t_2}\times\\ \eta(0,t_1,t_2,0,0)
\bigg\{ \frac{ \log(N(\q))t_1t_2}{2} {\rm{res}}_{z=0} \left( \frac{ N(\q)^{\frac{z}{2}}(z+t_1)(z+t_2)}{z^3}\right)\\
-t_1t_2{\rm{res}}_{z=0}\left( \frac{ N(\q)^{\frac{z}{2}}(z+t_1)(z+t_2)}{z^4} \right) \\
+(t_1+t_2){\rm{res}}_{z=0}\left( \frac{N(\q)^{\frac{z}{2}}(z+t_1)(z+t_2)}{z^3}\right)\bigg\}
\frac{dt_1dt_2}{t_1t_2}.
\nonumber
\end{multline}
Tous calculs faits, cela donne (en r\'eunissant les premier et troisi\`eme termes):
\begin{multline}
\mathcal{M}_2^{diag}(\q)=\\ \frac{2L(1/2,\pi_\infty)^2N(\q)^{1/2}}{(2i\pi)^2{\rm{res}}_1(\zeta_F)^2}\int_{(\frac{1}{2})\times (\frac{1}{2})}
\frac{M^{t_1+t_2} \PM(t_1)\PM(t_2)}{t_1+t_2}\eta(0,t_1,t_2,0,0)\\
\bigg\{\left( \frac{\log(N(\q))^2}{8}t_1t_2+\frac{\log(N(\q))}{2}(t_1+t_2)+1\right)\left(\frac{\log(N(\q))}{2}t_1t_2+t_1+t_2\right)
\\-\left( \frac{\log(N(\q))^3}{8.3!}t_1t_2+\frac{\log(N(\q))^2}{2^2.2}(t_1+t_2)+\frac{\log(N(\q))}{2}\right)t_1t_2\bigg\}\frac{dt_1dt_2}{t_1t_2}
\\+\mathcal{O}\left(\frac{N(\q)^{1/2}}{\log(N(\q))}\right).\nonumber
\end{multline}
Avec les valeurs des int\'egrales ci-dessous d\'ej\`a rencontr\'ees lors de la section \ref{Moment2} ($n$ est un entier positif):
\begin{eqnarray}{}
\frac{1}{(2i\pi)^2}\iint_{(1),(1)}M^{t_1+t_2}\PM(t_1)\PM(t_2)\frac{t_1t_2}{t_1+t_2}dt_1dt_2 &=& \frac{1}{\log(M)^3}\int_0^1P''(u)^2du\nonumber\\
\frac{1}{2i\pi}\int_{(1)}M^{t}\PM(t)t^n\frac{dt}{t} &=& (\log(M))^{-n}P^{(n)}(1)\nonumber
\end{eqnarray}
on trouve, apr\`es des calculs \'el\'ementaires, l'expression
\begin{multline}
\mathcal{M}_2^{diag}(\q)=\frac{2L(1/2,\pi_\infty)^2\zeta_F(2)^2N(\q)^{\frac{1}{2}}}{{\rm{res}}_1(\zeta_F)^2}\times\\
\bigg( \frac{1}{3\Delta^3}\int_0^1P''(u)^2du+\frac{1}{\Delta^2}P'(1)^2+ \frac{2}{\Delta}P'(1)P(1)+P(1)^2\bigg)
+\mathcal{O}\left(\frac{N(\q)^{1/2}}{\log(N(\q))}\right).
\end{multline}

\subsection{Optimisation}
En faisant tendre $\Delta$ vers 1, on trouve donc que
$$\liminf_{N(\q)\to\infty}\sum^h_{\pi\in\Pi_\q^{\g{k}}}\mathbbmss{1}_{\Lambda'(1/2,\pi)\neq 0}\geq
\frac{ \big(P(1)+P'(1)\big)^2}
{2\left(\frac{1}{3}\int_0^1P''(t)^2dt+P'(1)^2+2P'(1)P(1)+P(1)^2\right)}.$$
Il s'agit d'optimiser cette fonction sur tous les polyn\^omes $P$ satisfaisant \`a $P(0)=P'(0)=0$, ce qui est un exercice 
facile d'analyse fonctionnelle. D'abord, l'existence d'une fonction maximisante est garantie par le
\begin{lemme}
Soit $\mathfrak{H}$ un espace de Hilbert, $\mathscr{L}$ une forme lin\'eaire continue et $\mathscr{Q}$ une forme quadratique 
continue et d\'efinie-positive. Alors la fonction
\begin{displaymath}
\varphi:  \begin{array}{ll}\mathfrak{H}\setminus \{0\}\longrightarrow  \R \\ h \longmapsto  \frac{\mathscr{L}(h)^2}{\mathscr{Q}(h)+\mathscr{L}(h)^2}
\end{array}
\end{displaymath}
est born\'ee et atteint ses bornes.
\end{lemme}
\preuve: Visiblement, $0\leq \varphi\leq 1$. On a $\varphi(\lambda h)=\varphi(h)$ pour tout $\lambda$ non nul: ainsi 
$\displaystyle{\sup_{\mathfrak{H}\setminus \{0\}}\varphi=\sup_{\mathbb{S}}\varphi}$ avec $\mathbb{S}=\{h\in\mathfrak{H},||h||=1\}$. 
Soit $h_n\in\mathbb{S}$ une suite maximisante (telle que $\varphi(h_n)\to\sup\varphi$). \`A extraction pr\`es, on peut supposer que 
$h_n$ tend faiblement vers $h$, avec $||h||\leq 1$. Par convexit\'e $\mathscr{Q}(h)\leq \liminf \mathscr{Q}(h_n)$, et donc 
puisque $\mathscr{L}(h_n)\to\mathscr{L}(h)$, $\varphi(h)\geq \lim \varphi(h_n)$. $\blacksquare$\\

On se place donc dans l'espace de Sobolev $$H^2(0,1)=\{f\in L^2(0,1), f',f''\in L^2(0,1)\}$$ muni de sa norme hilbertienne $$\|f\|^2=\|f\|^2+\|f'\|^2+\|f''\|^2$$ et soit $$\mathfrak{H}=
\{f\in H^2(0,1),f(0)=f'(0)=0\}\subset H^2(0,1)$$ muni de la topologie induite. Le lemme pr\'ec\'edent s'applique avec $\mathscr{L}(f)=f(1)+f'(1)$ et 
$\mathscr{Q}(f)=\int_0^1f''(t)^2dt$: il existe donc une $f\in\mathfrak{H}$ maximisant notre fonctionnelle. Calculons maintenant $f$: on va voir que c'est 
en fait un polyn\^ome.

Remarquons d'abord l'\'equivalence:
$$f \textrm{ maximise } \frac{\mathscr{L}(.)^2}{\mathscr{Q}(.)+\mathscr{L}(.)^2} \Longleftrightarrow f \textrm{ minimise }\frac{ \mathscr{Q}(.)}{\mathscr{L}(.)^2 }.$$
Puis, $f$  doit annuler la diff\'rentielle de $\frac{\mathscr{Q}(.)}{\mathscr{L}(.)^2}$, soit:
$$\forall h\in\mathfrak{H}, \left(f(1)+f'(1)\right)^2\int_0^1f''h''- (f(1)+f'(1))(h(1)+h'(1))\int_0^1f''^2 =0$$
on doit avoir $f(1)+f'(1)\neq 0$, sans quoi $\mathscr{L}(f)=0$,  soit:
$$\forall h\in\mathfrak{H}, \left(f(1)+f'(1)\right)\int_0^1f''h''- (h(1)+h'(1))\int_0^1f''^2 =0$$
une int\'egration par parties donne:
\begin{multline}
\forall h\in\mathfrak{H},\left( (f(1)+f'(1))f''(1)-\int_0^1f''^2\right)h'(1) \\=\int_0^1\left( (f(1)+f'(1))f^{(3)}(t)
+\int_0^1f''^2\right)h'(t)dt     \nonumber
\end{multline}
Il est bien connu qu'une telle \'egalit\'e ne peut qu'\^etre triviale, i.e. qu'elle implique les deux \'equations:
\begin{displaymath}
\left\{\begin{array}{cc} (f(1)+f'(1))f^{(3)}(t)+\int_0^1f''^2=0\\(f(1)+f'(1))f''(1)-\int_0^1f''^2=0\end{array}\right.
\end{displaymath}
La premi\`ere dit que $f$ est un polyn\^ome de degr\'e 3, donc de la forme $aX^3+X^2$ (par homog\'en\'eit\'e de la fonctionnelle) 
et en ajoutant les deux \'egalit\'es, retirant le cas $f'(1)+f(1)=0$ encore une fois, on trouve que $a=-1/6$ convient. Le 
polyn\^ome $X^2-X^3/6$ est donc le polyn\^ome optimal recherch\'e, pour lequel
$$\frac{ \big(P(1)+P'(1)\big)^2}
{2\left(\frac{1}{3}\int_0^1P''(t)^2dt+P'(1)^2+2P'(1)P(1)+P(1)^2\right)}=\frac{7}{16}$$
ce qui \'etait annonc\'e.

\section{De la moyenne harmonique \`a la moyenne naturelle}\label{Naturel}
Le but de cette section est d'expliquer la preuve des th\'eor\`emes \ref{1} et \ref{2}. Nous ne traiterons que le premier, car 
le second est analogue. Il serait int\'eressant de le d\'eduire directement du th\'eor\`eme \ref{nonannulation}, mais cette 
approche n\'ecessiterait des r\'esultats de nature probabiliste encore inconnus: nous allons donc en fait reprendre l'\'etude 
des moments amollis, en incorporant dans l'amollisseur le poids harmonique, ce qui a comme cons\'equence de compliquer 
consid\'erablement les calculs (et c'est pour cela que l'on a commenc\'e par les moments harmoniques), mais on s'apercevra 
que les modifications ne portent que sur des termes d'ordre 2, laissant ainsi les termes principaux inchang\'es \`a une constante absolue pr\`es. Cette section 
a \'et\'e incorpor\'ee car c'est la non-annulation en moyenne naturelle qui int\'eresse le plus les g\'eom\`etres (pour toute \'etude 
du rang des vari\'et\'es ab\'eliennes modulaires).

Avec les notations pr\'ec\'edentes, nous noterons les moments naturels amollis pour $i=1,2$ (toujours pour $\q$ maximal):
$$\Mo_i(\q)=\frac{\sum_{\pi\in\Pi_\q^{\g{k}}}\Lambda(1/2,\pi)^i\M(\pi)^i}{|\Pi_\q^{\g{k}}|}$$
et donc
$$\Mo_i(\q)=\frac{\sum_{\pi\in\Pi_\q^{\g{k}}}^h\omega_\pi^{-1}\Lambda(1/2,\pi)^i\M(\pi)^i}
{\sum_{\pi\in\Pi_\q^{\g{k}}}^h\omega_\pi^{-1}}$$
avec, on le rappelle, $\omega_\pi=\frac{\Gamma(\g{k}-\g{1})}{(4\pi)^{\g{k}-\g{1}}|\disc|^{1/2}||\varphi_\pi||_{\mathcal{H}_\q^{\g{k}}}^2}$. 
Se souvenant que 
\begin{multline}
||\varphi_\pi||_{\mathcal{H}_\q^{\g{k}}}^2=\int_{Z(\adele)\gl(F)\backslash \gl(\adele)/K_0(\q)}|\varphi_\pi
(g)|^2dg\\=(N(\q)+1)\int_{Z(\adele)\gl(F)\backslash \gl(\adele)}|\varphi_\pi(g)|^2dg.\nonumber
\end{multline}
On a donc: $$\Mo_i(\q)=\frac{\sum_{\pi\in\Pi_\q^{\g{k}}}^h||\varphi_\pi||^2\Lambda(1/2,\pi)^i\M(\pi)^i}
{\sum_{\pi\in\Pi_\q^{\g{k}}}^h||\varphi_\pi||^2}$$
On va donc chercher une expression de $$||\varphi_\pi||^2=\int_{Z(\adele)\gl(F)\backslash \gl(\adele)}|\varphi_\pi(g)|^2dg$$ 
en termes des coefficients de Fourier de $\pi$: ceci est connu pour $F=\Q$, c'est la formule de Shimura. 

Ce qui suit est divis\'e en deux parties: la premi\`ere est d\'evou\'ee \`a des calculs locaux pour le carr\'e sym\'etrique et 
la convolution de Rankin-Selberg, n\'ecessaires aux \'enonc\'e et preuve de la formule de Shimura. On en d\'eduit au 
passage une estimation du nombre de formes de conducteur $\q$ (qui s'adapterait d'ailleurs pour $\q$ sans facteurs carr\'es). 
Ensuite, on revient \`a notre probl\`eme sur les moments naturels, et on y prouve les estimations n\'ecessaires.

\subsection{Le carr\'e sym\'etrique et fonctions $L$ de Rankin-Selberg}
Nous supposons ici que $F$ est totalement r\'eel, et $\q$ id\'eal maximal; cependant bon nombre des rappels qui suivent 
sont valables dans une plus grande g\'en\'eralit\'e.\\
Soit $\pi_1,\pi_2$ deux repr\'esentations automorphes paraboliques, que nous supposerons de caract\`ere central trivial (l\`a 
encore parce que nous ne nous int\'eressons qu'\`a ce cas). Gr\^ace aux travaux de Jacquet (cf. le livre de Jacquet \cite{J}, et les notes de Cogdell \cite{C}), on sait d\'efinir 
la fonction $L$ de Rankin-Selberg de $\pi_1,\pi_2$, not\'ee $L(s,\pi_1\times\pi_2)$, par un produit eul\'erien convergent 
pour $\Re(s)>1$, et admettant un pole en $s=1$ si et seulement si $\pi_1\cong \check{\pi}_2$ (contragr\'ediente de $\pi_2$). 
Depuis peu, on sait gr\^ace \`a Ramakrishnan \cite{R} que cette fonction $L$ est celle d'une repr\'esentation automorphe de ${\rm{GL}}(4)$, 
not\'ee $\pi_1\boxtimes\pi_2$. 

Pour $j\in\{1,2\}$ et chaque place finie $\p$
notons $\alpha_{\pi,j}(\p)$ les param\`etres de Langlands, i.e. tels que:
$$L(s,\pi)=(1-\alpha_{\pi,1}(\p)N(\p)^{-s})^{-1}(1-\alpha_{\pi,2}(\p)N(\p)^{-s})^{-1}$$
alors, en toute place non ramifi\'ee, on a:
$$L(s,\pi_1\times\pi_2)=\prod_{i,j}(1-\alpha_{\pi_2,i}(\p)\alpha_{\pi_1,j}(\p)N(\p)^{-s})^{-1}$$
Une propri\'et\'e importante de la convolution de Rankin-Selberg est d'avoir une repr\'esentation int\'egrale \`a l'aide 
de s\'eries d'Eisenstein. Rappelons-en la d\'efinition en laissant le lecteur consulter l'article \cite{JZ} pour plus 
de d\'etails. 

Soit donc $\Phi=\otimes_v\Phi_v\in\mathcal{S}(\adele^2)$ une fonction de la classe de Schwartz (la partie archim\'edienne 
$\Phi_\infty$ est une fonction de la classe de Schwartz au sens classique, la partie finie $\Phi_f$ est localement constante, 
valant $\mathbbmss{1}_{\mathcal{O}_v^2}$ en presque toute place). On d\'efinit
$$f(g,\Phi,s)=|{\rm{det}}(g)|_{\adele}^s\int_{\adele^\times}\Phi\left( (0,t)g\right) |t|^{2s}d^\times t$$
et la s\'erie d'Eisentein correspondante, absolument convergente pour $\Re(s)>1$:
\begin{multline}
E(g,\Phi,s)=\sum_{\gamma\in P(F)\backslash \gl(F)}f(\gamma g,\Phi,s)=|{\rm{det}}g|_{\adele}^s\int_{F^\times \setminus \adele^{\times}}
\sum_{\xi\in F^2\setminus \{0\}}\Phi(\xi tg)|t|^{2s}d^\times t\nonumber
\end{multline}
$E$ satisfait \`a une \'equation fonctionnelle, et se prolonge m\'eromorphiquement \`a tout le plan complexe, mais ce qui 
est important pour nous c'est qu'elle admet un p\^ole en 1 et en 0, avec un r\'esidu donn\'e par:
\begin{eqnarray}
{\rm{res}}_{s=1}E(g,\Phi,s)=\frac{\hat{\Phi}(0)}{2}=\frac{1}{2}\int_{\adele^2}\Phi(x,y)dxdy.
\end{eqnarray}
Ceci assure donc que pour toute forme automorphe parabolique $\varphi$
\begin{eqnarray}
||\varphi||^2=\frac{2}{\hat{\Phi}(0)}{\rm{res}}_{s=1}\left(\int_{Z(\adele)\gl(F)\backslash\gl(\adele)}|\varphi(g)|^2E(g,\Phi,s)dg\right).
\end{eqnarray}
D'autre part, quand $\varphi$ est un vecteur sp\'ecial, l'int\'egrale de droite est li\'ee \`a la fonction $L$ de 
Rankin-Selberg $L(s,\pi\times\pi)$: ceci est r\'esum\'e dans les lemmes suivants.

\begin{lemme}
Soient $\pi_1,\pi_2$ repr\'esentations automorphes paraboliques de caract\`ere central trivial, et $\mathcal{W}(\pi_j,\psi)$ 
leur mod\`ele de Whittaker. Pour $j=1,2$ soit $\varphi_{\pi_j}$ leur vecteur sp\'ecial (et $W_{\pi_j}$ leur fonction de Whittaker 
respective). On a alors, pour $\Re(s)>1$ et $\epsilon:=\left(\begin{array}{cc}-1 & 0\\0 & 1\\ \end{array}\right) $:
\begin{multline}
\int_{Z(\adele)\gl(F)\backslash\gl(\adele)}\varphi_{\pi_1}(g)\varphi_{\pi_2}(g)E(g,\Phi,s)dg\\=\int_{Z(\adele)N(\adele)\backslash \gl(\adele)}
W_{\pi_1}(g)W_{\pi_2}\left(\epsilon g\right)f(g,\Phi,s)dg.
\end{multline}
\end{lemme}
C'est ce qu'on nomme ``m\'ethode'' de Rankin-Selberg (ou unfolding). Cf \cite{Bu}, proposition (3.8.2). 

\begin{lemme}\label{Phi}
Avec les m\^emes hypoth\`eses, en supposant que le caract\`ere additif local $\psi_v$ est non ramifi\'e: si $\pi_1,\pi_2$ 
sont non ramifi\'ees en $v$, en prenant $\Phi_v=\mathbbmss{1}_{\mathcal{O}_v^2}$ on a:
\begin{eqnarray}
\int_{Z(F_v)N(F_v)\backslash\gl(F_v)}W_{\pi_{1,v}}(g_v)W_{\pi_{2,v}}(\epsilon g_v) f_v(g_v,\Phi_v,s)dg_v=L(s,\pi_{1,v}\times\pi_{2,v}).\nonumber
\end{eqnarray}
Si $v\leftrightarrow \q$ et $\pi_{v}$ a pour conducteur $\q\mathcal{O}_v$ (auquel cas $\pi_v$ est sp\'eciale), le choix 
$\Phi_v=\mathbbmss{1}_{\q\mathcal{O}_v\times\mathcal{O}_v^\times}$ donne:
\begin{eqnarray}
\int_{Z(F_v)N(F_v)\backslash\gl(F_v)}W_{\pi_{v}}(g_v)W_{\pi_{v}}(\epsilon g_v) f_v(g_v,\Phi_v,s)dg_v
=\frac{L(s,\pi_v\times\pi_v)}{N(\q)+1}.\nonumber
\end{eqnarray}
\end{lemme}
\preuve : La prem\`ere \'egalit\'e est connue: voir par exemple \cite{Bu}, chapitre 3.8, ou \cite{Zh}, proposition 2.5.1. Pour la 
seconde, on part de l'\'egalit\'e:
$$\int_{F_v^\times} W_{\pi_{v}}\left( \left(\begin{array}{cc}a & 0\\0 & 1\end{array}\right) \right)
W_{\pi_{v}}\left( \left(\begin{array}{cc}-a & 0\\0 & 1\end{array}\right) \right)|a|^{s-1}d^\times a=L(s,\pi_{v}\times\pi_{v})$$
qui tient au fait que $\pi$ a un conducteur en $v$ d'exposant 1, et $\psi_v$ non ramifi\'e: voir \cite{popa2}. On constate 
alors qu'avec le choix $\Phi_v=\mathbbmss{1}_{\q\mathcal{O}_v\times\mathcal{O}_v^\times}$, on a 
\begin{multline}
f_v\left( \left(\begin{array}{cc}a_v & 0\\0 & 1\end{array}\right) k_v ,\Phi_v,s\right):=
|a_v|_v^{s}\int_{F_v^\times}\Phi_v\left((0,t_v)\left(\begin{array}{cc}a_v & 0\\0 & 1\end{array}\right) k_v \right)|t_v|_v^{2s}d^\times t_v
\\= |a_v|_v^{s} \mathbbmss{1}_{K_0(\q\mathcal{O}_v)}(k_v)\nonumber
\end{multline}
En outre la d\'ecomposition d'Iwasawa $\gl(F_v)=Z(F_v)N(F_v)A(F_v)\gl(\mathcal{O}_v)$ assure l'\'egalit\'e: 
$$\int_{Z(F_v)N(F_v)\backslash\gl(F_v)}W(g_v)dg_v=\int_{F_v^\times\times\gl(\mathcal{O}_v)}W\left( \left(\begin{array}{cc}a_v & 0\\0 & 1\end{array}\right) 
k_v\right)\frac{d^\times a_v}{|a_v|_v}dk_v$$
pour toute fonction int\'egrable $W$. On a donc:
\begin{multline}
\int_{Z(F_v)N(F_v)\backslash\gl(F_v)}W_{\pi_{v}}(g_v)W_{\pi_{v}}(\epsilon g_v) f_v(g_v,\Phi_v,s)dg_v
\\=\int_{F_v^\times\times\gl(\mathcal{O}_v)} W_{\pi_{v}}\left( \left(\begin{array}{cc}a & 0\\0 & 1\end{array}\right)k_v \right)
W_{\pi_{v}}\left( \left(\begin{array}{cc}-a & 0\\0 & 1\end{array}\right)k_v \right)\times\\|a|^{s}
\mathbbmss{1}_{K_0(\q\mathcal{O}_v)}(k_v)\frac{d^\times a}{|a_v|_v}dk_v
=L(s,\pi_v\times\pi_v)\times{\rm{vol}}(K_0(\q\mathcal{O}_v))\nonumber
\end{multline}
et ${\rm{vol}}(K_0(\q\mathcal{O}_v))=[\gl(\mathcal{O}_v):K_0(\q\mathcal{O}_v)]^{-1}=(N(\q)+1)^{-1}$. $\blacksquare$\\

\rem: Il existe aussi aux places infinies une fonction $\Phi_\infty$ qui a les m\^emes propri\'et\'es. Comme dans ce qui suit, 
le type \`a l'infini est fix\'e, la valeur de cette fonction ne nous importe pas.\\

En conclusion, avec le choix global $\Phi_{\g{k},\q}=\otimes_v\Phi_v$, nous avons montr\'e pour $\Re(s)>1$:

\begin{multline}
|\disc|^{s-1/2}\frac{L(s,\pi\times\pi)}{N(\q)+1}=\int_{Z(\adele)\gl(F)\backslash\gl(\adele)}\varphi_{\pi}(g)^2E(g,\Phi_{\g{k},\q},s)dg\\=\int_{Z(\adele)N(\adele)\backslash \gl(\adele)}
W_{\pi}(g)W_{\pi}\left(\epsilon g\right)f(g,\Phi_{\g{k},\q},s)dg\\
=\int_{\adele^\times\times K} W_\pi\left(\left(\begin{array}{cc}a & 0\\0 & 1\end{array}\right) k\right)
W_\pi\left(\left(\begin{array}{cc}-a & 0\\0 & 1\end{array}\right)k \right)|a|^s\mathbbmss{1}_{K_0(\q)}(k)\frac{d^\times a}{|a|}dk\\
=\int_{\adele^\times} W_\pi\left(\left(\begin{array}{cc}a & 0\\0 & 1\end{array}\right)\right)
\overline{W_\pi\left(\left(\begin{array}{cc}-a & 0\\0 & 1\end{array}\right)\right)}|a|^s\frac{d^\times a}{|a|}
{\textrm{  ($W_\pi\left(\left(\begin{array}{cc}a & 0\\0 & 1\end{array}\right)\right)\in\R$) }}\\
=\int_{Z(\adele)\gl(F)\backslash\gl(\adele)}|\varphi_{\pi}(g)|^2E(g,\Phi_{\g{k},\q},s)dg\nonumber
\end{multline}
(le facteur $|\disc|^{s-1/2}$
vient de ce que le caract\`ere additif global est ramifi\'e en la diff\'erente), ce qui assure que: 
\begin{eqnarray}
\frac{\widehat{\Phi}_{\g{k},\q}(0)}{2}||\varphi_\pi||^2=\frac{|\disc|^{1/2}}{N(\q)+1}{\rm{res}}_{s=1}\{L(s,\pi\times\pi)\}.
\end{eqnarray}
D'autre part, gr\^ ace \`a Gelbart-Jacquet \cite{GJ}, on sait qu'il existe une repr\'esentation automorphe sur ${\rm{GL}}_3$, not\'ee 
${\rm{sym}}^2\pi$, telle que 
$$L(s,\pi\times\pi)=L(s,1)L(s,{\rm{sym}}^2\pi)$$
$L(s,1)$ \'etant la fonction z\^eta de Dedekind compl\'et\'ee par les facteurs \`a l'infini. La formule de Shimura 
s'\'enonce ainsi, avec $\Phi_{\g{k},\q}$ la fonction d\'ecrite ci-dessus:
\begin{proposition}
Soit $\pi$ une forme modulaire de poids $\g{k}\in 2\N^d_{\geq 1}$ et conducteur $\q$, de caract\`ere central trivial,  $\varphi_\pi$ son vecteur 
nouveau. On a la relation:
\begin{eqnarray}
||\varphi_\pi||^2=\frac{2|\disc|^{1/2}}{\widehat{\Phi}_{\g{k},\q}(0)(N(\q)+1)}L(1,{\rm{sym}}^2\pi)\res_{s=1}L(s,1).
\end{eqnarray}
\end{proposition}
Dans le cas d'une forme modulaire de Hilbert de poids $\g{k}$ et conducteur $\q$, la fonction $\Phi_{\g{k},\q}=\Phi_\infty\otimes\Phi_f$ 
ne d\'epend aux places archim\'ediennes que du poids et aux places finies que du conducteur: cette fonction est donc la m\^eme 
pour toute forme $\pi\in\Pi_\q^{\g{k}}$, et on note $\Phi_\infty=\Phi_{\g{k}},\Phi_f=\Phi_{\q}$. En r\'esum\'e, pour $i=1,2$:
\begin{eqnarray}
\Mo_i(\q)=\frac{\sum_{\pi\in\Pi_\q^{\g{k}}}^h L(1,{\rm{sym}}^2\pi)\Lambda(1/2,\pi)^i\M(\pi)^i}
{\sum_{\pi\in\Pi_\q^{\g{k}}}^h L(1,{\rm{sym}}^2\pi)}.
\end{eqnarray}
D\'ecrivons plus pr\'ecis\'ement la fonction $L$ du carr\'e sym\'etrique d'une forme modulaire de Hilbert de poids 
$\g{k}\in\N^d$ et de conducteur $\cond=\q$. On a $L(s,{\rm{sym}}^2\pi)=L(s,{\rm{sym}}^2\pi_\infty)
L(s,{\rm{sym}}^2\pi_f)$ avec 
$$L(s,{\rm{sym}}^2\pi_\infty)=\pi^{-3ds/2} \Gamma\left( \frac{s+1}{2} \right) ^d\Gamma\left( \frac{\g{s}+\g{k}}{2} \right) 
\Gamma\left( \frac{\g{s}+\g{k}-\g{1}}{2} \right) $$
la partie finie est une s\'erie de Dirichlet convergeant pour $\Re(s)>1$:
$$L(s,{\rm{sym}}^2\pi_f)=\sum_{\n\subset \entier}\rho_\pi(\n)N(\n)^{-s}$$
avec $$\rho_\pi(\n)=\sum_{\substack{\E^2\f=\n\\(\E,\q)=\entier}}\lambda_\pi(\f^2).$$
D'autre part, outre le produit eul\'erien dont on ne se servira pas, on a une \'equation fonctionnelle. En fait, 
$\Lambda(s,{\rm{sym}}^2\pi):=N(\cond)^sL(s,{\rm{sym}}^2\pi)$ v\'erifie:
$$\Lambda(s,{\rm{sym}}^2\pi)=\Lambda(1-s,{\rm{sym}}^2\pi).$$
La valeur explicite de $L(1,{\rm{sym}}^2\pi)$ se calcule par la m\^eme m\'ethode que $L(1/2,\pi)$. L'\'equation fonctionnelle donne ici:
\begin{multline}\label{val1}
2i\pi L(1,{\rm{sym}}^2\pi)=\sum_{\n\subset \entier}\frac{\rho_\pi(\n)}{N(\n)}\int_{(\delta)} \left(\frac{N(\q)}{N(\n)}\right)^s
L(s+1,{\rm{sym}}\pi_\infty)\frac{ds}{s}+ \\ \sum_{\n\subset \entier}\frac{\rho_\pi(\n)}{N(\n)}\int_{(\delta)} \left(\frac{N(\q)}
{N(\n)}\right)^sL(s+1,{\rm{sym}}\pi_\infty)\frac{ds}{s+1}
\end{multline}
pour tout $\delta>0$ fix\'e. Ceci permet de voir que $L(1,{\rm{sym}}^2\pi)\ll_\varepsilon N(\q)^\varepsilon$ pour tout $\varepsilon>0$ (voir \cite{M}, page 26, (1.24) et \cite{Mo}, th\'eor\`eme 4). On en d\'eduit une formule asymptotique pour le 
nombre de repr\'esentations modulaires de Hilbert de poids $\g{k}$ et conducteur $\q$ premier:

\begin{proposition}\label{nombre}
Pour $\g{k}\in 2\N^d_{\geq 1}$ fix\'e, quand $N(\q)$ tend vers l'infini parmi les id\'eaux maximaux de $\entier$, on a:
\begin{eqnarray}\label{s}
\sum_{\pi\in\Pi_\q^{\g{k}}}^hL(1,{\rm{sym}}^2\pi)\,\,
_{ \overrightarrow{N(\q)\to\infty} }\,\,
\zeta_F(2)L(1,{\rm{sym}}^2\pi_\infty).
\end{eqnarray} 
Par cons\'equent, $\Phi_{\g{k}}$ \'etant la partie archim\'edienne de $\Phi$ (cf la remarque suivant le lemme \ref{Phi}), 
ne d\'ependant que de $\g{k}$:
\begin{eqnarray}
|\Pi_\q^{\g{k}}|\Equi{N(\q)}{\infty}\frac{2(4\pi)^{\g{k}-\g{1}}|\disc|^2{\rm{res}}_{s=1}L(s,1)L(1,{\rm{sym}}^2\pi_\infty)\zeta_F(2)}
{\widehat{\Phi}_{\g{k}}(0)\Gamma(\g{k}-\g{1})}N(\q).
\end{eqnarray}
\end{proposition}
\preuve : Admettant la premi\`ere assertion un instant, on \'ecrit:
\begin{multline}
|\Pi_\q^{\g{k}}|=\sum_{\pi\in\Pi_\q^{\g{k}}}\omega_\pi\omega_\pi^{-1}=\sum^h_{\pi\in\Pi_\q^{\g{k}}}\omega_\pi^{-1}=
\frac{(4\pi)^{\g{k}-\g{1}}|\disc|^{1/2}}{\Gamma(\g{k}-\g{1})}\sum^h_{\pi\in\Pi_\q^{\g{k}}}||\varphi_\pi||_{\mathcal{H}_\q^{\g{k}}}^2
\\=\frac{(4\pi)^{\g{k}-\g{1}}|\disc|^{1/2}}{\Gamma(\g{k}-\g{1})}\sum^h_{\pi\in\Pi_\q^{\g{k}}}||\varphi_\pi||^2(N(\q)+1)
\\=\frac{2(4\pi)^{\g{k}-\g{1}}|\disc|{\rm{res}}_{s=1}L(s,1)}{\Gamma(\g{k}-\g{1})\widehat{\Phi}_{\g{k},\q}(0)(N(\q)+1)}\sum^h_{\pi\in\Pi_\q^{\g{k}}}L(1,{\rm{sym}}^2\pi)(N(\q)+1).
\nonumber\end{multline}
La valeur $\widehat{\Phi}_{\g{k},\q}(0)$ se calcule ais\'ement, puisque $\Phi_{\g{k},\q}=\Phi_{\g{k}}\otimes\bigotimes_{v<\infty}\Phi_v$:
\begin{multline}
\widehat{\Phi}_{\g{k},\q}(0)=\int_{\adele^2}\Phi(x,y)dxdy=\int_{F_\infty^2}\Phi_{\g{k}}(x_\infty,y_\infty)dx_\infty dy_\infty\times\\
\prod_{v<\infty}\int_{F_v^2}\Phi_v(x_v,y_v)dx_vdy_v
\end{multline}
si $v\neq \q, \Phi_v=\mathbbmss{1}_{\mathcal{O}_v^2}$ soit $\widehat{\Phi}_v(0)={\rm{vol}}(\mathcal{O}_v)^2=
|N(\mathfrak{D}_{F_v})|^{-1}$ (avec la normalisation de Tate). \\
Si $v=\q$, $\Phi_v=\mathbbmss{1}_{\q\mathcal{O}_v\times\mathcal{O}_v^\times}$ 
soit: $$\widehat{\Phi}_v(0)=\int_{\q\mathcal{O}_v\times\mathcal{O}_v^\times}dxdy=N(\mathfrak{D}_{F_v})^{-1}N(\q)^{-1}(1-N(\q)^{-1}).$$ 
Puisque $|\disc|=\prod_{v<\infty}N(\mathfrak{D}_{F_v})$, on a bien:
$$|\Pi_\q^{\g{k}}|\sim \frac{2(4\pi)^{\g{k}-\g{1}}|\disc|^2{\rm{res}}_{s=1}L(s,1)}
{\Gamma(\g{k}-\g{1})\widehat{\Phi}_{\g{k}}(0)N(\q)^{-1}(1-N(\q)^{-1})}\lim_{N(\q)\to\infty}\sum_{\pi\in\Pi_\q^{\g{k}}}^hL(1,{\rm{sym}}^2\pi).
$$
Commen\c cons donc \`a prouver (\ref{s}).
Pour cela on utilise la formule (\ref{val1}) pour $L(1,{\rm{sym}}^2\pi)$, pour tout $\delta>0$ (petit):
\begin{multline}
2i\pi\sum_{\pi\in\Pi_\q^{\g{k}}}^h L(1,{\rm{sym}}^2\pi)=
\sum_{\n\subset \entier}\sum_{\pi\in\Pi_\q^{\g{k}}}^h \frac{\rho_\pi(\n)}{N(\n)}\int_{(\delta)} \left(\frac{N(\q)}{N(\n)}\right)^s
L(s+1,{\rm{sym}}^2\pi_\infty)\frac{ds}{s}+ \\ \sum_{\n\subset \entier}\sum_{\pi\in\Pi_\q^{\g{k}}}^h \frac{\rho_\pi(\n)}{N(\n)}\int_{(\delta)} \left(\frac{N(\q)}
{N(\n)}\right)^sL(s+1,{\rm{sym}}^2\pi_\infty)\frac{ds}{s+1}\nonumber
\end{multline}

\begin{multline}
=\sum_{\E,\f \subset \entier}\sum_{\pi\in\Pi_\q^{\g{k}}}^h \frac{\lambda_\pi(\f^2)\chi_q(\E)}{N(\E^2\f)}\int_{(\delta)} \left(\frac{N(\q)}{N(\E^2\f)}\right)^s
L(s+1,{\rm{sym}}^2\pi_\infty)\frac{ds}{s}+ \\ \sum_{\E,\f\subset \entier}\sum_{\pi\in\Pi_\q^{\g{k}}}^h \frac{\lambda_\pi(\f^2)\chi_\q(\E)}{N(\E^2\f)}\int_{(\delta)} \left(\frac{N(\q)}
{N(\E^2\f)}\right)^sL(s+1,{\rm{sym}}^2\pi_\infty)\frac{ds}{s+1}.\nonumber
\end{multline}
On rappelle que $\chi_\q(\E)=0$ si $\q|\E$, 1 sinon. Afin d'appliquer la formule de Petersson, on \'ecrit toujours
 $\f=\xi\A$, $\A$ parcourant un ensemble de repr\'esentants de $\mathscr{C}\ell^+(F)$, et $\xi\in(\A^{-1})^{\gg 0}/\unite$. 
Comme $\A$ ne parcourant qu'un ensemble fini, on peut consid\'erer que $\A$ est fixe dans ce qui suit. La proposition \ref{petersson} 
dit:
\begin{multline}
\sum_{ \pi\in\Pi_\q^{\g{k}}}^h
\lambda_\pi(\xi^2\A^2)=\mathbbmss{1}_{\xi\A=\entier}\\
+\frac{C}{|\disc|^{1/2}} \sum_{\substack{\overline{\cc}^2=\overline{\A}^2\\
c\in\cc^{-1}\q\setminus \{0\}\\ 
\varepsilon\in\mathcal{O}_F^{\times +}/\mathcal{O}_F^{\times 2}}}
\frac{\kl(\varepsilon\xi^2,\A^2;1,\entier;c,\cc)}{N(c\cc)}J_{\g{k-1}}\left(4\pi\frac{\sqrt{
\g{\varepsilon\xi^2}[\A^2\cc^{-2}]}}{|\g{c}|}\right)\\
-\sum_{\varphi}^h\lambda(\xi^2,\A^2\diff^{-1},\varphi)\overline{\lambda(1,\diff^{-1},\varphi)}
\end{multline}
o\`u $\varphi$ parcourt dans la derni\`ere somme une base orthonormale de formes anciennes. $\xi\A=\entier\Leftrightarrow 
\A=\entier {\textrm{ et }}\xi=1$. On incorpore cette derni\`ere \'egalit\'e dans chacun des deux blocs successifs donnant 
$\sum_\pi^h L(1,{\rm{sym}}^2\pi)$, en respectant l'ordre ``terme diagonal-terme Kloosterman'':

\begin{multline}
2i\pi\sum_{\pi\in\Pi_\q^{\g{k}}}^h L(1,{\rm{sym}}^2\pi)
=\sum_{\E\subset \entier} \frac{\chi_q(\E)}{N(\E^2)}
\int_{(\delta)} \left(\frac{N(\q)}{N(\E^2)}\right)^s
L(s+1,{\rm{sym}}^2\pi_\infty)\frac{ds}{s}+ \\
\frac{C}{|\disc|^{1/2}}\sum_{\substack{\E \subset \entier\\ \bar{\A}\in\mathscr{C}\ell^+(F)\\ \xi\in(\A^{-1})^{\gg 0}}}\frac{\chi_q(\E)}{N(\E^2\xi\A)}
\int_{(\delta)} \left(\frac{N(\q)}{N(\E^2\xi\A)}\right)^s
L(s+1,{\rm{sym}}^2\pi_\infty)\frac{ds}{s}\times\\
 \sum_{\substack{\overline{\cc}^2=\overline{\A}^2\\
c\in\cc^{-1}\q\setminus \{0\}\\ 
\varepsilon\in\mathcal{O}_F^{\times +}/\mathcal{O}_F^{\times 2}}}
\frac{\kl(\varepsilon\xi^2,\A^2;1,\entier;c,\cc)}{N(c\cc)}J_{\g{k-1}}\left(4\pi\frac{\sqrt{
\g{\varepsilon\xi^2}[\A^2\cc^{-2}]}}{|\g{c}|}\right) \nonumber 
\end{multline}

\begin{multline}
+ \sum_{\E\subset \entier} \frac{\chi_\q(\E)}{N(\E^2)}\int_{(\delta)} \left(\frac{N(\q)}
{N(\E^2)}\right)^sL(s+1,{\rm{sym}}^2\pi_\infty)\frac{ds}{s+1}+\\
\frac{C}{|\disc|^{1/2}}\sum_{\substack{\E \subset \entier\\\bar{\A}\in\mathscr{C}\ell^+(F)\\ \xi\in(\A^{-1})^{\gg 0}}}\frac{\chi_q(\E)}{N(\E^2\xi\A)}
\int_{(\delta)} \left(\frac{N(\q)}{N(\E^2\xi\A)}\right)^s
L(s+1,{\rm{sym}}^2\pi_\infty)\frac{ds}{s+1}\times\\
 \sum_{\substack{\overline{\cc}^2=\overline{\A}^2\\
c\in\cc^{-1}\q\setminus \{0\}\\ 
\varepsilon\in\mathcal{O}_F^{\times +}/\mathcal{O}_F^{\times 2}}}
\frac{\kl(\varepsilon\xi^2,\A^2;1,\entier;c,\cc)}{N(c\cc)}J_{\g{k-1}}\left(4\pi\frac{\sqrt{
\g{\varepsilon\xi^2}[\A^2\cc^{-2}]}}{|\g{c}|}\right) \\
-(\textrm{Formes Anciennes}).
\nonumber
\end{multline}
$\bullet $ Une analyse semblable \`a celle de \cite{Va1}, lemme 2.3 donne:
$$\int_{(\delta)} y^{-s}
L(s+1,{\rm{sym}}^2\pi_\infty)\frac{ds}{s}\ll \exp(-Cy^{1/3d})$$
(l'exposant $1/3d$ tient au nombre de facteurs $\Gamma$ apparaissant dans $L(s+1,{\rm{sym}}^2\pi_\infty)$). Ceci permet 
d'\'evaluer les termes Kloosterman: on coupe la somme en $\E,\xi,\A$ en sommes sur $N(\E\xi\A)\leq N(\q)^{1+\eta}$ et 
$N(\E\xi\A)\geq N(\q)^{1+\eta}$ ($\eta>0$ fix\'e). Avec la borne de Weil, et l'estimation de l'int\'egrale, on a:

\begin{multline}
\sum_{N(\E\xi\A)\geq N(\q)^{1+\eta}}\frac{\chi_q(\E)}{N(\E^2\xi\A)}
\int_{(\delta)} \left(\frac{N(\q)}{N(\E^2\xi\A)}\right)^s
L(s+1,{\rm{sym}}^2\pi_\infty)\frac{ds}{s}\times\\
 \sum_{\substack{\overline{\cc}^2=\overline{\A}^2\\
c\in\cc^{-1}\q\setminus \{0\}\\ 
\varepsilon\in\mathcal{O}_F^{\times +}/\mathcal{O}_F^{\times 2}}}
\frac{\kl(\varepsilon\xi^2,\A^2;1,\entier;c,\cc)}{N(c\cc)}J_{\g{k-1}}\left(4\pi\frac{\sqrt{
\g{\varepsilon\xi^2}[\A^2\cc^{-2}]}}{|\g{c}|}\right)  \\
\ll \sum_{N(\E\xi\A)\geq N(\q)^{1+\eta}}\frac{\chi_q(\E)}{N(\E^2\xi\A)}
\exp\left(-C\left(\frac{N(\E^2\xi\A)}{N(\q)} \right)^{1/3d}\right)\times\\
\sum_{\substack{\overline{\cc}^2=\overline{\A}^2\\
c\in\cc^{-1}\q\setminus \{0\}\\ 
\varepsilon\in\mathcal{O}_F^{\times +}/\mathcal{O}_F^{\times 2}}}
\frac{\tau(c\cc)}{N(c\cc)^{1/2}}J_{\g{k-1}}\left(4\pi\frac{\sqrt{
\g{\varepsilon\xi^2}[\A^2\cc^{-2}]}}{|\g{c}|}\right)  \ll \exp (-CN(\q)^{\eta/3d}).\nonumber
\end{multline}

Pour la partie en $N(\E\xi\A)\leq N(\q)^{1+\eta}$:

\begin{multline}
\sum_{N(\E\xi\A)\leq N(\q)^{1+\eta}}\frac{\chi_q(\E)}{N(\E^2\xi\A)}
\int_{(\delta)} \left(\frac{N(\q)}{N(\E^2\xi\A)}\right)^s
L(s+1,{\rm{sym}}^2\pi_\infty)\frac{ds}{s}\times\\
 \sum_{\substack{\overline{\cc}^2=\overline{\A}^2\\
c\in\cc^{-1}\q\setminus \{0\}\\ 
\varepsilon\in\mathcal{O}_F^{\times +}/\mathcal{O}_F^{\times 2}}}
\frac{\kl(\varepsilon\xi^2,\A^2;1,\entier;c,\cc)}{N(c\cc)}J_{\g{k-1}}\left(4\pi\frac{\sqrt{
\g{\varepsilon\xi^2}[\A^2\cc^{-2}]}}{|\g{c}|}\right)\\ \ll  
\sum_{N(\E\xi\A)\leq N(\q)^{1+\eta}}\frac{\chi_q(\E)}{N(\E^2\xi\A)}
\int_{(\delta)} \left(\frac{N(\q)}{N(\E^2\xi\A)}\right)^s
L(s+1,{\rm{sym}}^2\pi_\infty)\frac{ds}{s}\times\\
 \sum_{\substack{\overline{\cc}^2=\overline{\A}^2\\
c\in\cc^{-1}\q\setminus \{0\}\\ 
\varepsilon\in\mathcal{O}_F^{\times +}/\mathcal{O}_F^{\times 2}}}
\frac{\tau(c\cc)}{N(c\cc)^{1/2}}J_{\g{k-1}}\left(4\pi\frac{\sqrt{
\g{\varepsilon\xi^2}[\A^2\cc^{-2}]}}{|\g{c}|}\right)  \\ \ll
\sum_{N(\E\xi\A)\leq N(\q)^{1+\eta}}\frac{\chi_q(\E)}{N(\E^2\xi\A)}
\int_{(\delta)} \left(\frac{N(\q)}{N(\E^2\xi\A)}\right)^s
L(s+1,{\rm{sym}}^2\pi_\infty)\frac{ds}{s}\times\\
N(\xi\A)N(\q)^{-3/2}\ll_{\delta,\eta}N(\q)^{-1/2}\nonumber
\end{multline}
la derni\`ere majoration \'etant grossi\`ere. L'autre terme de Kloosterman se traite mot pour mot de la m\^eme mani\`ere.\\
$\bullet$ Les termes diagonaux. D'abord:
\begin{multline}
\sum_{\E\subset \entier} \frac{\chi_q(\E)}{N(\E^2)}
\int_{(\delta)} \left(\frac{N(\q)}{N(\E^2)}\right)^s
L(s+1,{\rm{sym}}^2\pi_\infty)\frac{ds}{s}\\=\int_{(\delta)}N(\q)^s\zeta_F^{(\q)}(2+2s)L(s+1,{\rm{sym}}^2\pi_\infty)\frac{ds}{s}
\\={\rm{res}}_{s=0}\left( N(\q)^s\zeta_F^{(\q)}(2+2s)L(s+1,{\rm{sym}}^2\pi_\infty)s^{-1}\right)\\
+\int_{(-1/4)}N(\q)^s\zeta_F^{(\q)}(2+2s)L(s+1,{\rm{sym}}^2\pi_\infty)\frac{ds}{s}\\
=\zeta_F^{(\q)}(2)L(1,{\rm{sym}}^2\pi_\infty)+\mathcal{O}(N(\q)^{-1/4})\nonumber
\end{multline}
qui tend bien vers $\zeta_F(2)L(1,{\rm{sym}}^2\pi_\infty)$. L'autre:

\begin{multline}
\sum_{\E\subset \entier} \frac{\chi_q(\E)}{N(\E^2)}
\int_{(\delta)} \left(\frac{N(\q)}{N(\E^2)}\right)^s
L(s+1,{\rm{sym}}^2\pi_\infty)\frac{ds}{s+1}\\=\int_{(\delta)}N(\q)^s\zeta_F^{(\q)}(2+2s)L(s+1,{\rm{sym}}^2\pi_\infty)\frac{ds}{s+1}
\\={\rm{res}}_{s=-1/2}\left( N(\q)^s\zeta_F^{(\q)}(2+2s)L(s+1,{\rm{sym}}^2\pi_\infty)(s+1)^{-1}\right)\\
+\int_{(-3/4)}N(\q)^s\zeta_F^{(\q)}(2+2s)L(s+1,{\rm{sym}}^2\pi_\infty)\frac{ds}{s+1}\ll N(\q)^{-1/2}.
\end{multline}
$\bullet$ Le terme en les formes anciennes est trait\'e en appendice pour plus de clart\'e: voir la section \ref{app2}.    $\blacksquare$

\subsection{Les moments naturels}
Comme on l'a montr\'e dans la section pr\'ec\'edente, l'expression des deux premiers moments est de la forme, pour $i=1,2$:
$$\Mo_i(\q)=\frac{\sum_{\pi\in\Pi_\q^{\g{k}}}^h L(1,{\rm{sym}}^2\pi)\Lambda(1/2,\pi)^i\M(\pi)^i}
{\sum_{\pi\in\Pi_\q^{\g{k}}}^h L(1,{\rm{sym}}^2\pi)}.$$
La limite du d\'enominateur vient d'\^etre calcul\'ee: Il reste donc l'\'etude des num\'erateurs des moments, et pour prouver 
le th\'eor\`eme \ref{1} (et le th\'eor\`eme \ref{2}!), on doit g\'en\'eraliser les estimations (\ref{moment1}) et 
(\ref{moment2}), et c'est le contenu pr\'evisible de la 

\begin{proposition}
Quand $N(\q)\to\infty$, parmi les id\'eaux maximaux de $\entier$, $\g{k}\in 2\N^d_{\geq 1}$ fix\'e, on a:
\begin{multline}\label{momentnat1}
\sum_{\pi\in\Pi_\q^{\g{k}}}^h L(1,{\rm{sym}}^2\pi)\Lambda(1/2,\pi) \M(\pi)\\ =\frac{P'(1)\zeta_F(2)^2L(1/2,\pi_\infty)L(1,{\rm{sym}}^2\pi_\infty)}
{\res_{s=1}(\zeta_F)}\times \frac{2N(\q)^{1/4}}{\Delta\log(N(\q))}
\Big(1+\mathcal{O}\Big(\frac{1}{\log(N(\q))}\Big)\Big)
\end{multline}
\begin{multline}\label{momentnat2}
\sum_{\pi\in\Pi_\q^{\g{k}}}^h L(1,{\rm{sym}}^2\pi)\Lambda(1/2,\pi)^2\M(\pi)^2 = \frac{\zeta_F(2)^3L(1/2,\pi_\infty)^2L(1,{\rm{sym}}^2\pi_\infty)}
{\res_{s=1}(\zeta_F)^2}\\ \times\frac{8N(\q)^{1/2}}{\log(N(\q))^2}
\Big( \frac{||P''||_{L^2(0,1)}^2}
{\Delta^3}+\frac{P'(1)^2}{\Delta^2}+\mathcal{O}\Big(\frac{1}
{\log(N(\q))}\Big)\Big).
\end{multline}
\end{proposition}
L'id\'ee, sugg\'er\'ee par Iwaniec et Sarnak \cite{IS} (et mise en \oe uvre par Kowalski et Michel \cite {KM} pour $F=\Q$), est la suivante: on incorpore 
$L(1,{\rm{sym}}^2\pi)$ dans l'amollisseur $\M(\pi)$, obtenant ainsi un nouvel amollisseur $\widetilde{\M}(\pi)$. L'expression 
d\'ej\`a donn\'ee du carr\'e sym\'etrique est en terme des coefficients de la fonction $L$ de $\pi$, et ainsi on obtient 
une \'ecriture permettant le recours \`a la formule de Petersson. Une difficult\'e s'ajoute \`a la complexit\'e 
des termes produits: sans autre argument, la s\'erie d\'efinissant $L(1,{\rm{sym}}^2\pi)$ est a priori de longueur $N(\q)$ 
et par cons\'equent le nouvel amollisseur est trop long (et les techniques utilis\'ees ne permettent pas de 
montrer que le terme Kloosterman est n\'egligeable face \`a la diagonale). Kowalski et Michel ont r\'esolu cette difficult\'e: si  la s\'erie en question est \emph{individuellement} de longueur $N(\q)$, une s\'erie de longueur 
$N(\q)^\varepsilon$ (quel que soit $\varepsilon>0$) est suffisante \emph{en moyenne} (sur $\pi$). Tant que $\Delta+\varepsilon<1$, le terme Kloosterman reste n\'egligeable, ainsi que la contribution des formes anciennes. Pour r\'ealiser cela, on va suivre une d\'emarche 
de \cite{Ro}, qui passe par une \'etude des z\'eros de Siegel de $L(s,{\rm{sym}}^2\pi_f)$.\\

Soit $\eta>0$ un r\'eel (petit) fix\'e. On consid\`ere:
\begin{eqnarray}\nonumber
\Pi_\q^{\g{k}}(\eta)=\{\pi\in\Pi_\q^{\g{k}}|L(s,{\rm{sym}}^2\pi_f)\neq 0, \forall \Re(s)>1-\eta, |\Im(s)|\leq \log(N(\q)^3\}.
\end{eqnarray}
Nous montrerons en appendice (section \ref{densite}) que
$$|\Pi_\q^{\g{k}}\setminus \Pi_\q^{\g{k}}(\eta)| \ll_{\eta} N(\q)^{b\eta}$$
ce qui signifie qu'il y a tr\`es peu de carr\'es sym\'etriques dont la fonction $L$ s'annule pr\`es de 1. Cette propri\'et\'e a \'et\'e prouv\'ee par Kowalski et Michel, dans un contexte tr\`es g\'en\'eral, pour des repr\'esentations automorphes sur $\Q$. La preuve repose sur une in\'egalit\'e de grand crible; on l'inclut par souci de compl\'etude.

Pour $\alpha>0$ fix\'e, on pose:
\begin{multline}
\mathscr{D}_\alpha({\rm{sym}}^2\pi_f)=\frac{1}{2i\pi}\int_{(1)}L(s+1,{\rm{sym}}^2\pi_f)\Gamma(s)N(\q)^{\alpha s}ds
\\=\sum_{\n\subset\entier}\frac{\rho_\pi(\n)}{N(\n)}\exp\left( -\frac{N(\n)}{N(\q)^\alpha}\right).
\end{multline}

Soit $\mathscr{C}(\eta,\q)$ le chemin polygonal reliant $1-i\infty,1-i\log(N(\q)^2,-\frac{\eta}{2}-i\log(N(\q))^2,
-\frac{\eta}{2}+i\log(N(\q))^2,1+i\log(N(\q))^2,1+i\infty$. Le th\'eor\`eme des r\'esidus assure que:
$$L(1,{\rm{sym}}^2\pi_f)=\mathscr{D}_\alpha(\sym\pi_f)+\frac{1}{2i\pi}\int_{\mathscr{C}(\eta,\q)}L(s+1,\sym \pi_f)N(\q)^{\alpha s}\Gamma(s)ds.$$

Il est alors clair que $\Dir$ est de longueur $N(\q)^\alpha$: on a donc la s\'erie courte que l'on souhaitait. D'autre part, le lemme suivant assure que l'autre terme n'affecte pas les moments:

\begin{lemme}
Quand $\q$ parcourt les id\'eaux maximaux de $\entier$, $\g{k}\in 2\N^d_{\geq 1}$ fix\'e, on a pour tout $\varepsilon>0$ et $j\in\{1,2\}$:
\begin{multline}
\sum_{\pi\in\Pi_\q^{\g{k}}}^h \left(\int_{\mathscr{C}(\eta,\q)}L(s+1,\sym \pi_f)N(\q)^{\alpha s}\Gamma(s)ds\right) \Lambda(1/2,\pi)^j
\M(\pi)^j \\ \ll_{\eta,\varepsilon}N(\q)^{\frac{j}{4}+\varepsilon-\frac{\eta\alpha}{2}}.
\end{multline}
\end{lemme}
\preuve : On \'ecrit $\sum_{\pi\in\Pi_\q^{\g{k}}}^h=\sum_{\pi\in\Pi_\q^{\g{k}}(\eta)}^h+\sum_{\pi\in\Pi_\q^{\g{k}}\setminus \Pi_\q^{\g{k}}(\eta)}^h$, et on traite s\'epar\'ement les deux termes:\\
$\bullet$ D'apr\`es l'in\'egalit\'e de convexit\'e, $\Lambda(1/2,\pi)\ll_{\g{k},\varepsilon}N(\q)^{\frac{1}{2}+\varepsilon}$ (voir par exemple \cite{M}, (1.22), et l'article de Venkatesh \cite{V2} qui r\'esout le probl\`eme de sous-convexit\'e dans ce cas -- et bien d'autres! ) et de fa\c con \'el\'ementaire $\M(\pi)\ll_P M^{1/2}$. En outre, sur le contour $\mathscr{C}(\eta,\q)$, $L(s+1,\sym \pi_f)\ll_{\g{k},\varepsilon}N(\q)^{\varepsilon}$ pour tout $\varepsilon>0$. Soit:
\begin{multline}
\sum_{\pi\in\Pi_\q^{\g{k}}\setminus \Pi_\q^{\g{k}}(\eta)}^h\left(\int_{\mathscr{C}(\eta,\q)}L(s+1,\sym \pi_f)N(\q)^{\alpha s}\Gamma(s)ds\right) \Lambda(1/2,\pi)^j \M(\pi)^j \\ \ll N(\q)^{\varepsilon+\alpha+j/2}M^{j/2}\sum_{\pi\in\Pi_\q^{\g{k}}\setminus \Pi_\q^{\g{k}}(\eta)}^h 1\ll N(\q)^{\varepsilon+\alpha+j/2+b\eta-1}M^{j/2} 
\end{multline}
et donc, pour $\eta$ et $\alpha$ assez petits, l'in\'egalit\'e voulue est vraie ($M=N(\q)^{\Delta/2}$ avec $0<\Delta<1$).\\
$\bullet$ Pour $\pi\in\Pi_\q^{\g{k}}(\eta)$, on a (cf \cite{Ro}, lemme 3):
$$\int_{\mathscr{C}(\eta,\q)}L(s+1,\sym \pi_f)N(\q)^{\alpha s}\Gamma(s)ds\ll N(\q)^{\varepsilon-\alpha\eta/2}$$
donc par Cauchy-Schwarz:
\begin{multline}
\left|\sum_{\pi\in\Pi_\q^{\g{k}}(\eta)}^h\left(\int_{\mathscr{C}(\eta,\q)}L(s+1,\sym \pi_f)N(\q)^{\alpha s}\Gamma(s)ds\right) \Lambda(1/2,\pi)\M(\pi)\right| \\ \leq  N(\q)^{\varepsilon-\alpha\eta/2} 
\left( \sum_{\pi\in\Pi_\q^{\g{k}}}^h\Lambda(1/2,\pi)^2\M(\pi)^2\right)^{1/2}\ll N(\q)^{\varepsilon-\alpha\eta/2+1/4}\nonumber
\end{multline}
la derni\`ere majoration venant de l'estimation du second moment (\ref{moment2}). \\
De m\^eme:

\begin{multline}
\left|\sum_{\pi\in\Pi_\q^{\g{k}}(\eta)}^h\left(\int_{\mathscr{C}(\eta,\q)}L(s+1,\sym \pi_f)N(\q)^{\alpha s}\Gamma(s)ds\right) \Lambda(1/2,\pi)^2\M(\pi)^2\right|\\ \ll N(\q)^{\varepsilon-\alpha\eta/2}\sum_{\pi\in\Pi_\q^{\g{k}}}^h\Lambda(1/2,\pi)^2\M(\pi)^2\ll N(\q)^{\varepsilon-\alpha\eta/2+1/2}.\nonumber
\end{multline}
Ceci ach\`eve la preuve du lemme. $\blacksquare$\\
Il reste donc \`a \'etudier:
$$\sum_{\pi\in\Pi_\q^{\g{k}}}^h \Dir \Lambda(1/2,\pi)^j\M(\pi)^j$$
pour $j=1,2$, en esp\'erant trouver les \'equivalents (\ref{momentnat1}) et (\ref{momentnat2}).

\subsubsection{Le premier moment naturel}
On montre ici
\begin{multline}
\sum_{\pi\in\Pi_\q^{\g{k}}}^h \Dir \Lambda(1/2,\pi)\M(\pi)\\
=\frac{P'(1)\zeta_F(2)^2L(1/2,\pi_\infty)}
{\res_{s=1}(\zeta_F)}\times \frac{2N(\q)^{1/4}}{\Delta\log(N(\q))}
\Big(1+\mathcal{O}\Big(\frac{1}{\log(N(\q))}\Big)\Big)
\end{multline}
ce qui est bien, apr\`es multiplication par $L(1,\sym\pi_\infty)$, l'estimation (\ref{momentnat1}) (en effet, $\Dir$ approxime la partie finie $L(1,\sym\pi_f)$). Pour cela, on pose 
$$\widetilde{\M}(\pi)=\Dir\M(\pi)$$
La propri\'et\'e de multiplicativit\'e des coefficients $\lambda_\pi$ donne:
$$\widetilde{\M}(\pi)=\sum_{\substack{(\E,\q)=1\\N(\m)\leq M\\\f\subset\entier}}
\sum_{\D|(\m,\f^2)}\frac{\exp(-N(\E^2\f)/N(\q)^\alpha)}
{N(\E^2\f)}\frac{\mu(\m)P\left(\frac{\log(M/N(\m))}{\log(M)}\right)}{\psi(\m)N(\m)^{1/2}}\lambda_\pi(\m\f^2\D^{-2})$$
ce qui se r\'e\'ecrit:
$$\widetilde{\M}(\pi)=\sum_{\D}
\sum_{\substack{ (\E,\q)=1\\N(\m\D)\leq M\\ \f\D=\diamond  }}
\frac{\exp(-N(\E^2)\sqrt{N(\f\D)}/N(\q)^\alpha)}{N(\E)^2\sqrt{N(\f\D})}
\frac{\mu(\m\D)P\left(\frac{\log(M/N(\m\D))}{\log(M)}\right)}{\psi(\m\D)N(\m\D)^{1/2}}\lambda_\pi(\m\f)$$
la relation $\f\D=\diamond$ signifiant qu'il existe un id\'eal $\A$ tel que $\f\D=\A^2$. \\
Avec la valeur d\'ej\`a vue en section \ref{Moment1} (\'equation (\ref{valeur1})):

$$\Lambda(1/2,\pi)=(1+\varepsilon_\pi)N(\q)^{1/4}\sum_{\n\subset \entier}
\frac{\lambda_\pi(\n)}{\sqrt{N(\n)}}F(N(\n)/N(\q)^{1/2})$$
on trouve donc:
\begin{multline}
\sum_{\pi\in\Pi_\q^{\g{k}}}^h \Dir \Lambda(1/2,\pi)\M(\pi)=\\N(\q)^{1/4}\sum_{\n\subset \entier}
\sum_{\D}\sum_{\substack{ (\E,\q)=1\\N(\m\D)\leq M\\ \f\D=\diamond  }}
\frac{F(N(\n)/N(\q)^{1/2})}{\sqrt{N(\n)}}\frac{\exp(-N(\E^2)\sqrt{N(\f\D)}/N(\q)^\alpha)}{N(\E)^2\sqrt{N(\f\D})}\times\\
\frac{\mu(\m\D)P\left(\frac{\log(M/N(\m\D))}{\log(M)}\right)}{\psi(\m\D)N(\m\D)^{1/2}}
\sum_{\pi\in\Pi_\q^{\g{k}}}^h(1+\varepsilon_\pi)\lambda_\pi(\m\f)\lambda_\pi(\n).\nonumber
\end{multline}
Comme dans la section \ref{Moment1}, on utilise la relation $$\varepsilon_\pi=-i^{\g{k}}\lambda_\pi(\q)N(\q)^{1/2}$$ et on trouve:
\begin{multline}
\sum_{\pi\in\Pi_\q^{\g{k}}}^h \Dir \Lambda(1/2,\pi)\M(\pi)=\\N(\q)^{1/4}\sum_{\n\subset \entier}
\sum_{\D}\sum_{\substack{ (\E,\q)=1\\N(\m\D)\leq M\\ \f\D=\diamond  }}
\frac{F(N(\n)/N(\q)^{1/2})}{\sqrt{N(\n)}}\frac{\exp(-N(\E^2)\sqrt{N(\f\D)}/N(\q)^\alpha)}{N(\E)^2\sqrt{N(\f\D})}\times\\
\frac{\mu(\m\D)P\left(\frac{\log(M/N(\m\D))}{\log(M)}\right)}{\psi(\m\D)N(\m\D)^{1/2}}
\sum_{\pi\in\Pi_\q^{\g{k}}}^h\lambda_\pi(\m\f)\lambda_\pi(\n)
\end{multline}
\begin{multline}
-i^{\g{k}}N(\q)^{3/4}\sum_{\n\subset \entier}
\sum_{\D}\sum_{\substack{ (\E,\q)=1\\N(\m\D)\leq M\\ \f\D=\diamond  }}
\frac{F(N(\n)/N(\q)^{1/2})}{\sqrt{N(\n)}}\times\\\frac{\exp(-N(\E^2)\sqrt{N(\f\D)}/N(\q)^\alpha)}{N(\E)^2\sqrt{N(\f\D})}
\frac{\mu(\m\D)P\left(\frac{\log(M/N(\m\D))}{\log(M)}\right)}{\psi(\m\D)N(\m\D)^{1/2}}
\sum_{\pi\in\Pi_\q^{\g{k}}}^h\lambda_\pi(\m\f)\lambda_\pi(\n\q).
\nonumber
\end{multline}

Apr\`es les param\'etrisations habituelles des id\'eaux $\m,\f,\n$ (cf. les sections pr\'ec\'edentes), on applique la formule de Petersson, qui donne deux termes diagonaux, deux termes Kloosterman, 
et deux termes en formes anciennes. \\

Il est imm\'ediat de voir que, pour peu que $2\alpha+\Delta<1$, ce qu'il est loisible de supposer, les termes Kloosterman et ceux venant des formes anciennes se majorent comme on l'avait fait, en $N(\q)^{1/4-\beta}$, avec un $\beta>0$, d\'ependant cette fois aussi de $\alpha$; cela r\'esulte d'un calcul p\'enible, mais n'apportant rien de nouveau, analogue au premier moment harmonique. En fait, la perturbation engendr\'ee par le terme $\Dir$ se ma\^itrise bien car c'est une s\'erie courte, pr\'ecis\'ement:
\begin{multline}
\Dir= \sum_{\substack{(\E,\q)=1\\\f\subset\entier\\N(\E^2\f)\leq N(\q)^{2\alpha}}}
\frac{\exp(-N(\E^2\f)/N(\q)^\alpha)}
{N(\E^2\f)}\lambda_\pi(\f^2)+\mathcal{O}(\exp(-N(\q)^\alpha))\nonumber
\end{multline}
On s'int\'eresse donc au terme diagonal qui a la forme $D_1+D_2$; pour $D_1$ la diagonale est \og $\n=\m\f$\fg{}, soit:
\begin{multline}
D_1=N(\q)^{1/4}\sum_{\D}\sum_{\substack{ (\E,\q)=1\\N(\m\D)\leq M\\ \f\D=\diamond  }}
\frac{F(N(\m\f)/N(\q)^{1/2})}{\sqrt{N(\m\f)}}\frac{\exp(-N(\E^2)\sqrt{N(\f\D)}/N(\q)^\alpha)}{N(\E)^2\sqrt{N(\f\D})}\times\\
\frac{\mu(\m\D)P\left(\frac{\log(M/N(\m\D))}{\log(M)}\right)}{\psi(\m\D)N(\m\D)^{1/2}}
\end{multline}
et pour $D_2$, la diagonale est \og $\n\q=\m\f$\fg soit $\q|\f$ et $\n=\m\f\q^{-1}$:

\begin{multline}
D_2=-i^{\g{k}}N(\q)^{3/4}\sum_{\D}\sum_{\substack{ (\E,\q)=1\\N(\m\D)\leq M\\ \f\D=\diamond \\ \q|\f }}
\frac{F(N(\m\f)/N(\q)^{3/2})}{\sqrt{N(\m\f\q^{-1})}}\times\\\frac{\exp(-N(\E^2)\sqrt{N(\f\D)}/N(\q)^\alpha)}{N(\E)^2\sqrt{N(\f\D})}\frac{\mu(\m\D)P\left(\frac{\log(M/N(\m\D))}{\log(M)}\right)}{\psi(\m\D)N(\m\D)^{1/2}}
\end{multline}
mais comme $\q|\f$, en posant $\f=\q\f'$, on voit que:
$$\sum_{\substack{ (\E,\q)=1\\ \q\f'\D=\diamond}}\frac{\exp(-N(\E^2)\sqrt{N(\q\f'\D)}/N(\q)^\alpha)}{N(\E)^2\sqrt{N(\q\f'})}\ll 
\exp(-N(\q)^{1/2-\alpha})$$
d'o\`u $D_2\ll N(\q)^A\exp(-N(\q)^{1/2-\alpha})$ pour un $A$ positif, et donc $D_2$ ne contribuera pas.\\

Reste \`a traiter $D_1$, qui donnera le terme principal. Avec les expressions de $F,P$ et de l'exponentielle comme transform\'ees de Mellin, on trouve:
$$D_1=\frac{N(\q)^{\frac{1}{4}}}{(2i\pi)^3}\iiint_{(1),(1),(1)}N(\q)^{\frac{s}{2}+\alpha u}M^tL(s+\frac{1}{2},\pi_\infty)\PM(t)\Gamma(u)F(s,t,u)\frac{ds}{s}\frac{dt}{t}du$$

avec 
\begin{multline}
F(s,t,u)=\sum_{\substack{ (\E,\q)=1\\\m, \D \\ \f\D=\diamond  }}\frac{\mu(\m\D)}
{\psi(\m\D)N(\m)^{1+s+t}N(\D)^{1+t+\frac{u}{2}}N(\E)^{2+2u}N(\f)^{1+s+\frac{u}{2}}}.
\end{multline}
Cette fonction poss\`ede les propri\'et\'es habituelles:

\begin{lemme}
Pour $\Re(s),\Re(t),\Re(u)>0$, on a le d\'eveloppement eul\'erien:
$$F(s,t,u)=\zeta_F^{(\q)}(2+2u)\zeta_F(2+2s+u)\prod_{\p}\left( 1-\frac{N(\p)^{-(1+s+t)}}{\psi(\p)}-\frac{N(\p)^{-(2+s+t+u)}}
{\psi(\p)}\right).$$
On peut \'ecrire:
$$F(s,t,u)=\zeta_F^{(\q)}(2+2u)\zeta_F(2+2s+u)\zeta_F(1+s+t)^{-1}\eta(s,t,u)$$
avec $\eta$ fonction holomorphe sur $\{\Re(s+t)>-1,\Re(s+t+u)>-1\}$, et  $\eta(0,0,0)=1$.
\end{lemme}

Avec ce lemme, le calcul du terme principal de $D_1$ suit la m\^eme d\'emarche que pour le premier moment harmonique. On 
change d'abord de ligne d'int\'egration en $s$ de $\Re(s)=1$ \`a  $\Re(s)=-1$, donc tant que $\Delta+2\alpha<1$ on a :
\begin{multline}
D_1=\frac{N(\q)^{\frac{1}{4}}L(1/2,\pi_\infty)}{(2i\pi)^2}\iint_{(1),(1)}N(\q)^{\alpha u}M^t\PM(t)\Gamma(u)F(0,t,u)\frac{dt}{t}du\\ +\mathcal{O}(N(\q)^{1/4+\frac{\Delta+2\alpha-1}{2}}).\nonumber
\end{multline}

Vu que $F(0,t,u)=\zeta_F^{(\q)}(2+2u)\zeta_F(2+u)\zeta_F(1+t)^{-1}\eta(0,t,u)$, on fait le changement de ligne en $t$ de 
$\Re(t)=1$ \`a $\Re(t)=-1/2$, soit:

\begin{multline}
D_1=\frac{N(\q)^{1/4}L\left(\frac{1}{2},\pi_\infty\right) P'(1)}{{\textrm{res}}_1(\zeta_F)\log(M)}\times\\
\frac{1}{2i\pi}\int_{(1)}\Gamma(u)N(\q)^{\alpha u}\zeta_F^{(\q)}(2+2u)\zeta_F(2+u)\eta(0,0,u)du 
\bigg(1+\mathcal{O}\Big(1/\log(N(\q))\Big)\bigg)\\+\frac{N(\q)^{\frac{1}{4}}L(1/2,\pi_\infty)}{(2i\pi)^2}\iint_{(-1/2),(1)}N(\q)^{\alpha u}M^t\PM(t)\Gamma(u)F(0,t,u)\frac{dt}{t}du.\nonumber
\end{multline}
La toute derni\`ere int\'egrale est un $\mathcal{O}(N(\q)^{1/4-\Delta/4+\alpha})$. En outre, en changeant la ligne en $u$ de 
$\Re(u)=1$ en $\Re(u)=-1/2$:
$$\frac{1}{2i\pi}\int_{(1)}\Gamma(u)N(\q)^{\alpha u}\zeta_F^{(\q)}(2+2u)\zeta_F(2+u)\eta(0,0,u)du =\zeta_F(2)^2+\mathcal{O}
(N(\q)^{-\alpha/2}).$$
Au final, ceci donne:
$$D_1=\frac{N(\q)^{1/4}L\left(\frac{1}{2},\pi_\infty\right)\zeta_F(2)^2 P'(1)}{{\textrm{res}}_1(\zeta_F)\log(M)}\times\\
\bigg(1+\mathcal{O}\Big(1/\log(N(\q))\Big)\bigg)$$
ce que nous voulions!

\subsubsection{Le deuxi\`eme moment naturel}

Pour achever cette section, il nous reste \`a montrer:
\begin{multline}
\sum_{\pi\in\Pi_\q^{\g{k}}}^h \Dir \Lambda(1/2,\pi)^2\M(\pi)^2
=\frac{\zeta_F(2)^3L(1/2,\pi_\infty)^2}
{\res_{s=1}(\zeta_F)^2}\\ \times\frac{8N(\q)^{1/2}}{\log(N(\q))^2}
\Big( \frac{||P''||_{L^2(0,1)}^2}
{\Delta^3}+\frac{P'(1)^2}{\Delta^2}+\mathcal{O}\Big(\frac{1}
{\log(N(\q))}\Big)\Big).
\end{multline}

Comme pr\'ec\'edemment, on incorpore $\Dir$ dans l'amollisseur, ce qui donne:

\begin{multline}
\widetilde{\M}_2(\pi):=\M(\pi)^2\Dir=\sum_{\D\subset\entier}\sum_{\substack{N(\m_1\D)\leq M\\N(\m_2\D)\leq M}}
\sum_{\substack{\f\subset\entier\\(\E,\q)=1}}\frac{\exp\left( -\frac{N(\E^2\f)}{N(\q)^\alpha}\right)}{N(\E^2\f)}\times\\
\frac{\mu(\m_1\D)P\left(\frac{\log(M/N(\m_1\D))}{\log(M)}\right)}{\psi(\m_1\D)N(\m_1\D)^{1/2}}\cdot
\frac{\mu(\m_2\D)P\left(\frac{\log(M/N(\m_2\D))}{\log(M)}\right)}{\psi(\m_2\D)N(\m_2\D)^{1/2}}
\lambda_\pi(\f^2)\lambda_\pi(\m_1\m_2).
\end{multline}
Par commodit\'e, on pose $\m=\m_1\m_2$:
\begin{multline}
\widetilde{\M}_2(\pi)=\sum_{\D\subset\entier}\sum_{\m\subset\entier}
\sum_{\substack{\f\subset\entier\\(\E,\q)=1}}\lambda_\pi(\f^2)\lambda_\pi(\m)
\frac{\exp\left( -\frac{N(\E^2\f)}{N(\q)^\alpha}\right)}{N(\D)N(\E^2\f)}\times\\
\frac{1}{N(\m)^{1/2}} \left(\sum_{\substack{\m_1\m_2=\m\\N(\m_1\D)\leq M\\N(\m_2\D)\leq M}}
\frac{\mu(\m_1\D)P\left(\frac{\log(M/N(\m_1\D))}{\log(M)}\right)}{\psi(\m_1\D)}\cdot
\frac{\mu(\m_2\D)P\left(\frac{\log(M/N(\m_2\D))}{\log(M)}\right)}{\psi(\m_2\D)}\right).\nonumber
\end{multline}

Puisque $N(\m)\leq M^2<N(\q)$, $\m$ est premier avec $\q$ donc:
$$\lambda_\pi(\f^2)\lambda_\pi(\m)=\sum_{\de|(\f^2,\m)}\lambda_\pi(\f^2\m\de^{-2})$$
donc, en faisant le changement de variables $\m\leftarrow \m\de^{-1},\f\leftarrow \f^2\de^{-1}$, il vient:

\begin{multline}
\widetilde{\M}_2(\pi)=\sum_{\D\subset\entier}\sum_{\de\subset\entier}
\sum_{\m\subset\entier}
\sum_{\substack{\f\subset\entier\\\f\de=\diamond \\(\E,\q)=1}}\lambda_\pi(\f\m)
\frac{\exp\left( -\frac{N(\E)^2 \sqrt{N(\f\de)}}{N(\q)^\alpha}\right)}{N(\D)N(\de)N(\E)^2\sqrt{N(\f)}}\times\\
\frac{1}{\sqrt{N(\m)}} \left(\sum_{\substack{\m_1\m_2=\m\de\\N(\m_1\D)\leq M\\N(\m_2\D)\leq M}}
\frac{\mu(\m_1\D)P\left(\frac{\log(M/N(\m_1\D))}{\log(M)}\right)}{\psi(\m_1\D)}\cdot
\frac{\mu(\m_2\D)P\left(\frac{\log(M/N(\m_2\D))}{\log(M)}\right)}{\psi(\m_2\D)}\right)\nonumber
\end{multline}
ce qui est une forme utilisable en vue d'appliquer la formule des traces de Petersson:

\begin{multline}
\sum_{\pi\in\Pi_\q^{\g{k}}}^h \Dir \Lambda(1/2,\pi)^2\M(\pi)^2=
2N(\q)^{1/2}\sum_{\substack{\n,\m\\ \D,\de}}
\sum_{\substack{\f\subset\entier\\\f\de=\diamond \\(\E,\q)=1}}
\frac{G\left(\frac{N(\n)}{N(\q)}\right)}{\sqrt{N(\n)}}\tau(\n) \times\\
\frac{1}{\sqrt{N(\m)}} \left(\sum_{\substack{\m_1\m_2=\m\de\\N(\m_1\D)\leq M\\N(\m_2\D)\leq M}}
\frac{\mu(\m_1\D)P\left(\frac{\log(M/N(\m_1\D))}{\log(M)}\right)}{\psi(\m_1\D)}\cdot
\frac{\mu(\m_2\D)P\left(\frac{\log(M/N(\m_2\D))}{\log(M)}\right)}{\psi(\m_2\D)}\right)\\
\times\frac{\exp\left( -\frac{N(\E)^2 \sqrt{N(\f\de)}}{N(\q)^\alpha}\right)}{N(\D)N(\de)N(\E)^2\sqrt{N(\f)}}
\sum_{\pi\in\Pi_\q^{\g{k}}}^h \lambda_\pi(\n)\lambda_\pi(\f\m).\nonumber
\end{multline}

Le terme en sommes de Kloosterman se traite comme lors du deuxi\`eme moment harmonique, ainsi que celui issu des formes anciennes. La diagonale est ici \og $\n=\m\f$\fg:

\begin{multline}
D(\q)=2N(\q)^{1/2}\sum_{\substack{\m\\ \D,\de}}
\sum_{\substack{\f\subset\entier\\\f\de=\diamond \\(\E,\q)=1}}
\frac{G\left(\frac{N(\m\f)}{N(\q)}\right)}{\sqrt{N(\m\f)}}\tau(\m\f)
\cdot\frac{\exp\left( -\frac{N(\E)^2 \sqrt{N(\f\de)}}{N(\q)^\alpha}\right)}{N(\D)N(\de)N(\E)^2\sqrt{N(\f)}}\times
\\
\frac{1}{\sqrt{N(\m)}} \left(\sum_{\substack{\m_1\m_2=\m\de\\N(\m_1\D)\leq M\\N(\m_2\D)\leq M}}
\frac{\mu(\m_1\D)P\left(\frac{\log(M/N(\m_1\D))}{\log(M)}\right)}{\psi(\m_1\D)}\cdot
\frac{\mu(\m_2\D)P\left(\frac{\log(M/N(\m_2\D))}{\log(M)}\right)}{\psi(\m_2\D)}\right).\nonumber
\end{multline}

La description int\'egrale des fonctions intervenant implique:

\begin{multline}
D(\q)=\frac{2N(\q)^{1/2}}{(2i\pi)^4}\iiiint_{(1),(1),(1),(1)}N(\q)^{s+\alpha u}M^{t_1+t_2}L(s+1/2,\pi_\infty)^2\times\\
\zeta_F^{(q)}(1+2s) \PM(t_1)\PM(t_2)\Gamma(u)F(s,t_1,t_2,u)\frac{ds}{s}\frac{dt_1}{t_1}\frac{dt_2}{t_2}du
\end{multline}
avec la fonction $F$, s\'erie absolument convergente dans la r\'egion $\{(s,t_1,t_2,u)\in \C^4; \Re(s),\Re(t_1),\Re(t_2),\Re(u)>1\}$ valant alors:
\begin{multline}
\sum_{\substack{\m\\ \D,\de}}
\sum_{\substack{\f\subset\entier\\\f\de=\diamond \\(\E,\q)=1}}
\frac{\tau(\m\f)}{N(\m)^{1+s}N(\f)^{1+s+\frac{u}{2}}N(\E)^{2+2u}N(\D)^{1+t_1+t_2}N(\de)^{1+\frac{u}{2}}}\times\\
\left(\sum_{\m_1\m_2=\m\de}\frac{\mu(\m_1\D)\mu(\m_2\D)}{\psi(\m_1\D)\psi(\m_2\D)}N(\m_1)^{-t_1}N(\m_2)^{-t_2}\right).
\end{multline}
 On v\'erifie que l'on a la 
factorisation en produit infini, convergent sur le m\^eme domaine:
$$F(s,t_1,t_2,u)=\prod_{\p}F_\p(s,t_1,t_2,u)$$ avec:

\begin{multline}
F_\p(s,t_1,t_2,u)=\sum_{\substack{m,d,D\geq 0\\(\q,\p^e)=1 }}
\sum_{\substack{f\geq 0\\ f+D {\textrm{ pair }} }}\tau(\p^{m+f})
\bigg(\sum_{\substack{m_1+m_2\\=m+D}}\frac{\mu(\p^{m_1+d})\mu(\p^{m_2+d})}{\psi(\p^{m_1+d})\psi(\p^{m_2+d})}\times\\
N(\p)^{-m_1t_1}N(\p)^{-m_2t_2}\bigg)
  N(\p)^{-e(2+2u)-m(1+s)-f(1+s+\frac{u}{2})-d(1+t_1+t_2)-D(1+\frac{u}{2})}. \nonumber
\end{multline}
On peut calculer explicitement $F_\p$. Les 
valeurs permises pour $(d,m,D)$ sont (0,0,0), (1,0,0), (0,1,0), (0,0,1), (0,1,1), (0,2,0), (0,0,2). C'est dans cet ordre que l'on
 pr\'esente les facteurs produits par ces conditions et la somme sur $f$, not\'es entre parenth\`eses; le terme produit par la somme \og libre\fg en $e$ 
est tout \`a gauche ($\chi_\q(\p)$ vaut 1 si $\q\neq \p$, 0 sinon):

\begin{multline}
F_\p(s,t_1,t_2,u)=(1-\chi_\q(\p)N(\p)^{-2-2u})^{-1}\times\\
\bigg\{  
\left( \frac{1}{1-N(\p)^{-2(1+s+u/2)} }+\frac{2N(\p)^{-2(1+s+u/2)}}{(1-N(\p)^{-2(1+s+u/2)})^2} \right)\\
+\frac{N(\p)^{-(1+t_1+t_2)}}{\psi(\p)^2}\left( \frac{1}{1-N(\p)^{-2(1+s+u/2)} }+\frac{2N(\p)^{-2(1+s+u/2)}}{(1-N(\p)^{-2(1+s+u/2)})^2} \right)\\
-\frac{N(\p)^{-(1+s+t_1)}+N(\p)^{-(1+s+t_2)}}{\psi(\p)}\times \\\left( \frac{2}{1-N(\p)^{-2(1+s+u/2)} }+\frac{2N(\p)^{-2(1+s+u/2)}}{(1-N(\p)^{-2(1+s+u/2)})^2} \right)\\
-\frac{N(\p)^{-(1+u/2+t_1)}+N(\p)^{-(1+u/2+t_2)}}{\psi(\p)}\times\\ 
\left( \frac{N(\p)^{-(1+s+u/2)}}{1-N(\p)^{-2(1+s+u/2)} }+\frac{N(\p)^{-(1+s+u/2)}(1+N(\p)^{-2(1+s+u/2)})}{(1-N(\p)^{-2(1+s+u/2)})^2} \right)\\
+\frac{N(\p)^{-(2+s+u/2+t_1+t_2)}}{\psi(\p)^2} \times\\ \left( \frac{2N(\p)^{-(1+s+u/2)}}{1-N(\p)^{-2(1+s+u/2)} }+\frac{N(\p)^{-(1+s+u/2)}(1+N(\p)^{-2(1+s+u/2)})}{(1-N(\p)^{-2(1+s+u/2)})^2} \right)\\
+\frac{N(\p)^{-(2+s+u/2+t_1+t_2)}}{\psi(\p)^2} \left( \frac{3}{1-N(\p)^{-2(1+s+u/2)} }+\frac{2N(\p)^{-2(1+s+u/2)}}{(1-N(\p)^{-2(1+s+u/2)})^2} \right)\\
+\frac{N(\p)^{-(2+s+u/2+t_1+t_2)}}{\psi(\p)^2} \left( \frac{1}{1-N(\p)^{-2(1+s+u/2)} }+\frac{2N(\p)^{-2(1+s+u/2)}}{(1-N(\p)^{-2(1+s+u/2)})^2} \right)
\bigg\}.\nonumber
\end{multline}

Il est clair que l'on peut mettre en facteur $(1-N(\p)^{-2(1+s+u/2)})^{-1}$. Les termes d'ordre un sont  
$N(\p)^{-(1+t_1+t_2)}$, $-2N(\p)^{-(1+s+t_1)}$ et $-2N(\p)^{-(1+s+t_2)}$. Leur mise en facteur forc\'ee permet d'\'elargir 
le domaine de convergence:

\begin{lemme}
On peut \'ecrire 
$$F(s,t_1,t_2,u)=\frac{\zeta_F^{(\q)}(2+2u)\zeta_F(2+2s+u)\zeta_F(1+t_1+t_2)}{\zeta_F(1+t_1+s)^2\zeta_F(1+t_2+s)^2}\eta(s,t_1,t_2,u)$$
la fonction $\eta$ \'etant donn\'ee par un produit eul\'erien convergent (au moins) sur le domaine 
$\Re(s+u/2)>-1/2$, $\Re(t_1+t_2)>-1/2$, $\Re(s+t_1)>-1/2$, $\Re(s+t_2)>-1/2$, avec $\eta(0,0,0,0)=\zeta_F(2)$.
\end{lemme}
Cette affirmation r\'esulte d'un calcul \'el\'ementaire. Elle permet d'effectuer les changements de droite d'int\'egration 
ad\'equats. On part donc de:

\begin{multline}
D(\q)=\frac{2N(\q)^{1/2}}{(2i\pi)^4}\iiiint_{(1),(1),(1),(1)}N(\q)^{s+\alpha u}M^{t_1+t_2}L(s+1/2,\pi_\infty)^2
\times\\\zeta_F^{(q)}(1+2s)
\PM(t_1)\PM(t_2)\Gamma(u)
F(s,t_1,t_2,u)\frac{ds}{s}\frac{dt_1}{t_1}\frac{dt_2}{t_2}du.\nonumber
\end{multline}
On peut d\'ej\`a d\'eplacer toutes les lignes d'int\'egration jusqu'\`a $(1/2)$. Puis, on change de ligne en $s$:
\begin{multline}
D(\q)=\frac{2N(\q)^{1/2}}{(2i\pi)^3}\iiint_{(\frac{1}{2}),(\frac{1}{2}),(\frac{1}{2})}N(\q)^{\alpha u}M^{t_1+t_2}
\PM(t_1)\PM(t_2)\Gamma(u)\times\\
{\rm{res}}_{s=0}\left( N(\q)^sL\left(\frac{1}{2},\pi_\infty\right)F(s,t_1,t_2,u)\right)
\frac{dt_1}{t_1}\frac{dt_2}{t_2}du\\
+\frac{2N(\q)^{1/2}}{(2i\pi)^4}\iiiint_{(-\frac{1}{2}),(\frac{1}{2}),(\frac{1}{2})}N(\q)^{s+\alpha u}M^{t_1+t_2}L(s+1/2,\pi_\infty)^2
\times\\\zeta_F^{(q)}(1+2s)
\PM(t_1)\PM(t_2)\Gamma(u)
F(s,t_1,t_2,u)
\frac{ds}{s}\frac{dt_1}{t_1}\frac{dt_2}{t_2}du
\nonumber
\end{multline}
le second terme est un $\mathcal{O}(N(\q)^{1/2-\delta})$ (pour un $\delta>0$), si on a pris la peine de choisir $\alpha<1-\Delta$. 
Pour \'evaluer le premier, on continue de proc\'eder comme \`a la section \ref{Moment2}: il faut calculer le r\'esidu, ou 
du moins la partie participant au terme dominant. Celui-ci s'\'ecrit:

\begin{multline}
 \frac{{\textrm{res}}_{1}(\zeta_F^{(\q)})}{2}\log(N(\q))F(0,t_1,t_2)L\left(\frac{1}{2},\pi_\infty\right)^2   +
{\textrm{res}}_{1}(\zeta_F^{(\q)})L\left(\frac{1}{2},\pi_\infty\right)^2  \times\\ \frac{\partial F}{\partial s}(0,t_1,t_2)+\ldots\nonumber
\end{multline}
\og $\ldots$\fg d\'esignant des termes ayant un effet de second ordre. Plus pr\'ecis\'ement, il s'agit de:

\begin{multline}
\frac{ {\textrm{res}}_{1}(\zeta_F^{(\q)})}{2}\log(N(\q))\frac{\zeta_F^{(\q)}(2+2u)\zeta_F(2+u)\zeta_F(1+t_1+t_2)}{\zeta_F(1+t_1)^2\zeta_F(1+t_2)^2}\eta(0,t_1,t_2,u)
\times\\ L\left(\frac{1}{2},\pi_\infty\right)^2   +
{\textrm{res}}_{1}(\zeta_F^{(\q)})L\left(\frac{1}{2},\pi_\infty\right)^2 \zeta_F^{(\q)}(2+2u) \zeta_F(1+t_1+t_2)  
\zeta_F(2+u)\times\\ \eta(0,t_1,t_2,u) \frac{\partial}{\partial s}
\left(  \zeta_F(1+t_1+s)^{-2}\zeta_F(1+t_2+s)^{-2} \right)|_{s=0}.
\end{multline}
La s\'eparation des variables $u$ et $(t_1,t_2)$ permet de finir le calcul en s'appuyant sur la section \ref{Moment2}. On a 
en effet:
$$D(\q)=D_1(\q)+D_2(\q)+\mathcal{O}(N(\q)^{1/2-\eta})$$ ($\eta>0$) avec 

\begin{multline}
D_1(\q)=\frac{{\textrm{res}}_{1}(\zeta_F^{(\q)})N(\q)^{1/2}\log((N(\q))L\left(\frac{1}{2},\pi_\infty\right)^2}{(2i\pi)^3}
\int_{(\frac{1}{2})} N(\q)^{\alpha u}\Gamma(u)\times\\\zeta_F^{(\q)}(2+2u)\zeta_F(2+u)
\bigg(\iint_{(\frac{1}{2}),(\frac{1}{2})} N(\q)^{\alpha u}M^{t_1+t_2}
\PM(t_1)\PM(t_2) \times\\ \frac{ \zeta_F(1+t_1+t_2)}{\zeta_F(1+t_1)^2\zeta_F(1+t_2)^2}
\eta(0,t_1,t_2,u)
\frac{dt_1}{t_1}\frac{dt_2}{t_2} \bigg)du\nonumber
\end{multline}

\begin{multline}
D_2(\q)=\frac{2{\textrm{res}}_{1}(\zeta_F^{(\q)})N(\q)^{1/2}L\left(\frac{1}{2},\pi_\infty\right)^2}{(2i\pi)^3}
\int_{(\frac{1}{2})} N(\q)^{\alpha u}\Gamma(u)\times\\\zeta_F^{(\q)}(2+2u)\zeta_F(2+u)
\bigg(\iint_{(\frac{1}{2}),(\frac{1}{2})} N(\q)^{\alpha u}M^{t_1+t_2}
\PM(t_1)\PM(t_2) \times\\ \zeta_F(1+t_1+t_2)\frac{\partial}{\partial s}
\left(\zeta_F(1+t_1)^{-2}\zeta_F(1+t_2)^{-2}\right)|_{s=0}
\eta(0,t_1,t_2,u)
\frac{dt_1}{t_1}\frac{dt_2}{t_2} \bigg)du.\nonumber
\end{multline}

Les int\'egrales internes ont \'et\'e calcul\'ees en section \ref{Moment2}.
\begin{multline}
D_1(\q)=\frac{ N(\q)^{1/2}\log((N(\q))L\left(\frac{1}{2},\pi_\infty\right)^2}{(2i\pi){\textrm{res}}_{1}(\zeta_F^{(\q)})^2}
\int_{(\frac{1}{2})} N(\q)^{\alpha u}\Gamma(u)\times\\\zeta_F^{(\q)}(2+2u)\zeta_F(2+u)
\eta(0,0,0,u)
du\times \left(\log(M)^{-3}\int_0^1P''(t)^2dt+\mathcal{O}(\log(M)^{-4})\right)
\nonumber
\end{multline}

\begin{multline}
D_2(\q)=\frac{2 N(\q)^{1/2}L\left(\frac{1}{2},\pi_\infty\right)^2}{(2i\pi){\textrm{res}}_{1}(\zeta_F^{(\q)})^2}
\int_{(\frac{1}{2})} N(\q)^{\alpha u}\Gamma(u)\times\\\zeta_F^{(\q)}(2+2u)\zeta_F(2+u)\eta(0,0,0,u)
du\times\left( \log(M)^{-2}P'(1)^2+\mathcal{O}(\log(M)^{-3})\right).\nonumber
\end{multline}

Il ne reste plus qu'\`a calculer l'int\'egrale en $u$:

\begin{multline}
\frac{1}{2i\pi}\int_{(\frac{1}{2})} N(\q)^{\alpha u}\Gamma(u)\zeta_F^{(\q)}(2+2u)\zeta_F(2+u)\eta(0,0,0,u)
du\\={\rm{res}}_{s=0}\left(N(\q)^{\alpha u}\Gamma(u)\zeta_F^{(\q)}(2+2u)\zeta_F(2+u)\eta(0,0,0,u)\right)\\+
\frac{1}{2i\pi}\int_{(-\frac{1}{2})} N(\q)^{\alpha u}\Gamma(u)\zeta_F^{(\q)}(2+2u)\zeta_F(2+u)\eta(0,0,0,u)
du\\=\zeta_F(2)^{3}+\mathcal{O}(N(\q)^{-\alpha/2})
\end{multline}

et cela ach\`eve le calcul.\\

\rem Nous n'avons fait les calculs des moments naturels que pour $\Lambda(1/2,\pi)$. Il est clair cependant que ceux-ci 
s'\'etendent au cas de la d\'eriv\'ee. C'est ainsi que l'on conclut la preuve des th\'eor\`emes 
\ref{1} et \ref{2}.\\

\rem {\sc{ Une approche probabiliste}}\\Le lecteur aura remarqu\'e qu'en fait, nous n'avons pas d\'eduit les valeurs des moments naturels de celles des 
moments harmoniques: nous avons refait les calculs, et aurions pu m\^eme les pr\'esenter directement (en perdant certainement 
un peu de lisibilit\'e). Il est donc logique de se demander s'il existe une preuve naturelle, donnant les moments naturels 
gr\^ace aux moments harmoniques.

Pour donner une id\'ee de r\'eponse, consid\'erons l'ensemble fini $\Pi_\q^{\g{k}}$, muni de la probabilit\'e (asymptotique)
$\int Xd\mathbb{P}=\sum_{\pi\in\Pi_\q^{\g{k}}}\omega_\pi X(\pi)$. On 
peut consid\'erer les variables al\'eatoires suivantes:
\begin{displaymath}
X_\q:\begin{array}{ll}\Pi_\q^{\g{k}}\longrightarrow \R\\\pi\longmapsto L(1/2,\pi)\M(\pi)\end{array}
\end{displaymath}
et
\begin{displaymath}
Y_\q:\begin{array}{ll}\Pi_\q^{\g{k}}\longrightarrow \R\\\pi\longmapsto L(1,{\rm{sym}}^2\pi)\end{array}.
\end{displaymath}
Les travaux de Royer (notamment \cite{Ro}) montrent que $Y_\q$ est faiblement convergente (quand $F=\Q$). Faisons 
l'hypoth\`ese que ce soit aussi le cas de $X_\q$ (nous avons montr\'e la convergence des deux premiers moments, donc 
cette hypoth\`ese est raisonnable vis-\`a-vis de ce qui pr\'ec\`ede). On a envie de poser:\\

{\sc{Conjecture}}: les limites de $X_\q$ et $Y_\q$ sont ind\'ependantes. \\

Supposons la conjecture vraie. On \'ecrit, en se rappelant de l'\'etude faite sur le carr\'e sym\'etrique:
\begin{multline}
\frac{ \sum_{\pi\in\Pi_\q^{\g{k}}} L(1/2,\pi)\M(\pi)}{|\Pi_\q^{\g{k}}|}
=\frac{\sum_{\pi\in\Pi_\q^{\g{k}}}^h L(1,{\rm{sym}}^2\pi) L(1/2,\pi)\M(\pi)}{\sum_{\pi\in\Pi_\q^{\g{k}}}^h L(1,{\rm{sym}}^2\pi)}
(1+{\rm{o}}(1))\\=\frac{ \sum_{\pi\in\Pi_\q^{\g{k}}}^h L(1,{\rm{sym}}^2\pi)  \sum_{\pi\in\Pi_\q^{\g{k}}}^h  L(1/2,\pi)\M(\pi)}
{\sum_{\pi\in\Pi_\q^{\g{k}}}^h L(1,{\rm{sym}}^2\pi)}(1+{\rm{o}}(1))\\= \sum_{\pi\in\Pi_\q^{\g{k}}}^h  L(1/2,\pi)\M(\pi)
(1+{\textrm{o}}(1)).
\end{multline}
Notons que par l'\'etude des deux premiers moments amollis, que nous avons faite \og \`a la main\fg, nous avons montr\'e que 
les variables $X_\q$ et $Y_\q$ sont d\'ecorr\'el\'ees \`a l'infini, ce qui donne une raison de croire en la conjecture. 
Malheureusement, la m\'ethode des moments semble limit\'ee, \`a cause de la complexit\'e des calculs d'une part, mais aussi 
parce que d'autres termes dominants font leur apparition dans les moments d'ordre sup\'erieur \`a trois, qui ne sont 
pas issus de la diagonale: c'est d\'ej\`a le cas du moment d'ordre quatre. Il faudrait donc une autre approche pour esp\'erer une 
preuve g\'en\'erale de la conjecture, qui devrait \^etre formul\'ee dans un cadre automorphe plus g\'en\'eral.

\section{Appendice}\label{Appendice}
Nous allons ici expliquer le traitement des termes issus des formes anciennes, laiss\'es  de c\^ot\'e lors des sections 
\ref{Moment1}, \ref{Moment2} et \ref{Derivee}. Cet appendice s'appuie essentiellement sur l'article \cite{ILS}, qui construit dans un 
contexte plus g\'en\'eral, sur $\Q$, une base orthogonale des formes anciennes. Avec la famille d'amollisseurs utilis\'ee 
ici, cette base particuli\`ere permet de majorer la contribution des formes anciennes aux deux premiers moments amollis, et 
d'assurer qu'elles n'affectent pas l'\'equivalent trouv\'e. \\

On suppose ici que $\q$ est premier (non nul) et $\g{k}\in 2\Z^d_{\geq 1}$.\\
Rappelons que l'on a la d\'ecomposition suivante (cf (\ref{newold})):
$$\mathcal{H}_\q^{\g{k}}=\underbrace{\bigoplus_{\pi\in\Pi_\q^{\g{k}}}\pi^{K_\infty K_0(\q)}}_{\mathcal{H}_\q^{\g{k}}(new)}\oplus
\underbrace{\bigoplus_{\pi\in\Pi_{\entier}^{\g{k}}} \pi^{K_\infty K_0(\q)}}_{\mathcal{H}_\q^{\g{k}}(old)}$$
et on supposera dor\'enavant que $\mathcal{H}_\q^{\g{k}}(old)$ est non vide (i.e. qu'il existe des formes non ramifi\'ees). 
Le th\'eor\`eme de Casselman assure que pour toute $\pi$ non ramifi\'ee, ${\rm{dim}}_\C\pi^{K_\infty K_0(\q)}=2$. Si $\varphi_\pi$ d\'esigne un vecteur 
nouveau, il fait \'evidemment partie de cet espace . Soit $\varpi_\q$ une uniformisante de la place $\q$, et notons aussi 
$\varpi_\q={\rm{id}}(\q)$ l'id\`ele valant $\varpi_\q$ \`a la place $\q$, 1 ailleurs: alors $\varphi_\pi\left(g\left(
\begin{array}{cc}\varpi_\q^{-1}&0\\0&1\end{array}\right)\right)$ est aussi $K_0(\q)$-invariant. Comme ces deux fonctions  
sont lin\'eairement ind\'ependantes, elles forment une base de $\mathcal{H}_\q^{\g{k}}(old)$. Pour en d\'eduire une base 
orthogonale, on a besoin du

\begin{lemme}Soit $X(\q)=Z(\adele)\gl(F)\backslash \gl(\adele)/K_0(\q)$.
Soit $\pi\in\Pi_{\entier}^{\g{k}}$ et $\varphi_\pi$ un \'el\'ement sp\'ecial de $\pi$. On a:
\begin{eqnarray} \int_{X(\q)}|\varphi_\pi(g)|^2dg=
\int_{X(\q)}\left|\varphi_\pi\left(g\left(
\begin{array}{cc}\varpi_\q^{-1}&0\\0&1\end{array}\right)\right)\right|^2dg\nonumber
\end{eqnarray}
\begin{eqnarray}\int_{X(\q)} \varphi_\pi(g)\overline{\varphi_\pi\left(g\left(
\begin{array}{cc}\varpi_\q^{-1}&0\\0&1\end{array}\right)\right)}dg=\frac{N(\q)^{1/2}\lambda_\pi(\q)}{N(\q)+1}
\int_{X(\q)}|\varphi_\pi(g)|^2dg.\nonumber\end{eqnarray}
\end{lemme}
{\sc{Preuve}}: Montrons la deuxi\`eme \'egalit\'e, la premi\`ere \'etant cons\'equence \'evidente du fait que $dg$ est issue 
d'une mesure de Haar. On a d'une part, avec $K_f=\prod_{v<\infty} \gl(\mathcal{O}_{F_v})$:
\begin{multline}
\int_{Z(\adele)\gl(F)\backslash \gl(\adele)/K_f}\varphi_\pi(g)\overline{\varphi_\pi\left(g\left(
\begin{array}{cc}\varpi_\q^{-1}&0\\0&1\end{array}\right)\right)}dg\\=\int_{Z(\adele)\gl(F)\backslash \gl(\adele)/K_0(\q)}
[K:K_0(\q)]^{-1}
\varphi_\pi(g)\overline{\varphi_\pi\left(g\left(
\begin{array}{cc}\varpi_\q^{-1}&0\\0&1\end{array}\right)\right)}dg\nonumber
\end{multline}
d'autre part:
\begin{multline}
\int_{Z(\adele)\gl(F)\backslash \gl(\adele)/K_f}\varphi_\pi(g)\overline{\varphi_\pi\left(g\left(
\begin{array}{cc}\varpi_\q^{-1}&0\\0&1\end{array}\right)\right)}dg\\=\int_{Z(\adele)\gl(F)\backslash \gl(\adele)/K_f}
\int_{K_f}
\varphi_\pi(gk^{-1})\overline{\varphi_\pi\left(g\left(
\begin{array}{cc}\varpi_\q^{-1}&0\\0&1\end{array}\right)\right)}dkdg\\
=\int_{Z(\adele)\gl(F)\backslash \gl(\adele)/K_f}\varphi_\pi(g)\int_{K_f}\overline{\varphi_\pi\left(gk\left(
\begin{array}{cc}\varpi_\q^{-1}&0\\0&1\end{array}\right)\right)}dkdg\\
=\int_{Z(\adele)\gl(F)\backslash \gl(\adele)/K_f}\varphi_\pi(g)\times\frac{1}{N(\q)+1}\overline{{\rm{T}}_\q(\varphi_\pi)(g)}
dg\nonumber
\end{multline}
${\rm{T}}_\q$ d\'esignant l'op\'erateur de Hecke en $\q$ (tel que d\'efini dans \cite{G}, (3.15) page 48); on sait de plus 
(m\^eme r\'ef\'erence, (4.16) page 72) que 
$${\rm{T}}_\q(\varphi)=N(\q)^{1/2}\lambda_\pi(\q)$$
et ceci donne le r\'esultat. $\blacksquare$\\

\rem : ce lemme correspond au lemme 2.4 de \cite{ILS}, et en donne une preuve ad\'elique (dans un cas plus simple, auquel le 
lemme 2.4 pourrait se ramener).\\

Maintenant que l'on a g\'en\'eralis\'e ce dont on avait besoin, on peut, en suivant \cite{ILS}, poser:
\begin{multline}
\psi_\pi(g):=\left(\frac{N(\q)}{\rho_\pi(\q)}\right)^{1/2}\sum_{\cc\D=\q}\frac{\mu(\cc)\lambda_\pi(\cc)}{N(\cc)\psi(\cc)}
N(\D)^{-1/2}\varphi_\pi\left(g\left(\begin{array}{cc}{\rm{id}}(\D)^{-1}&0\\0&1\end{array}\right)\right)
\end{multline}
avec $$\rho_\pi(\cc)=\prod_{\p|\cc}\left(1-N(\p)\left(\frac{\lambda_\pi(\p)}{N(\p)+1}\right)^2\right).$$
Gr\^ace au lemme pr\'ec\'edent, $\{\varphi_\pi,\psi_\pi\}_{\pi\in \Pi_{\entier}^{\g{k}}}$ forme une base orthogonale de $\mathcal{H}_\q^{\g{k}}(old)$, telle 
que $||\varphi_\pi||_{\mathcal{H}_\q^{\g{k}}}= ||\psi_\pi||_{\mathcal{H}_\q^{\g{k}}}$. \\

Soit $L(s,{\rm{sym}}^2\pi)$ le carr\'e sym\'etrique de $\pi$. On a un d\'eveloppement eul\'erien 
$L(s,{\rm{sym}}^2\pi)=\prod_vL(s,{\rm{sym}}^2\pi_v)$ index\'e par les 
places finies et infinies, et on d\'eduit de l'expression des facteurs locaux 
que $$(1+N(\q)^{-1})\rho_\pi(\q)=(1-N(\q)^{-2})^{-1}L(1,{\rm{sym}}^2\pi_\q)^{-1}$$ (cf \cite{ILS}: (2.50) et 
section 3). \\

Les coefficients de Fourier de $\psi_\pi$ sont donc donn\'es par:
\begin{multline}
\lambda(1,\n\diff^{-1},\psi_\pi)=\Big(N(\q)(1-N(\q)^{-2})(1+N(\q)^{-1})L(1,{\rm{sym}}^2\pi_\q)\Big)^{1/2}\times\\
\bigg(\frac{-1}{N(\q)+1}\lambda_\pi(\q)\lambda_\pi(\n)+\lambda_\pi(\n\q^{-1})\mathbbmss{1}_{\q|\n}\bigg).
\end{multline}

\subsection{La contribution des formes anciennes aux premiers et deuxi\`emes moments}
La contribution des formes anciennes au premier moment s'\'ecrit:
\begin{multline}
M_1(\q,old)=N(\q)^{1/4}\sum_{\substack{\n\subset \entier\\ N(\m)\leq M}}\frac{F(N(\n)/N(\q)^{1/2})}{N(\n)^{1/2}}\cdot 
\frac{\mu(\m)P\left(\frac{
\log(M/N(\m))}{\log(M)}\right)}{\psi(\m)N(\m)^{1/2}}\times\\
\sum_{f}^h\lambda(1,\n\diff^{-1},f)\lambda(1,\m\diff^{-1},f)\\
-N(\q)^{3/4}\sum_{\substack{\n\subset \entier\\ N(\m)\leq M}}\frac{F(N(\n)/N(\q)^{1/2})}{N(\n)^{1/2}}\cdot 
\frac{\mu(\m)P\left(\frac{
\log(M/N(\m))}{\log(M)}\right)}{\psi(\m)N(\m)^{1/2}}\times\\
\sum_{f}^h\lambda(1,\n\q\diff^{-1},f)\lambda(1,\m\diff^{-1},f)
\end{multline}
$\{f\}$ d\'esignant la base orthogonale form\'ee par $\{\varphi_\pi,\psi_\pi\}_{\pi\in\Pi_{\entier}^{\g{k}}}$, et $\sum^h$ 
 d\'esignant la somme normalis\'ee pour le produit scalaire de $\mathcal{H}_\q^{\g{k}}$.\\
\indent La partie de cette somme form\'ee en ne gardant que les termes en $\varphi_\pi$ est facile \`a traiter. En effet, 
on a 
\begin{eqnarray}
||\varphi_\pi||_{\mathcal{H}_\q^{\g{k}}}^2=[K_f:K_0(\q)]||\varphi_\pi||_{\mathcal{H}_{\entier}^{\g{k}}}^2=
(N(\q)+1)||\varphi_\pi||_{\mathcal{H}_{\entier}^{\g{k}}}^2
\end{eqnarray}
on peut par exemple alors utiliser la formule de Petersson pour $\Pi_{\entier}^{\g{k}}$, puisqu'on a la bonne normalisation 
harmonique, et on peut majorer la contribution des termes en $\varphi_\pi$ $M_1(\q,\varphi_\pi)$ par:
$$M_1(\q,\varphi_\pi)\ll\frac{ N(\q)^{3/4}N(\q)^{1/4}M^{1/2}}{N(\q)}\ll N(\q)^{1/4-\eta}$$
avec $\eta>0$.\\

Reste \`a traiter la contribution des $\psi_\pi$, not\'ee $M_1(\q,\psi_\pi)$. On \'ecrit:
\begin{multline}
M_1(\q,\psi_\pi)=N(\q)^{1/4}\sum_{\substack{\n\subset \entier\\ N(\m)\leq M}}\frac{F(N(\n)/N(\q)^{1/2})}{N(\n)^{1/2}}\cdot 
\frac{\mu(\m)P\left(\frac{
\log(M/N(\m))}{\log(M)}\right)}{\psi(\m)N(\m)^{1/2}}\times\\
\sum_{\pi}^h\lambda(1,\n\diff^{-1},\psi_\pi)\lambda(1,\m\diff^{-1},\psi_\pi)\\
-N(\q)^{3/4}\sum_{\substack{\n\subset \entier\\ N(\m)\leq M}}\frac{F(N(\n)/N(\q)^{1/2})}{N(\n)^{1/2}}\cdot 
\frac{\mu(\m)P\left(\frac{
\log(M/N(\m))}{\log(M)}\right)}{\psi(\m)N(\m)^{1/2}}\times\\
\sum_{\pi}^h\lambda(1,\n\q\diff^{-1},\psi_\pi)\lambda(1,\m\diff^{-1},\psi_\pi)
\end{multline}
la somme en $\pi$ portant bien entendu sur $\Pi_{\entier}^{\g{k}}$. On va traiter les deux termes apparaissant dans cette 
expression s\'epar\'ement, en notant : $$M_1(\q,\psi_\pi)=A(\q)-B(\q)$$ dans l'expression pr\'ec\'edente. Notons pour commencer 
que, puisque $N(\m)\leq M<N(\q)$, on a 
\begin{multline}\lambda(1,\m\diff^{-1},\psi_\pi)= \Big(N(\q)(1-N(\q)^{-2})(1+N(\q)^{-1})L(1,
{\rm{sym}}^2\pi_\q)\Big)^{1/2}\times\\
\bigg(\frac{-1}{N(\q)+1}\lambda_\pi(\q)\lambda_\pi(\m)\bigg)\nonumber
\end{multline}
en outre les termes $(1-N(\q)^{-2})(1+N(\q)^{-1})$ n'influencent pas les expressions en question, asymptotiquement en $N(\q)$, 
de sorte qu'on note toujours $A(\q),B(\q)$ les expressions simplifi\'ees suivantes:
\begin{multline}
A(\q)=N(\q)^{1/4}\sum_{\substack{\n\subset \entier\\ N(\m)\leq M}}\frac{F(N(\n)/N(\q)^{1/2})}{N(\n)^{1/2}}\cdot 
\frac{\mu(\m)P\left(\frac{
\log(M/N(\m))}{\log(M)}\right)}{\psi(\m)N(\m)^{1/2}}\times\\
\sum_{\pi}^h\Big(N(\q)L(1,{\rm{sym}}^2\pi_\q)\Big)
\bigg(\frac{-1}{N(\q)+1}\lambda_\pi(\q)\lambda_\pi(\n)+\lambda_\pi(\n\q^{-1})\mathbbmss{1}_{\q|\n}\bigg)\times\\
\bigg(\frac{-1}{N(\q)+1}\lambda_\pi(\q)\lambda_\pi(\m)\bigg)
\end{multline}
\begin{multline}
B(\q)=N(\q)^{3/4}\sum_{\substack{\n\subset \entier\\ N(\m)\leq M}}\frac{F(N(\n)/N(\q)^{1/2})}{N(\n)^{1/2}}\cdot 
\frac{\mu(\m)P\left(\frac{
\log(M/N(\m))}{\log(M)}\right)}{\psi(\m)N(\m)^{1/2}}\times\\
\sum_{\pi}^h\Big(N(\q)L(1,{\rm{sym}}^2\pi_\q)\Big)
\bigg(\frac{-1}{N(\q)+1}\lambda_\pi(\q)\lambda_\pi(\n\q)+\lambda_\pi(\n)\bigg)\times\\
\bigg(\frac{-1}{N(\q)+1}\lambda_\pi(\q)\lambda_\pi(\m)\bigg).
\end{multline}
Pour majorer ces expressions, on note d'abord que $||\psi_\pi||_{\mathcal{H}_\q^{\g{k}}}^2=||\varphi_\pi||
_{\mathcal{H}_\q^{\g{k}}}^2=(N(\q)+1)||\psi_\pi||_{\mathcal{H}_{\entier}^{\g{k}}}^2$; en outre, on a la factorisation 
suivante du facteur local du carr\'e sym\'etrique, en termes des param\`etres de Langlands de $\pi$:
\begin{multline}
L(1,{\rm{sym}}^2\pi_\q)=(1-\alpha_{\pi,1}(\q)^2N(\q)^{-1})^{-1}(1-\alpha_{\pi,1}(\q)\alpha_{\pi,2}(\q)N(\q)^{-1})^{-1}\times\\
(1-\alpha_{\pi,2}(\q)^2N(\q)^{-1})^{-1}\nonumber
\end{multline}
 et ceci implique que 
$$L(1,{\rm{sym}}^2\pi_\q)\ll 1$$
Traitons d'abord le terme $A(\q)$:
\begin{multline}
A(\q)=N(\q)^{1/4}\sum_{\substack{\n\subset \entier\\ N(\m)\leq M}}\frac{F(N(\n)/N(\q)^{1/2})}{N(\n)^{1/2}}\cdot 
\frac{\mu(\m)P\left(\frac{
\log(M/N(\m))}{\log(M)}\right)}{\psi(\m)N(\m)^{1/2}}\times\\
\sum_{\pi}\frac{\Gamma(\g{k-1})N(\q)L(1,{\rm{sym}}^2\pi_\q)}{|\disc|^{1/2}(4\pi)^{\g{k-1}}(N(\q)+1)||\varphi_\pi||
_{\mathcal{H}_{\entier}^{\g{k}}}^2}\times\\
\bigg(\frac{-1}{N(\q)+1}\lambda_\pi(\q)\lambda_\pi(\n)+\lambda_\pi(\n\q^{-1})\mathbbmss{1}_{\q|\n}\bigg)
\bigg(\frac{-1}{N(\q)+1}\lambda_\pi(\q)\lambda_\pi(\m)\bigg)
\end{multline}
soit, avec les estimations ci-dessus:
\begin{multline}
A(\q)\ll_F N(\q)^{1/4}\sum_{\substack{\n\subset \entier\\ N(\m)\leq M}}\left|\frac{F(N(\n)/N(\q)^{1/2})}{N(\n)^{1/2}}\cdot 
\frac{\mu(\m)P\left(\frac{
\log(M/N(\m))}{\log(M)}\right)}{\psi(\m)N(\m)^{1/2}}\right|\times\\
\sum_{\pi}\frac{1}{||\varphi_\pi||_{\mathcal{H}_{\entier}^{\g{k}}}^2}
\bigg(\frac{1}{(N(\q)+1)^2}\lambda_\pi(\q)^2|\lambda_\pi(\n)\lambda_\pi(\m)|\bigg)\\
+N(\q)^{1/4}\sum_{\substack{\q| \n\\ N(\m)\leq M}}\left|\frac{F(N(\n)/N(\q)^{1/2})}{N(\n)^{1/2}}\cdot 
\frac{\mu(\m)P\left(\frac{
\log(M/N(\m))}{\log(M)}\right)}{\psi(\m)N(\m)^{1/2}}\right|\times\\
\sum_{\pi}\frac{1}{||\varphi_\pi||_{\mathcal{H}_{\entier}^{\g{k}}}^2}
\bigg(\frac{1}{N(\q)+1}|\lambda_\pi(\q)\lambda_\pi(\n\q^{-1})\lambda_\pi(\m)|\bigg).\nonumber
\end{multline}
On r\'eindexe d'abord la deuxi\`eme somme en $\n$ en posant $\n'=\n\q^{-1}$: il appara\^it alors un facteur $F((2\pi)^dN(\n')N(\q))$ et en 
d\'enominateur un facteur $N(\q)^{1/2}$.  
On majore brutalement les coefficients $\lambda_\pi(\cdot)$ par Ramanujan, soit:
$$\lambda_\pi(\n)\ll_{\varepsilon} N(\n)^\varepsilon$$ (les bornes de Kim-Shahidi prouv\'ees dans \cite{KS}, $\lambda_\pi(\n)\ll_{\varepsilon} N(\n)^{1/9+\varepsilon}$, suffisent ici), en se 
rappelant qu'il n'y a qu'un nombre fini, \emph{ind\'ependant} de $\q$, de formes non ramifi\'ees:
\begin{multline}
A(\q)\ll_F N(\q)^{1/4}\sum_{\n\subset \entier}\left|\frac{F(N(\n)/N(\q)^{1/2})N(\n)^{\varepsilon}}{N(\n)^{1/2}}\right|\times\\ 
\sum_{N(\m)\leq M}\left|\frac{\mu(\m)P\left(\frac{
\log(M/N(\m))}{\log(M)}\right)N(\m)^\varepsilon}{\psi(\m)N(\m)^{1/2}}\right|\times
\frac{N(\q)^{2\varepsilon}}{(N(\q)+1)^2}\\
+N(\q)^{1/4}\sum_{ \n'\subset\entier}\left|\frac{F((2\pi)^d N(\n')N(\q)^{1/2})N(\n')^\varepsilon}{N(\n'\q)^{1/2}}\right|\times\\ \sum_{\N(\m)\leq M}
\left|\frac{\mu(\m)P\left(\frac{
\log(M/N(\m))}{\log(M)}\right)N(\m)^\varepsilon}{\psi(\m)N(\m)^{1/2}}\right|\times
\frac{N(\q)^\varepsilon}{N(\q)+1}.\nonumber
\end{multline}
On utilise 
alors que $F(y)\ll_A y^{-A}$ pour tout $A>0$, quand $y\to\infty$, et que $F(y)\ll 1$ pour $0<y<1$:
\begin{multline}
A(\q)\ll_{A} N(\q)^{1/4}N(\q)^{1/4+\varepsilon}M^{1/2+\varepsilon}(N(\q)+1)^{-2}\\+N(\q)^{1/4}\times\frac{N(\q)^{-A}}{N(\q)^{1/2}}\times M^{1/2+\varepsilon}(1+N(\q))^{-1}
\ll_{A,F} N(\q)^{-1}
\end{multline}
alors qu'on avait seulement besoin d'un majorant en $N(\q)^{1/4-\eta}$.\\

Pour $B(\q)$ on fait la m\^eme chose (sauf que les bornes de Kim-Shahidi sont ici insuffisantes):
\begin{multline}
B(\q)=N(\q)^{3/4}\sum_{\substack{\n\subset \entier\\ N(\m)\leq M}}\frac{F(N(\n)/N(\q)^{1/2})}{N(\n)^{1/2}}\cdot 
\frac{\mu(\m)P\left(\frac{
\log(M/N(\m))}{\log(M)}\right)}{\psi(\m)N(\m)^{1/2}}\times\\
\sum_{\pi}\frac{\Gamma(\g{k-1})N(\q)L(1,{\rm{sym}}^2\pi_\q)}{|\disc|^{1/2}(4\pi)^{\g{k-1}}(N(\q)+1)||\varphi_\pi||
_{\mathcal{H}_{\entier}^{\g{k}}}^2}
\bigg(\frac{-1}{N(\q)+1}\lambda_\pi(\q)\lambda_\pi(\n\q)+\lambda_\pi(\n)\bigg)\times\\
\bigg(\frac{-1}{N(\q)+1}\lambda_\pi(\q)\lambda_\pi(\m)\bigg)
\end{multline}
la m\^eme proc\'edure donne:
\begin{multline}
B(\q)\ll_F N(\q)^{3/4}N(\q)^{1/4+4\varepsilon}M^{1/2+\varepsilon}(N(\q)+1)^{-2}\\+N(\q)^{3/4}N(\q)^{1/4+2\varepsilon}M^{1/2+\varepsilon}(N(\q)+1)^{-1}\ll_F N(\q)^{\Delta/4}
\end{multline}
ce qui est bien une majoration en $N(\q)^{1/4-\eta}$, comme voulu.\\

L'\'etude de la contribution des formes anciennes au second moment se fait selon les m\^emes lignes, et ne pr\'esente pas 
de difficult\'e de nouvel ordre: nous ne l'incluons donc pas au pr\'esent travail. Il en va de m\^eme pour la contribution des formes anciennes dans les premier et second moments de la d\'eriv\'ee.

\subsection{La contribution des formes anciennes au nombre de formes modulaires de Hilbert: fin de la proposition \ref{nombre} }\label{app2}
 Les param\`etres des sommes sont 
$\E\subset \entier, \bar{\A}\in\mathscr{C}\ell^+(F),\xi\in(\A^{-1})^{\gg 0}/\unite$ ($\A$ parcourt toujours un ensemble fini).
Il s'agit de :
\begin{multline}
\sum_{\E,\bar{\A},\xi}
\sum_{\varphi}^h \frac{\lambda(\xi^2,\A^2\diff^{-1},\varphi)\overline{\lambda(1,\diff^{-1},\varphi)}\chi_q(\E)}{N(\E^2\xi\A)}\times\\
\int_{(\delta)} \left(\frac{N(\q)}{N(\E^2\xi\A)}\right)^s
L(s+1,{\rm{sym}}^2\pi_\infty)\frac{ds}{s}+\\
\sum_{\E,\bar{\A},\xi}
\sum_{\varphi}^h \frac{\lambda(\xi^2,\A^2\diff^{-1},\varphi)\overline{\lambda(1,\diff^{-1},\varphi)}\chi_q(\E)}{N(\E^2\xi\A)}\times\\
\int_{(\delta)} \left(\frac{N(\q)}{N(\E^2\xi\A)}\right)^s
L(s+1,{\rm{sym}}^2\pi_\infty)\frac{ds}{s+1}.\nonumber
\end{multline}
La base orthogonale des formes anciennes choisie est  $\{\varphi\}=\{\varphi_\pi,\psi_\pi\}_{\pi\in\Pi_{\entier}^{\g{k}}}$. 
On rappelle que $\sum^h_\varphi x_\varphi=\sum_\varphi\omega_\varphi x_\varphi$ avec $$\omega_\varphi=
\frac{ \Gamma(\g{k}-\g{1}) }{(4\pi)^{\g{k}-\g{1}}|\disc|^{1/2}||\varphi||^2_{\mathcal{H}_\q^{\g{k}}}}$$
et $||\varphi||^2_{\mathcal{H}_\q^{\g{k}}}=(N(\q)+1)||\varphi||^2$. La contribution des formes $\{\varphi_\pi\}_ {\pi\in\Pi
_{\entier}^{\g{k}}}$ est imm\'ediate (par exemple apr\`es application de la formule de Peterson \`a $\Pi_{\entier}^{\g{k}}$):

\begin{multline}
\sum_{\E,\bar{\A},\xi}
\sum_{\pi\in\Pi
_{\entier}^{\g{k}}}^h \frac{\lambda(\xi^2,\A^2\diff^{-1},\varphi_\pi)\overline{\lambda(1,\diff^{-1},\varphi_\pi)}\chi_q(\E)}{N(\E^2\xi\A)}\times\\
\int_{(\delta)} \left(\frac{N(\q)}{N(\E^2\xi\A)}\right)^s
L(s+1,{\rm{sym}}^2\pi_\infty)\frac{ds}{s}+\\
\sum_{\E,\bar{\A},\xi}
\sum_{\pi\in\Pi
_{\entier}^{\g{k}}}^h \frac{\lambda(\xi^2,\A^2\diff^{-1},\varphi_\pi)\overline{\lambda(1,\diff^{-1},\varphi_\pi)}\chi_q(\E)}{N(\E^2\xi\A)}\times\\
\int_{(\delta)} \left(\frac{N(\q)}{N(\E^2\xi\A)}\right)^s
L(s+1,{\rm{sym}}^2\pi_\infty)\frac{ds}{s+1}
\ll_{\delta}N(\q)^{-1+\delta}.\nonumber
\end{multline}

En ce qui concerne la contribution des $\psi_\pi$, on a :
$$\lambda(\xi^2,\A\diff^{-1},\psi_\pi)=\sqrt{N(\q)L(1,{\rm{sym}}^2\pi_\q)} \left( \frac{-1}{N(\q)+1}\lambda_\pi(\xi\A)+
\lambda_\pi(\xi\A\q^{-1})\mathbbmss{1}_{\q|\xi\A}\right)$$
$$\lambda(1,\diff^{-1},\psi_\pi)=-\sqrt{N(\q)L(1,{\rm{sym}}^2\pi_\q)} \frac{\lambda_\pi(\q)}{N(\q)+1}$$
et donc:
\begin{multline}
\sum_{\E,\bar{\A},\xi}
\sum_{\pi\in\Pi_{\entier}^{\g{k}}}^h 
\frac{\lambda(\xi^2,\A^2\diff^{-1},\psi_\pi)\overline{\lambda(1,\diff^{-1},\psi_\pi)}\chi_q(\E)}{N(\E^2\xi\A)}\times\\
\int_{(\delta)} \left(\frac{N(\q)}{N(\E^2\xi\A)}\right)^s
L(s+1,{\rm{sym}}^2\pi_\infty)\frac{ds}{s}+\\
\sum_{\E,\bar{\A},\xi}
\sum_{\pi\in\Pi
_{\entier}^{\g{k}}}^h \frac{\lambda(\xi^2,\A^2\diff^{-1},\psi_\pi)\overline{\lambda(1,\diff^{-1},\psi_\pi)}\chi_q(\E)}{N(\E^2\xi\A)}\times\\
\int_{(\delta)} \left(\frac{N(\q)}{N(\E^2\xi\A)}\right)^s
L(s+1,{\rm{sym}}^2\pi_\infty)\frac{ds}{s+1}\nonumber
\end{multline}

\begin{multline}
=\sum_{\E,\bar{\A},\xi}\frac{\chi_q(\E)}{N(\E^2\xi\A)}
\int_{(\delta)} \left(\frac{N(\q)}{N(\E^2\xi\A)}\right)^s
L(s+1,{\rm{sym}}^2\pi_\infty)\frac{ds}{s}\times \\ \sum_{\pi\in\Pi_{\entier}^{\g{k}}}
\frac{\Gamma(\g{k}-\g{1})N(\q)L(1,{\rm{sym}}^2\pi_\q)}{(4\pi)^{\g{k}-\g{1}}|\disc|^{1/2}(N(\q)+1)||\psi_\pi||^2}
 \left( \frac{-\lambda_\pi(\xi\A) }{N(\q)+1} +
\lambda_\pi(\xi\A\q^{-1})\mathbbmss{1}_{\q|\xi\A}\right) 
\frac{\lambda_\pi(\q)}{N(\q)+1}\\
+\sum_{\E,\bar{\A},\xi}\frac{\chi_q(\E)}{N(\E^2\xi\A)}
\int_{(\delta)} \left(\frac{N(\q)}{N(\E^2\xi\A)}\right)^s
L(s+1,{\rm{sym}}^2\pi_\infty)\frac{ds}{s+1}\times \\ \sum_{\pi\in\Pi_{\entier}^{\g{k}}}
\frac{\Gamma(\g{k}-\g{1})N(\q)L(1,{\rm{sym}}^2\pi_\q)}{(4\pi)^{\g{k}-\g{1}}|\disc|^{1/2}(N(\q)+1)||\psi_\pi||^2}
 \left( \frac{-\lambda_\pi(\xi\A) }{N(\q)+1} +
\lambda_\pi(\xi\A\q^{-1})\mathbbmss{1}_{\q|\xi\A}\right) 
\frac{\lambda_\pi(\q)}{N(\q)+1}.\nonumber
\end{multline}
On termine en utilisant la borne (grossi\`ere, mais suffisante)
$$\lambda_\pi(\n)\ll N(\n)^{1/4}$$
cons\'equence en fait de la formule de Petersson, la majoration  $L(1,{\rm{sym}}^2\pi_\q)\ll 1$ et donc la contribution des $\psi_\pi$ est domin\'ee par:

\begin{multline}
\sum_{\E,\bar{\A},\xi}\frac{\chi_q(\E)}{N(\E^2\xi\A)}
\int_{(\delta)} \left(\frac{N(\q)}{N(\E^2\xi\A)}\right)^s
L(s+1,{\rm{sym}}^2\pi_\infty)\frac{ds}{s}\times \\ \sum_{\pi\in\Pi_{\entier}^{\g{k}}}
\frac{\Gamma(\g{k}-\g{1})N(\q)L(1,{\rm{sym}}^2\pi_\q)}{(4\pi)^{\g{k}-\g{1}}|\disc|^{1/2}(N(\q)+1)||\psi_\pi||^2}\times\\
 \left( \frac{N(\xi\A)^{1/2}N(\q)^{1/4} }{N(\q)+1} +
N(\xi\A)^{1/2}N(\q)^{-1/4}\mathbbmss{1}_{\q|\xi\A}\right) 
\frac{N(\q)^{1/4}}{N(\q)+1}\\
+\sum_{\E,\bar{\A},\xi}\frac{\chi_q(\E)}{N(\E^2\xi\A)}
\int_{(\delta)} \left(\frac{N(\q)}{N(\E^2\xi\A)}\right)^s
L(s+1,{\rm{sym}}^2\pi_\infty)\frac{ds}{s+1}\times \\ \sum_{\pi\in\Pi_{\entier}^{\g{k}}}
\frac{\Gamma(\g{k}-\g{1})N(\q)L(1,{\rm{sym}}^2\pi_\q)}{(4\pi)^{\g{k}-\g{1}}|\disc|^{1/2}(N(\q)+1)||\psi_\pi||^2}\times\\
 \left( \frac{N(\xi\A)^{1/2}N(\q)^{1/4} }{N(\q)+1} +
N(\xi\A)^{1/2}N(\q)^{-1/4}\mathbbmss{1}_{\q|\xi\A}\right) 
\frac{N(\q)^{1/4}}{N(\q)+1}\\
\ll N(\q)^{-1} \sum_{\E,\bar{\A},\xi}\frac{\chi_q(\E)}{N(\E^2)N(\xi\A)^{1/2}}
\int_{(\delta)} \left(\frac{N(\q)}{N(\E^2\xi\A)}\right)^s
L(s+1,{\rm{sym}}^2\pi_\infty)(\frac{ds}{s}+\frac{ds}{s+1}).\nonumber
\end{multline}
En d\'epla\c cant la ligne d'int\'egration de $(\delta)$ en $(1/2+\delta)$, pout $\delta>0$, de sorte \`a ce que la 
somme en $\xi,\A$ soit absolument convergente, on trouve
\begin{multline}
\sum_{\E,\bar{\A},\xi}\frac{\chi_q(\E)}{N(\E^2)N(\xi\A)^{1/2}}
\int_{(1/2+\delta)} \left(\frac{N(\q)}{N(\E^2\xi\A)}\right)^s
L(s+1,{\rm{sym}}^2\pi_\infty)(\frac{ds}{s}+\frac{ds}{s+1})\\
\ll N(\q)^{\frac{1}{2}+\delta}\nonumber
\end{multline}
et donc une contribution n\'egligeable. Ainsi s'ach\`eve la preuve de la proposition \ref{nombre}.

\subsection{Un th\'eor\`eme de densit\'e}\label{densite}
Soit $\eta>0$ un r\'eel (petit) fix\'e, $T>0$,
\begin{eqnarray}
R(\eta, T)=\{ s\in\C;\, 1-\eta<\Re(s)<1, |\Im(s)|\leq T\}
\end{eqnarray}
 et:
\begin{eqnarray}
\Pi_\q^{\g{k}}(\eta,T)=\{\pi\in\Pi_\q^{\g{k}};\,L(s,{\rm{sym}}^2\pi_f)\neq 0, \, \forall s \in R(\eta,T) \}.
\end{eqnarray}
 Le but de cette courte section est d'\'etablir une borne du type:
\begin{eqnarray}\label{dens}
|\Pi_\q^{\g{k}}\setminus \Pi_\q^{\g{k}}(\eta,T)| \ll_{\eta} T^D N(\q)^{b\eta}.
\end{eqnarray}
$D,b>0$ \'etant des r\'eels que l'on ne cherchera pas \`a optimiser.

Ce type de majorations est appel\'e \og{} th\'eor\`eme de densit\'e \fg{}, et le travail de Kowalski et Michel \cite{KM2} sur ce probl\`eme (quand le niveau varie, commme c'est notre cas) montre que des familles tr\`es g\'en\'erales de formes automorphes sur $\Q$ satisfont \`a celles-ci; ce qui suit est inspir\'e de leur preuve, \`a ceci pr\`es que nous \'evitons d'utiliser la borne de Ramanujan, comme Luo \cite{luo2} dans son th\'eor\`eme de densit\'e sur le carr\'e sym\'etrique (quand les param\`etres \`a l'infini varient): il est d'ailleurs pr\'ecis\'e dans \cite{KM2} que les bornes
$$\lambda_\pi(\n)\ll N(\n)^{\delta}$$
 pour $\delta\in ]0,1/4[$ uniforme en $\pi$ , suffisent -- ce qui pour $\gl$ est connu. En outre, l'article \cite{LW} de Lau et Wu donne une simplification de la preuve de Luo dans le contexte des puissances sym\'etriques, quand le poids varie, inspir\'ee de \cite{KM2}.

Soit donc $\pi$ une forme modulaire de Hilbert, de poids $\g{k}$ et de conducteur $\q$ -- pour le moment quelconque. Soit ${\rm{sym}}^2\pi$ le carr\'e sym\'etrique de $\pi$, et pour $\Re(s)>1$:
$$L(s,{\rm{sym}}^2\pi_f)=\sum_{\n\subset \entier} \rho_\pi(\n) N(\n)^{-s}$$
sa fonction $L$ finie, qui se prolonge en une fonction analytique sur le plan complexe.

 On aura besoin d'un amollisseur, approximant la valeur de $L(s,{\rm{sym}}^2\pi_f)^{-1}$; la valeur explicite:
$$L(s,{\rm{sym}}^2\pi_f)^{-1}=\sum_{\n\subset \entier} \rho_\pi^{-}(\n) N(\n)^{-s}$$
avec -- quand le conducteur $\q$ de $\pi$ est sans facteurs carr\'es:
$$\rho_\pi^{-}(\n'\n'')=\rho_\pi^{-}(\n')\rho_\pi^{-}(\n'') $$
pour $\n',\n''\subset\entier$, $\n'|\q^\infty$ et $(\n'',\q)=1$ et:

\begin{displaymath}
\rho_\pi^{-}(\n)=\left\{ \begin{array}{ccc}
\mu(\n)N(\n)^{-1} & \textrm{ si $\n|\q^\infty$ }\\
\mu(\A\B\cc)\mu(\B)\lambda_\pi(\A^2\B^2) & \textrm{ si $\n=\A\B^2\cc^3$ et $(\A,\B,\cc)=(\n,\q)=1$} \\
0 & \textrm{ sinon}
\end{array}\right.
\end{displaymath}
peut \^etre utilis\'ee pour d\'efinir l'amollisseur n\'ecessaire-- Luo le fait avec succ\`es. On utilise alors la formule de Petersson pour \'etablir une in\'egalit\'e de grand crible pour le carr\'e sym\'etrique; cependant elle a l'inconv\'enient de g\'en\'erer des calculs compliqu\'es. Nous allons donc contourner cette difficult\'e, selon \cite{LW}, inspir\'es par \cite{KM2}.

La fonction $L$ du carr\'e sym\'etrique admet le d\'evloppement eul\'erien:
$$L(s,{\rm{sym}}^2\pi_f)=\prod_\p L(s,{\rm{sym}}^2\pi_\p)$$
avec 
$$L(s,{\rm{sym}}^2\pi_\p)=\prod_{1\leq j\leq 2} (1-\alpha_\pi(\p)^{2-2j} N(\p)^{-s})^{-1}=(1-\rho_\pi(\p)N(\p)^{-s}+\ldots)^{-1}$$
quand $\pi$ est non ramifi\'ee en $\p$, et $(\alpha_\pi(\p),\alpha_\pi(\p)^{-1})$ sont les param\`etres de Langlands en $\p$. Quand $\pi$ est ramifi\'ee en $\p$, on peut quand m\^eme \'ecrire sous forme d'inverse d'un polyn\^ome de degr\'e au plus 3 en $N(\p)^{-s}$:
$$L(s,{\rm{sym}}^2\pi_\p)=(1-\rho_\pi(\p)N(\p)^{-s}+\cdots)^{-1}$$
Par exemple, quand le conducteur $\q$ de $\pi$ est sans facteur carr\'e, et que $\p|\q$ on a:
$$L(s,{\rm{sym}}^2\pi_\p)=(1-N(\p)^{-1}N(\p)^{-s})^{-1}$$


Notons qu'avec les bornes de Kim-Shahidi ($|\alpha_\pi(\p)|\leq N(\p)^{\frac{1}{9}}$) , le facteur:
$$(1-\rho_\pi(\p)N(\p)^{-s})$$
ne s'annule pas quand $\Re(s)\geq \frac{3}{4}$ et $N(\p)\geq 7$, ce qui permet d'\'ecrire, avec $\n_0$ produit des id\'eaux maximaux de norme moindre que 7, et $\Re(s)>1$:
$$L(s,{\rm{sym}}^2\pi_f)^{-1}=G_\pi(s)\sum_{(\n,\n_0)=\entier}\mu(\n)\rho_\pi(\n)N(\n)^{-s}$$
avec:
$$G_\pi(s)=\prod_{N(\p)<7}L(s,{\rm{sym}}^2\pi_\p)^{-1}\times\prod_{N(\p)\geq 7}\Big( (1-\rho_\pi(\p)N(\p)^{-s})^{-1}\times L(s,{\rm{sym}}^2\pi_\p)^{-1}\Big)$$
qui converge absolument, car pour $N(\p)\geq 7$, le terme g\'en\'eral du produit eul\'erien est $1+\mathcal{O}\left(N(\p)^{\frac{2}{9}-2s}\right)$ -- et uniform\'ement par rapport \`a $\pi$ -- dans le domaine $\Re(s)>\frac{3}{4}$; de plus on a:
$$G_\pi(s)\ll 1$$
dans cette m\^eme r\'egion. En fait cela est vrai dans un domaine plus grand, mais le th\'eor\`eme de densit\'e n'est non-trivial que pour $s$ proche de $1$ de toute fa\c con.

On d\'efinit alors un amollisseur $M_\pi^X(s)$:
\begin{eqnarray}
M_\pi^X(s):=G_\pi(s)\sum_{\substack{(\n,\n_0)=\entier\\N(\n)\leq X}}\mu(\n)\rho_\pi(\n)N(\n)^{-s}
\end{eqnarray}
Notons tout de suite la relation:

\begin{multline}\nonumber
1-L(s,{\rm{sym}}^2\pi_f)M_\pi^X(s)=L(s,{\rm{sym}}^2\pi_f)\cdot G_\pi(s)\sum_{\substack{(\n,\n_0)=\entier\\N(\n)\geq X}}\mu(\n)\rho_\pi(\n)N(\n)^{-s}
\end{multline}

Soient $Y>0$, $\sigma>0$, on a:
$$\exp(-1/Y)=\frac{1}{2i\pi}\int_{(\sigma)}\Gamma(w)Y^w dw$$
soit, pour tout $\rho\in\C$ tel que $\Re(\rho+\sigma)\geq 3/4$:
\begin{multline}\nonumber
\exp(-1/y)= \frac{1}{2i\pi}\int_{(\sigma)}\Big(1-L(\rho+w,{\rm{sym}}^2\pi_f)M_\pi^X(\rho+w)\Big) \Gamma(w)Y^w dw\\+
\frac{1}{2i\pi}\int_{(\sigma)}L(\rho+w,{\rm{sym}}^2\pi_f)M_\pi^X(\rho+w) \Gamma(w)Y^w dw
\end{multline}

Choisissons $\sigma=\sigma_1=1-\Re(\rho)+\varepsilon_1$ 
 supposons que $\rho\in R(\eta,T)$ est un z\'ero de $L(.,{\rm{sym}}^2\pi_f)$ -- \emph{que l'on devrait noter $ \rho=\rho_\pi$, car on fait un tel choix pour chaque forme $\pi$ quand il existe}. On peut alors d\'eplacer la ligne d'int\'egration de la deuxi\`eme int\'egrale de $(\sigma)$ \`a $(\sigma_2)$ avec $\sigma_2=3/4-\Re(\rho)<0$ d\`es que $\eta < 1/4$; comme $L(\rho+w,{\rm{sym}}^2\pi_f)$ s'annule en $w=0$, le changement de droite d'int\'egration n'introduit pas de nouveau terme, soit:
\begin{multline}\nonumber
\exp(-1/Y)= \frac{1}{2i\pi}\int_{(\sigma_1)}\Big(1-L(\rho+w,{\rm{sym}}^2\pi_f)M_\pi^X(\rho+w)\Big) \Gamma(w)Y^w dw\\+
\frac{1}{2i\pi}\int_{(\sigma_2)}L(\rho+w,{\rm{sym}}^2\pi_f)M_\pi^X(\rho+w) \Gamma(w)Y^w dw
\end{multline}

On \'ecrit $$X=N(\q)^A\, , \, Y=N(\q)^B, $$ les valeurs de $A,B>0$ seront fix\'ees plus loin; mais d\'ej\`a on peut dire qu'\`a cause de la croissance polyn\^omiale de $M_\pi^X$ et $L(.,{\rm{sym}}^2\pi_f)$ (bornes de convexit\'e), et de la d\'ecroissance exponentielle de la fonction $\Gamma$:
\begin{multline}\nonumber
\int_{ |\Im(w)|\geq \log(N(\q))^2}L(\rho+w,{\rm{sym}}^2\pi_f)M_\pi^X(\rho+w) \Gamma(w)Y^w dw \\ \ll T^D N(\q)^C \exp(-\log(N(\q))^2)
\end{multline}
pour un certain $D,C>0$, et de m\^eme pour la premi\`ere int\'egrale. On en d\'eduit donc, avec $\rho$ un tel z\'ero de $L(.,{\rm{sym}}^2\pi_f)$ et  $I(\q)=[-\log(N(\q))^2,\log(N(\q))^2]$:

\begin{multline}\nonumber
1\ll  T^D\left|\int_{I(\q)}\Big(1-L(\rho+\sigma_1+it,{\rm{sym}}^2\pi_f)M_\pi^X(\rho+\sigma_1+it)\Big) \Gamma(\sigma_1+it)Y^{\sigma_1+it} dt\right|\\+
T^D\left|\int_{I(\q)}L(\rho+\sigma_2+it,{\rm{sym}}^2\pi_f)M_\pi^X(\rho+\sigma_2+it) \Gamma(\sigma_2+it)Y^{\sigma_2+it} dt\right|\\
\ll T^DY^{2\eta}\log(N(\q))^2\int_{I(\q)}\Big|1-L(\rho+\sigma_1+it,{\rm{sym}}^2\pi_f)M_\pi^X(\rho+\sigma_1+it)\Big|^2  dt\\
+T^DY^{\eta-1/4}\int_{I(\q)}\Big|L(\rho+\sigma_2+it,{\rm{sym}}^2\pi_f)M_\pi^X(\rho+\sigma_2+it)\Big| dt\\
\end{multline}

On en d\'eduit donc, par positivit\'e, en assignant n'importe quelle valeur \`a $ \rho$ pour les $ \pi$ n'ayant pas un z\'ero dans $R(\eta,T)$:
\begin{multline}
\Pi_\q^{\g{k}}(\eta,T)\ll \\ T^DY^{2\eta}\log(N(\q))^2\int_{I(\q)}\sum_{\pi\in\Pi_\q^{\g{k}}}\Big|1-L(\rho+\sigma_1+it,{\rm{sym}}^2\pi_f)M_\pi^X(\rho+\sigma_1+it)\Big|^2  dt\\
+T^DY^{\eta-1/4}\int_{I(\q)}\sum_{\pi\in\Pi_\q^{\g{k}}}\Big|L(\rho+\sigma_2+it,{\rm{sym}}^2\pi_f)M_\pi^X(\rho+\sigma_2+it)\Big| dt\\
=I(1)+I(2)
\end{multline}
\noindent $\bullet$ La deuxi\`eme int\'egrale $I(2)$ est la plus facile \`a estimer: en effet, 
$$M_\pi^X(\rho+\sigma_2+it)\ll \sum_{N(\n)\leq X}N(\n)^{1/9-3/4}\ll X^{1/2}$$
et l'in\'egalit\'e de convexit\'e donne pour $\Re(s)=3/4$:
$$L(s,{\rm{sym}}^2\pi_f)\ll |s|^AN(\q)^{1/4+\varepsilon} .$$
On en d\'eduit:
\begin{eqnarray}
I(2)\ll_{\varepsilon}T^D Y^{\eta-1/4}X^{1/2}N(\q)^{1/4+\varepsilon} 
\end{eqnarray}

\noindent $\bullet$ En ce qui concerne $I(1)$, on rappelle que:

\begin{multline}\nonumber
\Big|1-L(\rho+\sigma_1+it,{\rm{sym}}^2\pi_f)M_\pi^X(\rho+\sigma_1+it)\Big|^2 \\ =|L(\rho+\sigma_1+it,{\rm{sym}}^2\pi_f)|^2\cdot |G_\pi(\rho+\sigma_1+it)|^2\Big|\sum_{\substack{(\n,\n_0)=\entier\\N(\n)\geq X}}\mu(\n)\rho_\pi(\n)N(\n)^{-\rho-\sigma_1-it}\Big|^2
\end{multline}

Avec la borne $L(1,{\rm{sym}}^2\pi_f)\ll_{\varepsilon} N(\q)^\varepsilon $, on en d\'eduit:

\begin{multline}\nonumber
\Big|1-L(\rho+\sigma_1+it,{\rm{sym}}^2\pi_f)M_\pi^X(\rho+\sigma_1+it)\Big|^2 \\ \ll_{\varepsilon}N(\q)^\varepsilon \Big|\sum_{\substack{(\n,\n_0)=\entier\\X \leq N(\n)\leq \exp(\log(N(\q))^2)}}\mu(\n)\rho_\pi(\n)N(\n)^{-\rho-\sigma_1-it}\Big|^2 + N(\q)^{\varepsilon-1}
\end{multline}
ce qui donne:

\begin{multline}\nonumber
I(1)  \ll_{\varepsilon}T^D Y^{2\eta}N(\q)^\varepsilon\times\\ \int_{I(\q)}\sum_{\pi\in\Pi_\q^{\g{k}}} \Big|\sum_{\substack{(\n,\n_0)=\entier\\X \leq N(\n)\leq \exp(\log(N(\q))^2)}}\mu(\n)\rho_\pi(\n)N(\n)^{-\rho-\sigma_1-it}\Big|^2 dt
\end{multline}

C'est ainsi qu'une in\'egalit\'e de grand crible nous permet d'achever la preuve. Essentiellement, l'article de Duke-Kowalski \cite{DK} permet de traiter les familles finies de formes automorphes $ \mathcal{F}$ sur $ {\rm GL}(n) / \Q$, dont le conducteur analytique tend vers l'infini quand $ |\mathcal{F}|$ tend vers l'infini, sous Ramanujan (en fait sous $ \mathcal{H}_\delta$: \og $|\alpha_{i,\pi(\p)}|\ll N(\p)^{1/4-\delta} $\fg{} pour n'importe quel $\delta >0$). Le cas du carr\'e sym\'etrique pose le l\'eger probl\`eme de multiplicit\'e, r\'esolu dans l'appendice  de \cite{DK} par Ramakrishnan par des arguments galoisiens, puis purement automorphes dans \cite{R}, ce qui permet d'\'etendre le r\'esultat \`a toutes les familles de formes sur $ {\rm GL}(2)$, holomorphes ou de Maass.\\

En d'autres termes, on doit trouver une majoration de 
$$\sum_{\pi\in\Pi_\q^{\g{k}}} \Big|\sum_{\substack{(\n,\n_0)=\entier\\ N \leq N(\n)\leq 2N}}\rho_\pi(\n) x_\n\Big|^2  $$
meilleure que la triviale $\displaystyle{|\Pi_\q^{\g{k}}| T \sum_{ N \leq N(\n)\leq 2N}|x_\n|^2 }$. 
Parce qu'un op\'erateur a la m\^eme norme hilbertienne que son adjoint, on peut majorer de fa\c con \'equivalente:
$$\sum_{\substack{(\n,\n_0)=\entier\\ N \leq N(\n)\leq 2N}} \Big| \sum_{\pi\in\Pi_\q^{\g{k}}}\rho_\pi(\n) x_\pi \Big|^2 .$$
En d\'eveloppant le carr\'e, on est conduit \`a:
$$ \sum_{\pi\in\Pi_\q^{\g{k}}}\sum_{\pi'\in\Pi_\q^{\g{k}}}x_\pi\overline{x_{\pi'}}\sum_{\substack{(\n,\n_0)=\entier\\ N \leq N(\n)\leq 2N}}\rho_\pi(\n)\overline{\rho_{\pi'}(\n)} $$
Soit $\Phi$ une fonction lisse, \`a support contenu dans $[N-1,2N+1] $, de sorte que: 
\begin{multline}\nonumber
\sum_{\substack{(\n,\n_0)=\entier\\ N \leq N(\n)\leq 2N}}\rho_\pi(\n)\overline{\rho_{\pi'}(\n)}=
\sum_{\substack{(\n,\n_0)=\entier\\ N \leq N(\n)\leq 2N}}\rho_\pi(\n)\overline{\rho_{\pi'(\n)}}\Phi(\n)\\
=\frac{1}{2i\pi}\int_{(2)}\widetilde{\Phi}(s)\sum_{(\n,\n_0)=\entier}\rho_\pi(\n)\overline{\rho_{\pi'}(\n)}N(\n)^{-s}ds.
\end{multline}

La somme interne est presque une fonction $L$ de Rankin-Selberg -- modulo les places ramifi\'ees. Plus pr\'ecis\'ement, pour $\Re(s)> 1 $:

$$ \sum_{(\n,\n_0)=\entier}\rho_\pi(\n)\overline{\rho_{\pi'}(\n)}N(\n)^{-s}= H(s,{\rm{sym}}^2\pi\otimes {\rm{sym}}^2\pi')L(s,{\rm{sym}}^2\pi\otimes {\rm{sym}}^2\pi')$$
$H(s,{\rm{sym}}^2\pi\otimes {\rm{sym}}^2\pi') $ \'etant un produit eul\'erien fini, portant sur les places de ramifications de $ \pi, \pi'$. Les bornes de Kim-Shahidi (c'est ici qu'intervient l'hypoth\`ese $ \mathcal{H}_\delta$) assurent que 
\begin{itemize}
\item La fonction $s\mapsto H(s,{\rm{sym}}^2\pi\otimes {\rm{sym}}^2\pi') $ converge sur le demi-plan $ \Re(s) \geq 1/2$.
\item On a la majoration uniforme $H(s,{\rm{sym}}^2\pi\otimes {\rm{sym}}^2\pi')\ll_\varepsilon N({\rm q}_{{\rm{sym}}^2\pi\otimes {\rm{sym}}^2\pi'})^\varepsilon $.
\end{itemize}
En d'autres termes: $H(\cdot,\cdot)$ est ignorable! Bushnell et Henniart ont donn\'e une majoration utile du conducteur de paires:
$$N({\rm q}_{\pi\otimes \pi'})\ll N({\rm q}_\pi)^2 N({\rm q}_{\pi'})^2 $$
on a donc:
\begin{multline}
\sum_{\substack{(\n,\n_0)=\entier\\ N \leq N(\n)\leq 2N}}\rho_\pi(\n)\overline{\rho_{\pi'}(\n)}\\=
\widetilde{\Phi}(1)H(1,{\rm{sym}}^2\pi\otimes {\rm{sym}}^2\pi'){\rm res}_{s=1}\Big( L(s,{\rm{sym}}^2\pi\otimes {\rm{sym}}^2\pi')\Big) \\ +
\frac{1}{2i\pi}\int_{(1/2)}\widetilde{\Phi}(s)H(s,{\rm{sym}}^2\pi\otimes {\rm{sym}}^2\pi')L(s,{\rm{sym}}^2\pi\otimes {\rm{sym}}^2\pi')ds\\
\ll N(\q)^\varepsilon\widetilde{\Phi}(1){\rm res}_{s=1}\Big( L(s,{\rm{sym}}^2\pi\otimes {\rm{sym}}^2\pi')\Big)\\
+N(\q)^\varepsilon\int_{(1/2)}|\widetilde{\Phi}(s)|\cdot|L(s,{\rm{sym}}^2\pi\otimes {\rm{sym}}^2\pi')|ds
\end{multline}

La th\'eorie de Rankin-Selberg dit que, \'etant donn\'e deux formes cuspidales $\Pi,\Pi'$,  $ L(s,\Pi\otimes \Pi')$ a un p\^ole en  $1+it$ si et seulement si $\Pi'\cong \Pi^{\vee}\otimes |{\rm det}|^{it} $. Dans le cas pr\'esent, $L(s,{\rm{sym}}^2\pi\otimes {\rm{sym}}^2\pi')$ ne peut avoir un p\^ole qu'en $ s=1$, et ceci si et seulement si $ {\rm{sym}}^2\pi'\cong ({\rm{sym}}^2\pi)^{\vee}$. 

Le travail de Ramakrishnan (\cite{R}, Theorem 4.1.2, Corollary 4.1.3) nous enseigne que cela implique:
\begin{itemize}
\item $(\pi')^{\vee}\cong\pi\otimes \chi $, pour $\chi $ caract\`ere id\'elique unitaire (quadratique si les caract\`eres centraux de $ \pi,\pi'$ co\"incident).
\item de plus, si les conducteurs de $ \pi, \pi'$ sont sans facteurs carr\'es, on a m\^eme $ (\pi')^{\vee}\cong \pi$.
\end{itemize}

Puisque la famille de formes qui nous int\'eresse est de conducteur l'id\'eal premier $\q $, c'est le second cas qui s'applique; il est cependant clair que l'on peut compter les caract\`eres quadratiques satisfaisant \`a la premi\`ere condition, quel que soit le niveau $ \n$.

On arrive donc \`a:
\begin{multline}
\sum_{\substack{(\n,\n_0)=\entier\\ N \leq N(\n)\leq 2N}} \Big| \sum_{\pi\in\Pi_\q^{\g{k}}}\rho_\pi(\n) x_\pi \Big|^2 \\
\ll_\varepsilon N(\q)^\varepsilon\widetilde{\Phi}(1)\sum_{\pi\in\Pi_\q^{\g{k}}}{\rm res}_{s=1}\Big( L(s,{\rm{sym}}^2\pi\otimes {\rm{sym}}^2\pi)\Big)|x_\pi|^2\\
+N(\q)^\varepsilon\int_{-\infty}^{+\infty}\Big|\widetilde{\Phi}\Big(\frac{1}{2}+it\Big)\Big|\sum_{\pi,\pi'\in\Pi_\q^{\g{k}}}\Big|L\Big(\frac{1}{2}+it,{\rm{sym}}^2\pi\otimes {\rm{sym}}^2\pi'\Big)\Big|x_{\pi}\overline{x_{\pi'}}dt\\
\ll N(\q)^\varepsilon \Big( N\sum_{\pi}|x_\pi|^2 + N^{1/2}N(\q)\sum_{\pi,\pi'}x_\pi\overline{x_{\pi'}}  \Big)
\end{multline}
par l'in\'egalit\'e de convexit\'e sur les fonctions $L$ de Rankin-Selberg: 
$$L\Big(\frac{1}{2}+it,\Pi\otimes \Pi'\Big)\ll_\varepsilon N({\rm q}_{\pi\otimes \pi'})^{1/4+\varepsilon}.$$
Enfin, l'in\'egalit\'e de Cauchy-Schwarz donne:

$$\sum_{\substack{(\n,\n_0)=\entier\\ N \leq N(\n)\leq 2N}} \Big| \sum_{\pi\in\Pi_\q^{\g{k}}}\rho_\pi(\n) x_\pi \Big|^2 
\ll_\varepsilon N(\q)^\varepsilon \Big( N + N^{1/2}N(\q)^2  \Big)
 \sum_{\pi}|x_\pi|^2.$$
 
 De fa\c con \'equivalente, on a montr\'e l'in\'egalit\'e de grand crible:
 $$\sum_{\pi\in\Pi_\q^{\g{k}}} \Big|\sum_{\substack{(\n,\n_0)=\entier\\ N \leq N(\n)\leq 2N}}\rho_\pi(\n) x_\n\Big|^2\ll_\varepsilon  N(\q)^\varepsilon \Big( N + N^{1/2}N(\q)^2  \Big)
 \sum_{N\leq \n\leq 2N}|x_\n|^2. $$
 
 Revenons maintenant \`a la majoration de l'int\'egrale $ I(1)$. La norme de l'id\'eal $ \n$ varie entre $X=N(\q)^A$ et $ \exp(\log(N(\q))^2)$; il faut d'abord d\'ecouper cet intervalle en intervalles dyadiques disjoints, en nombre moindre que $ \log(N(\q)^2)$, soit:
 $$[X, \exp(\log(N(\q))^2)]=\bigsqcup_{N {\textrm { dyadique}}} [N,2N]   $$
 puis, en majorant brutalement:
 \begin{multline}\nonumber
 \sum_{\pi\in\Pi_\q^{\g{k}}} \Big|\sum_{X \leq N(\n)\leq \exp(\log(N(\q))^2)}\mu(\n)\rho_\pi(\n)N(\n)^{-\rho-\sigma_1-it}\Big|^2\\ \leq 
 \sum_{\pi\in\Pi_\q^{\g{k}}} \log(N(\q)^2)\sup_{N}\Big|\sum_{N \leq N(\n)\leq 2N}\mu(\n)\rho_\pi(\n)N(\n)^{-\rho-\sigma_1-it}\Big|^2  
 \end{multline}
On se \og d\'ebarrasse\fg{} ainsi du terme en $\exp(\log(N(\q))^2)$ en appliquant l'in\'egalit\'e de grand crible:

\begin{multline}\nonumber
I(1)  \ll_{\varepsilon}T^DY^{2\eta} N(\q)^\varepsilon\times\\ \int_{I(\q)}\sum_{\pi\in\Pi_\q^{\g{k}}} \Big|\sum_{\substack{(\n,\n_0)=\entier\\X \leq N(\n)\leq \exp(\log(N(\q))^2)}}\mu(\n)\rho_\pi(\n)N(\n)^{-\rho-\sigma_1-it}\Big|^2 dt\\
\ll_\varepsilon T^DY^{2\eta}N(\q)^\varepsilon {\rm vol}(I(\q)) \sup_{N\geq X} \Big\{   \big(N+N^{1/2}N(\q)\big)\sum_{ N(\n)\geq N} \mu(\n)^2N(\n)^{-2-2\varepsilon_1}\Big\}
\\ \ll T^D Y^{2\eta}
\end{multline}

pour peu qu'on choisisse $\varepsilon<\varepsilon_1 $  ($ \varepsilon_1$ \'etant choisi, on le rappelle, tel que $\sigma_1=1-\Re(\rho)+\varepsilon_1$ donne la ligne d'int\'egration de la premi\`ere int\'egrale).
 
 Ayant \`a disposition d\'esormais des estimations pour les deux int\'egrales majorant le cardinal de $\Pi_\q^{\g{k}}(\eta,T)$, on peut donc conclure que
$$\Pi_\q^{\g{k}}(\eta,T)\ll T^D\Big( Y^{\eta-1/4}X^{1/2}N(\q)^{1/4+\varepsilon} + Y^{2\eta}\Big).$$
Reste \`a choisir $ X=N(\q)^A\, , \, Y=N(\q)^B$. Tout ce dont on a besoin est de maintenir:
$$ (\eta -\frac{1}{4})B+\frac{A}{2}+\frac{1}{4}+\varepsilon<0$$ 
ce qui est clairement possible, pour $B $ assez grand, tant que $\eta<\eta_0<1/4$ -- et en fait, en utilisant Ramanujan, il on a seulement besoin que $\eta_0<1/2$. Cela ach\`eve donc la preuve du th\'eor\`eme de densit\'e (\ref{dens}).


Stanford University, Department of Mathematics, building 380, Stanford, California 94305, USA.\\
E-mail: \href{mailto:trotabas@math.stanford.edu}{\texttt{trotabas@math.stanford.edu}}

\end{document}